\documentclass{acta_lim}
\usepackage{amssymb,pifont}
\usepackage{amsmath}
\usepackage{eucal}
\usepackage{mathrsfs} 
\usepackage{mathtools}
\usepackage{multicol}
\usepackage[shortlabels]{enumitem}
\usepackage{url,adjustbox}
\usepackage[normalem]{ulem}
\usepackage{mathscinet}
\usepackage[all]{xy}
\usepackage{stackrel}

\usepackage[cmyk]{xcolor}

\usepackage[longnamesfirst]{natbib}
\usepackage{har2nat}
\setcitestyle{comma,aysep={}}
\shortcites{Demaine,Lorentz,Etingof}

\usepackage{newtxtext}
\usepackage{newtxmath}
\DeclareSymbolFont{CMlargesymbols}{OMX}{cmex}{m}{n}
\DeclareMathDelimiter{(}{\mathopen} {operators}{"28}{CMlargesymbols}{"00}
\DeclareMathDelimiter{)}{\mathclose}{operators}{"29}{CMlargesymbols}{"01}

\pagerange{\pageref{firstpage}--\pageref{lastpage}}
\allowdisplaybreaks
\doi{21000076}
\usepackage[UKenglish]{babel}

\usepackage[twoside,
paperheight=247mm,
paperwidth=174mm,
marginratio={3:5,2:3},
left=18mm,
top=20mm]{geometry}
\usepackage[hang]{subfigure}

\newcommand*{\ie}{{\it i.e.}}
\newcommand*{\tie}{{that~is}}
\newcommand*{\eg}{{\it e.g.}}
\numberwithin{figure}{section}
\numberwithin{table}{section}
\newcommand{\rme}{{\mathrm{e}}}
\newcommand{\rmi}{{\mathrm{i}}}
\newcommand{\dd}{{\mathrm{d}}}
\newcommand{\mypara}[1]{\medskip\noindent{\bf #1.}}
\newcommand{\hspz}{\hspace{-6pt}}
\newcommand{\zza}{\hspace{4.7pt}}
\newcommand{\zzb}{\hspace{3.6pt}}
\newcommand{\zzc}{\hspace{5.5pt}}

\makeatletter
\newsavebox{\@brx}
\newcommand{\llangle}[1][]{\savebox{\@brx}{\(\m@th{#1\langle}\)}%
  \mathopen{\copy\@brx\kern-0.5\wd\@brx\usebox{\@brx}}}
\newcommand{\rrangle}[1][]{\savebox{\@brx}{\(\m@th{#1\rangle}\)}%
  \mathclose{\copy\@brx\kern-0.5\wd\@brx\usebox{\@brx}}}
\makeatother

\DeclareFontFamily{U}{wncy}{}
\DeclareFontShape{U}{wncy}{m}{n}{<->wncyr10}{}
\DeclareSymbolFont{mcy}{U}{wncy}{m}{n}
\DeclareMathSymbol{\Sh}{\mathord}{mcy}{"58} 

\usepackage{tikz} 
\usetikzlibrary{perspective}
\usetikzlibrary{3d}
\usetikzlibrary{cd}

\usepackage{tkz-graph}
\usepackage{pgfplots}
\usepackage{tikzlings}

\title[Tensors in computations]{Tensors in computations}
\author[L.-H.~Lim]{%
Lek-Heng Lim\\
Computational and Applied Mathematics Initiative,\\
University of Chicago,
Chicago, IL~60637, USA\\
E-mail: {lekheng@uchicago.edu}
}

\newtheorem{theorem}{Theorem}[section]
\newtheorem{definition}[theorem]{Definition}
\newtheorem{example}[theorem]{Example}

\DeclareMathOperator{\sgn}{sgn}
\DeclareMathOperator{\spn}{span}
\DeclareMathOperator{\im}{im}
\DeclareMathOperator{\tr}{tr}
\DeclareMathOperator{\diag}{diag}
\DeclareMathOperator{\GA}{GA}
\DeclareMathOperator{\GL}{GL}
\DeclareMathOperator{\SL}{SL}
\DeclareMathOperator{\SOr}{SO}
\DeclareMathOperator{\SU}{SU}
\DeclareMathOperator{\Un}{U}
\DeclareMathOperator{\Or}{O}
\DeclareMathOperator{\Sp}{Sp}
\DeclareMathOperator{\SE}{SE}
\DeclareMathOperator{\E}{E}
\DeclareMathOperator{\Diff}{Diff}
\DeclareMathOperator{\Iso}{Iso}
\DeclareMathOperator{\Symp}{Symp}
\DeclareMathOperator{\Fn}{\mathcal{F}}
\DeclareMathOperator{\Lin}{\mathcal{L}}
\DeclareMathOperator{\Mult}{\mathcal{M}}
\DeclareMathOperator{\Bd}{\mathcal{B}}
\DeclareMathOperator{\Ten}{\mathsf{T}}
\DeclareMathOperator{\Det}{Det}
\DeclareMathOperator{\rank}{rank}
\DeclareMathOperator{\mrank}{\mu rank}
\DeclareMathOperator{\brank}{\overline{\rank}}
\DeclareMathOperator{\nnz}{nnz}

\DeclareMathOperator{\seg}{\sigma_{\otimes}}
\DeclareMathOperator{\segt}{\sigma_{\mathbb{T}}}
\DeclareMathOperator{\vect}{vec}

\newcommand{\hatotimes}{\mathbin{\widehat{\otimes}}}

\newcommand{\tp}{{\scriptscriptstyle\mathsf{T}}}
\newcommand{\F}{{\scriptscriptstyle\mathsf{F}}}
\newcommand{\G}{{\scriptscriptstyle\mathsf{G}}}
\newcommand{\Beta}{\mathrm{B}}
\newcommand{\Mu}{\mathrm{M}}

\DeclareSymbolFont{bbold}{U}{bbold}{m}{n}
\DeclareSymbolFontAlphabet{\mathbbold}{bbold}
\newcommand{\bbone}{\mathbbold{1}}

\newcommand{\p}{{\scriptscriptstyle+}}
\newcommand{\pp}{{\scriptscriptstyle++}}
\newcommand{\ppp}{{\scriptscriptstyle+++}}
\newcommand{\mh}{{\scriptscriptstyle-1/2}}
\newcommand{\ph}{{\scriptscriptstyle1/2}}
\newcommand{\mrt}{{\scriptscriptstyle-1/\sqrt{3}}}
\newcommand{\prt}{{\scriptscriptstyle1/\sqrt{3}}}

\pgfkeys{/thing/.unknown/.code=\relax}

\begin{document}

\label{firstpage}
\maketitle

\begin{abstract}
The notion of a tensor captures three great ideas: equivariance,
multilinearity, separability. But trying to be three things at once makes
the notion difficult to understand. We will explain tensors in an
accessible and elementary way through the lens of linear algebra and
numerical linear algebra, elucidated with examples from computational and
applied mathematics.
\end{abstract}

\tableofcontents

\section{Introduction}\label{sec:intro}

We have two goals in this article: the first is to explain in detail and in the simplest possible terms what a tensor is; the second is to discuss the main ways in which tensors play a role in computations. The two goals are interwoven: what defines a tensor is also what makes it useful in computations, so it is important to gain a genuine understanding of tensors. We will take the reader through the three common definitions of a tensor: as an object that satisfies certain transformation rules, as a multilinear map, and as an element of a tensor product of vector spaces. We will explain the motivations behind these definitions, how one definition leads to the next, and how they all fit together. All three definitions are useful in computational mathematics but in different ways; we will intersperse our discussions of each definition with considerations of how it is employed in computations, using the latter as impetus for the former.

Tensors arise in computations in one of three ways: equivariance under coordinate changes, multilinear maps and separation of variables, each corresponding to one of the aforementioned definitions.

\mypara{Separability} The last of the three is the most prevalent, arising in various forms in problems, solutions and algorithms. One may exploit separable structures in ordinary differential equations, in integral equations, in Hamiltonians, in objective functions, \emph{etc.} The separation-of-variables technique may be applied to solve partial differential and integral equations,  or even finite difference and integro-differential equations. Separability plays an indispensable role in Greengard and Rokhlin's fast multipole method, Grover's quantum search algorithm, Hartree--Fock approximations of various stripes, Hochbaum and Shanthikumar's non-linear separable convex optimization, Smolyak's quadrature, and beyond. It also underlies tensor product construction of various objects -- bases, frames, function spaces, kernels, multiresolution analyses, operators and quadratures among them.

\mypara{Multilinearity} Many common operations, ranging from multiplying two complex numbers  or matrices to the convolution or bilinear Hilbert transform of functions, are bilinear operators. Trilinear functionals arise in self-concordance and numerical stability, and in fast integer and fast matrix multiplications. The last leads us to the matrix multiplication exponent, which characterizes the computational complexity of not just matrix multiplication but also inversion, determinant, null basis, linear systems, LU/QR/eigenvalue/Hessenberg decompositions, and more. Higher-order multilinear maps arise in the form of multidimensional Fourier/Laplace/Z/discrete cosine transforms and as cryptographic multilinear maps. Even if a multivariate function is not multilinear, its derivatives always are, and thus Taylor and multipole series are expansions in multilinear maps.

\mypara{Equivariance}  Equivariance is an idea as old as tensors and has been implicitly used throughout algorithms in numerical linear algebra. Every time we transform the scalars, vectors or matrices in a problem from one form to another via a sequence of Givens or Jacobi rotations, Householder reflectors or Gauss transforms, or through a Krylov or Fourier or wavelet basis, we are employing equivariance in the form of various $0$-, $1$- or $2$-tensor transformation rules. Beyond numerical linear algebra, the equivariance of tools from interior point methods to deep neural networks is an important reason for their success. The equivariance of Newton's method allows it to solve arbitrarily ill-conditioned optimization problems in theory and highly ill-conditioned ones in practice. \textsf{AlphaFold~2} conquered the protein structure prediction problem with an equivariant neural network.
\medskip

 Incidentally, out of the so-called `top ten algorithms of the twentieth century' \cite{top10}, fast multipole \cite{Board}, fast Fourier transform \cite{Rockmore}, Krylov subspaces \cite{VDV}, Francis's QR \cite{Parlett} and matrix decompositions \cite{Stewart}
all draw from these three aforementioned tensorial ideas in one way or another. Nevertheless, by `computations' we do not mean just algorithms, although these will constitute a large part of our article; the title of our article also includes the use of tensors as a tool in the analysis of algorithms (\eg\ self-concordance in polynomial-time convergence), providing intractability guarantees (\eg\ cryptographic multilinear maps), reducing complexity of models (\eg\  equivariant neural networks), quantifying computational complexity (\eg\ exponent of matrix multiplication) and in yet other ways. It would be prudent to add that while tensors are an essential ingredient in the aforementioned computational tools, they are rarely the \emph{only} ingredient: it is usually in combination with other concepts and techniques -- calculus of variations, Gauss quadrature, multiresolution analysis, the power method, reproducing kernel Hilbert spaces, singular value decomposition, {\em etc.} -- that they become potent tools.

The article is written with accessibility and simplicity in mind. Our exposition assumes only standard undergraduate linear algebra: vector spaces, linear maps, dual spaces, change of basis, {\em etc.} Knowledge of numerical linear algebra is a plus since this article mainly targets computational mathematicians. 
As physicists have played an outsize role in the development of tensors, it is inevitable that motivations for certain aspects, explanations for certain definitions, {\em etc.}, are best understood from a physicist's perspective and to this end we will include some discussions to provide context. When discussing applications or the physics origin of certain ideas, it is inevitable that we have to assume slightly more background, but we strive to be self-contained and limit ourselves to the most pertinent basic ideas. In particular, it is not our objective to provide a comprehensive survey of all things tensorial. We focus on the foundational and fundamental; results from current research make an appearance only if they illuminate these basic aspects. We have avoided a formal textbook-style  treatment and opted for a more casual exposition that hopefully makes for easier reading (and  in part to keep open the possibility of a future book project).

\subsection{Overview}

There are essentially three ways to define a tensor, reflecting the chronological evolution of the notion through the last 140 years or so:
\begin{dingautolist}{192} 
\item\label{st:tensor1} as a multi-indexed object that \emph{satisfies certain transformation rules},
\item\label{st:tensor2} as a multilinear map,
\item\label{st:tensor3} as an element of a tensor product of vector spaces.
\end{dingautolist}
The key to the coordinate-dependent definition in \ref{st:tensor1} is the emphasized part: the transformation rules are what define a tensor when one chooses to view it as a multi-indexed object, akin to the group laws in the definition of a group. 
The more modern coordinate-free definitions \ref{st:tensor2} and \ref{st:tensor3} automatically encode these transformation rules. The multi-indexed object, which could be a hypermatrix, a polynomial, a differential form, {\em etc.}, is then a coordinate representation of either \ref{st:tensor2} or \ref{st:tensor3} with respect to a choice of bases, and the corresponding transformation rule is a change-of-basis theorem.

Take the following case familiar to anyone who has studied linear algebra. Let $\mathbb{V}$ and $\mathbb{W}$ be finite-dimensional vector spaces over $\mathbb{R}$. A linear map $f \colon  \mathbb{V} \to \mathbb{W}$ is an order-$2$ tensor of covariant order $1$ and contravariant order $1$; this is the description according to definition~\ref{st:tensor2}. We may construct a new vector space $\mathbb{V}^* \otimes \mathbb{W}$, where $\mathbb{V}^*$ denotes the dual space of $\mathbb{V}$, and show that $f$ corresponds to a unique element $T \in \mathbb{V}^* \otimes \mathbb{W}$; this is the description according to definition~\ref{st:tensor3}. With a choice of bases $\mathscr{B}_\mathbb{V}$ and $\mathscr{B}_\mathbb{W}$, $f$ or $T$ may be represented as a matrix $A \in \mathbb{R}^{m \times n}$, where $m = \dim \mathbb{V}$ and $n = \dim \mathbb{W}$. Note, however, that the matrix $A$ is not unique and depends on our choice of bases, and different  bases $\mathscr{B}'_\mathbb{V}$ and $\mathscr{B}'_\mathbb{W}$ would give a different matrix representation $A' \in \mathbb{R}^{m \times n}$. The change-of-basis theorem states that the two matrices are related via the transformation rule\footnote{Equivalently, $A = XA'Y^{-1}$. In numerical linear algebra, there is a tendency to view this as matrix decomposition as opposed to change of basis.} $A' = X^{-1}A Y$,  where $X$ and $Y$ are the change-of-basis matrices on $\mathbb{W}$ and $\mathbb{V}$ respectively. Definition~\ref{st:tensor1} is essentially an attempt to define a linear map using its change-of-basis theorem  --  possible but  awkward.  The reason for such an awkward definition is one of historical necessity: definition~\ref{st:tensor1} had come before any of the modern notions we now take for granted in linear algebra  --  vector space, dual space, linear map, bases, {\em etc.}  --  were invented.

We stress that if one chooses to work with definition~\ref{st:tensor1}, then it is the transformation rule/change-of-basis theorem and not the multi-indexed object that defines a tensor. For example, depending on whether the transformation rule takes $A \in \mathbb{R}^{n \times n}$ to $XAY^{-1}$, $XAY^\tp$, $X^{-1} A Y^{-\tp}$, $XAX^{-1}$, $XAX^\tp$ or $X^{-1} AX^{-\tp}$, we obtain different tensors with entirely different properties. Also, while we did not elaborate on the change-of-basis matrices $X$ and $Y$, they play an important role in the transformation rule. If $\mathbb{V}$ and $\mathbb{W}$ are vector spaces without any additional structures, then $X$ and $Y$ are just required to be invertible; but if  $\mathbb{V}$ and $\mathbb{W}$ are, say, norm or inner product or symplectic vector spaces, then $X$ and $Y$ would need to preserve these structures too. More importantly, every notion we define, every property we study for any tensor  --  rank, determinant, norm, product, eigenvalue, eigenvector, singular value, singular vector, positive definiteness, linear system, least-squares problems, eigenvalue problems, {\em etc.}  --  must conform to the respective transformation rule. This is a point that is often lost; it is not uncommon to find mentions of `tensor such and such' in recent literature that makes no sense for a tensor.

As we will see, the two preceding paragraphs extend in a straightforward way to order-$d$ tensors (henceforth $d$-tensors) of contravariant order $p$ and covariant order $d-p$ for any integers $d \ge p \ge 0$, of which a linear map corresponds to the case $p=1$, $d=2$.

While discussing tensors, we will also discuss their role in computations. The most salient applications are often variations of the familiar \emph{separation-of-variables} technique that one encounters when solving ordinary and partial differential equations, integral equations or even integro-differential equations. Here the relevant perspective of a tensor is that in definition~\ref{st:tensor3}; we will see that $L^2(X_1) \otimes L^2(X_2) \otimes \dots \otimes L^2(X_d) = L^2(X_1 \times X_2 \times \dots \times X_d)$, \tie, multivariate functions are tensors in an appropriate sense. Take $d = 3$ for simplicity. Separation of variables, in its most basic manifestation, exploits functions $f \colon  X \times Y \times Z\to \mathbb{R}$ that take a multiplicatively separable form
\begin{equation}\label{eq:sep-1}
f(x,y,z) = \varphi (x)\psi(y) \theta (z),
\end{equation}
or, equivalently,
\begin{equation}\label{eq:sep0}
f = \varphi \otimes \psi \otimes \theta.
\end{equation}
For real-valued functions, this defines the tensor product $\otimes$, \tie, \eqref{eq:sep0} simply means \eqref{eq:sep-1}. A tensor taking the form in \eqref{eq:sep0} is called \emph{decomposable} or \emph{pure} or \emph{rank-one} (the last only if it is non-zero).
This deceptively simple version is already an essential ingredient in some of the most important  algorithms, such as Greengard and Rokhlin's fast multipole method, Grover's quantum search algorithm, Hochbaum and Shanthikumar's non-linear separable convex optimization, and more.

A slightly more involved version assumes that $f$ is a sum of separable functions,
\[
f(x,y,z) = \sum_{i=1}^r \varphi_i(x) \psi_i(y) \theta_i(z),
\]
or, equivalently,
\begin{equation}\label{eq:t0}
f = \sum_{i=1}^r \varphi_i \otimes \psi_i \otimes \theta_i,
\end{equation}
and this underlies various fast algorithms for evaluating bilinear operators, which is a $3$-tensor in the sense of definition~\ref{st:tensor2}. If one allows for a slight extension of \ref{st:tensor2} to modules (essentially a vector space where the field of scalars is a ring like $\mathbb{Z}$), then various fast integer multiplication algorithms such as those of  Karatsuba, Toom and Cook, Sch\"onhage and Strassen, and F\"urer may also be viewed in the light of \ref{st:tensor2}. But even with the standard definition of allowing only scalars from $\mathbb{R}$ or $\mathbb{C}$, algorithms that exploit \eqref{eq:t0} include Gauss's algorithm for complex multiplication, Strassen's algorithm for fast matrix multiplication/inversion and various algorithms for structured matrix--vector products (\eg\ Toeplitz, Hankel, circulant matrices). These algorithms show that a combination of both perspectives in definitions~\ref{st:tensor2} and \ref{st:tensor3} can be fruitful.

An even more involved version of the separation-of-variables technique allows more complicated structures on the indices such as
\[
f(x,y,z) = \sum_{i,j,k=1}^{p,q,r} \varphi_{ij}(x) \psi_{jk}(y) \theta_{ki}(z),
\]
or, equivalently,
\begin{equation}\label{eq:mps0}
f = \sum_{i,j,k=1}^{p,q,r}  \varphi_{ij} \otimes \psi_{jk} \otimes \theta_{ki}.
\end{equation}
This is the celebrated \emph{matrix product state} in the tensor network literature. Techniques such as DMRG simplify computations of eigenfunctions of Schr\"odinger operators by imposing such structures on the desired eigenfunction. Note that \eqref{eq:mps0}, like \eqref{eq:t0}, is also a decomposition into a sum of separable functions but the indices are captured by the following graph: 
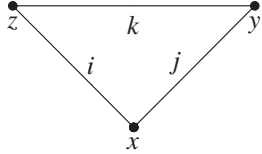
\[
\begin{tikzpicture}[scale=0.8]
\filldraw
    (1.75,0.75) node[align=center, below] {}
    (0,0) circle (2pt) node[align=center, below] {\small{$x$}}
     -- (0.7,0.7) circle (0pt) node[align=center, above] {\small{$j$}}
 -- (2,2) circle (2pt)  node[align=center, below] {\small{$y$}}
  -- (0,2) circle (0pt) node[align=center, below] {\small{$k$}}
 -- (-2,2) circle (2pt)  node[align=center, below] {\small{$z$}}  
  -- (-0.7,0.7) circle (0pt) node[align=center, above] {\small{$i$}}
 -- (0,0) circle (2pt) node[align=center, below]{}% {$\varphi$}
-- cycle;
\end{tikzpicture}
\]
More generally, for any undirected graph $G$ on $d$ vertices, we may decompose a $d$-variate function into a sum of separable functions in a similar manner whereby the indices are summed over the edges of $G$ and the different variables correspond to the different vertices.  Such representations of a tensor are often called \emph{tensor network} states and the graph $G$ a tensor network. They often make an effective ansatz
for analytical and numerical solutions of PDEs arising from quantum chemistry.

The oldest definition of a tensor, \ie\ definition~\ref{st:tensor1}, encapsulates a fundamental notion, namely that of \emph{equivariance} (and also \emph{invariance}, as a special case) under coordinate transformations, which continues to be relevant in modern applications. For example, the affine invariance of Newton's method and self-concordant barrier functions plays an important role in interior point methods for convex optimization \cite{Nesterov}; exploiting equivariance in deep convolutional neural networks has led to groundbreaking performance on real-world tasks such as classifying images in the \texttt{CIFAR10} database \cite{Welling} and, more recently, predicting the shape of
protein molecules in the \texttt{CASP14} protein-folding competition \label{pg:alpha}\cite{Nature,Science}.\footnote{Research papers that describe \textsf{AlphaFold~2} are still unavailable at the time of writing but the fact that it uses an equivariant neural network \cite{Fuchs} may be found in Jumper's slides at \href{https://predictioncenter.org/casp14/doc/presentations/2020_12_01_TS_predictor_AlphaFold2.pdf}{https://predictioncenter.org/casp14/doc/presentations/ 2020\_12\_01\_TS\_predictor\_AlphaFold2.pdf}.} More generally, we will see that the notion of equivariance has guided the design of algorithms and methodologies through the ages, and is implicit in many classical algorithms in numerical linear algebra.

The three notions we raised  --  separability, multilinearity, equivariance  --  largely account for the usefulness of tensors in computations, but there are also various measures of the `size' of a tensor, notably tensor ranks and tensor norms, that play an indispensable auxiliary role. Of these, two notions of rank and three notions of norm stand out for their widespread applicability. Given the decomposition \eqref{eq:tr0}, a natural way to define rank is evidently
\begin{equation}\label{eq:tr0}
\rank(f) \coloneqq \min \biggl\{ r \colon   f= \sum_{i=1}^r \varphi_i \otimes \psi_i \otimes \theta_i \biggr\},
\end{equation}
and a natural way to define a corresponding norm is
\begin{equation}\label{eq:nn0}
\lVert f \rVert_\nu \coloneqq \inf \biggl\{ \sum_{i=1}^r \ \lVert\varphi_i\rVert \lVert \psi_i\rVert \lVert \theta_i \rVert \colon   f = \sum_{i=1}^r \varphi \otimes \psi_i \otimes \theta_i \biggr\}.
\end{equation}
If we have an inner product, then its dual norm is
\begin{equation}\label{eq:sn0}
\lVert f \rVert_{\sigma} \coloneqq \sup \dfrac{\lvert \langle f, \varphi \otimes \psi \otimes \theta \rangle\rvert}{ \lVert\varphi\rVert \lVert \psi\rVert \lVert \theta \rVert}.
\end{equation}
Definitions \eqref{eq:tr0}, \eqref{eq:nn0} and \eqref{eq:sn0} are the \emph{tensor rank}, \emph{nuclear norm} and \emph{spectral norm} of $f$ respectively. They generalize to tensors of arbitrary order $d$ in the obvious way, and for $d=2$ reduce to rank, nuclear norm  and spectral norm of linear operators and matrices. The norms in \eqref{eq:nn0} and \eqref{eq:sn0} are in fact the two standard ways to turn a tensor product of Banach spaces into a Banach space in functional analysis, and in this context they are called projective and injective norms respectively.

From a computational perspective,  we will see in Sections~\ref{sec:mult} and \ref{sec:tenprod} that in many instances there is a tensor, usually of order three, at the heart of a computational problem, and  issues regarding complexity, stability, approximability, {\em etc.}, of that problem can often be translated to finding or bounding \eqref{eq:tr0}, \eqref{eq:nn0} and \eqref{eq:sn0}.  One issue with \eqref{eq:tr0}, \eqref{eq:nn0} and \eqref{eq:sn0} is that they are all NP-hard to compute when $d \ge 3$. As computationally tractable alternatives, we also consider the \emph{multilinear rank}
\[
\mrank(f) \coloneqq \min \biggl\{ (p,q,r) \colon    \sum_{i,j,k=1}^{p,q,r} \varphi_i \otimes \psi_j \otimes \theta_k \biggr\},
\]
and the \emph{Frobenius norm} (also called the \emph{Hilbert--Schmidt norm})
\[
\lVert f \rVert_\F =\sqrt{\langle f, f\rangle}.
\]
While there are many other notions in linear algebra such as eigenvalues and eigenvectors, singular values and singular vectors, determinants, characteristic polynomials, {\em etc.}, that do generalize to higher orders $d \ge 3$, their value in computation is somewhat limited relative to ranks and norms. As such we will only mention these other notions in passing.

\subsection{Goals and scope}

We view our article in part as an elementary alternative to the heavier treatments in \citet{BCS} and \citet{Land1,Land2,Land3}, all excellent treatises on tensors with substantial discussions of their role in computations. Nevertheless they require a heavy dose of algebra and algebraic geometry, or at least a willingness to learn these topics. They also tend to gloss over analytic aspects such as tensor norms. There are also excellent treatments such as those of \citet{Cheney,Defant,Diestel} and \citet{Ryan} that cover the analytic perspective, but they invariably focus on tensors in infinite-dimensional spaces, and with that comes the many complications that can be avoided in a finite-dimensional setting. None of the aforementioned references are easy reading, being aimed more at specialists than beginners, with little discussion of elementary questions.

Our article takes a middle road, giving equal footing to both algebraic and analytic aspects, \tie, we discuss tensors over vector spaces (algebraic) as well as tensors over norm or inner product spaces (analytic), and we explain why they are different and how they are related. As we target readers whose main interests are in one way or another related to computations --  numerical linear algebra, numerical PDEs, optimization, scientific computing, theoretical computer science, machine learning, {\em etc.}\ -- and such a typical target reader would tend to be far more conversant with analysis than with algebra, the manner in which we approach algebraic and analytic topics is calibrated accordingly. We devote considerable length to explaining and motivating algebraic notions such as modules or commutative diagrams, but tend to gloss over analytic ones such as distributions or compact operators. We assume some passing familiarity with computational topics such as quadrature or von~Neumann stability analysis but none with `purer' ones such as formal power series or Littlewood--Richardson coefficients. We also draw almost all our examples from topics close to the heart of computational mathematicians: kernel SVMs, Krylov subspaces, Hamilton equations, multipole expansions, perturbation theory, quantum chemistry and physics, semidefinite programming, wavelets, {\em etc.} As a result almost all our examples tend to have an analytic bent, although there are also exceptions such as cryptographic multilinear maps or equivariant neural networks, which are more algebraic.

Our article is intended to be completely elementary. Modern studies of tensors largely fall into four main areas of mathematics: algebraic geometry, differential geometry, representation theory and functional analysis. We have avoided the first three almost entirely and have only touched upon the most rudimentary aspects of the last. The goal is to show that even without bringing in highbrow subjects, there is still plenty that could be said about tensors; all we really need is standard undergraduate linear algebra and some multivariate calculus -- vector spaces, linear maps, change of basis, direct sums and products, dual spaces, derivatives as linear maps, {\em etc.} Among this minimal list of prerequisites, one would find that dual spaces tend to be a somewhat common blind spot of our target readership, primarily because they tend to work over inner product spaces where duality is usually a non-issue. Unfortunately, duality cannot be avoided if one hopes to gain any reasonably complete understanding of tensors, as tensors defined over inner product spaces represent a relatively small corner of the subject. A working familiarity with dual spaces is a must: how a dual basis behaves under a change of basis, how to define the transpose and adjoint of abstract linear maps, how taking the dual changes direct sums into direct products, {\em etc.} Fortunately these are all straightforward, and we will remind readers of the relevant facts as and when they are needed. Although we will discuss tensors over infinite-dimensional spaces, modules, bundles, operator spaces, {\em etc.}, we only do so when they are relevant to computational applications. When discussing the various properties and notions pertaining to tensors, we put them against a backdrop of familiar analogues in linear algebra and numerical linear algebra, by either drawing a parallel or a contrast between the multilinear and the linear.

Our approach towards explaining the roles of tensors in computations is to present them as examples alongside our discussions of the three standard definitions of tensors on page~\pageref{st:tensor1}. The hope is that our reader would gain both an appreciation of these three definitions, each capturing a different feature of a tensor, as well as how each is useful in computations. While there are a very small number of inevitable overlaps with \citet{BCS} and \citet{Land1,Land2,Land3} such as fast matrix multiplication, we try to offer a different take on these topics. Also, most of the examples in our article cannot be found in~them. 

Modern textbook treatments of tensors in algebra and geometry tend to sweep definition~\ref{st:tensor1} under the rug of the coordinate-free approaches in definitions~\ref{st:tensor2} and~\ref{st:tensor3}, preferring not to belabour the coordinate-dependent point of view in definition~\ref{st:tensor1}. For example, there is no trace of it in the widely used textbook on algebra by \citet{Lang} and that on differential geometry by \citet{Lee}. Definition~\ref{st:tensor1}  tends to be found only in much older treatments \cite{Borisenko,Brand,Hay,Lovelock,McConnell,Michal,Schouten,Spain,Synge,Wrede} that are often more confusing than illuminating to a modern reader. This has resulted in a certain amount of misinformation. As a public service, we will devote significant space to definition~\ref{st:tensor1} so as to make it completely clear. For us, belabouring the point is the point. By and large we  will discuss definitions~\ref{st:tensor1}, \ref{st:tensor2} and \ref{st:tensor3} in chronological order, as we think there is some  value in seeing how the concept of tensors evolved through time. Nevertheless we will occasionally get ahead of ourselves to point out how a later definition would clarify an earlier one or revisit an earlier definition with the hindsight of a later one.

There is no lack of reasonably elementary, concrete treatments of tensors in physics texts \cite{AMR,Martin,MTW,Wald} but the ultimate goal in these books is invariably \emph{tensor fields} and not tensors. The pervasive use of Einstein's summation convention,\footnote{No doubt a very convenient shorthand, and its avoidance of the summation symbol $\sum$ a big  typesetting advantage in the pre-\TeX\ era.} the use of both upper and lower indices (sometimes also left and right, \ie\ indices on all four corners)  in these books can also be confusing, especially when one needs to switch between Einstein's conventions and the usual conventions in linear algebra. Do the raised indices refer to the inverse of the matrix or are they there simply because they needed to be summed over? How does one tell the inverse apart from the inverse transpose of a matrix? The answers depend on the source, somewhat defeating the purpose of a convention. In this article, tensor fields are mentioned only in passing and we do not use Einstein's summation convention. Instead, we frame our discussions of tensors in notation standard in linear algebra and numerical linear algebra or slight extensions of these and, by so doing, hope that our readers would find it easier to see where tensors fit in.

Our article is not about algorithms for solving various problems related to higher-order tensors. Nevertheless, we  take the opportunity to make a few points in Section~\ref{sec:comp} about the prospects of such algorithms. The bottom line is that we do not think such problems, even in low-order ($d =3$ or $4$) low-dimensional ($n < 10$) instances, have efficient or even just provably correct algorithms. This is unfortunate as there are various `tensor models' for problems in areas from quantum chemistry to data analysis that are based on the optimistic assumption that there exist such algorithms. 
It is fairly common to see claims of extraordinary things that follow from such `tensor models',  but upon closer examination one would invariably find an NP-hard problem such~as
\[
\min_{\varphi,\psi,\theta}\; \lVert f - \varphi \otimes \psi \otimes \theta \rVert
\]
embedded in the model as a special case. If one could efficiently solve this problem, then the range of earth-shattering things that would follow is limitless.\footnote{In the words of a leading authority  \cite[pp.~ix and 11]{fortnow}, `society as we know it would change dramatically, with immediate, amazing advances in medicine, science, and entertainment and the automation of nearly every human task,' and `the world will change in ways that will make the Internet seem like a footnote in history.'} Therefore, these claims of extraordinariness are hardly surprising given such an enormous caveat; they are simply consequences of the fact that if one could efficiently solve \emph{any} NP-hard problems, one could efficiently solve all NP-complete problems.

Nearly all the materials in this article are classical. It could have been written twenty years ago.  We would not have been able to mention the resolution of the salmon conjecture or discuss applications such as \textsf{AlphaFold}, the matrix multiplication exponent and complexity of integer multiplication would be slightly higher, and some of the books cited would be in their earlier editions. But 95\%, if not more, of the content would have remained intact. We limit our discussion in this article to results that are by-and-large rigorous, usually in the mathematical sense of having a proof but occasionally in the sense of having extensive predictions consistent with results of scientific experiments (like the effectiveness of DMRG). While this article contains no original research,  we would like to think that it offers abundant new insights on existing topics: our treatment of multipole expansions in Example~\ref{eg:multpole}, separation of variables in Examples~\ref{eg:sepabs}--\ref{eg:ide}, stress tensors in Example~\ref{eg:stress}, our interpretation of the universal factorization property in Section~\ref{sec:tensor3c}, discussions of the various forms of higher-order derivatives in Examples~\ref{eg:hod}, \ref{eg:hog}, \ref{eg:linhod}, the way we have presented and motivated the tensor transformation rules in Section~\ref{sec:trans}, {\em etc.}, have never before appeared elsewhere to the best of our knowledge. In fact, we have made it a point to not reproduce anything verbatim from existing literature; even standard materials are given a fresh take.

\subsection{Notation and convention}

Almost all results, unless otherwise noted, will hold for tensors over both $\mathbb{R}$ and $\mathbb{C}$ alike. In order to not have to state `$\mathbb{R}$ or $\mathbb{C}$' at every turn, we will denote our field as $\mathbb{R}$ with the implicit understanding that all discussions remain true over $\mathbb{C}$ unless otherwise stated. We will adopt numerical linear algebra convention and regard  vectors in $\mathbb{R}^n$ as column vectors, \ie\ $\mathbb{R}^n \equiv \mathbb{R}^{n \times 1}$.  Row vectors in $\mathbb{R}^{1 \times n}$ will always be denoted as $v^\tp$ with $v \in \mathbb{R}^n$. To save space, we also adopt the convention that an $n$-tuple delimited with parentheses always denotes a column vector, \tie,
\[
(a_1,\ldots, a_n) = \begin{bmatrix}a_1 \\ \vdots \\ a_n \end{bmatrix} \in \mathbb{R}^n.
\]
All bases in this article will be ordered bases. We write $\GL(n)$ for the general linear group, \tie, the set of $n \times n$ invertible matrices and $\Or(n)$ for the orthogonal group of $n \times n$ orthogonal matrices. Other matrix groups will be introduced later. For any $X \in \GL(n)$, we write
\[
X^{-\tp} \coloneqq (X^{-1})^\tp = (X^\tp)^{-1}.
\]

The vector space of polynomials in variables $x_1,\dots, x_n$ with real coefficients will be denoted by $\mathbb{R}[x_1,\dots,x_n]$. For Banach spaces $\mathbb{V}$ and $\mathbb{W}$, we write $\Bd(\mathbb{V};\mathbb{W})$ for the set of all bounded linear operators $\Phi \colon \mathbb{V} \to \mathbb{W}$. We abbreviate $\Bd(\mathbb{V};\mathbb{V})$ as $\Bd(\mathbb{V})$.  A word of caution may be in order. We use the term `separable' in two completely different and unrelated senses with little risk of confusion: a Banach or Hilbert space is separable if it has a countable dense subset whereas a function is separable if it takes the form in \eqref{eq:sep-1}.

We restrict some of our discussions to $3$-tensors for a few reasons: (i) the sacrifice in generality brings about an enormous reduction in notational clutter, (ii) the generalization to $d$-tensors for $d >3$ is straightforward once we present the $d=3$ case, (iii) $d =3$ is the first unfamiliar case for most readers given that $d = 0,1,2$ are already well treated in linear algebra, and (iv) $d=3$ is the most important case for many of our applications to computational mathematics.

\section{Tensors via transformation rules}\label{sec:trans}

We will begin with definition~\ref{st:tensor1}, which first appeared in a book on crystallography   by the physicist Woldemar \citet{Voigt}:
\begin{quote}\label{Voigt}
An abstract entity represented by an array of components
that are functions of coordinates such that, under a transformation of
coordinates, the new components are related to the transformation and to the
original components in a definite way.
\end{quote}
While the idea of a tensor had appeared before in works of Cauchy and Riemann \cite{Conrad,Yau} and  also around the same time in \citet{Ricci}, this is believed to be the earliest appearance of the word `tensor'. Although we will refer to the quote above as `Voigt's definition', it is not a direct translation of a specific sentence from \citet{Voigt} but a paraphrase  attributed to \citet{Voigt} in the {\em Oxford English Dictionary}. Voigt's definition is essentially the one adopted in all early textbooks on tensors such as \citet{Brand,Hay,Lovelock,McConnell,Michal,Schouten,Spain,Synge} and \citet{Wrede}.\footnote{These were originally published in the 1940s--1960s but all have been reprinted by Dover Publications and several remain in print.} This is not an easy definition to work with and likely contributed to the reputation of tensors being a tough subject to master. Famously, Einstein struggled with tensors \cite{Einstein} and the definition he had to dabble with would invariably have been definition~\ref{st:tensor1}.   Nevertheless, we should be reminded that linear algebra as we know it today was an obscure art in its infancy in the 1900s when Einstein was learning about tensors. Those trying to learn about tensors in modern times enjoy the benefit of a hundred years of pedagogical progress. By building upon concepts such as vector spaces, linear transformations, change of basis, {\em etc.}, that we take for granted today but were not readily accessible a century ago, the task of explaining tensors is significantly simplified.

In retrospect, the main issue with definition~\ref{st:tensor1} is that it is a `physicist's definition', \tie, it attempts to define a quantity by describing the change-of-coordinates rules that the quantity must satisfy without specifying the quantity itself. This approach may be entirely natural in physics where one is interested in questions such as `Is stress a tensor?' or `Is electromagnetic field strength a tensor?', \tie, the definition is always applied in a way where some physical quantity such as stress takes the place of the unspecified quantity, but it makes for an awkward definition in mathematics. The modern definitions~\ref{st:tensor2} and~\ref{st:tensor3} remedy this by stating unequivocally what this unspecified quantity is.

\subsection{Transformation rules illustrated with linear algebra}\label{sec:transrules}

The defining property of a tensor in definition~\ref{st:tensor1} is the `transformation rules' alluded to in the second and third lines of Voigt's definition, \tie, how the coordinates transform under change of basis. In modern parlance, these transformation rules express the notion of \emph{equivariance} under a linear action of a matrix group, an idea that has proved to be important even in trendy applications such as deep convolutional and attention networks. There is no need to go to higher order tensors to appreciate them; we will see that all along we have been implicitly working with these transformation rules in linear algebra and numerical linear algebra. Instead of stating them in their most general forms right away, we will allow readers to discover these transformation rules for themselves via a few simple familiar examples; the meaning of the terminology appearing in italics may be surmised from the examples and will be properly defined in due course.

\mypara{Eigenvalues and eigenvectors}\label{eg:eig} An eigenvalue $\lambda \in \mathbb{C}$ and a corresponding eigenvector $0 \ne v \in \mathbb{C}^n$ of a matrix $A \in\mathbb{C}^{n \times n}$ satisfy the eigenvalue/eigenvector equation $Av =\lambda v$. For any invertible $X \in \mathbb{C}^{n \times n}$, 
\begin{equation}\label{eq:ev}
(XAX^{-1}) Xv = \lambda Xv.
\end{equation}
The transformation rules here are $A \mapsto XAX^{-1}$, $v \mapsto Xv$ and $\lambda \mapsto \lambda$. An eigenvalue is an \emph{invariant} $0$-tensor and an eigenvector is a \emph{contravariant} $1$-tensor of a \emph{mixed} $2$-tensor, \tie, eigenvalues of $A$ and $XAX^{-1}$ are always the same whereas an eigenvector $v$ of $A$ transforms into an eigenvector $Xv$ of $XAX^{-1}$. This says that eigenvalues and eigenvectors are defined for mixed $2$-tensors.

\mypara{Matrix product} For matrices $A \in \mathbb{C}^{m \times n}$ and $B \in \mathbb{C}^{n \times p}$, we have
\begin{equation}\label{eq:matprod}
(XAY^{-1})(YBZ^{-1}) = X(AB)Z^{-1}
\end{equation}
for any invertible $X \in \mathbb{C}^{m \times m}$, $Y \in \mathbb{C}^{n \times n}$, $Z \in \mathbb{C}^{p \times p}$. The transformation rule is again of the form $A \mapsto XAY^{-1}$. This says that the standard matrix--matrix product is defined on mixed $2$-tensors.

\mypara{Positive definiteness} A matrix $A \in \mathbb{R}^{n \times n}$ is positive definite if and only if
\[
XAX^\tp \quad\text{or}\quad X^{-\tp}AX^{-1}
\]
is positive definite for any invertible $X \in \mathbb{R}^{n \times n}$. This says that positive definiteness is a property of covariant $2$-tensors that have transformation rule $X \mapsto XAX^\tp$ or of \emph{contravariant} $2$-tensors that have transformation rule $X \mapsto X^{-\tp}AX^{-1}$.

\mypara{Singular values and singular vectors} A singular value $\sigma \in \mathbb{R}$ and its corresponding left singular vector $0 \ne u \in \mathbb{R}^m$ and right singular vector $0 \ne v \in \mathbb{R}^n$ of a matrix $A \in\mathbb{R}^{m \times n}$ satisfy the singular value/singular vector equation
\[
\biggl\{
\begin{aligned}
Av &=\sigma u,\\
A^\tp u & =  \sigma v.
\end{aligned}
\]
For any \emph{orthogonal} matrices $X \in \mathbb{R}^{m \times m}$,  $Y \in \mathbb{R}^{n \times n}$, we have
\[
\biggl\{
\begin{aligned}
(XAY^\tp) Yv &=\sigma Xu,\\
(XAY^\tp)^\tp Xu & =  \sigma Yv.
\end{aligned}
\]
Singular values are invariant and singular vectors are contravariant as they transform as $\sigma \mapsto \sigma$, $u \mapsto Xu$, $v \mapsto Yv$ for a matrix that transforms as $A \mapsto XAY^\tp$. This tells us that singular values are \emph{Cartesian} $0$-tensors, left and right singular vectors are contravariant Cartesian $1$-tensors, defined on contravariant Cartesian $2$-tensors. `Cartesian' means that these transformations use orthogonal matrices instead of invertible ones.

\mypara{Linear equations} Let $A \in \mathbb{R}^{m \times n}$ and $b \in \mathbb{R}^m$. Clearly  $Av = b$ has a solution if and only if
\[
(XAY^{-1}) Yv = Xb
\]
has a solution for any invertible  $X \in \mathbb{R}^{m \times m}$ and $Y \in \mathbb{R}^{n \times n}$. If $m = n$, then we may instead want to consider
\[
(XAX^{-1}) Xv = Xb,
\]
and if furthermore $A$ is symmetric, then
\[
(XAX^\tp) X^{-\tp} v  = Xb ,
\]
where $X \in \mathbb{R}^{n \times n}$  is either invertible or orthogonal (for the latter $X^{-\tp} = X$). We will examine this ambiguity in transformation rules later. The solution in the last case is a covariant $1$-tensor, \tie, it  transforms as $v \mapsto X^{-\tp} v$ with an invertible~$X$.

\mypara{Ordinary and total least-squares} Let $A \in \mathbb{R}^{m \times n}$ and $b \in \mathbb{R}^m$. The ordinary least-squares problem is given by
\begin{equation}\label{eq:ls}
\min_{v\in \mathbb{R}^n}\; \lVert Av - b\rVert^2 = \min_{v\in \mathbb{R}^n}\; \lVert (XAY^{-1}) Yv - X b \rVert^2
\end{equation}
for any orthogonal $X \in \mathbb{R}^{m \times m}$ and invertible $Y \in \mathbb{R}^{n \times n}$. Unsurprisingly, the normal equation $A^\tp A v = A^\tp b$ has the same property:
\[
(XAY^{-1})^\tp (XAY^{-1} ) Yv = (XAY^{-1})^\tp Xb.
\]
The total least-squares problem is defined by
\begin{align*}
&\min \{ \lVert E \rVert^2 + \lVert r \rVert^2 \colon  (A+E)v = b+r \}\\*
&\quad =\min \{ \lVert XEY^\tp \rVert^2 + \lVert Xr \rVert^2 \colon  ( XAY^\tp+ XE Y^\tp) Yv = Xb+X r \}
\end{align*}
for any orthogonal $X \in \mathbb{R}^{m \times m}$ and orthogonal $Y \in \mathbb{R}^{n \times n}$. Here the minimization is over $E \in \mathbb{R}^{m\times n}$, $r \in \mathbb{R}^m$, and the constraint is interpreted as the linear system being consistent. Both ordinary and total least-squares are defined on tensors  --  minimum values transform as invariant $0$-tensors, $v,b,r$ as contravariant $1$-tensors and $A,E$ as mixed $2$-tensors  --  but the $2$-tensors involved are different as $Y$ is not required to be orthogonal in \eqref{eq:ls}.

\mypara{Rank, determinant, norm}\label{eg:rank} Let $A \in \mathbb{R}^{m \times n}$. Then
\begin{align*}
\rank(XAY^{-1}) = \rank(A), \quad
\det(XAY^{-1}) = \det(A), \quad
\lVert XAY^{-1} \rVert = \lVert A \rVert ,
\end{align*}
where $X$ and $Y$ are, respectively, invertible, special linear or orthogonal matrices. Here $\lVert \, \cdot \, \rVert$ may denote either the spectral, nuclear or Frobenius norm  and the determinant is regarded as identically zero whenever $m \ne n$. Rank, determinant and norm are defined on mixed $2$-tensors, special linear mixed $2$-tensors and Cartesian mixed $2$-tensors respectively.
\medskip

The point that we hope to make with these familiar examples is the following. The most fundamental and important concepts, equations, properties and problems in linear algebra and numerical linear algebra  --  notions that extend far beyond linear algebra into other areas of mathematics, science and engineering  --  all conform to tensor transformation rules. These transformation rules are not merely accessories to the definition of a tensor: they are the very crux of it and are what make tensors useful.

On the one hand, it is intuitively clear what the transformation rules are in the examples on pages~\pageref{eg:eig}--\pageref{eg:rank}: They are transformations that preserve either the \emph{form} of an equation (\eg\ if we write $A' = XAX^{-1}$, $v' = Xv$, $\lambda' = \lambda$, then \eqref{eq:ev} becomes $A' v' = \lambda' v'$), or the \emph{value} of a quantity such as rank/determinant/norm, or a \emph{property} such as positive definiteness. In these examples, the `multi-indexed object' in definition~\ref{st:tensor1} can be a matrix like $A$ or $E$, a vector like $u,v,b,r$ or a scalar like $\lambda$; the coordinates of the matrix or vector, \ie\ $a_{ij}$ and $v_i$, are sometimes called the `components of the tensor'. The matrices  $X, Y$ play a different role in these transformation rules and should be distinguished from the multi-indexed objects; for reasons that will be explained in Section~\ref{sec:multmaps}, we call them change-of-basis matrices.

On the other hand, these examples also show why the transformation rules can be confusing:  they are ambiguous in multiple ways.
\begin{enumerate}[\upshape (a)]
\setlength\itemsep{3pt}
\item\label{it:ambiguity1} The transformation rules can take several different forms. For example, $1$-tensors may  transform as
\[
v' = Xv\quad\text{or}\quad v' = X^{-\tp} v,
\]
$2$-tensors may transform as
\[
A' = XAY^{-1},\quad A' = XAY^\tp,\quad A' = XAX^{-1}, \quad A' = XAX^\tp,
\]
or yet other possibilities we have not discussed.

\item\label{it:ambiguity2} The change-of-basis matrices in these transformation rules may also take several different forms, most commonly invertible or orthogonal. In the examples, the change-of-basis matrices may be a single matrix\footnote{For those unfamiliar with this matrix group notation, it will be defined in \eqref{eq:groups}.}
\[
X \in \GL(n), \; \SL(n), \; \Or(n),
\]
and a pair of them
\[
(X,Y) \in \GL(m) \times \GL(n), \; \SL(m) \times \SL(n), \; \Or(m) \times \Or(n),
\]
or, as we saw in the case of ordinary least-squares \eqref{eq:ls}, $(X,Y) \in \Or(m) \times \GL(n)$. There are yet other possibilities for change-of-basis matrices we have not discussed, such as Lorentz, symplectic, unitary, {\em etc.}

\item\label{it:ambiguity3} Yet another ambiguity on top of \ref{it:ambiguity1} is that the roles of $v$ and $v'$ or $A$ and $A'$ are sometimes swapped and the transformation rules written~as
\[
v' = X^{-1}v,\quad v' = X^\tp v
\]
or
\[
A' = X^{-1}AY,\quad A' = X^{-1}AY^{-\tp},\quad A' = X^{-1}AX, \quad A' = X^{-1}AX^{-\tp}
\]
respectively. Note that the transformation rules here and those in \ref{it:ambiguity1} are all but identical: the only difference is in whether we label the multi-indexed object on the left or that on the right with a prime.
\end{enumerate}
We may partly resolve the ambiguity in  \ref{it:ambiguity1} by introducing covariance and contravariance type: tensors with transformation rules of the form $A' = XAX^\tp$ are covariant $2$-tensors or $(0,2)$-tensors; those of the form $A' = XAX^{-1}$ are mixed $2$-tensors or $(1,1)$-tensors; those of the form $A' = X^{-\tp}AX^{-1}$ are contravariant $2$-tensors or $(2,0)$-tensors. Invariance, covariance and contravariance are all special cases of \emph{equivariance} that we will discuss later. Nevertheless, we are still unable to distinguish between $A' = XAX^\tp$ and $A' = XAY^\tp$: both are legitimate transformation rules for covariant $2$-tensors.

These ambiguities \ref{it:ambiguity1}, \ref{it:ambiguity2}, \ref{it:ambiguity3} are the result of a  defect in definition~\ref{st:tensor1}. One ought to be asking: What is the tensor in definition~\ref{st:tensor1}? The answer is that it is actually missing from the definition. The multi-indexed object $A$ \emph{represents} the tensor whereas the transformation rules on $A$ \emph{defines} the tensor but the tensor itself has been left unspecified. This is a key reason why definition~\ref{st:tensor1} has been so confusing to early learners. Fortunately, it is easily  remedied by definition~\ref{st:tensor2} or \ref{st:tensor3}. We should, however, bear in mind that definition~\ref{st:tensor1} predated modern notions of vector spaces and linear maps, which are necessary for definitions~\ref{st:tensor2} and~\ref{st:tensor3}.

When we introduce definitions~\ref{st:tensor2} or~\ref{st:tensor3}, we will see that the ambiguity in \ref{it:ambiguity1}  is just due to different transformation rules for different tensors. Getting ahead of ourselves, for vector spaces $\mathbb{V}$ and $\mathbb{W}$, the transformation rule $A' = XAX^{-1}$ applies to tensors in $\mathbb{V} \otimes \mathbb{V}^*$ whereas $A' = XAY^{-1}$ applies to those in $\mathbb{V} \otimes \mathbb{W}^*$; $A' = XAX^\tp$ and $A' = X^{-\tp}AX^{-1}$ apply to tensors in $\mathbb{V} \otimes \mathbb{V}$ and $\mathbb{V}^* \otimes \mathbb{V}^*$ respectively. The matrices $A$ and $A'$ are representations of a tensor with respect to two different bases, and the ambiguity in \ref{it:ambiguity3} is just a result of which basis is regarded as the `old' basis and which as the `new' one.

The ambiguity in \ref{it:ambiguity2}  is also easily resolved with definitions~\ref{st:tensor2} or~\ref{st:tensor3}. The matrices $X$ and $Y$ are change-of-basis matrices on $\mathbb{V}$ and $\mathbb{W}$ and are always invertible, but they also have to preserve additional structures on $\mathbb{V}$ and $\mathbb{W}$ such as inner products or norms. For example,
singular values and singular vectors are only defined for inner product spaces and so we require the matrices $X$ and $Y$ to preserve inner products; for the Euclidean inner product, this simply means that $X$ and $Y$ are orthogonal matrices. Tensors with transformation rules involving orthogonal matrices are sometimes called `Cartesian tensors'.

These examples on pages~\pageref{eg:eig}--\pageref{eg:rank} illustrate why the transformation rules in definition~\ref{st:tensor1} are as crucial in mathematics as they are in physics. Unlike in physics, we do not use the transformation rules to check whether a physical quantity such as stress or strain is a tensor; instead we use them to ascertain whether an equation (\eg\ an eigenvalue/eigenvector equation), a property (\eg\  positive definiteness), an operation (\eg\ a matrix product) or a mathematical quantity (\eg\ rank) is defined on a tensor. In other words, when working with tensors, it is not a matter of simply writing down an expression involving multi-indexed quantities, because if the expression does not satisfy the transformation rules, it is undefined on a tensor. For example, multiplying two matrices via the  Hadamard (also known as Schur) product
\begin{equation}\label{eq:hp}
\begin{bmatrix}
a_{11} & a_{12} \\
a_{21} & a_{22}
\end{bmatrix} \circ
\begin{bmatrix}
b_{11} & b_{12} \\
b_{21} & b_{22}
\end{bmatrix}
=
\begin{bmatrix}
a_{11} b_{11}  &  a_{12} b_{12}\\
a_{21} b_{21} &  a_{22} b_{22}
\end{bmatrix}
\end{equation}
seems a lot more obvious than the standard matrix product
\begin{equation}\label{eq:mm}
\begin{bmatrix}
a_{11} & a_{12} \\
a_{21} & a_{22}
\end{bmatrix}
\begin{bmatrix}
b_{11} & b_{12} \\
b_{21} & b_{22}
\end{bmatrix}
=
\begin{bmatrix}
a_{11} b_{11} + a_{12} b_{21} & a_{11} b_{12} + a_{12} b_{22}\\
a_{21} b_{11} + a_{22} b_{21} & a_{21} b_{12} + a_{22} b_{22}
\end{bmatrix}\!.
\end{equation}
However, \eqref{eq:mm} defines a product on tensors whereas \eqref{eq:hp} does not, the reason being that the standard matrix product satisfies the transformation rules in \eqref{eq:matprod} whereas the Hadamard product does not. While there may be occasional scenarios where the Hadamard product could be useful, the standard matrix product is far more prevalent in all areas of mathematics, science and engineering. This is not limited to matrix  products; as we can see from the list of examples, the most useful and important notions we encounter in linear algebra are invariably tensorial notions that satisfy various transformation rules for $0$-, $1$- and $2$-tensors. We will have more to say about these issues in the next few sections. 

\subsection{Transformation rules for $d$-tensors}\label{sec:transdtensor}

We  now state the transformation rules in their most general form, laying out the full details of definition~\ref{st:tensor1}. Let $A = [a_{j_1 \cdots j_d } ] \in \mathbb{R}^{n_1 \times \dots \times n_d}$ be a $d$-dimensional matrix or \emph{hypermatrix}. Readers who require a rigorous definition may take it as a real-valued function $A \colon \{1,\ldots,n_1\} \times \cdots \times \{1,\ldots,n_d\} \to \mathbb{R}$, 
  a perspective we will discuss in Example~\ref{eg:hyp}, but those who are fine with a `$d$-dimensional matrix' may picture it as such. For matrices, \ie\ the usual two-dimensional matrices, $X = [x_{i_1j_1}] \in \mathbb{R}^{m_1 \times n_1}, Y = [y_{i_2j_2}] \in \mathbb{R}^{m_2 \times n_2}, \ldots, Z=[z_{i_dj_d}] \in \mathbb{R}^{m_d \times n_d}$, we define
\begin{equation}\label{eq:mmm}
(X,Y,\ldots, Z) \cdot A = B ,
\end{equation}
where $B = [b_{i_1 \cdots i_d}] \in  \mathbb{R}^{m_1 \times \dots \times m_d}$ is given by
\[
b_{i_1 \cdots i_d} = \sum_{j_1 = 1}^{n_1}\sum_{j_2 = 1}^{n_2}\dots \sum_{j_d = 1}^{n_d} x_{i_1 j_1} y_{i_2 j_2} \cdots z_{i_d j_d} a_{j_1 \cdots j_d}
\]
for $i_1 =1,\ldots,n_1,\; i_2 =1,\ldots,n_2, \;\ldots,\; i_d =1,\ldots,n_d$. We call the operation \eqref{eq:mmm} multilinear matrix multiplication; as we will discuss after stating the tensor transformation rules, the notation \eqref{eq:mmm} is chosen to be consistent with standard notation for group action. For $d=1$, it reduces to $Xa = b$ for $a \in \mathbb{R}^n$, $b\in \mathbb{R}^m$, and for
$d =2$, it reduces to
\[
(X,Y) \cdot A = XAY^\tp.
\]
For now, we will let the hypermatrix $A \in  \mathbb{R}^{n_1 \times \dots \times n_d}$ be the multi-indexed object in definition~\ref{st:tensor1} to keep things simple; we will later see that this multi-indexed object does not need to be a hypermatrix.

Let $X_1 \in \GL(n_1), X_2 \in \GL(n_2), \ldots, X_d \in \GL(n_d)$. The covariant $d$-tensor transformation rule is
\begin{equation}\label{eq:co1}
A' = (X_1^\tp,X_2^\tp,\ldots,X_d^\tp) \cdot A.
\end{equation}
The contravariant $d$-tensor transformation rule is
\begin{equation}\label{eq:contra1}
A' = (X_1^{-1},X_2^{-1},\ldots,X_d^{-1}) \cdot A.
\end{equation}
The mixed $d$-tensor transformation rule is\footnote{We state it in this form for simplicity. For example, we do not distinguish between the $3$-tensor transformation rules $A' = (X^{-1},Y^\tp, Z^\tp) \cdot A$, $A' = (X^\tp,Y^{-1}, Z^\tp) \cdot A$ and $A' = (X^\tp,Y^\tp, Z^{-1}) \cdot A$. All three are of  type $(1,2)$.}
\begin{equation}\label{eq:mixed1}
A' = (X_1^{-1},\ldots,X_p^{-1},X_{p+1}^\tp,\ldots,X_d^\tp) \cdot A.
\end{equation}
We say that the transformation rule in \eqref{eq:mixed1} is of \emph{type} or valence $(p,d-p)$ or, more verbosely, of contravariant order $p$ and covariant order $d - p$. As such, covariance is synonymous with type $(0,d)$ and contravariance with type $(d,0)$. 

For $d=1$, the hypermatrix is just $a \in \mathbb{R}^n$. If it transforms as $a' = X^{-1}a$, then it is the coordinate representation of a contravariant $1$-tensor or \emph{contravariant vector}; if it transforms as $a' = X^\tp a$, then it is the coordinate representation of a covariant $1$-tensor or \emph{covariant vector}. For $d=2$, the hypermatrix is just a matrix $A \in \mathbb{R}^{m \times n}$; writing $X_1 = X$, $X_2 = Y$, the transformations in \eqref{eq:co1}, \eqref{eq:contra1}, \eqref{eq:mixed1} become $A' = X^\tp AY$, $A' = X^{-1}AY^{-\tp}$, $A' = X^{-1}AY$, which are the transformation rules for a covariant, contravariant, mixed $2$-tensor respectively. These look different from the transformation rules in the examples on pages~\pageref{eg:eig}--\pageref{eg:rank} as a result of the ambiguity \ref{it:ambiguity3} on page~\pageref{it:ambiguity3}, which we elaborate below.

Clearly, the equalities in the middle and right columns below only differ in which side we label with a prime but are otherwise identical:\label{transrulesummary1}
\begin{alignat*}{5}
&\text{contravariant $1$-tensor} &a' &= X^{-1} a, &a' &= X  a,\\*
&\text{covariant $1$-tensor} &a' &= X^\tp a, &a' &= X^{-\tp} a,\\
&\text{contravariant $2$-tensor}\quad &A' &= X^{-1} A Y^{-\tp},\quad &A' &= X A Y^\tp,\\
&\text{covariant $2$-tensor} &A' &= X^\tp A Y, &A' &= X^{-\tp} A Y^{-1},\\*
&\text{mixed $2$-tensor} &A' &= X^{-1}A Y, &A' &= XA Y^{-1}.
\end{alignat*}
We will see in Section~\ref{sec:multmaps} after introducing definition~\ref{st:tensor2} that these transformation rules come  from the change-of-basis theorems for vectors, linear functionals, dyads,  bilinear functionals and linear operators, respectively, with $X$ and $Y$ the change-of-basis matrices. The versions in the middle and right column differ in terms of how we write our change-of-basis theorem. For example, take the standard change-of-basis theorem in a vector space. Do we write $a_{\text{new}} = X^{-1} a_{\text{old}}$ or $a_{\text{old}} = X  a_{\text{new}}$?

Observe, however, that there is no repetition in either column: the transformation rule, whether in the middle or right column, uniquely identifies the type of tensor. It is not possible to confuse, say, a contravariant $2$-tensor with a covariant $2$-tensor just because we use the transformation rule in the middle column to describe one and the  transformation rule in the right column to describe the other.

The version in the middle column is consistent with  the names `covariance' and `contravariance', which are  based on whether the hypermatrix `co-varies', \ie\ transforms with the same $X$, or `contra-varies', \ie\ transforms with its inverse $X^{-1}$. This is why we have stated our transformation rules in \eqref{eq:co1}, \eqref{eq:contra1}, \eqref{eq:mixed1} to be consistent with those in the middle column. But there are also occasions, as in the examples on pages~\pageref{eg:eig}--\pageref{eg:rank}, when it is more natural to express the transformation rules as  those in the right column. 

When $n_1 = n_2 = \dots = n_d=n$, the transformation rules in \eqref{eq:co1}, \eqref{eq:contra1}, \eqref{eq:mixed1} may take on a different form with a single change-of-basis matrix $X \in \GL(n)$ as opposed to $d$ of them. In this case the hypermatrix $A \in \mathbb{R}^{n \times \dots \times n}$ is hypercubical, \ie\ the higher-order equivalent of a square matrix, and the covariant tensor transformation rule is
\begin{equation}\label{eq:co2}
A' = (X^\tp,X^\tp,\ldots,X^\tp) \cdot A,
\end{equation}
the contravariant tensor transformation rule is
\begin{equation}\label{eq:contra2}
A' = (X^{-1},X^{-1},\ldots,X^{-1}) \cdot A,
\end{equation}
and the mixed tensor  transformation rule is
\begin{equation}\label{eq:mixed2}
A' = (X^{-1},\ldots,X^{-1},X^\tp,\ldots,X^\tp) \cdot A.
\end{equation}
Again, we will see in Section~\ref{sec:multmaps} after introducing definition~\ref{st:tensor2} that the difference between, say, \eqref{eq:co1} and  \eqref{eq:co2} is that the former expresses the change-of-basis theorem for a multilinear map $f \colon  \mathbb{V}_1 \times \dots \times \mathbb{V}_d \to \mathbb{R}$ whereas the latter is for a multilinear map $f \colon  \mathbb{V} \times \dots \times \mathbb{V} \to \mathbb{R}$.

For $d=2$, we have $A \in \mathbb{R}^{n \times n}$, and each transformation rule  again takes two  different forms that are identical in substance:\label{transrulesummary2}
\begin{alignat*}{5}
&\text{covariant $2$-tensor} &A' &= X^\tp A X, &A' &= X^{-\tp} A X^{-1},\\*
&\text{contravariant $2$-tensor}\quad &A' &= X^{-1} A X^{-\tp},\quad &A' &= X A X^\tp,\\*
&\text{mixed $2$-tensor} &A' &= X^{-1}A X, &A' &= XA X^{-1}.
\end{alignat*}
Note that either form uniquely identifies the tensor type.

\begin{example}[covariance versus contravariance]\hspace{-7pt}%
The three transformation rules for $2$-tensors are $A' = XAX^\tp$, $A' = X^{-\tp}AX^{-1}$, $A' = XAX^{-1}$. We know from linear algebra that there is a vast difference between the first two (congruence) and the last one (similarity). For example, eigenvalues and eigenvectors are defined for mixed $2$-tensors by virtue of \eqref{eq:ev} but undefined for covariant or contravariant $2$-tensors since the eigenvalue/eigenvector equation $Av = \lambda v$ is incompatible with the first two transformation rules. While both the contravariant and covariant transformation rules describe congruence of matrices, the difference between them is best seen to be the difference between the quadratic form and the second-order partial differential operator:
\[
\sum_{i=1}^n \sum_{j =1}^n a_{ij} v_i v_j \quad\text{and}\quad
\sum_{i=1}^n \sum_{j =1}^n a_{ij} \dfrac{\partial^2}{\partial v_i \partial v_j}.
\]
$A' = X^{-\tp}AX^{-1}$  and $A' = XAX^\tp$, respectively, describe how their coefficient matrices transform under a change of coordinates $v' = Xv$.
\end{example}

As we have mentioned in Section~\ref{sec:transrules}, definition~\ref{st:tensor1} leaves the tensor unspecified. In physics, this is to serve as a placeholder for a physical quantity under discussion such as deformation or rotation or strain. In mathematics, if we insist on not bringing in definitions~\ref{st:tensor2} or~\ref{st:tensor3}, then the way to use definition~\ref{st:tensor1} is best illustrated with an example.

\begin{example}[transformation rule determines tensor type]
We say  the eigen\-value/eigenvector equation $Av = \lambda v$ conforms to the transformation rules since $(XAX^{-1}) Xv = \lambda Xv$, and thus $A$ transforms like a mixed $2$-tensor $A' = XAX^{-1}$, the eigenvector $v$ transforms like a covariant $1$-tensor $v' = Xv$ and the eigenvalue $\lambda$ transforms like an invariant $0$-tensor. It is only with context (here we have eigenvalues/eigenvectors) that it makes sense to speak of $A$, $v$, $\lambda$ as tensors. Is $A =\begin{bsmallmatrix}1&2&3\\ 2 &3&4 \\3 &4 &5\end{bsmallmatrix}$ a tensor? It is a tensor if it represents, say, stress (in which case $A$ transforms as a contravariant $2$-tensor) or, say, we are interested in its eigenvalues and eigenvectors (in which case $A$ transforms as a mixed $2$-tensor). Stripped of such contexts, the question as to whether a matrix or hypermatrix is a tensor is not meaningful.
\end{example}

We have stressed that the transformation rules define the tensor but one might think that this only applies to the type of the tensor, \ie\ whether it is covariant or contravariant or mixed. This is not the case: it defines \emph{everything} about the tensor.

\begin{example}[transformation rule determines tensor order]\label{eg:order}\hspace{-7pt}%
Does $A =\begin{bsmallmatrix}1&2&3\\ 2 &3&4 \\3 &4 &5\end{bsmallmatrix}$ represent a $1$-tensor or a $2$-tensor? One might think that since a matrix is a doubly indexed object, it must represent a $2$-tensor. Again this is a fallacy. If the transformation rules we apply to a matrix $A \in \mathbb{R}^{m \times n}$ are of the forms
\begin{equation}\label{eq:mat1tensor}
A' = XA = [Xa_1,\ldots, Xa_n],\quad  A' = X^{-\tp}A = [X^{-\tp}a_1,\ldots, X^{-\tp}a_n],
\end{equation}
where $a_1,\ldots,a_n \in \mathbb{R}^m$ are the columns of $A$ and $X \in \GL(m)$, then $A$ represents a contravariant or covariant $1$-tensor respectively, \tie, the \emph{order} of a tensor is also determined by the transformation rules. In numerical linear algebra, algorithms such as Gaussian elimination and Householder or Givens QR are all based on transformations of a matrix as a $1$-tensor as in \eqref{eq:mat1tensor}; we will have more to say about this in Section~\ref{sec:equivariance}. In equivariant neural networks, it does not matter that filters are  represented as five-dimensional arrays \cite[Section~7.1]{Welling} or that inputs are a matrix of vectors of `irreducible fragments' \cite[Definition~1]{Risi2}: the transformation rules involved are ultimately just those of contravariant $1$-tensors and mixed $2$-tensors, as we will see in Example~\ref{eg:neural}.
\end{example}

The change-of-basis matrices $X_1,\ldots,X_d$ in \eqref{eq:co1}, \eqref{eq:contra1}, \eqref{eq:mixed1}  and $X$ in \eqref{eq:co2}, \eqref{eq:contra2}, \eqref{eq:mixed2}  do not need to be just invertible matrices, \ie\ elements of $\GL(n)$: they could instead come from any classical groups, most commonly the orthogonal group $\Or(n)$ or unitary group $\Un(n)$, possibly with unit-determinant constraints, \ie\ $\SL(n)$, $\SOr(n)$, $\SU(n)$. The reason is that as the name implies, these are change-of-basis matrices for the vector spaces in definitions~\ref{st:tensor2} or~\ref{st:tensor3}, and if those vector spaces are  equipped with additional structures, we expect the change-of-basis matrices to preserve those structures.

\begin{example}[Cartesian and Lorentzian tensors]\label{eg:cartesian}\hspz%
For the tensor transformation rule over $\mathbb{R}^4$ equipped with the Euclidean inner product
\[
\langle x,y\rangle = x_0 y_0 + x_1 y_1 + x_2 y_2 + x_3 y_3,
\]
we want the change-of-basis matrices to be from the orthogonal group $\Or(4)$ or special orthogonal group $\SOr(4)$, but if $\mathbb{R}^4$ is equipped instead with the Lorentzian scalar product in relativity,
\[
\langle x,y\rangle = x_0 y_0 - x_1 y_1 - x_2 y_2 - x_3 y_3,
\]
then the Lorentz group $\Or(1,3)$ or proper Lorentz group $\SOr(1,3)$ or restricted Lorentz group $\SOr^\p(1,3)$ is more appropriate. Names like Cartesian tensors or Lorentzian tensors are sometimes used to indicate whether the transformation rules involve orthogonal or Lorentz groups but sometimes physicists would assume that these are self-evident from the context and are left unspecified, adding to the confusion for the uninitiated. The well-known \emph{electromagnetic field tensor} or Faraday tensor, usually represented in matrix form by
\[
F = \begin{bmatrix}
     0     & -E_x/c & -E_y/c & -E_z/c \\[2pt]
     E_x/c &  0     & -B_z   &  B_y    \\[2pt]
     E_y/c &  B_z   &  0     & -B_x   \\[2pt]
     E_z/c & -B_y   &  B_x   &  0
  \end{bmatrix}\!,
\]
where we have written $(x_0, x_1, x_2, x_3 ) = (ct,x,y,z)$,
is a contravariant Lorentz $2$-tensor as it satisfies a transformation rule of the form $F' = XFX^\tp$ with a change-of-basis matrix $X \in \Or(1,3)$.
\end{example}

For easy reference, here is a list of the most common matrix Lie groups that the change-of-basis matrices $X_1,\ldots,X_d$ or $X$ may belong to:
\begin{align}\label{eq:groups}
\GL(n) &= \{ X \in \mathbb{R}^{n \times n} \colon  \det(X) \ne 0\}, \notag \\*
\SL(n) &= \{ X \in \mathbb{R}^{n \times n} \colon  \det(X) = 1\}, \notag\\
\Or(n) &= \{ X \in \mathbb{R}^{n \times n} \colon  X^\tp X = I\}, \notag\\
\SOr(n) &= \{ X \in \mathbb{R}^{n \times n} \colon  X^\tp X = I,\; \det(X)=1\}, \notag\\
\Un(n) &= \{ X \in \mathbb{C}^{n \times n} \colon  X^* X = I\}, \notag\\
\SU(n) &= \{ X \in \mathbb{C}^{n \times n} \colon  X^* X = I,\; \det(X)=1\}, \notag\\
\Or(p,q) &= \{ X \in \mathbb{R}^{n \times n} \colon  X^\tp I_{p,q} X = I_{p,q}\}, \notag\\
\SOr(p,q) &= \{ X \in \mathbb{R}^{n \times n} \colon  X^\tp I_{p,q} X = I_{p,q},\; \det(X)=1\}, \notag\\
\Sp(2n,\mathbb{R}) &= \{ X \in \mathbb{R}^{2n \times 2n} \colon  X^\tp J X = J\}, \notag\\*
\Sp(2n) &= \{ X \in \mathbb{C}^{2n \times 2n} \colon  X^\tp J X = J,\; X^* X = I\}.
\end{align}
In the above, $\mathbb{R}$ may be replaced by $\mathbb{C}$ to obtain corresponding groups of complex matrices, $p$ and $q$ are positive integers with $p + q = n$, and $I =I_n$ denotes the $n \times n$ identity matrix. The matrices $I_{p,q} \in \mathbb{C}^{n \times n} $ and $J \in \mathbb{C}^{2n \times 2n}  $, respectively, are
\[
I_{p,q} \coloneqq \begin{bmatrix} I_p & 0 \\ 0 & - I_q \end{bmatrix}\!, \quad J \coloneqq \begin{bmatrix}  0 & I \\ I & 0 \end{bmatrix}\!.
\]
There are some further possibilities that deserve a mention here, if only for cultural reasons. While not as common as the classical Lie groups above, depending on the application at hand, the groups of invertible lower/upper or unit lower/upper triangular matrices, Heisenberg groups, {\em etc.}, may all serve as the group of change-of-basis matrices. Another three particularly important examples are the general affine, Euclidean and special Euclidean groups
\begin{equation}\label{eq:SE}
\begin{aligned}
\GA(n) &= \biggl\{ \begin{bmatrix} X & y \\ 0 &  1 \end{bmatrix} \in \mathbb{R}^{(n+1) \times (n+1)} \colon  X \in \GL(n), \; y \in \mathbb{R}^n \biggr\},\\
\E(n) &= \biggl\{ \begin{bmatrix} X & y \\ 0 &  1 \end{bmatrix} \in \mathbb{R}^{(n+1) \times (n+1)} \colon  X \in \Or(n), \; y \in \mathbb{R}^n \biggr\},\\
\SE(n) &= \biggl\{ \begin{bmatrix} X & y \\ 0 &  1 \end{bmatrix} \in \mathbb{R}^{(n+1) \times (n+1)} \colon  X \in \SOr(n), \; y \in \mathbb{R}^n \biggr\},
\end{aligned}
\end{equation}
which encode a linear transformation by $X$ and translation by $y$ in $\mathbb{R}^n$. Note that even though the matrices involved are $(n+1) \times (n+1)$, they act on vectors $v\in\mathbb{R}^n$ embedded as $(v,0) \in \mathbb{R}^{n+1}$. As we will see in Example~\ref{eg:neural}, the last two groups play a special role in equivariant neural networks for image classification ($n =2$) and protein structure predictions ($n=3$).

The above discussions extend to an infinite-dimensional setting in two different ways. In general relativity, the tensors are replaced by tensor fields and the change-of-basis matrices are replaced by change-of-coordinates matrices whose entries are functions of coordinates; these are called \emph{diffeomorphisms} and they form infinite-dimensional Lie groups.  In quantum mechanics, the tensors over finite-dimensional vector spaces are replaced by tensors over infinite-dimensional Hilbert spaces and the change-of-basis matrices are replaced by \emph{unitary operators}, which may be regarded as infinite-dimensional matrices (see Example~\ref{eg:hyp}). In this case the tensor transformation rules are particularly important as they describe the evolution of the quantum system. We will have more to say about these in Examples~\ref{eg:tenfield2} and \ref{eg:grover}.

We are now ready to state a slightly modernized version of definition~\ref{st:tensor1}.

\begin{definition}[tensors via transformation rules]\label{def:tensor1}\hspz%
A hypermatrix $A \in \mathbb{R}^{n_1 \times \dots \times n_d}$  \emph{represents} a $d$-tensor if it satisfies one of the transformation rules in \eqref{eq:co1}--\eqref{eq:mixed2} for some matrix groups in \eqref{eq:groups}.
\end{definition}

Notice that this only defines what it means for a hypermatrix to represent a tensor: it does not fix the main drawback of definition~\ref{st:tensor1}, \ie\ the lack of an object to serve as the tensor. Definition~\ref{st:tensor1} cannot be easily turned into a standalone definition and is best understood in conjunction with either definition~\ref{st:tensor2} or \ref{st:tensor3}; these transformation rules then become change-of-basis theorems for a multilinear map or an element in a tensor product of vector spaces.  Nevertheless, definition~\ref{st:tensor1} has its value: the transformation rules encode the important notion of equivariance that we will discuss next over Sections~\ref{sec:imtransrules} and \ref{sec:equivariance}.

There are some variants of the tensor transformation rules that come up in various areas of physics. We state them below for completeness and as examples of hypermatrices that do not represent tensors.
  
\begin{example}[pseudotensors and tensor densities]\label{eg:pseudo}
  For a change-of-basis matrix $X \in \GL(n)$ and a hypermatrix $A \in \mathbb{R}^{n \times \dots \times n}$,  the pseudotensor transformation rule \cite[Section~3.7.2]{Borisenko}, the tensor density transformation rule \cite[Section~3.8]{Plebanskieu} and the pseudotensor density transformation rule  are, respectively,
\begin{equation}\label{eq:pmixed2}
\begin{aligned}
A' &= \sgn(\det X)\, (X^{-1},\ldots,X^{-1},X^\tp,\ldots,X^\tp) \cdot A,\\
A' &= (\det X)^p\, (X^{-1},\ldots,X^{-1},X^\tp,\ldots,X^\tp) \cdot A,\\
A' &= \sgn(\det X) (\det X)^p\, (X^{-1},\ldots,X^{-1},X^\tp,\ldots,X^\tp) \cdot A.
\end{aligned}
\end{equation}
Evidently they differ from the tensor transformation rule \eqref{eq:mixed2} by a scalar factor. The power $p \in \mathbb{Z}$ is called the weight with $p =0$ giving us \eqref{eq:mixed2}.
An alternative name for a pseudotensor is an `axial tensor' and in this case a true tensor is called a `polar tensor' \cite[Section~4]{Hartmann}. 
Note that the transformation rules in \eqref{eq:pmixed2} are variations of \eqref{eq:mixed2}. While it is conceivable to define similar variations of \eqref{eq:mixed1} with $\det (X_1\cdots X_d) = \det X_1\cdots \det X_d$ in place of $\det X$, we are unable to find a reference for this and {would} therefore refrain from stating it formally.
\end{example}

\paragraph{A word about notation.} Up to this point we have introduced only one notion that is not from linear algebra, namely \eqref{eq:mmm} for multilinear matrix multiplication; this notion dates back to at least \citet{Hitch1} and the notation is standard too. We denote it so that \eqref{eq:co1}, \eqref{eq:contra1}, \eqref{eq:mixed1} reduce to the standard notation $g \cdot x$ for a group element $g \in G$ acting on an element $x$ of a $G$-set. Here our group $G$ is a product of other groups $G = G_1 \times \dots \times G_d$, and so an element takes the form $g = (g_1,\ldots,g_d)$. We made a conscious decision to cast everything in this article in terms of \emph{left} action so as not to introduce another potential source of confusion. We could have introduced a \emph{right} multilinear matrix multiplication of  $A \in \mathbb{R}^{n_1 \times \dots \times n_d}$ by matrices  $X  \in \mathbb{R}^{n_1 \times m_1}, Y  \in \mathbb{R}^{n_2 \times m_2}, \ldots, Z \in \mathbb{R}^{n_d \times m_d}$;  note that the dimensions of these matrices are now the transposes of those in  \eqref{eq:mmm},
\[
A \cdot (X,Y,\ldots, Z) = B ,
\]
with $B \in  \mathbb{R}^{m_1 \times \dots \times m_d}$ given by
\[
b_{i_1 \cdots i_d} = \sum_{j_1 = 1}^{n_1}\sum_{j_2 = 1}^{n_2}\dots \sum_{j_d = 1}^{n_d}  a_{j_1 \cdots j_d} x_{i_1 j_1} y_{i_2 j_2} \cdots z_{i_d j_d}.
\]
This would have allowed us to denote \eqref{eq:co1} and \eqref{eq:mixed1}  without transposes, with the latter in a two-sided product. Nevertheless, we do not think it is worth the trouble.

\subsection{Importance of the transformation rules}\label{sec:imtransrules}

When we say that a group is a set with an operation satisfying certain axioms, these axioms play a crucial role and cannot be disregarded. It is the same with tensors. The importance of the transformation rules in definition~\ref{st:tensor1} should be abundantly clear at this point.  In general, a formula in terms of the entries of the hypermatrix is undefined for the tensor since that formula is unlikely to conform to the stringent transformation rules; this mistake is unfortunately commonplace. Here we will discuss some erroneous ideas that result from treating a tensor as no more than the hypermatrix that represents it,  oblivious of the transformation rules it must satisfy.

\begin{example}[`multiplying higher-order tensors']\label{eg:homult}
  There have been many recent attempts at finding a formula for `multiplying higher-order tensors' that purports to extend the standard matrix--matrix product \eqref{eq:mm}. These proposed formulas are justified with a demonstration that they are associative, distributive, have a multiplicative identity, {\em etc.}, just like matrix multiplication. This is misguided. We need go no further than $2$-tensors to see the problem: both \eqref{eq:hp} and \eqref{eq:mm} are associative, distributive and have multiplicative identities, but only the formula in \eqref{eq:mm} defines a product of $2$-tensors because it conforms to the transformation rules \eqref{eq:matprod}. In fact no other formula aside from the standard matrix--matrix product defines a product for $2$-tensors, a consequence of the first fundamental theorem of invariant theory; see \citet[Chapter~9]{Procesi} or \citet[Section~5.3]{Goodman}.
  This theorem also implies that  there is no formula that will yield a product for any odd-ordered tensors. In particular, one may stop looking for a formula to multiply two $n \times n \times n$ matrices to get a third $n \times n \times n$ hypermatrix that defines a product for $3$-tensors, because such a formula does not exist. Incidentally,  the transformation rules also apply to adding tensors: we may not add two matrices representing two tensors of different types because they satisfy different transformation rules; for example, if $A \in \mathbb{R}^{n \times n}$ represents a mixed $2$-tensor and $B \in \mathbb{R}^{n \times n}$ a covariant $2$-tensor, then $A + B$ is rarely if ever meaningful.
\end{example}

\begin{example}[`identity tensors']
The identity matrix $I \in \mathbb{R}^{3 \times 3}$ is of course
\begin{equation}\label{eq:id2}
I = \sum_{i=1}^3 e_i \otimes e_i  \in \mathbb{R}^{3 \times 3},
\end{equation}
where $e_1,e_2,e_3 \in \mathbb{R}^3$ are the standard basis vectors.
What should be the extension of this notion to $d$-tensors? Take $d = 3$ for
 illustration. It would appear that
\begin{equation}\label{eq:id3}
A =\sum_{i=1}^3 e_i \otimes e_i \otimes e_i  \in \mathbb{R}^{3 \times 3 \times 3}
\end{equation}
is the obvious generalization but, as with the Hadamard product, obviousness does not necessarily imply correctness when it comes to matters tensorial. One needs to check the transformation rules. Note that \eqref{eq:id2} is independent of the choice of orthonormal basis: we obtain the same matrix with any orthonormal basis $q_1,q_2,q_3 \in \mathbb{R}^3$, a consequence of
\[
(Q,Q) \cdot I = Q I Q^\tp = I
\]
for any $Q \in \Or(3)$, \tie, the identity matrix is well-defined as a Cartesian tensor. On the other hand the hypermatrix $A$ in \eqref{eq:id3} does not have this property. One may show that up to scalar multiples, $M = I$ is the unique matrix satisfying
\begin{equation}\label{eq:isotropic2}
(Q,Q) \cdot M = M
\end{equation}
for any $Q \in \Or(3)$, a property known as \emph{isotropic}. An isotropic $3$-tensor would then be one that satisfies
\begin{equation}\label{eq:isotropic3}
(Q,Q,Q) \cdot T = T.
\end{equation}
Up to scalar multiples, \eqref{eq:isotropic3} has a unique solution given by the hypermatrix 
\[
J =\sum_{i=1}^3\sum_{j=1}^3\sum_{k=1}^3 \varepsilon_{ijk} e_i \otimes e_j \otimes e_k \in \mathbb{R}^{3 \times 3 \times 3},
\]
where $\varepsilon_{ijk}$ is the Levi-Civita symbol
\[
\varepsilon_{ijk}   =
\begin{cases}
+1 & \text{if }(i,j,k)=(1,2,3),(2,3,1),(3,1,2),\\
-1 & \text{if }(i,j,k)=(1,3,2),(2,1,3),(3,2,1),\\
0 & \text{if }i=j, j=k,k=i,
\end{cases}
\]
which is evidently quite different from \eqref{eq:id3}.  More generally, $n$-dimensional isotropic $d$-tensors, \ie\ tensors with the same hypermatrix representation regardless of choice of orthonormal bases in $\mathbb{R}^n$, have been explicitly determined for $n=3$ and $ d \le 8$ in \citet{KF},  and  studied for arbitrary values of $n$ and $d$ in \citet{WeylBook}. From a tensorial perspective, isotropic tensors, not a `hypermatrix with ones on the diagonal' like \eqref{eq:id3}, extend the notion of an identity matrix.
\end{example}

\begin{example}[hyperdeterminants]\hspz%
The simplest order-$3$ analogue for the deter\-minant, called the Cayley hyperdeterminant for $A = [a_{ijk}] \in \mathbb{R}^{2 \times 2 \times 2}$, is given~by
\begin{align*}
\Det(A) & = 
a_{000}^{2} a_{111}^{2}+a_{001}^{2}a_{110}^{2}+a_{010}^{2}
a_{101}^{2}+a_{011}^{2}a_{100}^{2}\\*
 & \quad -2(a_{000}a_{001}a_{110}a_{111}+a_{000}a_{010}a_{101}a_{111}+a_{000}
a_{011}a_{100}a_{111}\\
 & \quad +a_{001}a_{010}a_{101}a_{110}+a_{001}a_{011}a_{110}a_{100}+a_{010}
a_{011}a_{101}a_{100})\\*
&\quad +4(a_{000}a_{011}a_{101}a_{110}+a_{001}a_{010}a_{100}a_{111}),
\end{align*}
which looks nothing like the usual expression for a matrix determinant. Just as the matrix determinant is preserved under a transformation $A' = (X,Y) \cdot A = XAY^\tp $ for $X, Y\in \SL(n)$, this is preserved under a transformation of the form $A' = (X,Y,Z) \cdot A$ for any $X, Y, Z \in \SL(2)$. This  notion of hyperdeterminant has been extended to any $A \in \mathbb{R}^{n_1 \times \dots \times n_d}$ with
\[
n_i - 1 \le \sum_{j \ne i} (n_j - 1), \quad i =1,\ldots, d,
\]
in \citet{GKZpaper}. 
\end{example}

As these three examples reveal, extending any common linear algebraic notion  to higher order will require awareness of the transformation rules. Determining the  higher-order tensorial analogues of notions such as linear equations or eigenvalues or determinant is not a matter of simply adding more indices to
\[
\sum_{j=1}^n a_{ij} x_j = b_i, \quad  \sum_{j=1}^n a_{ij} x_j = \lambda x_i, \quad \sum_{\sigma \in \mathfrak{S}_n} \sgn(\sigma)\prod_{i=1}^n a_{i \sigma(i)}.
\]

Of course, one might argue that anyone is free to concoct any formula for `tensor multiplication' or call any hypermatrix an `identity tensor' even if they are undefined on tensors: there is no right or wrong. But there is. They would be wrong in the same way adding fractions as $a/b+ c/d = (a+b)/(c+d)$ is wrong, or at least far less useful than the standard way of adding fractions. Fractions are real, and adding half a cake to a third of a cake does not give us two-fifths of a cake. Likewise tensors are real: we discovered the transformation rules for tensors just as we discovered the arithmetic rules for adding fractions, \tie, we did not invent them arbitrarily.

While we have limited our examples to mathematical ones, we will end with a note about the physics perspective. In physics, these transformation rules are just as if not more important. One does not even have to go to higher-order tensors to see this: a careful treatment of \emph{vectors} in physics already requires such an approach.

\begin{example}[tensor transformation rules in physics]\hspz%
As we saw at the beginning of Section~\ref{sec:trans}, the tensor transformation rules originated from physics as change-of-coordinates rules for certain physical quantities. Not all physical quantities satisfy tensor transformation rules; those that do are called `tensorial quantities' or simply `tensors' \cite[p.~20]{Voigt}. We will take the simplest physical quantity, \emph{displacement}, \ie\ the location $q$ of a point particle in space, to illustrate. Once we select our $x$-, $y$- and $z$-axes, space becomes $\mathbb{R}^3$ and $q$ gets coordinates $a=(a_1,a_2,a_3) \in \mathbb{R}^3$ with respect to these axes. If we change our axes to a new set of $x'$-, $y'$- and $z'$-axes -- note that nothing physical has changed, the point $q$ is still where it was, the only difference is in the directions we decided to pick as axes -- the coordinates must therefore change accordingly to some $a'=(a_1',a_2',a_3') \in \mathbb{R}^3$ to compensate for this change in the choice of axes (see Figure~\ref{fig:coord}). If the matrix $X \in \GL(n)$ transforms the $x$-, $y$- and $z$-axes to the $x'$-, $y'$- and $z'$-axes, then the coordinates must transform in the opposite way $a' = X^{-1}a$ to counteract this change so that on the whole everything stays the same. This is the transformation rule for a contravariant $1$-tensor on page~\pageref{transrulesummary1}, and we are of course just paraphrasing the change-of-basis theorem for a vector in a vector space. This same transformation rule applies to many other physical quantities: any time derivatives of displacement such as velocity $\dot{q}$, acceleration $\ddot{q}$, momentum $m\dot{q}$, force $m\ddot{q}$. So the most familiar physical quantities we encountered in introductory mechanics are contravariant $1$-tensors.
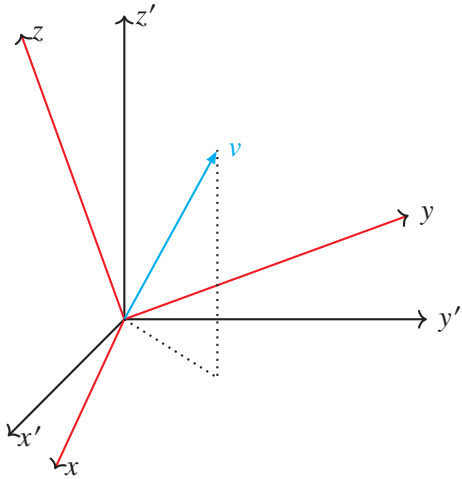
\begin{figure}[ht]
\centering
\begin{tikzpicture}[thick]
  % X',Y',Z' in black
  \draw[->] (0,0,0) -- (0,0,4)  node[right, text width=5em] {$x'$};
  \draw[->] (0,0,0) -- (4,0,0)  node[right, text width=5em] {$y'$};
  \draw[->] (0,0,0) -- (0,4,0)  node[right, text width=5em] {$z'$};
  
  % X,Y,Z in red
  \begin{scope}[rotate=20,draw=red]
      \draw[->] (0,0,0) -- (0,0,4)  node[right, text width=5em] {$x$};
      \draw[->] (0,0,0) -- (4,0,0)  node[right, text width=5em] {$y$};
      \draw[->] (0,0,0) -- (0,4,0)  node[right, text width=5em] {$z$};
  \end{scope}  

  % a small
  \draw[cyan, -latex] (0,0,0) -- (2,3,2)    node[right, text width=5em] (b) {$v$};
  
  \draw [dotted]
         (0,0,0) 
      -- (2,0,2) %node [below, midway] {$u$}
      -- (2,3,2);
\end{tikzpicture}
\parbox{230pt}{\caption{Linear transformation of coordinate axes.}}
\label{fig:coord}
\end{figure}

Indeed, contravariant $1$-tensors are how many physicists would regard vectors \cite{Weinreich}: a vector is an object represented by $a \in \mathbb{R}^n$ that satisfies the transformation rule $a' = X^{-1}a$, possibly with the additional requirement that the change-of-basis matrix $X$ be in $ \Or(n)$ \cite[Chapter~11]{Feynman1} or $X \in \Or(p,q)$ \cite[Chapter~4]{Rindler}.

The transformation rules perspective has proved to be very useful in physics. For example, special relativity is essentially the observation that the laws of physics are invariant under Lorentz transformations $X \in \Or(1,3)$ \cite{special}.  In fact, a study of the contravariant $1$-tensor transformation rule under the $\Or(1,3)$-analogue of Givens rotations,\label{Qadj}
\[
\adjustbox{max width=\textwidth}{$
  \begin{bmatrix} \cosh\theta & -\sinh\theta & 0 & 0\\ -\sinh\theta & \cosh\theta & 0 & 0\\ 0&0&1&0\\ 0&0&0&1\\ \end{bmatrix}\!, \ \ \begin{bmatrix} \cosh\theta & 0 & -\sinh\theta & 0\\ 0&1&0&0\\ -\sinh\theta &0& \cosh\theta &0\\ 0&0&0&1\\ \end{bmatrix}\!, \ \ \begin{bmatrix} \cosh\theta &0&0&-\sinh\theta\\ 0&1&0&0\\ 0&0&1&0\\ -\sinh\theta&0&0&\cosh\theta\\ \end{bmatrix}$,}
\]
is enough to derive most standard results of special relativity; see \citet[Section~6.9]{FIS} and \citet[Chapters~4 and 5]{Woodhouse}.

Covariant $1$-tensors are also commonplace in physics. If the coordinates of $q$ transform as $a' = X^{-1}a$, then the coordinates of its derivative $\partial/\partial q$ or $\nabla_q$ would transform as $a' = X^\tp a$. So  physical quantities that satisfy the covariant $1$-tensor transformation rule tend to have a contravariant $1$-tensor `in the denominator' such as conjugate momentum $p = \partial L/\partial \dot{q}$ or electric field $E =- \nabla \phi - \partial A/ \partial t$. In the former, $L$ is the Lagrangian, and we have the velocity $\dot{q}$  `in the denominator'; in the latter, $\phi$ and $A$, respectively, are the scalar and vector potentials, and the gradient $\nabla$ is a derivative with respect to spatial variables, and thus have displacement $q$ `in the denominator'.

More generally, what we said above applies to more complex physical quantities: when the axes -- also called the reference frame in physics -- are changed,  coordinates of physical quantities must change in a way that preserves the laws of physics; for tensorial quantities, this would be \eqref{eq:co1}--\eqref{eq:mixed2}, the `definite way' in Voigt's definition on page~\pageref{Voigt}.  This is the well-known maxim that the laws of physics should not depend on coordinates, a version of which the reader will find on page~\pageref{Thorne}. The reader will also find higher-order examples  in  Examples~\ref{eg:cartesian}, \ref{eg:tenfield} and \ref{eg:stress}.
\end{example}

\subsection{Tensors in computations I: equivariance}\label{sec:equivariance}

These transformation rules are no less important in computations. The most fundamental algorithms in numerical linear algebra are invariably based implicitly on one of these transformation rules combined with special properties of the change-of-basis matrices. Further adaptations for numerical stability and efficiency can sometimes obscure what is going on in these algorithms, but the basic principle is always to involve one of these transformation rules, or, in modern lingo, to take advantage of \emph{equivariance}.

\begin{example}[equivalence,\, similarity,\, congruence]\label{eg:canforms}\hspz%
The transformation rules\linebreak %%20210319
 for $2$-tensors have better-known names in linear algebra and numerical linear algebra.
\begin{itemize}
\setlength\itemsep{3pt}
\item Equivalence of matrices: $A' = XAY^{-1}$.

\item Similarity of matrices: $A' = XAX^{-1}$.

\item Congruence of matrices: $A' = XAX^\tp$.
\end{itemize}
The most common choices for $X$ and $Y$ are either orthogonal matrices or invertible ones (similarity and congruence are identical if $X$ is orthogonal). These transformation rules have  \emph{canonical forms}, and this is an important feature in linear algebra and numerical linear algebra alike, not withstanding the celebrated fact \cite{GW} that some of these cannot be computed in finite precision: Smith normal form for equivalence, singular value decomposition for orthogonal equivalence, Jordan and Weyr canonical forms for similarity \cite{Weyr}, Turnbull--Aitken and Hodge--Pedoe canonical forms for congruence \cite{DeTe}. 
One point to emphasize is that the transformation rules determine the canonical form for the tensor: it makes no sense to speak of Jordan form for a contravariant $2$-tensor or a Turnbull--Aitken form for a mixed $2$-tensor, even if both tensors are represented by exactly the same matrix $A \in \mathbb{R}^{n \times n}$.
\end{example}

At this juncture it is appropriate to highlight a key difference between $2$-tensors and higher-order ones:  while $2$-tensors have canonical forms, higher-order tensors in general do not \cite[Chapter~10]{Land1}. This is one of several reasons why we should not expect an extension of linear algebra or numerical linear algebra to $d$-dimensional hypermatrices in a manner that resembles the $d=2$ versions.

Before we continue, we will highlight two simple properties.
\begin{enumerate}[\upshape (i)]
\setlength\itemsep{3pt}
\item The change-of-basis matrices may be multiplied and inverted: if $X$ and $Y$ are orthogonal or invertible, then so is $XY$ and so is $X^{-1}$, \tie, the set of all change-of-basis matrices $\Or(n)$ or $\GL(n)$ forms a group.

\item The transformation rules may be composed: if we have $a' = X^{-\tp} a$ and $a'' = Y^{-\tp} a'$, then $a'' = (YX)^{-\tp}a$; if we have $A' = XAX^{-1}$ and $A'' = YA'Y^{-1}$, then $A'' = (YX)A(YX)^{-1}$, \tie, the transformation rule defines a group action.
\end{enumerate}
These innocuous observations say that to get a matrix $A$ into a desired form $B$, we may just work on a `small part' of the matrix $A$, \eg\  a $2 \times 2$-submatrix or a fragment of a column, by  applying a transformation  that affects that `small part'. We then repeat it on other parts to obtain a sequence of transformations:
\begin{align}\label{eq:sequence}
A &\to X_1A \to X_2(X_1A ) \to \dots \to B,\notag \\*
A &\to X_1^{-\tp}A \to X_2^{-\tp}(X_1^{-\tp} A ) \to \dots \to B,\notag \\
A &\to X_1AX_1^\tp \to X_2(X_1 A X_1^\tp )X_2^\tp \to \dots \to B,\notag \\*
A &\to X_1AX_1^{-1} \to X_2(X_1 A X_1^{-1})X_2^{-1} \to \dots \to B,\notag \\*
A &\to X_1AY_1^{-1} \to X_2(X_1 A Y_1^{-1})Y_2^{-1} \to \dots \to B,
\end{align}
and piece all change-of-basis matrices together to get the required $X$ as either $X_m X_{m-1} \dots X_1$ or its limit as $m \to \infty$ (likewise for $Y$). Algorithms for computing standard matrix decompositions such as LU, QR, EVD, SVD, Cholesky, Schur, {\em etc.}, all involve applying a sequence of such transformation rules \cite{GVL}. In numerical linear algebra, if $m$ is finite, \ie\ $X$ may be obtained in finitely many steps (in exact arithmetic), then the algorithm is called \emph{direct}, whereas if it requires $m \to \infty$, \ie\ $X$ may only be approximated with a limiting process, then it is called \emph{iterative}.  Furthermore, in numerical linear algebra, one tends to see such transformations as giving a \emph{matrix decomposition}, which may then be used to solve other problems involving $A$. This is sometimes called `the decompositional approach to matrix computation' \cite{Stewart}.

Designers of algorithms for matrix computations, even if they were not explicitly aware of these transformation rules and properties, were certainly observing them implicitly. For instance, it is rare to find algorithms that mix different transformation rules for different types of tensors, since what is incompatible for tensors tends to lead to meaningless results. Also, since eigenvalues are defined for mixed $2$-tensors but not for contravariant or covariant $2$-tensors, the transformation $A' = XAX^{-1}$ is pervasive in algorithms for eigenvalue decomposition but we rarely if ever find $A' = XAX^\tp$ (unless of course $X$ is orthogonal).

In numerical linear algebra, the use of the transformation rules in \eqref{eq:sequence} goes hand in hand with a salient property of the group of change-of-basis matrices.

\begin{example}[Givens rotations,\, Householder reflectors,\, Gauss transforms]\label{eg:GHG}
Recall that these are defined by 
\begin{align*}
G &= 
       {\begin{bmatrix}   1   & \cdots &    0   & \cdots &    0   & \cdots &    0   \\
                      \vdots & \ddots & \vdots &        & \vdots &        & \vdots \\
                         0   & \cdots & \cos \theta  & \cdots & -\sin \theta & \cdots &    0   \\
                      \vdots &        & \vdots & \ddots & \vdots &        & \vdots \\
                         0   & \cdots &  \sin \theta   & \cdots & \cos \theta & \cdots &    0   \\
                      \vdots &        & \vdots &        & \vdots & \ddots & \vdots \\
                         0   & \cdots &    0   & \cdots &    0   & \cdots &    1
       \end{bmatrix}}  \in \SOr(n), \\
H &= I -\dfrac{2vv^\tp}{v^\tp v} \in \Or(n), \quad
M = I -ve_i^\tp \in \GL(n) ,
\end{align*}
with the property that
the Givens rotation $a' = Ga$ is a rotation of $a$ in the $(i,j)$-plane by an angle $\theta$,
the Householder reflection
$a' = Ha$ is the reflection of $a$ in the hyperplane with normal $v/\lVert v \rVert$ and,  for a judiciously chosen $v$,
the Gauss transform $a' = Ma$ is in $\spn \{e_{i+1},\ldots,e_n\}$, \tie, it has $(i+1)$th to $n$th coordinates zero.\footnote{If $m = \alpha e_j$, then $MA$ adds an $\alpha$ multiple of the $j$th row to the $i$th row of $A$; so the Gauss transform includes the elementary matrices that perform this operation.}
Standard algorithms in numerical linear algebra such as Gaussian elimination for LU decomposition and Cholesky decomposition, Givens and Householder QR for QR decomposition, Francis's QR or Rutishauser's LR algorithms for EVD, Golub--Kahan bidiagonalization for SVD, {\em etc.}, all rely in part on applying a sequence of transformation rules as in \eqref{eq:sequence} with one of these matrices playing the role of the change-of-basis matrices. The reason this is possible is that:
\begin{itemize}
\setlength\itemsep{3pt}
\item any $X \in \SOr(n)$ is a product of Givens rotations,
\item any $X \in \Or(n)$ is a product of Householder reflectors,
\item any $X \in \GL(n)$ is a product of elementary matrices,
\item any unit lower triangular $X \in \GL(n)$ is a product of Gauss transforms.
\end{itemize}
In group-theoretic lingo, these matrices are \emph{generators} of the respective matrix Lie groups; in the last case, the set of all unit lower triangular matrices, \ie\ ones on the diagonal, is also a subgroup of $\GL(n)$.
\end{example}

Whether one seeks to solve a system of linear equations or find a least-squares solution or compute eigenvalue or singular value decompositions, the basic underlying principle in numerical linear algebra is  to transform the problem in such a way that the solution of the transformed problem is related to the original solution in a definite way; note that this is practically a paraphrase of Voigt's definition of a tensor on page~\pageref{Voigt}. Any attempt to give a comprehensive list of examples will simply result in our reproducing a large fraction of \citet{GVL}, so we will just give a familiar example viewed through the lens of tensor transformation rules.

\begin{example}[full-rank\, least-squares]\label{eg:qr}
As we saw in Section~\ref{sec:transrules}, the least-squares problem \eqref{eq:ls} satisfies the transformation rule of a mixed $2$-tensor $A' = XAY^{-1}$ with change-of-basis matrices $(X,Y) \in \Or(m) \times \GL(n)$. Suppose $\rank(A) = n$. Then, applying a sequence of covariant $1$-tensor transformation rules
\[
A \to Q_1^\tp A \to Q_2^\tp (Q_1^\tp A ) \to \dots \to Q^\tp A = \begin{bmatrix}
R\\
0
\end{bmatrix}
\]
given by the \emph{Householder QR algorithm}, we get
\[
A=Q
\begin{bmatrix}
R\\
0
\end{bmatrix}\!.
\]
As the minimum value is an invariant Cartesian $0$-tensor,
\begin{align*}
\min\, \lVert Av -b\rVert^2 &=\min\, \lVert Q^\tp
(Av-b)\rVert^2  =\min\left\lVert
\begin{bmatrix}
R\\
0
\end{bmatrix}
v-Q^\tp b\right\rVert^2\\*
&=\min\left\Vert
\begin{bmatrix}
R\\
0
\end{bmatrix}
v-
\begin{bmatrix}
c\\
d
\end{bmatrix}
\right\rVert^2=\min\, \lVert Rv-c\rVert^2 +\lVert d\rVert^2=\lVert d\rVert^2,
\end{align*}
where we have written
\[
Q^\tp b=
\begin{bmatrix}
c\\
d
\end{bmatrix}\!.
\]
In this case the solution of the transformed problem $Rv=c$ is in fact equal to that of the least-squares problem, and may be obtained through \emph{back-substitution}, \tie, another sequence of contravariant $1$-tensor transformation rules
\[
c \to Y_1^{-1}c \to Y_2^{-1}(Y_1^{-1} c ) \to \dots \to R^{-1}c = v,
\]
where the $Y_i$ are Gauss transforms. As noted above, the solution method reflects Voigt's definition: we transform the problem $\min \lVert Av -b\rVert^2$ into a form where the solution of the transformed problem $Rv=c$ is related to the original solution in a definite way.  Here we obtained the change-of-basis matrices $X = Q \in \Or(m)$ and $Y = R^{-1} \in \GL(n)$ via Householder QR and back-substitution respectively.
\end{example}

As we just saw, there are distinguished choices for the change-of-basis matrices that aid the solution of a numerical linear algebra problem. We will mention one of the most useful choices below.

\begin{example}[Krylov\, subspaces]\label{eg:Krylov}\hspz%
Suppose we have $A \in \mathbb{R}^{n \times n}$, with all  eigenvalues distinct and non-zero for simplicity. Take an arbitrary $b \in \mathbb{R}^n$. Then the matrix $K$ whose columns are \emph{Krylov basis} vectors of the form
\[
b, Ab, A^2 b, \ldots, A^{n-1} b
\]
is invertible, \ie\ $K \in \GL(n)$, and using it as the change-of-basis matrix gives us
\begin{equation}\label{eq:krylov}
A = K \begin{bmatrix}
0 & 0 & \cdots & 0 & -c_0 \\
1 & 0 & \cdots & 0 & -c_1 \\
0 & 1 & \cdots & 0 & -c_2 \\
\vdots & \vdots & \ddots & \vdots & \vdots \\
0 & 0 & \cdots & 1 & -c_{n-1}
\end{bmatrix} K^{-1}.
\end{equation}
This is a special case of the rational canonical form, a canonical form under similarity. The seemingly trivial observation \eqref{eq:krylov}, when combined with other techniques, becomes a powerful iterative method for a wide variety of computational tasks such as solving linear systems, least-squares, eigenvalue problems or evaluating various matrix functions \cite{VDV}. Readers unfamiliar with numerical linear algebra may find it odd that we do not use another obvious canonical form, one that makes the aforementioned problems trivial to solve, namely, the eigenvalue decomposition
\begin{equation}\label{eq:eigen}
A = X
\begin{bmatrix}
\lambda_1 & 0 & 0 & \cdots & 0 \\
0 & \lambda_2 & 0 &\cdots & 0 \\
0 & 0 & \lambda_3 & \cdots  & 0 \\
\vdots & \vdots & \vdots & \ddots & \vdots \\0 & 0 & 0 & \cdots & \lambda_n
\end{bmatrix}
X^{-1},
\end{equation}
where the change-of-basis matrix $X \in \GL(n)$ has columns given by the eigenvectors of $A$. The issue is that this is more difficult to compute than \eqref{eq:krylov}. In fact, to compute it, one way is to implicitly exploit \eqref{eq:krylov} and the relation between the two canonical forms:
\begin{equation*} 
\begin{bmatrix}
\lambda_1 & 0 & 0 & \cdots & 0 \\
0 & \lambda_2 & 0 &\cdots & 0 \\
0 & 0 & \lambda_3 & \cdots  & 0 \\
\vdots & \vdots & \vdots & \ddots & \vdots \\0 & 0 & 0 & \cdots & \lambda_n
    \end{bmatrix}=
V \begin{bmatrix}0 & 0 & \cdots & 0 & -c_0 \\
1 & 0 & \cdots & 0 & -c_1 \\
0 & 1 & \cdots & 0 & -c_2 \\
\vdots & \vdots & \ddots & \vdots & \vdots \\
0 & 0 & \cdots & 1 & -c_{n-1}
\end{bmatrix} V^{-1},
\end{equation*}
where
\begin{equation}\label{eq:vander}
V=
\begin{bmatrix}
1 & \lambda_1 & \lambda_1^2 & \dots & \lambda_1^{n-1}\\
1 & \lambda_2 & \lambda_2^2 & \dots & \lambda_2^{n-1}\\
1 & \lambda_3 & \lambda_3^2 & \dots & \lambda_3^{n-1}\\
\vdots & \vdots & \vdots & \ddots &\vdots \\
1 & \lambda_n & \lambda_n^2 & \dots & \lambda_n^{n-1}
\end{bmatrix}\!. 
\end{equation}
The pertinent point here is that the approach of solving a problem by finding an appropriate computational basis is also an instance of the $2$-tensor transformation rule. Aside from the Krylov basis, there are simpler examples such as diagonalizing a circulant matrix
\[
C =  \begin{bmatrix}
c_0     & c_{n-1} & \dots  & c_{2} & c_{1}  \\
c_{1} & c_0    & c_{n-1} &         & c_{2}  \\
\vdots  & c_{1}& c_0    & \ddots  & \vdots  \\
c_{n-2}  &        & \ddots & \ddots  & c_{n-1} \\
c_{n-1}  & c_{n-2} & \dots  & c_{1} & c_0   \\
\end{bmatrix}
\]
by expressing it in the Fourier basis, \ie\ with $\lambda_j^k = \rme^{2(j-1)k \pi \rmi/n}$ in \eqref{eq:vander}, and there are far more complex examples such as the wavelet approach to the fast multipole method of \citet{Beylkin}. This computes, in $O(n)$ operations and to arbitrary accuracy, a matrix--vector product $v \mapsto Av$ with a \emph{dense} $A \in \mathbb{R}^{n \times n}$ that is a finite-dimensional approximation of certain special integral transforms. The algorithm iteratively computes a special wavelet basis $X \in \GL(n)$ so that ultimately $X^{-1} A X =   B $ gives a banded matrix $B$ where $ v\mapsto Bv$ can be computed in time $O(n \log(1/\varepsilon))$ to $\varepsilon$-accuracy, and where $v \mapsto Xv$, $v \mapsto X^{-1}v$ are both computable in time $O(n)$. One may build upon this algorithm to obtain an $O(n \log^2 n \log \kappa(A))$ algorithm for the pseudoinverse \cite[Section~6]{Beylkin3} and potentially also to compute other matrix functions such as square root, exponential, sine and cosine, {\em etc.}\ \cite[Section~X]{Beylkin2}.
 We will discuss some aspects of this basis in Example~\ref{eg:wavelet}, which is constructed as the tensor product of multiresolution analyses.
\end{example} 

Our intention is to highlight the role of the tensor transformation rules in numerical linear algebra but we do not wish to overstate it. These rules are an important component of various algorithms but almost never the only one. Furthermore, far more work is required to guarantee  correctness in finite-precision arithmetic \cite{HighamBook}, although the tensor transformation rules can sometimes help with that too, as we will see next.

\begin{example}[affine invariance of Newton's method]\label{eg:Newton}
The discussion below is essentially the standard treatment in \citet{Boyd},
supplemented by commentary pointing out the tensor transformation rules. Consider the equality-constrained convex optimization problem
\begin{equation}\label{eq:eqconstopt}
\begin{tabular}{rll}
&minimize & $ f(v)$\\
&subject to & $Av = b$
\end{tabular}
\end{equation}
for a strongly convex $f\in C^2(\Omega)$ with $\beta I \preceq \nabla^2 f(v) \preceq \gamma I$.
The Newton step $\Delta v \in \mathbb{R}^n$ is defined as the solution to
\[
\begin{bmatrix}
\nabla^2 f(v) & A^\tp \\
A & 0
\end{bmatrix}
\begin{bmatrix}
\Delta v \\
\Delta \lambda
\end{bmatrix}
=
\begin{bmatrix}
-\nabla f(v) \\
0
\end{bmatrix}
\]
and the Newton decrement $\lambda(v) \in \mathbb{R}$ is defined as
\[
\lambda(v)^2 \coloneqq \nabla (v)^\tp \nabla f(v)^{-1} \nabla f(v).
\]
Let $X \in \GL(n)$ and suppose we perform a linear change of coordinates $X v' = v$. Then the Newton step in these new coordinates is given by
\[
\begin{bmatrix}
X^\tp \nabla^2 f(Xv) X & X^\tp A^\tp \\
A X & 0
\end{bmatrix}
\begin{bmatrix}
\Delta v' \\
\Delta \lambda'
\end{bmatrix}
=
\begin{bmatrix}
-X^\tp \nabla f(Xv) \\
0
\end{bmatrix}\!.
\]
We may check that $X \Delta v' = \Delta v$ and thus the iterates are related by $Xv'_k = v_k$ for all $k \in \mathbb{N}$ as long as we initialize with $Xv_0' = v_0$ \cite[Section~10.2.1]{Boyd}. Note that steepest descent satisfies no such property no matter which $1$-tensor transformation rule we use: $v' = Xv$, $v' = X^{-1}v$, $v' = X^\tp v$, or $v' = X^{-\tp}v$.
We also have that  $\lambda(Xv)^2 = \lambda(v)^2$, which is used in the stopping condition of Newton's method; thus the iterations stop at the same point  $Xv'_k = v_k$ when $\lambda(v'_k)^2 = \lambda(v_k)^2 \le 2 \varepsilon$ for a given $\varepsilon > 0$.  In summary, if we write $g(v') = f(Xv)$, then we have the following relations:
\begin{alignat*}{5}
&\text{coordinates} & &\text{contravariant $1$-tensor} & v' &= X^{-1} v,\\
&\text{gradient} & &\text{covariant $1$-tensor} &\nabla g(v') &= X^\tp \nabla f(Xv),\\
&\text{Hessian} & &\text{covariant $2$-tensor} \qquad\quad &\nabla^2 g(v') &= X^\tp \nabla^2 f(Xv) X,\\
&\text{Newton step} & &\text{contravariant $1$-tensor} &\Delta v' &= X^{-1} \Delta v,\\
&\text{Newton iterate} & &\text{contravariant $1$-tensor} &v'_k &= X^{-1}  v_k,\\
&\text{Newton decrement} &\quad &\text{invariant $0$-tensor} &\lambda(v'_k) &= \lambda(v_k).
\end{alignat*}
Strictly speaking, the gradient and Hessian are tensor fields, and we will explain the difference in Example~\ref{eg:tenfield}. We may extend the above discussion with the general affine group $\GA(n)$ in \eqref{eq:SE} in place of $\GL(n)$ to incorporate translation by a vector.

To see the implications on computations, it suffices to take, for simplicity, an unconstrained problem with $A = 0$ and $b= 0$ in \eqref{eq:eqconstopt}. Observe that the condition number of  $X^\tp \nabla^2 f(Xv) X $ can be scaled to any desired value in $[1,\infty)$ with an appropriate $X \in \GL(n)$. In other words, any Newton step is independent of the condition number of $\nabla^2 f(v) $ in exact arithmetic. In finite precision arithmetic, this manifests itself as insensitivity to condition number. Newton's method gives solutions of high accuracy for condition number as high as $\kappa \approx 10^{10}$ when steepest descent already stops working at $\kappa \approx 20$ \cite[Section~9.5.4]{Boyd}. Strictly speaking, $\kappa$ refers to the condition number of sublevel sets $\{ v \in \Omega \colon f(v) \le \alpha \})$. For any convex set $C \subseteq \mathbb{R}^n$, this is given by
\[
\kappa(C) \coloneqq \frac{\sup_{\lVert v \rVert =1} w(C,v)^2}{\inf_{\lVert v \rVert =1} w(C,v)^2}, \quad
w(C,v) \coloneqq \sup_{u \in C} u^\tp v - \inf_{u \in C} u^\tp v.
\]
However, if $\beta I \preceq \nabla^2 f(v) \preceq \gamma I$,
then  $\kappa(\{ v \in \Omega \colon f(v) \le \alpha \}) \le \gamma/\beta$ \cite[Section~9.1.2]{Boyd}, and so it is ultimately controlled by $\kappa (\nabla^2 f(v) )$. 
\end{example} 

Our last example shows that the tensor transformation rules are as important in the information sciences as they are  in the physical sciences.

\begin{example}[equivariant neural networks]\label{eg:neural}
  A feed-forward neural network is usually regarded as a function $f \colon \mathbb{R}^n \to \mathbb{R}^n$
  obtained by alternately composing affine maps $\alpha_i \colon \mathbb{R}^n \to \mathbb{R}^n $, $i=1,\dots,k$, with a non-linear function $\sigma \colon  \mathbb{R}^n \to \mathbb{R}^n$:
\[
\begin{tikzcd}
\mathbb{R}^n \arrow[urrrrrrd, bend left=25, "f = \alpha_k \circ \sigma \circ \alpha_{k -1} \circ \cdots \circ \sigma \circ \alpha_2 \circ \sigma \circ \alpha_1"]\arrow{r}{\alpha_1} & \mathbb{R}^n \arrow{r}{\sigma} & \mathbb{R}^n \arrow{r}{\alpha_2} & \mathbb{R}^n \arrow{r}{\sigma} & \cdots  \arrow{r}{\sigma} & \mathbb{R}^n \arrow{r}{\alpha_k} & \mathbb{R}^n.
\end{tikzcd}
\]
The \emph{depth}, also known as the number of \emph{layers}, is $k$ and the \emph{width}, also known as the number of \emph{neurons}, is $n$. We assume that our neural network has constant width throughout all layers. The non-linear function $\sigma$ is called \emph{activation} and we may assume that it is given by the ReLU function
$\sigma(t) \coloneqq \max(t,0)$ for $t \in \mathbb{R}$ and, by convention,
\begin{equation}\label{eq:ptwise}
\sigma(v) = (\sigma(v_1),\dots,\sigma(v_n)), \quad v = (v_1,\dots,v_n) \in \mathbb{R}^n.
\end{equation}
The affine function is defined by $\alpha_i(v) = A_iv + b_i$ for some $A_i \in \mathbb{R}^{n \times n}$ called the \emph{weight} matrix and some $b_i \in \mathbb{R}^n$ called the \emph{bias} vector.
We assume that $b_k = 0$ in the last layer.

Although convenient, it is somewhat misguided to be lumping the bias and weight together in an affine function. The biases $b_i$ are intended to serve as \emph{thresholds} for the activation function $\sigma$ \cite[Section~4.1.7]{Bishop} and should be part of it, detached from the weights $A_i$ that transform the input. If one would like to incorporate translations, one could do so with weights from a matrix group such as $\SE(n)$ in \eqref{eq:SE}. Hence a better but mathematically equivalent description of $f$ would be~as
\[
\begin{tikzcd}
\mathbb{R}^n \arrow[urrrrrrd, bend left=25, "f = A_k \sigma_{b_{k-1}}  A_{k -1}  \cdots  \sigma_{b_2} A_2  \sigma_{b_1}  A_1"]\arrow{r}{A_1} & \mathbb{R}^n \arrow{r}{\sigma_{b_1}} & \mathbb{R}^n \arrow{r}{A_2} & \mathbb{R}^n \arrow{r}{\sigma_{b_2}} & \cdots  \arrow{r}{\sigma_{b_{k-1}}} & \mathbb{R}^n \arrow{r}{A_k} & \mathbb{R}^n
\end{tikzcd}
\]
where we identify $A_i \in \mathbb{R}^{n \times n}$ with the linear operator $\mathbb{R}^n \to \mathbb{R}^n$, $v \mapsto A_i v$, and for any $b \in \mathbb{R}^n$ we define $\sigma_b \colon \mathbb{R}^n \to \mathbb{R}^n$ by $\sigma_b(v) = \sigma(v+b)\in \mathbb{R}^n$. We have dropped the composition symbol $\circ$ to avoid clutter and potential confusion with the Hadamard product \eqref{eq:hp}, and will continue to do so below. For a fixed  $\theta \in \mathbb{R}$,
\[
\sigma_\theta(t) = \begin{cases} t -\theta & t \ge \theta, \\ 0 & t < \theta, \end{cases}
\]
plays the role of a threshold for activation as was intended by \citet[p.~392]{Rosenblatt} and \citet[p.~120]{MP}.

A major computational issue with neural networks is the large number of unknown parameters, namely the $kn^2 + k(n-1)$ entries of the weights and biases, that have to be fitted with data, especially for deep neural networks where $k$ is  large.  Thus successful applications of neural networks require that we identify, based on the problem at hand, an appropriate low-dimensional subset of $\mathbb{R}^{n \times n}$ from which we will find our weights $A_1,\dots,A_k$. For instance, the very successful convolutional neural networks for image recognition \cite{AlexNet} relies on restricting $A_1,\dots,A_k$ to some block-Toeplitz--Toeplitz-block or BTTB matrices \cite[Section~13]{struct} determined by a very small number of parameters. It turns out that convolutional neural networks are a quintessential example of equivariant neural networks \cite{Welling}, and in fact every equivariant neural network may be regarded as a generalized convolutional neural network in an appropriate sense \cite{Risi1}. We will describe a simplified version that captures its essence and illustrates the tensor transformation rules.

Let $G \subseteq \mathbb{R}^{n \times n}$ be a matrix group.  A function $f \colon \mathbb{R}^n \to \mathbb{R}^n$ is said to be equivariant if it satisfies the condition that
\begin{equation}\label{eq:ENN0}
f(Xv) = X f(v)  \quad \text{for all } v\in \mathbb{R}^n,\; X \in G.
\end{equation}
An equivariant neural network is simply a feed-forward neural network $f \colon \mathbb{R}^n \to \mathbb{R}^n$ that satisfies \eqref{eq:ENN0}. The key to constructing an equivariant neural network is the trivial observation that
\begin{align*}
f(Xv) &= A_k \sigma_{b_{k-1}}  A_{k -1}  \cdots  \sigma_{b_2} A_2  \sigma_{b_1}  A_1 Xv\\*
&=X (X^{-1}A_k X) (X^{-1}\sigma_{b_{k-1}} X) (X^{-1} A_{k -1} X) \cdots \\*
&\qquad \cdots (X^{-1} \sigma_{b_2} X) (X^{-1} A_2 X)(X^{-1} \sigma_{b_1} X) (X^{-1} A_1 X) v\\*
&= XA_k' \sigma_{b_{k-1}}'  A_{k -1}'  \cdots  \sigma_{b_2}' A_2'  \sigma_{b_1}'  A_1' v
\end{align*}
and the last expression equals $X f(v)$ as long as we have
\begin{equation}\label{eq:ENN1}
A_i'  = X^{-1} A_i X = A_i, \quad \sigma_{b_i}' = X^{-1} \sigma_{b_i} X = \sigma_{b_i}
\end{equation}
for all $i=1,\dots,k$, and for all $X \in G$. The condition on the right is satisfied by any activation as long as it is applied coordinate-wise \cite[Section~6.2]{Welling} as we did in \eqref{eq:ptwise}. We will illustrate this below. The condition on the left, which says that $A_i$ is invariant under the mixed $2$-tensor transformation rule on page~\pageref{transrulesummary2}, is what is important. The condition limits the possible weights for $f$ to a (generally) much smaller subspace of matrices  that commute with all elements of $G$. However, finding this subspace (in fact a subalgebra) of \emph{intertwiners},
\begin{equation}\label{eq:intertwining}
\{A \in \mathbb{R}^{n \times n} \colon AX =XA \text{ for all } X \in G \},
\end{equation}
is usually where the challenge lies, although it is also a classical problem with numerous powerful tools  from group representation theory \cite{Brocker,Fulton,Mackey}.

There are many natural candidates for the group $G$ from the list in \eqref{eq:groups} and their subgroups, depending on the problem at hand. But we caution the reader that $G$ will generally be a very low-dimensional subset of  $\mathbb{R}^{n \times n}$. It will be pointless to pick, say, $G = \SOr(n)$ as the set in \eqref{eq:intertwining} will then be just $\{\lambda I \in \mathbb{R}^{n \times n}  \colon \lambda \in \mathbb{R} \}$, clearly too small to serve as meaningful weights for any neural network. Indeed, $G$ is usually chosen to be an image of a much lower-dimensional group $H$,
\[
G =  \rho(H) = \{ \rho(h) \in \GL(n) \colon h \in H\}
\]
for some group homomorphism $\rho \colon H \to \GL(n)$; here $\rho$ is called a \emph{representation}  of the group $H$. In image recognition applications, possibilities for $H$ are often discrete subgroups of $\SE(2)$ or $\E(2)$ such as the group of translations
\[
H = \biggl\{ \begin{bsmallmatrix} 1 & 0 & m_1 \\ 0 &  1 & m_2 \\  0 & 0 & 1\end{bsmallmatrix} \in \mathbb{R}^{3 \times 3} \colon  m_1,m_2 \in \mathbb{Z} \biggr\}
\]
for convolutional neural networks; the $p4$ group that augments translations with right-angle rotations,
\[
H = \biggl\{ \begin{bsmallmatrix} \cos(k\pi/2) & -\sin(k \pi/2) & m_1 \\  \sin(k\pi/2)  &   \cos(k\pi/2)  & m_2 \\  0 & 0 & 1\end{bsmallmatrix} \in \mathbb{R}^{3 \times 3} \colon k =0,1,2,3,\;  m_1,m_2 \in \mathbb{Z} \biggr\}
\]
in \citet[Section~4.2]{Welling}; the $p4m$ group that further augments $p4$ with reflections,
\[
H = \biggl\{ \begin{bsmallmatrix} (-1)^j \cos(k\pi/2) & (-1)^{j+1} \sin(k \pi/2) & m_1 \\  \sin(k\pi/2)  &   \cos(k\pi/2)  & m_2 \\  0 & 0 & 1\end{bsmallmatrix} \in \mathbb{R}^{3 \times 3} \colon \begin{multlined} k =0,1,2,3,\\  j=0,1,\;  m_1,m_2 \in \mathbb{Z} \end{multlined}\biggr\}
\]
in \citet[Section~4.3]{Welling}. Other possibilities for $H$ include the rotation group $\SOr(3)$ for 3D shape recognition \cite{Risi2}, the rigid motion group $\SE(3)$ for chemical {property} \cite{Fuchs} 
and protein structure (see page~\pageref{pg:alpha}) predictions and the Lorentz group $\SOr(1,3)$ for identifying top quarks in high-energy physics experiments \cite{Lorentz}, {\em etc.}

When people speak of $\SOr(3)$- or $\SE(3)$- or Lorentz-equivariant neural networks, they are referring to the group $H$ and not $G = \rho(H)$. A key step of these works is the construction of an appropriate representation $\rho$ for the problem at hand, or, equivalently, constructing a linear action of $H$ on $\mathbb{R}^n$. In these applications $\mathbb{R}^n$ should be regarded as the set of real-valued functions on a set $S$ of cardinality $n$, a perspective that we will introduce in Example~\ref{eg:hyp}. For concreteness, take an image recognition problem on a $60\,000$ collection of $28 \times 28$-pixel images of handwritten digits in greyscale levels $0,1,\dots,255$ \cite{mnist}. Then $S \subseteq \mathbb{Z}^2$ is the set of $n = 28^2 = 784$ pixel indices and an image is encoded as $v \in \mathbb{R}^{784}$ whose coordinates take values from $0$ (pitch black) to $255$ (pure white). Note that this is why the first three $H$ above are discrete: the elements $h \in H$ act on the pixel indices $S  \subseteq \mathbb{Z}^2$ instead of $\mathbb{R}^2$. These $60\,000$ images are then used to fit the neural network $f$, \ie\ to find the parameters $A_1,\dots,A_k \in \mathbb{R}^{784 \times 784}$ and $b_1,\dots,b_{k-1} \in \mathbb{R}^{784}$.

We end this example with a note on why a non-linear (and not even multilinear) function like ReLU can satisfy the covariant $2$-tensor transformation rule, which sounds incredible but is actually obvious once explained. Take a greyscale image $v \in \mathbb{R}^n$ drawn with a black outline (greyscale value $0$) but filled with varying shades of grey (greyscale values $1,\dots,255$) and consider the activation
\[
\sigma(t) = \begin{cases} 255 & t > 0, \\ 0 & t \le 0, \end{cases} \qquad t \in \mathbb{R},
\]
so that applying $\sigma$ to the image $v$ produces an image $\sigma(v) \in \mathbb{R}^n$ with all shadings removed, leaving just the black outline. Now take a $45^\circ$ rotation matrix $R \in \SOr(2)$ and let $X = \rho(R) \in \GL(n)$ be the corresponding matrix that rotates any image $v$ by $45^\circ$ to $Xv$:
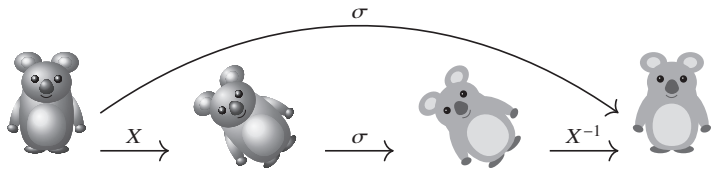
\[
\begin{tikzcd}
\tikz{\koala[3D,eye=gray,scale=0.6]} \arrow[urrrd, bend left=35, "\sigma"]\arrow{r}{X} & \tikz{\koala[3D,eye=gray,rotate=45,scale=0.6]} \arrow{r}{\sigma} & \tikz{\koala[rotate=45,contour=black,scale=0.6]} \arrow{r}{X^{-1}} & \tikz{\koala[contour=black,scale=0.6]}
\end{tikzcd}
\]
The bottom line is that $R$ acts on the indices of $v$ whereas $\sigma$ acts on the values of $v$ and the two actions are always independent, which is why $\sigma' = X^{-1}  \sigma X= \sigma$.

Our definition of equivariance in \eqref{eq:ENN0} is actually a simplification but suffices for our purposes. We state a formal definition here for completeness. A function $f \colon \mathbb{V} \to \mathbb{W}$ is said to be $H$-equivariant if there are two homomorphisms $\rho_1 \colon H \to \GL(\mathbb{V})$ and $\rho_2 \colon H \to \GL(\mathbb{W})$ so that
\[
f(\rho_1 (h) v) = \rho_2(h) f(v) \quad \text{for all } v\in \mathbb{V},\; h \in H.
\]
In \eqref{eq:ENN0}, we chose $\mathbb{V} = \mathbb{W} = \mathbb{R}^n$ and $\rho_1 = \rho_2$. The condition $f(Xv) = X f(v) $ could have been replaced by other tensor or pseudotensor transformation rules, say, $f(Xv) = \det(X)X^{-\tp} f(v) $, as long as they are given by homomorphisms. In fact, with slight modifications we may even allow for \emph{antihomomorphisms}, \ie\ $\rho(h_1 h_2) = \rho(h_2) \rho(h_1)$.
\end{example} 

There are many other areas in computations where the tensor transformation rules make an appearance. For example, both the primal and dual forms of a cone programming problem over a symmetric cone $\mathbb{K} \subseteq \mathbb{V}$ -- which include as special cases linear programming (LP), convex quadratic programming (QP), second-order cone programming (SOCP) and semidefinite programming (SDP) -- conform to the transformation rules for Cartesian $0$-, $1$- and $2$-tensors. However, the change-of-basis matrices would have to be replaced by a linear map from the \emph{orthogonal group of the cone} \cite{Hauser}:
\[
\Or(\mathbb{K}) \coloneqq \{ \varphi \colon  \mathbb{V} \to \mathbb{V} \colon  \varphi \text{ linear, invertible and } \varphi^* = \varphi^{-1} \}.
\]
Here the vector space $\mathbb{V}$ involved may not be $\mathbb{R}^n$  --  for SDP, $\mathbb{K} = \mathbb{S}_\pp^n$ and $\mathbb{V} = \mathbb{S}^n$, the space of $n\times n$ symmetric matrices  --  making the linear algebra notation we have been using to describe definition~\ref{st:tensor1} awkward and unnatural. Instead we should make provision to work with tensors over arbitrary vector spaces, for example the space of Toeplitz or Hankel or Toeplitz-plus-Hankel matrices, the space of polynomials or differential forms or differential operators, or, in the case of equivariant neural networks, the space of $L^2$-functions on homogeneous spaces \cite{Welling,Risi1}. This serves as another motivation for definitions~\ref{st:tensor2} and~\ref{st:tensor3}.

\subsection{Fallacy of `tensor $=$ hypermatrix'}

We conclude our discussion of definition~\ref{st:tensor1} with a few words about the fallacy of identifying a tensor with the hypermatrix that represents it. While we have alluded to this from time to time over the last few sections, there are three points that we have largely deferred until now.

Firstly, the `multi-indexed object' in definition~\ref{st:tensor1} is not necessarily a hypermatrix; they could be polynomials or exterior forms or non-commutative polynomials. For example, it is certainly more fruitful to represent symmetric tensors  as polynomials and alternating tensors as exterior forms. For one, we may take derivatives and integrals  or evaluate them at a point  --  operations that become awkward if we simply regard them as hypermatrices with symmetric or skew-symmetric entries.

Secondly, if one subscribes to this fallacy, then one would tend to miss important tensors hiding in plain sight. The object of essence in each of the following examples is a $3$-tensor, but one sees no triply indexed quantities anywhere.
\begin{enumerate}[\upshape (i)]
\setlength\itemsep{3pt}
\item\label{it:cplx} Multiplication of complex numbers:
\[
(a + ib)(c+ id) = (ac - bd) + i (bc+ad).
\]

\item\label{it:mm} Matrix--matrix products:
\[
\begin{bmatrix}
a_{11} & a_{12} \\
a_{21} & a_{22}
\end{bmatrix}
\begin{bmatrix}
b_{11} & b_{12} \\
b_{21} & b_{22}
\end{bmatrix}
=
\begin{bmatrix}
a_{11} b_{11} + a_{12} b_{21} & a_{11} b_{12} + a_{12} b_{22}\\
a_{21} b_{11} + a_{22} b_{21} & a_{21} b_{12} + a_{22} b_{22}
\end{bmatrix}\!.
\]

\item\label{it:GI} Grothendieck's inequality:
\[
\max_{\lVert x_i\rVert = \lVert y_j \rVert = 1}\Bigg| \sum_{i=1}^m\sum_{j=1}^n a_{ij} \langle x_i , y_j \rangle\Bigg|\le  K_\G\max_{\lvert \varepsilon_i \rvert = \lvert \delta_j \rvert = 1}\Bigg|\sum_{i=1}^m\sum_{j=1}^n a_{ij} \varepsilon_i \delta_j\Bigg|.
\]

\item\label{it:sep} Separable representations:
\begin{align*}
\sin(x + y + z) &  = \sin(x) \cos(y) \cos(z) + \cos(x) \cos(y) \sin(z) \\*
&\quad + \cos(x) \sin(y) \cos(z) - \sin(x) \sin(y) \sin(z).
\end{align*}
\end{enumerate}
Thirdly, the value of representing a $d$-tensor as a $d$-dimensional hypermatrix is often overrated.\label{pg:overrate}
\begin{enumerate}[\upshape (a)]
\setlength\itemsep{3pt}
\item \emph{It may not be meaningful.} For tensors in, say, $L^2(\mathbb{R}) \otimes L^2(\mathbb{R}) \otimes L^2(\mathbb{R})$, where the `indices' are continuous, the hypermatrix picture is not  meaningful. 

\item \emph{It may not be computable.} Given a bilinear operator, writing down its entries as a three-dimensional hypermatrix with respect to bases is in general a \#P-hard problem, as we will see in Example~\ref{eg:LRcoeff}.

\item \emph{It may not be useful.} In applications such as tensor networks, $d$ is likely large. While there may be some value in picturing a three-dimensional hypermatrix, it is no longer the case when $d \gg 3$.

\item \emph{It may not be possible.} There are applications such as Examples~\ref{eg:Gauss} and \ref{eg:int} that involve tensors defined over modules (vector spaces whose field of scalars is replaced by a ring). Unlike vector spaces, modules may not have bases and such tensors cannot be represented as hypermatrices.
\end{enumerate}
By far the biggest issue with thinking that a tensor is just a hypermatrix is that such an overly simplistic view disregards the transformation rules, which is the key to definition~\ref{st:tensor1} and  whose importance we have already discussed at length.  As we saw, $\begin{bsmallmatrix}1&2&3\\ 2 &3&4 \\3 &4 &5\end{bsmallmatrix}$ may represent a $1$-tensor or a $2$-tensor, it may represent a covariant or contravariant or mixed tensor, or it may not represent a tensor at all. What matters in definition~\ref{st:tensor1} is the transformation rules, not the multi-indexed object. If one thinks that a fraction is just a pair of numbers, adding them as $a/b + c/d = (a+c)/(b+d)$ is going to seem perfectly fine. It is the same with tensors.

\section{Tensors via multilinearity}\label{sec:mult}

The main deficiency of definition~\ref{st:tensor1}, namely that it leaves the tensor unspecified, is immediately remedied  by definition~\ref{st:tensor2}, which states unequivocally that the tensor is a multilinear map. The multi-indexed object in definition~\ref{st:tensor1} is then a coordinate representation of the multilinear map with respect to a choice of bases and the transformation rules become a change-of-basis theorem. By the 1980s, definition~\ref{st:tensor2} had by and large supplanted definition~\ref{st:tensor1} as the standard definition of a tensor in textbooks on mathematical methods for physicists \cite{AMR,Choquet,Hassani,Martin}, geometry \cite{Boothby,Helgason}, general relativity \cite{MTW,Wald} and gauge theory \cite{Bleecker}, although it dates back to at least \citet{Temple}.

Definition~\ref{st:tensor2} appeals to physicists for two reasons. The first is the well-known maxim that the laws of physics do not depend on coordinates, and thus the tensors  used to express these laws should not depend on coordinates either. This has been articulated in various forms in many places, but we quote a  modern version that appears as the very first sentence in Kip Thorne's lecture notes (now a book).

\mypara{Geometric principle}\label{Thorne} The laws of physics must all be expressible as geometric (coordinate-independent and reference-frame-independent) relationships between geometric objects (scalars, vectors, tensors, \dots) that represent physical entities \cite{Thorne}.
\medskip

Most well-known equations in physics express a relation between tensors. For example, Newton's second law $F = ma$ is a relation between two contravariant $1$-tensors -- force $F$ and acceleration $a$ -- and Einstein's field equation $G + \Lambda g =  \kappa T$ is a relation between three covariant $2$-tensors -- the Einstein tensor $G$, the metric tensor $g$ and the energy--momentum tensor $T$, with constants $\Lambda$ and $\kappa$. Since such laws of physics should not depend on coordinates, a tensor ought to be a coordinate-free object, and  definition~\ref{st:tensor2} meets this requirement: whether or not a map $f \colon  \mathbb{V}_1 \times \dots \times \mathbb{V}_d \to \mathbb{W}$ is multilinear does not depend on a choice of bases on $ \mathbb{V}_1, \ldots, \mathbb{V}_d,\mathbb{W} $.

For the same reason, we expect laws of physics expressed in terms of the Hadamard product to be rare if not non-existent, precisely because the Hadamard product is coordinate-dependent, as we saw at the end of Section~\ref{sec:imtransrules}.  Take two linear maps $f ,g \colon  \mathbb{V} \to \mathbb{V}$: it will not be possible to define the Hadamard product of $f$ and $g$ without choosing some basis on $\mathbb{V}$ and the value of the  product will depend on that chosen basis. So one advantage of using a coordinate-free definition of tensors such as definitions~\ref{st:tensor2} or~\ref{st:tensor3} is that it automatically avoids pitfalls such as performing coordinate-dependent operations that are undefined on tensors.

The second reason for favouring definition~\ref{st:tensor2} in physics is that the notion of multilinearity is central to the subject. We highlight another two maxims.

\mypara{Linearity principle} Almost any natural process is linear in small
amounts almost everywhere  \cite[p.~vii]{KM}.

\mypara{Multilinearity principle} If we keep
all but one factor constant, the varying factor obeys the linearity principle.
\medskip

It is probably fair to say that the world around us is understandable largely because of a combination of these two principles. The universal importance of linearity needs no elaboration; in this section we will discuss the importance of multilinearity, particularly in computations.

In mathematics, definition~\ref{st:tensor2} appeals because it is the simplest among the three definitions: a multilinear map is trivial to define for anyone who knows about linear maps. More importantly, definition~\ref{st:tensor2} allows us to define tensors over arbitrary vector spaces.\footnote{Note that we do not mean \emph{abstract} vector spaces. There are many \emph{concrete} vector spaces aside from $\mathbb{R}^n$ and $\mathbb{C}^n$ that are useful in applications and computations. Definitions~\ref{st:tensor2} and~\ref{st:tensor3} allow us to define tensors over these.} It is a misconception that in applications there is no need to discuss more general vector spaces because one can always get by with just two of them, $\mathbb{R}^n$ and $\mathbb{C}^n$. The fact of the matter is that other vector spaces often carry structures that are destroyed when one artificially identifies them with $\mathbb{R}^n$ or $\mathbb{C}^n$. Just drawing on the examples we have mentioned in Section~\ref{sec:trans}, semidefinite programming and equivariant neural networks already indicate why it is not a good idea to identify the vector space of symmetric $n\times n$ matrices with $\mathbb{R}^{n(n+1)/2}$ or an $m$-dimensional subspace of  real-valued functions $f \colon  \mathbb{Z}^2 \to \mathbb{R}$ with $\mathbb{R}^m$. We will elaborate on these later in the context of  definition~\ref{st:tensor2}. Nevertheless, we will also see that definition~\ref{st:tensor2} has one deficiency that will serve as an impetus for definition~\ref{st:tensor3}.

\subsection{Multilinear maps}\label{sec:multmaps}

Again we will begin from tensors of orders zero, one and two, drawing on familiar concepts from linear algebra -- vector spaces, linear maps, dual vector spaces -- before moving into less familiar territory with order-three tensors followed by the most general order-$d$ version.  In the following we let $\mathbb{U}$, $\mathbb{V}$, $\mathbb{W}$ be vector spaces over a field of scalars, assumed to be  $\mathbb{R}$ for simplicity but it may also be replaced by $\mathbb{C}$ or indeed any ring such as $\mathbb{Z}$ or $\mathbb{Z}/n\mathbb{Z}$, and we will see that such cases are also important in computations and applications.

Definition~\ref{st:tensor2} assumes that we are given at least one vector space $\mathbb{V}$.
A tensor of order zero is simply defined to be a \emph{scalar}, \ie\ an element of the field $\mathbb{R}$. A tensor of order one may either be a \emph{vector}, \ie\ an element of the vector space $v \in \mathbb{V}$, or a \emph{covector}, \ie\ an element of a dual vector space $\varphi \in \mathbb{V}^*$. A covector, also called a dual vector or linear functional, is a map $\varphi \colon  \mathbb{V} \to \mathbb{R}$ that satisfies
\begin{equation}\label{eq:lin}
\varphi( \lambda v + \lambda' v') = \lambda \varphi( v) + \lambda' \varphi(v')
\end{equation}
for all $v, v' \in \mathbb{V}$, $\lambda, \lambda' \in \mathbb{R}$. To distinguish between these two types of tensors, a vector is called a tensor of \emph{contravariant} order one and a covector a tensor of \emph{covariant} order one. To see how these relate to definition~\ref{st:tensor1}, we pick a basis $\mathscr{B} = \{v_1,\ldots,v_m\}$ of $\mathbb{V}$. This gives us a dual basis $\mathscr{B}^* = \{ v_1^*,\ldots,v_m^*\}$ on $\mathbb{V}^*$, \tie, $v_i^* \colon  \mathbb{V} \to \mathbb{R}$ is a linear functional satisfying
\[
v_i^*(v_j) = \begin{cases} 1 & i =j, \\ 0 & i \ne j, \end{cases}\quad j =1,\ldots,m.
\]
Any vector $v$ in $\mathbb{V}$ may be uniquely represented by a vector  in $\mathbb{R}^m$,
\begin{equation}\label{eq:vecrep}
\mathbb{V} \ni v =a_1 v_1 +\dots+a_m v_m \quad \longleftrightarrow \quad [v]_{\mathscr{B}} \coloneqq \begin{bmatrix} a_1 \\ \vdots \\ a_m \end{bmatrix} \in \mathbb{R}^m ,
\end{equation}
and any covector $\varphi$ in $\mathbb{V}^*$ may also be uniquely represented by a vector in $\mathbb{R}^m$,
\[
\mathbb{V}^* \ni \varphi = b_1 v_1^* +\dots+b_m v_m^* \quad \longleftrightarrow \quad [\varphi]_{\mathscr{B}^*} \coloneqq \begin{bmatrix} b_1 \\ \vdots \\ b_m \end{bmatrix} \in \mathbb{R}^m.
\]
The vector $ [v]_{\mathscr{B}} \in \mathbb{R}^m$ is the \emph{coordinate representation} of $v$ with respect to the basis $\mathscr{B}$; likewise $[\varphi]_{\mathscr{B}^*} \in \mathbb{R}^m$ is the coordinate representation of $\varphi$ with respect to $\mathscr{B}^*$. Here $ [v]_{\mathscr{B}}$ and $[\varphi]_{\mathscr{B}^*} $ are the multi-indexed (in this case single-indexed) object in definition~\ref{st:tensor1} but, as we have emphasized repeatedly, the multi-indexed object is not the crux of the definition. For example, take $m =3$ and $(3, -1, 2) \in \mathbb{R}^3$:  which vector does it represent? The answer is that it can represent any vector except the zero vector in $\mathbb{V}$ or any covector except the zero covector in $\mathbb{V}^*$, \tie, given any non-zero $v \in \mathbb{V}$, there exists a basis $\{v_1,v_2,v_3\}$ such that $v = 3v_1 -v_2 +2v_3$; ditto for non-zero covectors. Knowing that the multi-indexed object is $(3,-1,2)$ tells us absolutely nothing except that it is non-zero. One needs to know (i) the basis and (ii) the transformation rules, which are given by the change-of-basis theorem, stated below for easy reference.

Recall that if $\mathscr{C}=\{v'_1,\ldots, v'_m\}$ is another basis of $\mathbb{V}$ and $X \in \mathbb{R}^{m \times m}$ is such that
\begin{equation}\label{eq:cob1}
v'_j = \sum_{i=1}^m x_{ij} v_i, \quad j  =1,\ldots,m,
\end{equation}
then we must have $X \in \GL(m)$ and
\begin{equation}\label{eq:trans1tensors}
[v]_{\mathscr{C}} = X^{-1}[v]_{\mathscr{B}}, \quad  [\varphi]_{\mathscr{C}^*} = X^\tp  [\varphi]_{\mathscr{B}^*}
\end{equation}
for any $v \in \mathbb{V}$ and any $\varphi \in \mathbb{V}^*$. More precisely, if
\begin{align*}
a_1v_1 + \dots + a_m v_m &= a_1'v'_1 + \dots + a_m' v'_m, \\*
b_1v_1^* + \dots + b_m v_m^* &= b_1' v_1^{\prime*} + \dots + b_m' v_m^{\prime*},
\end{align*}
then the vectors of coefficients $a,a',b,b'\in \mathbb{R}^m$ must satisfy
\[
a' = X^{-1}a, \quad b' = X^\tp  b.
\]
We recover the transformation rules for contravariant and covariant $1$-tensors, respectively, in the middle column of the table on page~\pageref{transrulesummary1}, which we know are simply the change-of-basis theorems for vectors and covectors in linear algebra \cite{Berberian,FIS}. It is now clear why we have been calling the matrix $X$ a \emph{change-of-basis matrix} when discussing the transformation rules in Section~\ref{sec:trans}.

The most familiar $2$-tensor is a \emph{linear operator},\footnote{Following convention, linear and multilinear maps  are called \emph{functionals} if they are scalar-valued, \ie\ the codomain is $\mathbb{R}$, and \emph{operator} if they are vector-valued, \ie\ the codomain is a vector space of arbitrary dimension.}  \ie\ a map $\Phi \colon  \mathbb{U} \to \mathbb{V}$ that satisfies \eqref{eq:lin}. Let $\dim \mathbb{U} =n$ and $\dim \mathbb{V} = m$. Then, for any  bases $\mathscr{A} = \{u_1,\ldots,u_n\}$ on $\mathbb{U}$ and $\mathscr{B} = \{v_1,\ldots,v_m\}$ on $\mathbb{V}$, the linear operator $\Phi$ has a matrix representation with respect to these bases:
\begin{equation}\label{eq:matreplinop}
[\Phi]_{\mathscr{A},\mathscr{B}} = A \in \mathbb{R}^{m \times n} ,
\end{equation}
where the entries $a_{ij}$ are the coefficients in
\[
\Phi(u_j) = \sum_{i=1}^m a_{ij} v_i, \quad j =1,\ldots,n.
\]
Let  $\mathscr{A}'$ be another basis on $\mathbb{U}$ and $\mathscr{B}'$ another basis on $\mathbb{V}$. Let
\[
  [\Phi]_{\mathscr{A}',\mathscr{B}'} = A' \in \mathbb{R}^{m \times n}.
\]
The change-of-basis theorem for linear operators  tells us that if $X \in \GL(m)$ is a change-of-basis matrix on $\mathbb{V}$, \ie\ defined according to \eqref{eq:cob1}, and likewise $Y \in \GL(n)$ is a change-of-basis matrix on $\mathbb{U}$, then
\[
A' = X^{-1} A Y.
\]
For the special case $\mathbb{U} = \mathbb{V}$ with the same bases $\mathscr{A} =\mathscr{B}$ and $\mathscr{A}'=\mathscr{B}'$, we obtain
\[
A' = X^{-1} AX.
\]
Thus we have recovered the transformation rules for mixed $2$-tensors in the middle column of the tables on pages~\pageref{transrulesummary1} and \pageref{transrulesummary2}.

The next most common $2$-tensor is a \emph{bilinear functional}, \tie, a map $\beta \colon  \mathbb{U} \times \mathbb{V} \to \mathbb{R}$ that satisfies
\begin{equation}\label{eq:bl}
\begin{aligned}
\beta(\lambda u +\lambda' u', v) &=
\lambda\beta(u, v) +\lambda' \beta(u', v),\\
\beta(u,\lambda v +\lambda' v') &=
\lambda\beta(u, v) +\lambda' \beta(u, v')
\end{aligned}
\end{equation}
for all $u, u' \in \mathbb{U}$,  $v, v' \in \mathbb{V}$,  $\lambda, \lambda' \in \mathbb{R}$. In this case the matrix representation of $\beta$ with respect to bases $\mathscr{A} $ and $\mathscr{B}$ is particularly easy to describe: 
\[
[\beta]_{\mathscr{A},\mathscr{B}} = A \in \mathbb{R}^{m \times n}
\]
is simply given by
\[
a_{ij} = \beta(u_i,v_j), \quad i=1,\ldots,m, \; j =1,\ldots,n.
\]
We have the corresponding change-of-basis theorem for bilinear forms: if\label{pg:bilinear}
\[
[\beta]_{\mathscr{A}',\mathscr{B}'} = A' \in \mathbb{R}^{m \times n},
\]
then
\[
A' = X^\tp A Y,
\]
and when $\mathbb{U} = \mathbb{V}$, $\mathscr{A} =\mathscr{B}$, $\mathscr{A}'=\mathscr{B}'$, it reduces to
\[
A' = X^\tp AX.
\]
We have recovered the transformation rules for covariant $2$-tensors in the middle column of the tables on pages~\pageref{transrulesummary1} and \pageref{transrulesummary2}.

Everything up to this point is more or less standard in linear algebra \cite{Berberian,FIS}, summarized here for easy reference and to fix notation. The goal is to show that, at least in these familiar cases, the transformation rules in definition~\ref{st:tensor1} are the change-of-basis theorems of the multilinear maps in definition~\ref{st:tensor2}.

What about contravariant $2$-tensors? A straightforward way is simply to say that they are bilinear functionals on dual spaces $\beta \colon  \mathbb{U}^* \times \mathbb{V}^* \to \mathbb{R}$ and we will recover the transformation rules\label{pg:contravariant}
\[
A' = X^{-1} A Y^{-\tp}\quad \text{and}\quad A' = X^{-1} A X^{-\tp}
\]
via similar change-of-basis arguments. We would like to revive an older alternative, useful in engineering \cite{Chou,Tai} and physics \cite[Part I, Section 1.6]{Morse}, which has become somewhat obscure. Given vector spaces $\mathbb{U}$ and $\mathbb{V}$, we define the \emph{dyadic product} by
\[
(a_1 u_1 + \dots + a_n u_n)(b_1 v_1 + \dots + b_m v_m)
\coloneqq \sum_{i=1}^n \sum_{j=1}^m a_i b_j u_i v_j.
\]
Here the juxtaposed $u_i v_j$ is not given any further interpretation, and is simply taken to be an element (called a \emph{dyad}) in a new vector space of dimension $mn$ called the dyadic product of $\mathbb{U}$ and $\mathbb{V}$. While this appears to be a basis-dependent notion, it is actually not, and dyadics satisfy a contravariant $2$-tensor transformation rule.  This is a precursor to definition~\ref{st:tensor3}; we will see that once we insert a $\otimes$ between the vectors to get $u_i \otimes v_j$, the dyadic product of $\mathbb{U}$ and $\mathbb{V}$ is just the tensor product $\mathbb{U} \otimes \mathbb{V}$ that we will discuss in Section~\ref{sec:tenprod}.

Before moving to higher-order tensors, we highlight another impetus for definition~\ref{st:tensor3}. Note that there are many other maps that are also $2$-tensors according to definition~\ref{st:tensor2}. Aside from $\Phi \colon  \mathbb{U} \to \mathbb{V}$, we have other linear operators
\begin{equation}\label{eq:lin2}
\Phi_1 \colon  \mathbb{U}^* \to \mathbb{V}, \quad 
\Phi_2 \colon  \mathbb{U} \to \mathbb{V}^*, \quad
\Phi_3 \colon  \mathbb{U}^* \to \mathbb{V}^*,
\end{equation}
and, aside from $\beta \colon  \mathbb{U} \times \mathbb{V} \to \mathbb{R}$, other bilinear functionals
\begin{equation}\label{eq:blf2}
\beta_1 \colon  \mathbb{U}^* \times \mathbb{V} \to \mathbb{R}, \quad
\beta_2 \colon  \mathbb{U} \times \mathbb{V}^* \to \mathbb{R}, \quad
\beta_3 \colon  \mathbb{U}^* \times \mathbb{V}^* \to \mathbb{R},
\end{equation}
among other more esoteric maps we would not go into here. We begin to see an inkling of why definition~\ref{st:tensor2}, while certainly the simplest definition of a tensor, may not be the best one. After all, according to definition~\ref{st:tensor1}, there ought to be only three types of $2$-tensors: covariant, contravariant and mixed. If we choose bases on  $\mathbb{U}$ and  $\mathbb{V}$, we will find that a matrix representing the bilinear functional  $\beta \colon  \mathbb{U} \times \mathbb{V} \to \mathbb{R}$ and a matrix representing the linear operator $\Phi_2 \colon  \mathbb{U} \to \mathbb{V}^*$ satisfy exactly the same transformation rules, namely that of a covariant $2$-tensor $A' = X^\tp AY$. One benefit of  definition~\ref{st:tensor1} is that it classifies tensors according to their transformation rules, but this is lost in definition~\ref{st:tensor2}: there are many more types of multilinear maps than there are types of tensors. For example, the linear and bilinear maps above that are constructed out of two vector spaces $\mathbb{U}$ and $\mathbb{V}$ may be neatly categorized into three different types of $2$-tensors:
\begin{alignat*}{5}
&\text{covariant $2$-tensor} &\Phi_2 &\colon \mathbb{U} \to \mathbb{V}^*, &\beta &\colon \mathbb{U} \times \mathbb{V} \to \mathbb{R},\\*
&\text{contravariant 2-tensor}\quad  &\Phi_1 &\colon \mathbb{U}^* \to \mathbb{V}, &\beta_3 &\colon \mathbb{U}^* \times \mathbb{V}^* \to \mathbb{R},\\
&\text{mixed 2-tensor} &\Phi &\colon \mathbb{U} \to \mathbb{V}, &\beta_2 &\colon \mathbb{U} \times \mathbb{V}^* \to \mathbb{R},\\*
&&\Phi_3 &\colon \mathbb{U}^* \to \mathbb{V}^*, \quad&\beta_1 &\colon \mathbb{U}^* \times \mathbb{V} \to \mathbb{R}.
\end{alignat*}
Any two maps in the same category satisfy the same change-of-basis theorem, \tie, the matrices that represent them satisfy the same tensor transformation rule.

The number of possibilities grows as the order of the tensor increases. A $3$-tensor is most commonly a \emph{bilinear operator}, \ie\ a map $\Beta \colon  \mathbb{U} \times \mathbb{V} \to \mathbb{W}$ that satisfies the rules in \eqref{eq:bl} but is vector-valued. The next most common $3$-tensor is a trilinear functional, \ie\ a map $\tau \colon \mathbb{U} \times \mathbb{V} \times \mathbb{W} \to \mathbb{R}$ that satisfies
\begin{align}\label{eq:tl}
\tau(\lambda u +\lambda' u',v, w) &=
\lambda\tau(u, v,w) +\lambda' \tau(u', v, w), \notag \\*
\tau(u,\lambda v +\lambda' v', w) &=
\lambda\tau(u,v, w) +\lambda' \tau(u,v', w),\\*
\tau(u,v,\lambda w +\lambda' w') &=
\lambda\tau(u,v, w) +\lambda' \tau(u,v, w') \notag
\end{align}
for all $u, u' \in \mathbb{U}$,  $v, v' \in \mathbb{V}$,  $w, w' \in \mathbb{W}$, $\lambda, \lambda' \in \mathbb{R}$. However, there are also many other possibilities, including bilinear operators
\begin{equation}\label{eq:blo1}
\Beta \colon  \mathbb{U}^* \times \mathbb{V} \to \mathbb{W}, \;
\Beta \colon  \mathbb{U} \times \mathbb{V}^* \to \mathbb{W}, \ldots,
\Beta \colon  \mathbb{U}^* \times \mathbb{V}^* \to \mathbb{W}^*
\end{equation}
and trilinear functionals
\begin{equation}\label{eq:tlf1}
\tau \colon \mathbb{U}^* \times \mathbb{V} \times \mathbb{W} \to \mathbb{R}, \;
\tau \colon \mathbb{U} \times \mathbb{V}^* \times \mathbb{W} \to \mathbb{R}, \ldots,
\tau \colon \mathbb{U}^* \times \mathbb{V}^* \times \mathbb{W}^* \to \mathbb{R},
\end{equation}
among yet more complicated maps. For example, if $\Lin(\mathbb{U}; \mathbb{V})$ denotes the vector space of all linear maps from $\mathbb{U}$ to $\mathbb{V}$, then operator-valued linear operators or linear operators on operator spaces such as
\begin{equation}\label{eq:lo2}
 \Phi_1\colon  \mathbb{U} \to \Lin(\mathbb{V}; \mathbb{W}), \quad
\Phi_2 \colon  \Lin(\mathbb{U}; \mathbb{V}) \to \mathbb{W}
\end{equation}
or bilinear functionals such as
\begin{equation}\label{eq:blf3}
\beta_1 \colon  \mathbb{U} \times \Lin(\mathbb{V}; \mathbb{W}) \to \mathbb{R}, \quad
\beta_2 \colon  \Lin(\mathbb{U}; \mathbb{V}) \times \mathbb{W} \to \mathbb{R}
\end{equation}
are also $3$-tensors. We may substitute any of $\mathbb{U}$, $\mathbb{V}$, $\mathbb{W}$ with $\mathbb{U}^*$, $\mathbb{V}^*$, $\mathbb{W}^*$ in \eqref{eq:lo2} and \eqref{eq:blf3} to obtain yet other $3$-tensors.

Organizing this myriad of multilinear maps, which increases exponentially with order, is the main reason why we would ultimately want to adopt definition~\ref{st:tensor3}, so that we may classify the many different types of multilinear maps into a smaller number of tensors of different types.  For the treatments  that rely on definition~\ref{st:tensor2} -- a formal version of which will appear later as Definition~\ref{def:tensor2} --
such as \citet{AMR,Bleecker,Boothby,Choquet,Hassani,Helgason,Martin,MTW} and \citet{Wald}, this myriad of possibilities for a $d$-tensor is avoided by requiring that the codomain of the multilinear maps be $\mathbb{R}$ or $\mathbb{C}$, \tie, only multilinear \emph{functionals} are allowed; we will see why this approach makes for an awkward definition of tensors.

The change-of-basis theorems for order-$3$ tensors are no longer standard in linear algebra, but it is easy to extrapolate from the standard  change-of-basis theorems for vectors, linear functionals, linear operators and bilinear functionals we have discussed above.  Let $\mathbb{U},\mathbb{V},\mathbb{W}$ be vector spaces over $\mathbb{R}$ of finite dimensions $m,n,p$, and let
\begin{equation}\label{eq:basis0}
\mathscr{A} =\{u_1,\ldots,u_m\} , \quad \mathscr{B} =\{v_1,\ldots,v_n\}, \quad \mathscr{C} =\{w_1,\ldots,w_p\} 
\end{equation}
be arbitrary bases. We limit ourselves to bilinear operators $\Beta \colon  \mathbb{U} \times \mathbb{V} \to \mathbb{W}$ and trilinear functionals $\tau \colon  \mathbb{U} \times \mathbb{V} \times \mathbb{W} \to \mathbb{R}$, as other cases are easy to deduce from these.  Take any vectors $u \in \mathbb{U}$, $v \in \mathbb{V}$,  $w \in \mathbb{W}$ and express them in terms of linear combinations of basis vectors
\[
u = \sum_{i=1}^m a_i u_i, \quad v = \sum_{j=1}^n b_j v_j, \quad w = \sum_{k=1}^p c_k w_k.
\]
Since $\tau$ is trilinear, we get
\begin{align*}
\tau(u, v, w) = \sum_{i=1}^m  \sum_{j=1}^n \sum_{k=1}^p a_i b_j c_k \, \tau(u_i, v_j,  w_k),
\end{align*}
\tie, $\tau$ is completely determined by its values on $\mathscr{A} \times \mathscr{B} \times \mathscr{C}$ in \eqref{eq:basis0}. So if we let
\[
a_{ijk} \coloneqq \tau(u_i,v_j,w_k), \quad i=1,\ldots,m, \; j =1,\ldots,n,\; k =1,\ldots,p,
\]
then the hypermatrix representation of $\tau$ with respect to bases $\mathscr{A} $, $\mathscr{B}$ and $\mathscr{C}$ is
\begin{equation}\label{eq:tribasesrep}
[\tau]_{\mathscr{A},\mathscr{B},\mathscr{C}} = A \in \mathbb{R}^{m \times n \times p}.
\end{equation}
Note that this practically mirrors the discussion for the bilinear functional case on page~\pageref{pg:bilinear} and can be easily extended to any multilinear functional $\varphi \colon  \mathbb{V}_1 \times \dots \times \mathbb{V}_d \to \mathbb{R}$ to get a hypermatrix representation $A  \in \mathbb{R}^{n_1 \times \dots \times n_d}$ with respect to any bases $\mathscr{B}_i$ on $\mathbb{V}_i$, $i=1,\ldots,d$.

The argument for the bilinear operator $\Beta$ requires an additional step. By bilinearity, we get
\[
\Beta(u,v) = \sum_{i=1}^m  \sum_{j=1}^n \sum_{k=1}^p a_i b_j c_k \, \Beta(u_i, v_j),
\]
and since $\Beta(u_i, v_j) \in \mathbb{W}$, we may express it as a linear combination
\begin{equation}\label{eq:write1}
\Beta(u_i,v_j) = \sum_{i=1}^p a_{ijk} w_k,  \quad i=1,\ldots,m, \; j =1,\ldots,n,\; k =1,\ldots,p.
\end{equation}
By the fact that $\mathscr{C}$ is a basis, the above equation uniquely defines the values $a_{ijk}$ and the hypermatrix representation of $\Beta$ with respect to bases $\mathscr{A} $, $\mathscr{B}$ and $\mathscr{C}$ is
\begin{equation}\label{eq:bibasesrep}
[\Beta]_{\mathscr{A},\mathscr{B},\mathscr{C}} = A \in \mathbb{R}^{m \times n \times p}.
\end{equation}

The hypermatrices in \eqref{eq:tribasesrep} and \eqref{eq:bibasesrep}, even if they are identical, represent different types of $3$-tensors. Let $A' = [ a_{ijk}' ] \in \mathbb{R}^{m \times n \times p}$ be a hypermatrix representation of $\tau$ with respect to bases $\mathscr{A}',\mathscr{B}',\mathscr{C}'$ on $\mathbb{U}, \mathbb{V}, \mathbb{W}$, \ie\ $a_{ijk}' = \tau(u_i',v_j',w_k')$. Then it is straightforward to deduce that
\[
A' = (X^\tp, Y^\tp, Z^\tp) \cdot A
\]
with change-of-basis matrices $X\in \GL(m)$, $Y\in \GL(n)$, $Z \in \GL(p)$ similarly defined  as in \eqref{eq:cob1}. Hence, by \eqref{eq:co1}, trilinear functionals are covariant $3$-tensors. On the other hand,  had $A' \in \mathbb{R}^{m \times n \times p}$ been a hypermatrix representation of $\Beta$, then
\[
A' = (X^\tp, Y^\tp, Z^{-1}) \cdot A.
\]  
Hence, by \eqref{eq:mixed1}, bilinear operators are mixed $3$-tensors of covariant order $2$ and contravariant order $1$. The extension to other types of $3$-tensors in \eqref{eq:blo1}--\eqref{eq:blf3} and to arbitrary $d$-tensors may be carried out similarly.

For completeness, we state a formal definition of multilinear maps, if only to serve as a glossary of terms and notation and as a pretext for interesting examples.

\begin{definition}[tensors via multilinearity]\label{def:tensor2}
Let $\mathbb{V}_1,\ldots,\mathbb{V}_d$ and $\mathbb{W}$ be real vector spaces. A \emph{multilinear map}, or more precisely a $d$-linear map, is a map $\Phi \colon  \mathbb{V}_1 \times \dots \times \mathbb{V}_{d} \to \mathbb{W}$ that satisfies
\begin{align}\label{eq:multilinear}
\Phi(v,\ldots, \lambda v_k + \lambda' v_k', \ldots, v_d) 
=\lambda \Phi(v_1,\ldots, v_k, \ldots, v_d) + \lambda' \Phi(v_1,\ldots, v_k', \ldots, v_d) 
\end{align}
for all $v_1 \in \mathbb{V}_1, \ldots, v_k, v_k' \in \mathbb{V}_k, \ldots, v_d \in \mathbb{V}_d$, $\lambda, \lambda' \in \mathbb{R}$, and all $k =1,\ldots,  d$. The set of all such maps will be denoted $\Mult^d(\mathbb{V}_1,\ldots,\mathbb{V}_d; \mathbb{W})$.
\end{definition}

Note that $\Mult^d(\mathbb{V}_1,\ldots,\mathbb{V}_d; \mathbb{W})$ is itself a vector space:  a linear combination of two $k$-linear maps is again a $k$-linear map; in particular $\Mult^1(\mathbb{V};\mathbb{W}) = \Lin(\mathbb{V};\mathbb{W})$. If the vector spaces  $\mathbb{V}_1,\ldots,\mathbb{V}_d$ and $\mathbb{W}$ are endowed with norms, then
\begin{equation}\label{eq:specnorm1}
\lVert \Phi \rVert_\sigma \coloneqq \sup_{v_1,\ldots,v_d \ne 0} \dfrac{\lVert \Phi(v_1,\ldots,v_d) \rVert}{\lVert v_1 \rVert \cdots \lVert v_d \rVert} = \sup_{\lVert v_1 \rVert =\dots=\lVert v_d\rVert = 1} \lVert \Phi(v_1,\ldots,v_d) \rVert
\end{equation}
defines a norm on $\Mult^d(\mathbb{V}_1,\ldots,\mathbb{V}_d; \mathbb{W})$ that we will call a \emph{spectral norm}; we have slightly abused notation by using the same $\lVert\,\cdot\,\rVert$ to denote norms on different spaces.

In Section~\ref{sec:ticmult}, we will see various examples of why (higher-order) multilinear maps are useful in computations. In part to provide necessary background for those examples, we will review their role as higher derivatives in multivariate calculus.

\begin{example}[higher-order derivatives]\label{eg:hod}\hspace{-2pt}%
Let $\mathbb{V}$ and $\mathbb{W}$ be norm spaces and let $\Omega \subseteq \mathbb{V}$ be an open subset. For any function $F \colon  \Omega \to \mathbb{W}$, recall that the (total) \emph{derivative} at $v \in \Omega$ is a linear operator $DF(v) \colon  \mathbb{V} \to \mathbb{W}$ that satisfies
\begin{equation}\label{eq:firstD}
\lim_{h \to 0} \dfrac{\lVert F(v + h) - F(v) - [DF(v)](h) \rVert}{\lVert h \rVert} = 0,
\end{equation}
or, if there is no such linear operator, then $F$ is not differentiable at $v$. The definition may be recursively applied to obtain derivatives of arbitrary order, assuming that they exist on $\Omega$: since $DF(v) \in \Lin(\mathbb{V};\mathbb{W})$, we may apply the same definition to $DF \colon  \Omega \to \Lin(\mathbb{V};\mathbb{W})$ and get $D^2F(v) \colon  \mathbb{V} \to \Lin(\mathbb{V};\mathbb{W})$ as $D(DF)$, \tie,
\begin{equation}\label{eq:secondD}
\lim_{h \to 0} \dfrac{\lVert DF(v + h) - DF(v) - [D^2F(v)](h) \rVert}{\lVert h \rVert} = 0.
\end{equation}
Doing this recursively, we obtain
\[
DF(v) \in \Lin(\mathbb{V};\mathbb{W}), \ \ \;
D^2F(v) \in \Lin(\mathbb{V};\Lin(\mathbb{V};\mathbb{W})), \ \ \;
D^3F(v) \in \Lin(\mathbb{V};\Lin(\mathbb{V};\Lin(\mathbb{V};\mathbb{W}))),
\]
and so on. Avoiding such nested spaces of linear maps is a good reason to introduce multilinear maps. Note that we  have
\begin{equation}\label{eq:LinMult}
\Lin(\mathbb{V};\Mult^{d-1}(\mathbb{V},\ldots,\mathbb{V};\mathbb{W})) = \Mult^d(\mathbb{V},\ldots,\mathbb{V}; \mathbb{W})
\end{equation}
since, for a linear map $\Phi$ on $\mathbb{V}$ taking values in $\Mult^{d-1}(\mathbb{V},\ldots,\mathbb{V};\mathbb{W})$, 
\[
[\Phi(h)](h_1,\ldots,h_{d-1})
\]
must be linear in $h$ for any fixed $h_1,\ldots,h_{d-1}$ and $(d-1)$-linear in $h_1,\ldots,h_{d-1}$ for any fixed $h$, \tie, it is $d$-linear in the arguments $h,h_1,\ldots,h_{d-1}$. Thus we obtain the $d$th-order derivative of  $F \colon  \Omega \to \mathbb{W}$ at a point $v \in \Omega$ as a $d$-linear map,
\begin{equation}\label{eq:hod}
D^dF(v) \colon  \mathbb{V} \times \dots \times \mathbb{V}\to \mathbb{W}.
\end{equation}
With some additional arguments, we may show that as long as $D^dF(v)$ exists in an open neighbourhood of $v$, then it must be a \emph{symmetric} multilinear map in the sense that
\begin{equation}\label{eq:symmderivative}
[D^dF(v)](h_1,h_2,\ldots,h_d) = [D^dF(v)](h_{\sigma(1)},h_{\sigma(2)},\ldots,h_{\sigma(d)})
\end{equation}
for any permutation $\sigma \in \mathfrak{S}_d$, the permutation group on $d$ objects; and in addition we have Taylor's theorem,
\begin{align*}
F(v+h) &=F(v)+ [DF(v)](h)+  \dfrac{1}{2}[D^2F(v)](h,h) + \cdots\\*
&\qquad \cdots +\dfrac{1}{d!}[D^d F(v)](h,\ldots,h)+R(h),
\end{align*}
where the remainder term $\lVert R(h) \rVert/\lVert h \rVert^d \to 0$  as $h \to 0$ \cite[Chapter~XIII, Section~6]{Lang}.
Strictly speaking, the `$=$' in \eqref{eq:LinMult} should be `$\cong$', but elements in both spaces satisfy the same tensor transformation rules, so as tensors, the  `$=$' is perfectly justified;  making this statement without reference to the tensor transformation rules is a reason why we would ultimately want to bring in definition~\ref{st:tensor3}.
\end{example}

Following the convention in \citet{AMR,Bleecker,Boothby,Choquet,Hassani,Helgason,Martin,MTW} and \citet{Wald},  a $d$-tensor is then defined as a $d$-linear functional, \ie\ by setting $\mathbb{W} = \mathbb{R}$. The following is a formal version of definition~\ref{st:tensor2}.

\begin{definition}[tensors as multilinear functionals]\label{def:tensor2a}
Let $p \le d$ be non-negative integers and let $\mathbb{V}_1,\ldots,\mathbb{V}_d$ be vector spaces over $\mathbb{R}$. A \emph{tensor} of \emph{contravariant order} $p$ and \emph{covariant order} $d-p$ is a multilinear functional
\begin{equation}\label{eq:mixedmultfunc}
\varphi \colon  \mathbb{V}_1^*  \times \dots \times \mathbb{V}_p^*  \times \mathbb{V}_{p+1} \times \dots \times \mathbb{V}_d\to \mathbb{R}.
\end{equation}
We say that $\varphi$ is a tensor of \emph{type} $(p,d-p)$ and of \emph{order} $d$.
\end{definition}

This is a passable but awkward definition. For $d=2$ and $3$, the bilinear and trilinear functionals in \eqref{eq:blf2} and \eqref{eq:tlf1} are tensors, but the linear and bilinear operators in \eqref{eq:lin2} and \eqref{eq:blo1} strictly speaking are not; sometimes this means one needs to convert operators to functionals, for example by identifying a linear operator $\Phi \colon  \mathbb{V} \to \mathbb{W}$ as the bilinear functional defined by $\mathbb{V} \times \mathbb{W}^* \to \mathbb{R}$, $(v,\varphi) \mapsto \varphi(\Phi(v))$, before one may apply a result. Definition~\ref{def:tensor2a} is also peculiar considering that by far the most common $1$-, $2$- and $3$-tensors are vectors, linear operators and bilinear operators respectively, but the definition excludes them at the outset. So instead of simply speaking of $v \in \mathbb{V}$, one would need to regard it as a linear functional on the space of linear functionals, \ie\ $v^{**} \colon  \mathbb{V}^* \to \mathbb{R}$, in order to regard it as a $1$-tensor in the sense of Definition~\ref{def:tensor2a}. In fact, it is often more useful to do the reverse, \tie, given a $d$-linear functional, we prefer to convert it into a $(d-1)$-linear operator, and we will give an example.

\begin{example}[higher-order gradients]\label{eg:hog}
We consider the higher-order derivatives of  a \emph{real-valued} function $f \colon  \Omega \to \mathbb{R}$, where $\Omega \subseteq \mathbb{V}$ is an open set and where $\mathbb{V}$ is equipped with an inner product $\langle \,\cdot\,,\cdot\,\rangle$. This is a special case of Example~\ref{eg:hod} with $\mathbb{W} = \mathbb{R}$ and $\mathbb{V}$ an inner product space. Since
\[
D^df(v) \colon  \underbracket[0.5pt]{\mathbb{V} \times \dots \times \mathbb{V}}_{d \text{ copies}} \to \mathbb{R}
\]
is a multilinear \emph{functional}, for any fixed first $d-1$ arguments,
\[
[D^df(v)](h_1,\ldots,h_{d-1}, \cdot) \colon  \mathbb{V}\to \mathbb{R}
\]
is a linear functional and, by the Riesz representation theorem, there must be a vector, denoted
$[\nabla^df(v)](h_1,\ldots,h_{d-1}) \in \mathbb{V}$ as it  depends on $h_1,\ldots,h_{d-1}$, such that
\begin{equation}\label{eq:Dnabla}
[D^df(v)](h_1,\ldots,h_{d-1}, h_d) = \langle [\nabla^df(v)](h_1,\ldots,h_{d-1}), h_d \rangle.
\end{equation}
Since $D^df(v)$ is a $d$-linear functional,
\[
\nabla^df(v) \colon  \underbracket[0.5pt]{\mathbb{V} \times \dots \times \mathbb{V}}_{d-1 \text{ copies}} \to \mathbb{V}
\]
is a $(d-1)$-linear operator. The symmetry in \eqref{eq:symmderivative} shows that our argument does not depend on our having fixed the first $d-1$ arguments:  we could have fixed any $d-1$ arguments and still obtained the same $\nabla^df(v) $, which must itself be symmetric in its $d-1$ arguments. We will call $\nabla^df(v) $ the $d$th gradient of $f$; note that it depends on the choice of inner product. For $d=1$, the argument above is essentially how one would define a gradient $\nabla f \colon  \Omega \to \mathbb{V}$ in Riemannian geometry. For $d=2$, we see that $\nabla^2 f = D (\nabla f) $, which is  how a Hessian is defined in optimization \cite{Boyd,Renegar}. In fact, we may show more generally that, for $d \ge 2$,
\begin{equation}\label{eq:nablad}
\nabla^d f = D^{d-1} (\nabla f),
\end{equation}
which may be taken as an alternative definition of higher gradients. The utility of a multilinear map approach is familiar to anyone who has attempted to find the gradient and Hessian of functions such~as
\[
f (X) = \tr(X^{-1})\quad\text{or}\quad g(x,Y)= x^\tp Y^{-1}x,
\]
where $x \in \mathbb{R}^n$ and $X,Y \in \mathbb{S}^n_\pp$ \cite[Examples~3.4 and 3.46]{Boyd}. The  definitions in terms of partial derivatives,
\begin{align}\label{eq:dumb}
\nabla f = \begin{bmatrix} \dfrac{\partial f}{\partial x_1} \\ \vdots \\ \dfrac{\partial f}{\partial x_n}  \end{bmatrix}\!, \qquad \nabla^2f &= 
\begin{bmatrix}
\dfrac{\partial^2 f}{\partial x_1^2} &  \cdots &  \dfrac{\partial^2 f}{\partial x_1\partial x_n}  \\
\vdots  & \ddots & \vdots \\ 
\dfrac{\partial^2 f}{\partial x_n\partial x_1}  & \cdots & \dfrac{\partial f}{\partial x_n^2}
 \end{bmatrix}\!,\notag \\*
\nabla^d f &= \biggl[ \dfrac{\partial^d f}{\partial x_i \partial x_j \cdots \partial x_k} \biggr]_{i,j,\ldots,k=1}^n,
\end{align}
provide little insight when functions are defined on vector spaces other than $\mathbb{R}^n$. On the other hand, using $f (X) = \tr(X^{-1})$ for illustration, with \eqref{eq:firstD} and \eqref{eq:secondD}, we~get
\begin{alignat*}{3}
D f(X) \colon  \mathbb{S}^n &\to \mathbb{R},  &H &\mapsto -\tr( X^{-1}HX^{-1}),\\*
D^2 f(X) \colon  \mathbb{S}^n \times  \mathbb{S}^n &\to \mathbb{R},\quad  &(H_1,H_2) &\mapsto \tr(X^{-1} H_1 X^{-1}H_2 X^{-1}  +  X^{-1} H_2 X^{-1}H_1X^{-1}),
\end{alignat*}
and more generally
\[
\begin{aligned}
&D^d f(X) \colon  \mathbb{S}^n \times  \dots \times  \mathbb{S}^n \to \mathbb{R},  \\
&(H_1,\ldots,H_d) \mapsto 
(-1)^d\tr\Biggl[\sum_{\sigma \in \mathfrak{S}_d}X^{-1} H_{\sigma(1)} X^{-1}H_{\sigma(2)} X^{-1} \cdots X^{-1} H_{\sigma(d)} X^{-1}\Biggr].
\end{aligned} 
\]
By the cyclic invariance of trace,
\begin{align*}
[D f(X)]( H) &=  \tr[(-X^{-2})H],\\*
[D^2 f(X)] (H_1,H_2) &= \tr[(X^{-1} H_1 X^{-2}  + X^{-2} H_1 X^{-1})H_2],\\
[D^3 f(X)] (H_1,H_2,H_3) &=\tr[(-X^{-1} H_1 X^{-1} H_2 X^{-2}  - X^{-1} H_2 X^{-2} H_1 X^{-1} \\*
&\quad-X^{-2} H_1 X^{-1} H_2 X^{-1}-X^{-1} H_2 X^{-1} H_1 X^{-2}\\*
&\quad -X^{-1} H_1 X^{-2} H_2 X^{-1} -X^{-2} H_2 X^{-1} H_1 X^{-1})H_3],
\end{align*}
and so by \eqref{eq:Dnabla},
\begin{align*}
\nabla f(X) &= -X^{-2},\\*
[\nabla^2f(X)]( H )&= X^{-1} H X^{-2}  + X^{-2} H X^{-1},\\
[\nabla^3f(X)]( H_1,H_2 )&=-X^{-1} H_1 X^{-1} H_2 X^{-2}  - X^{-1} H_2 X^{-2} H_1 X^{-1} \\*
&\quad-X^{-2} H_1 X^{-1} H_2 X^{-1}-X^{-1} H_2 X^{-1} H_1 X^{-2}\\*
&\quad -X^{-1} H_1 X^{-2} H_2 X^{-1} -X^{-2} H_2 X^{-1} H_1 X^{-1},
\end{align*}
all without having to calculate a single partial derivative.  Example~\ref{eg:logdet} will provide a similar illustration of how the use of multilinear maps makes  such calculations routine.
\end{example}

Returning to Definition~\ref{def:tensor2a}, suppose we choose any bases $\mathscr{B}_1=\{u_1,\ldots,u_{n_1}\}$, $\mathscr{B}_2=\{u_2,\ldots,u_{n_2}\}\ldots,\mathscr{B}_d=\{w_1,\ldots,w_{n_d}\}$
   for the vector spaces $\mathbb{V}_1,\mathbb{V}_2,\ldots,\mathbb{V}_d$ and let $A \in \mathbb{R}^{n_1 \times n_2 \times \dots \times n_d}$ be given by
\[
a_{i j \cdots k} = \varphi(u_i,v_j,\ldots,w_k).
\]
Then a consideration of change of basis similar to the various cases we discussed for $d=1,2,3$ recovers the covariant $d$-tensor transformation rules \eqref{eq:co1} for
\[
\varphi \colon  \mathbb{V}_1 \times \dots \times \mathbb{V}_d \to \mathbb{R},
\]
the contravariant $d$-tensor transformation rules \eqref{eq:contra1} for
\[
\varphi \colon  \mathbb{V}_1^*  \times \dots \times \mathbb{V}_d^*  \to \mathbb{R},
\]
and the mixed $d$-tensor transformation rules \eqref{eq:mixed1} for the general case in \eqref{eq:mixedmultfunc}. Furthermore,
\[
[\varphi]_{\mathscr{B}_1,\ldots,\mathscr{B}_d} = A \in \mathbb{R}^{n_1 \times \dots \times n_d}
\]
is exactly the hypermatrix representation of $\varphi$ in Definition~\ref{def:tensor1}. The important special case where $\mathbb{V}_1=\dots=\mathbb{V}_d = \mathbb{V}$ and where we pick only a single basis $\mathscr{B}$ gives us the transformation rules in  \eqref{eq:co2},  \eqref{eq:contra2} and~\eqref{eq:mixed2}.

At this juncture it is appropriate to highlight two previously mentioned points (page~\pageref{pg:overrate}) regarding the feasibility and usefulness of representing a multilinear map as a hypermatrix.

\begin{example}[writing down a hypermatrix is \#P-hard]\label{eg:LRcoeff}
As we saw in \eqref{eq:bibasesrep}, given bases $\mathscr{A} $, $\mathscr{B}$ and $\mathscr{C}$, a bilinear operator $\Beta$ may be represented as a hypermatrix $A$. Writing down the entries $a_{ijk}$ as in \eqref{eq:write1} appears to be a straightforward process, but this is an illusion: the task is \#P-hard in general. Let $0 \le d_1 \le d_2 \le \dots \le d_n$ be integers. Define the generalized Vandermonde matrix
\[
V_{(d_1,\ldots,d_n)} (x) \coloneqq
\begin{bmatrix} 
x_1^{d_1} & x_2^{d_1} & \dots & x_n^{d_1}  \\
x_1^{d_2} & x_2^{d_2} & \dots & x_n^{d_2}  \\
\vdots & \vdots & \ddots & \vdots  \\
x_1^{d_{n-1}} & x_2^{d_{n-1}} & \dots & x_n^{d_{n-1}} \\
x_1^{d_n} & x_2^{d_n} & \dots & x_n^{d_n}
\end{bmatrix}\!,
\]
observing in particular that
\[
V_{(0,1,\ldots,n-1)} (x) =
\begin{bmatrix} 
1 & 1 & \dots & 1  \\
x_1 & x_2 & \dots & x_n \\
\vdots & \vdots & \ddots & \vdots  \\
x_1^{n-2} & x_2^{n-2} & \dots & x_n^{n-2} \\
x_1^{n-1} & x_2^{n-1} & \dots & x_n^{n-1}
\end{bmatrix} 
\]
is the usual Vandermonde matrix. Suppose $d_i \ge i$ for each $i=1,\ldots,n$; then it is not hard to show using the well-known formula $\det V_{(0,1,\ldots,n-1)} (x) = \prod_{i < j} (x_i - x_j)$ that $\det V_{(d_1,d_2,\ldots,d_n)} (x) $ is divisible  by $\det V_{(0,1,\ldots,n-1)} (x) $. So, for any integers $0 \le p_1 \le p_2 \le \dots \le p_n$,
\[
s_{(p_1,p_2,\ldots,p_n)} (x) \coloneqq \dfrac{\det V_{(p_1,p_2+1,\ldots,p_n+n-1)} (x) }{\det V_{(0,1,\ldots,n-1)} (x) }
\]
is a multivariate polynomial in the variables $x_1,\ldots,x_n$. These are symmetric polynomials, \ie\ homogeneous polynomials $s$ with
\[
s(x_1,x_2,\ldots,x_n) = s(x_{\sigma(1)},x_{\sigma(2)},\ldots,x_{\sigma(n)})
\]
for any $\sigma \in \mathfrak{S}_n$, the permutation group on $n$ objects. Let $\mathbb{U}$, $\mathbb{V}$, $\mathbb{W}$ be the vector spaces of symmetric polynomials of degrees $d$, $d'$ and $d+d'$ respectively, and let $\Beta \colon  \mathbb{U} \times \mathbb{V} \to \mathbb{W}$ be the bilinear operator given by  polynomial multiplication, \ie\ $\Beta(s(x), t(x)) = s(x)t(x)$ for any symmetric polynomials $s(x)$ of degree $d$ and $t(x)$ of degree $d'$. A well-known basis of the vector space of degree-$d$ symmetric polynomials is the \emph{Schur basis} given by
\[
\{ s_{(p_1,p_2,\ldots,p_n)} (x) \in \mathbb{U} \colon  p_1 \le p_2 \le \dots \le p_n \text{ is an integer partition of } d\}.
\]
Let  $\mathscr{A} $, $\mathscr{B}$, $\mathscr{C}$ be the respective Schur bases for  $\mathbb{U}$, $\mathbb{V}$, $\mathbb{W}$.  In this case the coefficients in \eqref{eq:write1} are called \emph{Littlewood--Richardson coefficients}, and determining their values is a \#P-complete problem\footnote{This means it is as intractable as evaluating the permanent of a matrix whose entries are zeros and ones \cite{Valiant}. \#P-complete problems include all NP-complete problems, for example, deciding whether a graph is $3$-colourable is NP-complete, but computing its chromatic number is \#P-complete.} \cite{Narayanan}. In other words, determining the hypermatrix representation $A$ of the bilinear operator $\Beta$ is \#P-hard. Littlewood--Richardson coefficients are not as esoteric as one might think but have significance in linear algebra and numerical linear algebra alike. Among other things they play a central role in the resolution of Horn's conjecture about the eigenvalues of a sum of Hermitian matrices \cite{Klyachko,Knutson}.
\end{example}

As we mentioned earlier, a side benefit of definition~\ref{st:tensor2} is that it allows us to work with arbitrary vector spaces. In the previous example,  $\mathbb{U}$, $\mathbb{V}$, $\mathbb{W}$ are vector spaces of degree-$d$ symmetric polynomials for various values of $d$; in the next one we will have $\mathbb{U}=\mathbb{V}=\mathbb{W}= \mathbb{S}^n$, the vector space of $n \times n$ symmetric matrices. Aside from showing that hypermatrix representations may be neither feasible nor useful, Examples~\ref{eg:LRcoeff} and \ref{eg:logdet} also show that there are usually good reasons to work intrinsically with whatever vector spaces one is given.

\begin{example}[higher gradients of log determinant]\label{eg:logdet}\hspz%
The key to all interior point methods is  a \emph{barrier function} that traps iterates within the feasible region of a convex program.
In semidefinite programming, the optimal barrier function with respect to iteration complexity is the log barrier function for the cone of positive definite matrices $\mathbb{S}_\pp^n$,
\[
f \colon  \mathbb{S}_\pp^n \to \mathbb{R}, \quad f(X) = -\log \det X.
\]
Using the characterization in Example~\ref{eg:hog} with the inner product $\mathbb{S}^n \times \mathbb{S}^n \to \mathbb{R}$, $(H_1,H_2) \mapsto \tr(H_1^\tp H_2)$, we may show that its gradient is given by
\begin{equation}\label{eq:grad}
\nabla f \colon  \mathbb{S}_\pp^n \to \mathbb{S}^n, \quad \nabla f(X) = -X^{-1},
\end{equation}
and its Hessian at any $X \in \mathbb{S}_\pp^n$  is the linear map
\begin{equation}\label{eq:Hess}
\nabla^2f(X) \colon   \mathbb{S}^n \to  \mathbb{S}^n, \quad H \mapsto X^{-1} H X^{-1},
\end{equation}
expressions we may also find in \citet{Boyd} and \citet{Renegar}. While we may choose bases on $\mathbb{S}^n$ and artificially write the gradient and Hessian of $f$ in the forms \eqref{eq:dumb}, interested readers may check that they are a horrid mess that obliterates all insights and advantages proffered by \eqref{eq:grad} and \eqref{eq:Hess}. Among other things, \eqref{eq:grad} and \eqref{eq:Hess} allow one to exploit specialized algorithms for matrix product and inversion.

The third-order gradient $\nabla^3 f$ also plays an important role as we need it to ascertain \emph{self-concordance}  in Example~\ref{eg:self}. By our discussion in  Example~\ref{eg:hog}, for any $X \in \mathbb{S}_\pp$, this is a bilinear operator
\[
\nabla^3 f(X) \colon  \mathbb{S}^n \times \mathbb{S}^n \to \mathbb{S}^n,
\]
and by \eqref{eq:nablad}, we may differentiate \eqref{eq:Hess} to get
\begin{equation}\label{eq:3rd}
[\nabla^3 f(X)](H_1, H_2) = -X^{-1} H_1 X^{-1} H_2 X^{-1} - X^{-1} H_2 X^{-1} H_1 X^{-1}.
\end{equation}
Repeatedly applying \eqref{eq:nablad} gives the $(d-1)$-linear operator
\begin{align*}
& \nabla^d f(X) \colon  \mathbb{S}^n \times \dots \times \mathbb{S}^n \to \mathbb{S}^n,\\*
 &(H_1, \ldots,H_{d-1})\mapsto (-1)^d \sum_{\sigma \in \mathfrak{S}_{d-1}} X^{-1} H_{\sigma(1)} X^{-1} H_{\sigma(2)} X^{-1} \cdots X^{-1} H_{\sigma(d-1)} X^{-1},
\end{align*}
as the $d$th gradient of $f$. As interested readers may again check for themselves, expressing this as a $d$-dimensional hypermatrix is even less illuminating than expressing \eqref{eq:grad} and \eqref{eq:Hess} as one- and two-dimensional hypermatrices. Multilinear maps are essential for discussing higher derivatives and gradients of multivariate functions.
\end{example}

In Examples~\ref{eg:hog} and \ref{eg:logdet}, the matrices in $\mathbb{S}^n$ are $1$-tensors even though they are doubly indexed objects. This should come as no surprise after Example~\ref{eg:order}: the order of a tensor is not determined by the number of indices. We will see this in Example~\ref{eg:Calderon} again, where we will encounter a hypermatrix with $3n+3$ indices which nonetheless represents a $3$-tensor.

We conclude this section with some infinite-dimensional examples. Much like linear operators, there is not much that one could say about multilinear operators over infinite-dimensional vector spaces that is purely algebraic. The topic becomes much more interesting when one brings in analytic notions by equipping the vector spaces with norms or inner products.

As is well known, it is easy to ascertain the continuity of linear operators between Banach spaces: $\Phi \colon \mathbb{V} \to \mathbb{W}$ is continuous if and only if it is bounded in the sense of $\lVert \Phi(v) \rVert \le c \lVert v \rVert$ for some constant $c > 0$ and for all $v \in \mathbb{V}$, \ie\ if and only if $\Phi \in \Bd(\mathbb{V}, \mathbb{W})$. We have slightly abused notation by not distinguishing the norms on different spaces and will continue to do so below. Almost exactly the same proof extends to multilinear operators on Banach spaces: $\Phi \colon \mathbb{V}_1 \times \dots \times \mathbb{V}_d \to \mathbb{W}$ is continuous if and only if it is bounded in the sense of
\[
\lVert \Phi(v_1,\dots,v_d) \rVert \le c \lVert v_1 \rVert \cdots  \lVert v_d \rVert
\]
for some constant $c > 0$ and for all $v_1 \in \mathbb{V}_1,\dots,v_d \in \mathbb{V}_d$ \cite[Chapter~IV, Section~1]{Lang2}, \ie\ if and only if its spectral norm as defined in \eqref{eq:specnorm1} is finite. So if $ \mathbb{V}_1, \dots , \mathbb{V}_d$ are finite-dimensional, then $\Phi$ is automatically continuous; in this case it does not matter whether $\mathbb{W}$ is finite-dimensional.  We will write $\Bd(\mathbb{V}_1,\dots,\mathbb{V}_d; \mathbb{W})$ for the set of all bounded/continuous multilinear operators. Then $\Bd(\mathbb{V}_1,\dots,\mathbb{V}_d; \mathbb{W})$ is itself a Banach space when equipped with the spectral norm.

\begin{example}[infinite-dimensional bilinear operators]\label{eg:Calderon}
  Possibly the simplest continuous bilinear operator is the infinite-dimensional
  generalization of matrix--vector product $\mathbb{R}^{m \times n} \times \mathbb{R}^n \to \mathbb{R}^m$, $(A,v) \mapsto Av$, but where we instead have
\[
\Mu \colon \Bd(\mathbb{V}, \mathbb{W}) \times \mathbb{V} \to \mathbb{W}, \quad (\Phi, v) \mapsto \Phi(v).
\]
If $\mathbb{V}$ and $\mathbb{W}$ are Banach spaces equipped with norms $\lVert\,\cdot\,\rVert$ and $\lVert\,\cdot\,\rVert'$ respectively, then $\Bd(\mathbb{V}, \mathbb{W})$ is a Banach space equipped with the operator norm
\[
\lVert \Phi \rVert'' = \sup_{v \ne 0} \dfrac{\lVert \Phi(v) \rVert'}{\lVert v \rVert}.
\]
It is straightforward to see that the spectral norm of $\Mu$ is then
\[
\lVert \Mu \rVert_\sigma =\sup_{\Phi \ne 0,\; v \ne 0} \dfrac{\lVert\Mu( \Phi, v) \rVert'}{\lVert \Phi \rVert'' \lVert v \rVert} =  1
\]
and thus it is continuous. 
Next we will look at an actual Banach space of functions. Another quintessential bilinear operator is  the convolution of two functions $f, g \in L^1(\mathbb{R}^n)$, defined by
\[
f \ast g (x) \coloneqq \int_{\mathbb{R}^n} f(x - y) g(y)\,  \dd y.
\]
The result is an $L^1$-function because $\lVert f\ast g\rVert_1 \le \lVert f \rVert_1\lVert g\rVert_1$. Thus we have a well-defined bilinear operator
\[
\Beta_* \colon  L^1(\mathbb{R}^n)\times L^1(\mathbb{R}^n)  \to L^1(\mathbb{R}^n), \quad (f,g) \mapsto f \ast g.
\]
Since $\lVert \Beta_*(f,g) \rVert_1 =  \lVert f\ast g\rVert_1 \le \lVert f \rVert_1 \lVert g\rVert_1$, we have
\[
\lVert \Beta_* \rVert_\sigma = \sup_{f\!,\, g \ne 0} \dfrac{\lVert f\ast g\rVert_1}{\lVert f \rVert_1 \lVert g\rVert_1} \le 1,
\]
and equality is attained by, say, choosing $f, g \in L^1(\mathbb{R}^n)$ to be non-negative and applying Fubini. So $\lVert \Beta_* \rVert_\sigma =1$ and in particular $\Beta_*$ is continuous.

For a continuous linear operator $\Phi \in \Bd(\mathbb{H};\mathbb{H})$ and a continuous bilinear operator $\Beta \in \Bd(\mathbb{H},\mathbb{H};\mathbb{H})$ on a separable Hilbert space $\mathbb{H}$ with a countable orthonormal basis $\mathscr{B} = \{e_i \colon i \in \mathbb{N}\} $, Parseval's identity $f = \sum_{i =1}^\infty \langle f, e_i\rangle e_i$ gives us the following orthogonal expansions:
\begin{align*}
\Phi(f) &= \sum_{i =1}^\infty \sum_{j =1}^\infty  \langle \Phi(e_i), e_j \rangle  \langle f, e_j \rangle e_j, \\*
\Beta(f,g) &= \sum_{i =1}^\infty\sum_{j=1}^\infty\sum_{k=1}^\infty \langle \Beta(e_i,e_j), e_k \rangle  \langle f, e_i \rangle \langle g, e_j \rangle e_k
\end{align*}
for any $f,g \in \mathbb{H}$. Convergence of these infinite series is, as usual, in the norm induced by the inner product, and it follows from convergence that the hypermatrices representing $\Phi$ and $\Beta$ with respect to $\mathscr{B}$ are $l^2$-summable, \tie,
\[
(\langle \Phi(e_i), e_j \rangle )_{i,j=1}^\infty \in l^2(\mathbb{N} \times \mathbb{N}), \quad
( \langle \Beta(e_i,e_j), e_k \rangle  )_{i,j,k=1}^\infty \in l^2(\mathbb{N} \times \mathbb{N} \times \mathbb{N})
\]
(we discuss infinite-dimensional hypermatrices formally in Example~\ref{eg:hyp}).
In fact the converse is also true. If we define a linear operator and a bilinear operator~by
\begin{align*}
\Phi(f) &\coloneqq \sum_{i =1}^\infty \sum_{j =1}^\infty   a_{ij}  \langle f, e_j \rangle e_j, \\*
\Beta(f,g) &\coloneqq \sum_{i =1}^\infty\sum_{j=1}^\infty\sum_{k=1}^\infty b_{ijk}  \langle f, e_i \rangle \langle g, e_j \rangle e_k
\end{align*}
for any $f,g \in \mathbb{H}$ and where
\[
(a_{ij})_{i,j=1}^\infty \in l^2(\mathbb{N} \times \mathbb{N}), \quad
( b_{ijk} )_{i,j,k=1}^\infty \in l^2(\mathbb{N} \times \mathbb{N} \times \mathbb{N}),
\]
then  $\Phi \in \Bd(\mathbb{H};\mathbb{H})$ and $\Beta \in \Bd(\mathbb{H},\mathbb{H};\mathbb{H})$. So the continuity of these operators is determined by the $l^2$-summability of the corresponding hypermatrices. What is perhaps surprising is that this also holds to various extents for Banach spaces with `almost orthogonal bases' and other growth conditions for the coefficient hypermatrices, as we will see below.

Going beyond elementary examples like $\Mu$ or $\Beta_*$ requires somewhat more background. The space of \emph{Schwartz functions} is
\[
S(\mathbb{R}^n) \coloneqq \biggl\{ f \in C^\infty(\mathbb{R}^n) \colon \sup_{x \in \mathbb{R}^n} \biggl\lvert x_1^{p_1} \cdots x_n^{p_n} \frac{\partial^{q_1 + \dots + q_n} f}{\partial x_1^{q_1} \cdots \partial x_n^{q_n}} \biggr\rvert < \infty , \; p_i,q_j \in \mathbb{N} \biggr\}.
\]
The defining condition essentially states that all its derivatives decay rapidly to zero faster than any negative powers. Examples include compactly supported smooth functions and functions of the form $p(x)\exp(-x^\tp A x)$, where $p$ is a polynomial and $A \in \mathbb{S}^n_\pp$. While $S(\mathbb{R}^n)$ is not a Banach space, it is a dense subset of many common Banach spaces such as $L^p(\mathbb{R}^n)$, $p \in [1,\infty)$. Its continuous dual space,
\[
S'(\mathbb{R}^n) \coloneqq \Bd(S(\mathbb{R}^n);\mathbb{R}),
\]
is the space of \emph{tempered distributions}. It follows from Schwartz's kernel theorem (see Example~\ref{eg:distributions}) that any continuous bilinear operator  $\Beta \colon S(\mathbb{R}^n) \times S(\mathbb{R}^n) \to S'(\mathbb{R}^n)$ must take the form
\[
\Beta(f,g)(x) = \int_{\mathbb{R}^n}\int_{\mathbb{R}^n} K(x,y,z)f(y)g(z) \D y \D z
\]
for some tempered distribution $K \in S'(\mathbb{R}^n \times \mathbb{R}^n \times \mathbb{R}^n)$ \cite{Grafakos}. It is a celebrated result of \citet{Frazier} that if $\psi \in S(\mathbb{R}^n)$ has Fourier transform $\widehat{\psi}(\xi)$ vanishing outside the annulus $\pi/4 \le \lvert \xi \rvert \le \pi$ and bounded away from zero on a smaller annulus $\pi/4+\varepsilon \le \lvert \xi \rvert \le \pi -\varepsilon$, then it is an \emph{almost orthogonal wavelet} for $L^p(\mathbb{R}^n)$, $p \in (1,\infty)$. This means that if we set
\[
\mathscr{B}_\psi \coloneqq \{\psi_{k,\nu} \colon  (k,\nu)\in \mathbb{Z} \times \mathbb{Z}^n \}, \quad \psi_{k,\nu}(x) \coloneqq 2^{kn/2}\psi(2^k x - \nu),
\]
then any $f \in L^p(\mathbb{R}^n)$ has an expansion
\[
f = \sum_{(k,\nu) \in \mathbb{Z}^{n+1}} \langle f, \psi_{k,\nu} \rangle \psi_{k,\nu}
\]
that converges in the $L^p$-norm, much like Parseval's identity for a Hilbert space, even though we do not have a Hilbert space and $\mathscr{B}_\psi $ is not an orthonormal basis (in the language of Example~\ref{eg:wavelet}, $\mathscr{B}_\psi$ is a tight wavelet frame with frame constant~$1$). Furthermore, if the matrix\footnote{Note that $A$ is a matrix in the sense of Example~\ref{eg:hyp}. Here the row and column indices take values in $\mathbb{Z}^{n+1}$. Likewise for the hypermatrix $B$ later.}
\[
A \colon \mathbb{Z}^{n+1} \times \mathbb{Z}^{n+1} \to \mathbb{R}
\]
satisfies an `almost diagonal' growth condition that essentially says that $a_{(i,\lambda), (j,\mu)}$ is small whenever $(i,\lambda)$ and $(j,\mu)$ are far apart, then defining
\[
\Phi_A(f) \coloneqq \sum_{(i,\lambda) \in \mathbb{Z}^{n+1}} \sum_{(j,\nu) \in \mathbb{Z}^{n+1}} a_{(i,\lambda), (j,\mu)} \langle f, \psi_{i,\lambda} \rangle \psi_{j,\mu}
\]
gives us a continuous linear operator $\Phi_A \colon L^p(\mathbb{R}^n) \to L^p(\mathbb{R}^n)$, $p \in (1,\infty)$. This extends to bilinear operators. If the $3$-hypermatrix
\[
B \colon \mathbb{Z}^{n+1} \times \mathbb{Z}^{n+1}  \times \mathbb{Z}^{n+1} \to \mathbb{R}
\]
satisfies an analogous `almost diagonal' growth condition with $b_{(i,\lambda), (j,\mu), (k,\nu)}$ small whenever $(i,\lambda)$, $(j,\mu)$, $(k,\nu)$ are far apart from each other, then defining
\[
\Phi_B(f,g) \coloneqq \sum_{(i,\lambda) \in \mathbb{Z}^{n+1}} \sum_{(j,\mu) \in \mathbb{Z}^{n+1}} \sum_{(k,\nu) \in \mathbb{Z}^{n+1}} b_{(i,\lambda), (j,\mu), (k,\nu)} \langle f, \psi_{i,\lambda} \rangle \langle g, \psi_{j,\mu} \rangle\psi_{k,\nu}
\]
gives us a continuous linear operator $\Phi_B \colon L^p(\mathbb{R}^n) \times L^q(\mathbb{R}^n)  \to L^r(\mathbb{R}^n) $ whenever
\[
\frac{1}{p} + \frac{1}{q} = \frac{1}{r}, \quad 1 < p,q,r < \infty.
\]
The linear result is due to \citet{Frazier} but our description is based on \citet[Theorem~A]{Grafakos}, which also contains the bilinear result \cite[Theorem~1]{Grafakos}. We refer the reader to these references for the exact statement of the `almost diagonal' growth conditions.

Aside from convolution, the best-known continuous bilinear operator is probably the bilinear Hilbert transform
\begin{equation}\label{eq:biHT}
\operatorname{H}(f,g)(x) \coloneqq \lim_{\varepsilon \to 0} \int_{\lvert y \rvert > \varepsilon} f(x+y)g(x-y) \, \frac{\dd y}{y},
\end{equation}
being a bilinear extension of the Hilbert transform
\[
\operatorname{H}(f)(x) \coloneqq \lim_{\varepsilon \to 0} \int_{\lvert y \rvert > \varepsilon} f(x-y) \, \frac{\dd y}{y}.
\]
The latter, according to \citet[p.~15]{Krantz} `is, without question, the most important operator in analysis.' While it is a standard result that the Hilbert transform is continuous as a linear operator $\operatorname{H}\colon L^p(\mathbb{R}) \to L^p(\mathbb{R})$, $p \in (1,\infty)$, and we even know the exact value of its operator/spectral norm \cite[Remark~5.1.8]{GrafakosBook},
\[
\lVert \operatorname{H} \rVert_\sigma =
\begin{cases}
\pi\tan\biggl(\dfrac{\pi}{2p}\biggr) & 1 < p \le 2, \\[9pt]
\pi\cot\biggl(\dfrac{\pi}{2p}\biggr) & 2 \le p < \infty,
\end{cases}
\]
the  continuity of its bilinear counterpart had been a long-standing open problem. It was resolved by \citet{Lacey1,Lacey2}, who showed that as a bilinear operator, $\operatorname{H} \colon L^p(\mathbb{R}) \times L^q(\mathbb{R})  \to L^r(\mathbb{R})$ is continuous whenever
\[
\frac{1}{p} + \frac{1}{q} = \frac{1}{r}, \quad 1 < p,q \le \infty, \quad \frac{2}{3} < r < \infty,
\]
\tie, there exists $c > 0$ such that
\[
\lVert \operatorname{H}(f,g) \rVert_r \le c \lVert f \rVert_p \lVert g \rVert_q
\]
for all $f \in L^p(\mathbb{R})$ and $g \in L^q(\mathbb{R})$. The special case $p =q =2$, $r =1$,  open for more than thirty years, was known as the Calder\'on conjecture.

The study of infinite-dimensional multilinear operators along the above lines has become a vast undertaking, sometimes called multilinear harmonic analysis \cite{Muscalu2}, with a multilinear Calder\'{o}n--Zygmund theory \cite{MCZT} and profound connections to wavelets \cite{MeyerCoifman} among its many cornerstones.
\end{example}

\subsection{Decompositions of multilinear maps}\label{sec:structmultmaps}

For simplicity, we will assume a finite-dimensional setting in this section. We will see that every multilinear map can be constructed out of linear functionals and vectors. This is a simple observation but often couched in abstract terms under headings like the `universal factorization property' or  `universal mapping property'. Its simplicity notwithstanding, this observation is a fruitful one that is the basis of the notion of tensor rank.

In the context of definition~\ref{st:tensor2}, $0$- and $1$-tensors are building blocks, and there is little more that we may say about them. The definition really begins at $d=2$. Let $\Phi \colon  \mathbb{U}  \to \mathbb{V}$ be a linear operator and let $\mathscr{B} = \{v_1,\ldots,v_n\}$ be a basis of $\mathbb{V}$. Then, for any $u\in \mathbb{U}$, we have
\begin{equation}\label{eq:ump0}
\Phi(u) = \sum_{j=1}^n a_j v_j
\end{equation}
for some $a_1,\ldots,a_n \in \mathbb{R}$. Clearly, if we change $u$ the coefficients $a_1,\ldots,a_n$ must in general change too, \tie, they are functions of $u$ and we ought to have written
\begin{equation}\label{eq:ump1}
\Phi(u) = \sum_{j=1}^n a_j(u) v_j
\end{equation}
to indicate this. What kind of function is $a_i \colon  \mathbb{U} \to \mathbb{R}$? It is easy to see if we take the linear functional $v_i^* \colon  \mathbb{V} \to \mathbb{R}^*$ from the dual basis $\mathscr{B}^* = \{v_1^*,\ldots,v_n^*\}$  and hit \eqref{eq:ump1} on the left with it; then
\[
v_i^*(\Phi(u)) = \sum_{j=1}^n a_i(u) v_i^*(v_j) = a_i(u)
\]
by linearity of $v_i^*$ and the fact that $v_i^*(v_j) = \delta_{ij}$. Since this holds for all $u \in \mathbb{U}$, we have
\[
v_i^* \circ \Phi = a_i.
\]
So $a_i \colon  \mathbb{U} \to \mathbb{R}$ is a linear functional as $v_i^*$ and $\Phi$ are both linear. Switching back to our usual notation of denoting linear functionals as $\varphi_i$ instead of $a_i$, we see that every linear operator $\Phi \colon  \mathbb{U}  \to \mathbb{V}$ takes the form
\[
\Phi(u) = \sum_{j=1}^n \varphi_j(u) v_j
\]
for some linear functionals $\varphi_1,\ldots,\varphi_n \in \mathbb{V}^*$, \tie, every linear operator is constructed out of linear functionals and vectors.  Now observe that in order to get \eqref{eq:ump1} we really did not need $v_1,\ldots,v_n$  to be a basis of $\mathbb{V}$: we just need $v_1,\ldots,v_r$ to be a basis of $\im(\Phi)$ with $r = \dim \im(\Phi) = \rank(\Phi)$. Since $r$ cannot be any smaller or else \eqref{eq:ump1} would not hold, we must have
\[
\rank(\Phi) =\min\biggl\{r \colon  \Phi(u) = \sum_{i=1}^{r} \varphi_i (u) v_i\biggr\}.
\]

The argument is similar for a bilinear functional $\beta \colon  \mathbb{U} \times \mathbb{V} \to \mathbb{R}$ but also different enough to warrant going over. For any $u\in \mathbb{U}$, $\beta(u, \cdot) \colon  \mathbb{V} \to \mathbb{R}$ is a linear functional, so we must have
\[
\beta(u, \cdot) =   \sum_{i=1}^n a_i v_i^*
\]
as $\mathscr{B}^* = \{v_1^*,\ldots,v_n^*\}$ is a basis for $\mathbb{V}^*$. Since $a_i$ depends on $u$, we write
\[
\beta(u, \cdot) =   \sum_{i=1}^n a_i(u) v_i^*.
\]
Evaluating at $v_j$ gives
\[
\beta(u, v_j) =   \sum_{i=1}^n a_i(u) v_i^*(v_j) = a_j(u)
\]
for any $u \in \mathbb{U}$. So $a_j \colon   \mathbb{U} \to \mathbb{R}$ is a linear functional. Since $v_i^* \colon  \mathbb{V} \to  \mathbb{R}$ is also a linear functional, switching to our usual notation for linear functionals, we conclude that a bilinear functional $\beta \colon  \mathbb{U} \times \mathbb{V} \to \mathbb{R}$ must take the form
\begin{equation}\label{eq:ump1a}
\beta(u,v) = \sum_{i=1}^n \varphi_i(u) \psi_i (v)
\end{equation}
for some $\varphi_1,\ldots,\varphi_n \in \mathbb{U}^*$ and $\psi_1,\ldots,\psi_n \in \mathbb{V}^*$,  \tie, every bilinear functional is constructed out of linear functionals. In addition,
\[
\rank(\beta) =\min\biggl\{r \colon  \beta(u,v) = \sum_{i=1}^r \varphi_i(u) \psi_i (v)\biggr\}
\]
defines a notation of rank for bilinear functionals.

For a bilinear operator $\Beta  \colon  \mathbb{U} \times \mathbb{V} \to \mathbb{W}$ and any $u\in \mathbb{U}$ and $v\in \mathbb{V}$, 
\[
\Beta(u,v) = \sum_{j=1}^p a_j w_j
\]
with $\mathscr{C} = \{w_1,\ldots,w_p\}$ a basis of $\mathbb{W}$, but the difference now is that $a_j $ depends on both $u$ and $v$ and so is a function $a_j \colon  \mathbb{U} \times \mathbb{V} \to \mathbb{R}$, \tie,
\[
\Beta(u,v) = \sum_{j=1}^p a_j(u,v) w_j.
\]
Hitting this equation on the left by $w_i^*$ gives
\[
w_i^*(\Beta(u,v)) = a_i(u,v),
\]
and since composing a bilinear operator and a linear functional $\mathbb{U} \times \mathbb{V} \xrightarrow{\Beta} \mathbb{W} \xrightarrow{w_i^*} \mathbb{R}$ gives a bilinear functional, we conclude that $a_i$ is a bilinear functional. Applying the result in the previous paragraph, $a_i$ must take the form in \eqref{eq:ump1a}, and changing to our notation for linear functionals we get
\begin{equation}\label{eq:ump3}
\Beta(u,v) = \sum_{i=1}^n\sum_{j=1}^p \varphi_i(u) \psi_i (v) w_j = \sum_{k=1}^r \varphi_k(u) \psi_k (v) w_k ,
\end{equation}
where the last step is simply relabelling the indices, noting that both are sums of terms of the form $\varphi(u)\psi(v) w$, with $r=np$. The smallest $r$, \tie,
\begin{equation}\label{eq:trankbilin}
\rank(\Beta) =\min\biggl\{r \colon  \Beta(u,v) = \sum_{i=1}^{r} \varphi_i (u) \psi_i (v) w_i\biggr\} ,
\end{equation}
is called the \emph{tensor rank} of $\Beta$, a notion that may be traced to \citet[equations~2 and $2_a$]{Hitch1} and will play a critical role in the next section.

The same line of argument may be repeated on a trilinear functional $\tau \colon  \mathbb{U} \times \mathbb{V} \times \mathbb{W} \to \mathbb{R} $ to show that they are just sums of products of linear functionals
\begin{equation}\label{eq:triprodlin}
\tau(u,v,w) = \sum_{i=1}^r \varphi_i(u)\psi_i(v)\theta_i(w)
\end{equation}
and with it a corresponding notion of tensor rank. More generally, any $d$-linear map $\Phi \colon  \mathbb{V}_1 \times \mathbb{V}_2 \times \dots  \times \mathbb{V}_d \to \mathbb{W}$ is built up of linear functionals and vectors,
\[
\Phi(v_1,\ldots,v_d) = \sum_{i=1}^r \varphi_i(v_1)\psi_i(v_2)\cdots \theta_i(v_d) w_i,
\]
where the last $w_i$ may be dropped if it is a $d$-linear functional with $\mathbb{W} = \mathbb{R}$. Consequently, we see that the `multilinearness' in any  multilinear map comes from that of
\[
\mathbb{R}^d \ni (x_1,x_2,\ldots,x_d) \mapsto x_1x_2\cdots x_d \in \mathbb{R}.
\]
Take \eqref{eq:triprodlin} as illustration: where does the `trilinearity' of $\tau$ come from? Say we look at the middle argument; then
\begin{align*}
\tau(u,\lambda v + \lambda'v',w) &= \sum_{i=1}^r \varphi_i(u)\psi_i(\lambda v + \lambda'v')\theta_i(w)\\*
&=\sum_{i=1}^r \varphi_i(u)[\lambda\psi_i( v) + \lambda'\psi(v')]\theta_i(w)\\
&=\lambda\biggl[ \sum_{i=1}^r \varphi_i(u)\psi_i(v)\theta_i(w)\biggr] + \lambda'\biggl[ \sum_{i=1}^r \varphi_i(u)\psi_i(v')\theta_i(w)\biggr]\\*
&=\lambda\tau(u,v,w) + \lambda' \tau(u,v',w).
\end{align*}
The reason why it is linear in the middle argument is simply a result of
\[
x(\lambda y + \lambda' y') z=\lambda x yz + \lambda'x y' z,
\]
which is in turn a result of the trilinearity of $(x,y,z) \mapsto xyz$. All `multilinearness' in all multilinear maps arises in this manner. In addition, this `multilinearness' may be `factored out' of a multilinear map, and  once we do that, whatever is left behind is a linear map.
This is essentially the simple idea behind the universal factorization property of tensor products that we will formally state and discuss in Section~\ref{sec:tensor3c}.

\subsection{Tensors in computations II: multilinearity}\label{sec:ticmult}

A basic utility of the multilinear map perspective of tensors is that it allows one to recognize $3$-tensors in situations like the examples~\ref{it:cplx}, \ref{it:mm}, \ref{it:GI} on page~\pageref{it:cplx}, and in turn allows one to apply tensorial notions such as tensor ranks and tensor norms to analyse them. The most common $3$-tensors are bilinear operators and this is the case in computations too, although we will see that trilinear functionals do also appear from time to time. When the bilinear operator is matrix multiplication, study of its tensor rank leads us to the fabled \emph{exponent of matrix multiplication} $\omega$. While most are aware that $\omega$ is important in the complexity of matrix--matrix products and matrix inversion, it actually goes far beyond. The complexity of computing $LU$ decompositions, symmetric eigenvalue decompositions, determinants and characteristic polynomials, bases for null spaces, matrix sparsification  --  note that none of these are bilinear operators  --  all have asymptotic complexity either exactly $\omega$ or bounded by it \cite[Chapter~16]{BCS}.  Determining the exact value of $\omega$ truly deserves the status of a Holy Grail problem in numerical linear algebra. Especially relevant to this article is the fact that $\omega$ is a tensorial notion.

Given three vector spaces $\mathbb{U}$, $\mathbb{V}$, $\mathbb{W}$, how could one construct a bilinear operator
\[
\Beta \colon  \mathbb{U} \times \mathbb{V} \to \mathbb{W}?
\]
Taking a leaf from Section~\ref{sec:structmultmaps}, a simple and natural way is to take a linear functional $\varphi \colon  \mathbb{U} \to \mathbb{R}$, a linear functional $\psi \colon  \mathbb{V} \to \mathbb{R}$, a vector $w\in \mathbb{W}$, and then define
\begin{equation}\label{eq:decompbilin}
\Beta(u,v) = \varphi (u) \psi(v) w
\end{equation}
for any $u \in \mathbb{U}$ and $v \in \mathbb{V}$; it is easy to check that this is bilinear. We will call a non-zero bilinear operator of such a form \emph{rank-one}.

An important point to note is that in evaluating $\Beta(u,v)$ in \eqref{eq:decompbilin}, only multiplication of variables matters, and for this $\Beta(u,v)$ requires exactly one multiplication. Take for example $\mathbb{U} = \mathbb{V} = \mathbb{W} = \mathbb{R}^3$ and $\varphi(u) = u_1 + 2u_2 + 3u_3$, $\psi(v) = 2v_1 + 3v_2 + 4v_3$, $w=(3,4,5)$; then
\[
\Beta(u,v) = \begin{bmatrix} 3(u_1 + 2u_2 + 3u_3)(2v_1 + 3v_2 + 4v_3) \\ 4(u_1 + 2u_2 + 3u_3)(2v_1 + 3v_2 + 4v_3) \\ 5(u_1 + 2u_2 + 3u_3)(2v_1 + 3v_2 + 4v_3) \end{bmatrix}\!.
\]
This appears to require far more than one multiplication, but multiplications such as $2u_2$ or $4v_3$ are all scalar multiplications, \tie, one of the factors is a constant, and these are discounted.

This is the notion of \emph{bilinear complexity} introduced in \citet{Strassen2} and is a very fruitful idea. Instead of trying to design algorithms that minimize all arithmetic operations at once, we focus on arguably the most expensive one, the multiplication of variables, as this is the part that cannot be hardwired or hardcoded. On the other hand, once we have fixed $\varphi$, $\psi$, $w$, the evaluation of these linear functionals can be implemented as specialized algorithms in hardware for maximum efficiency. One example is the discrete Fourier transform,
\begin{equation}\label{eq:dft}
\begin{bmatrix} x_0' \\ x_1' \\ x_2'\\ x_3' \\ \vdots \\ x_{n-1}' \end{bmatrix} =
\dfrac{1}{\sqrt{n}} \begin{bmatrix}
1&1&1&1&\cdots &1 \\
1&\omega&\omega^2&\omega^3&\cdots&\omega^{n-1} \\
1&\omega^2&\omega^4&\omega^6&\cdots&\omega^{2(n-1)}\\ 1&\omega^3&\omega^6&\omega^9&\cdots&\omega^{3(n-1)}\\
\vdots&\vdots&\vdots&\vdots&\ddots&\vdots\\
1&\omega^{n-1}&\omega^{2(n-1)}&\omega^{3(n-1)}&\cdots&\omega^{(n-1)(n-1)}
\end{bmatrix}\begin{bmatrix} x_0 \\ x_1 \\ x_2\\ x_3 \\ \vdots \\ x_{n-1} \end{bmatrix}\!. 
\end{equation}
Fast Fourier transform is a specialized algorithm that gives us the value of $x_k'$ by evaluating the linear functional $\varphi(x) = \sum_{j=0}^{n-1} \omega^{jk} x_j$ but not in the obvious way, and it is often built into software libraries or hardware. It is instructive to ask about the bilinear complexity of \eqref{eq:dft}: the answer is in fact zero, as it involves only scalar multiplications with constants $\omega^{jk}$.

By focusing only on variable multiplications, we may neatly characterize bilinear complexity in terms of its tensor rank in \eqref{eq:trankbilin}; in fact, bilinear complexity \emph{is} tensor rank. Furthermore, upon finding the algorithms that are optimal in terms of variable multiplications  --  these are almost never unique  --  we may then seek among them the ones that are optimal in terms of scalar multiplications (\ie\  sparsest), additions, numerical stability, communication cost, energy cost, {\em etc.} In fact, one may often bound the number of additions and scalar multiplications in terms of the number of variable multiplications. For example, if an algorithm takes $n^p$ variable multiplications, we may often show that it takes at most $cn^p$ additions  and scalar multiplications for some constant $c > 0$ as each variable multiplication in the algorithm is accompanied by at most $c$ other additions and scalar multiplications; as a result the algorithm is still $O(n^p)$ even when we count all arithmetic operations. This will be the case for all problems considered in Examples~\ref{eg:Strassen} and \ref{eg:beyond}, \tie, the operation counts therein include all arithmetic operations. Henceforth, by `multiplication' in the remainder of this section we will always mean `variable multiplication' unless specified otherwise.

In general, a bilinear operator would not be rank-one as in \eqref{eq:decompbilin}. Nevertheless, as we saw in \eqref{eq:ump3}, it will always be decomposable into a sum of rank-one terms
\begin{equation}\label{eq:decompbilin2}
\Beta(u,v) = \sum_{i=1}^{r} \varphi_i (u) \psi_i (v) w_i
\end{equation}
for some $r \in \mathbb{N}$ with the smallest possible $r$ given by the tensor rank of $\Beta$.
Note that any decomposition of the form in \eqref{eq:decompbilin2} gives us an explicit algorithm for computing $\Beta$ with $r$  multiplications, and thus  $\rank(\Beta)$ gives us the bilinear complexity or least number of  multiplications required to compute $\Beta$. This relation between tensor rank and evaluation of bilinear operators first appeared in \citet{Strassen3}.

As numerical computations go, there is no need to compute a quantity exactly. What if we just require the right-hand side (remember this gives an algorithm) of \eqref{eq:decompbilin2} to be an $\varepsilon$-approximation  of the left-hand side? This leads to the notion of \emph{border rank}:
\begin{equation}\label{eq:brankbilin}
\brank(\Beta) =\min\biggl\{r \colon  \Beta(u,v) = \lim_{\varepsilon \to 0^+ } \sum_{i=1}^{r} \varphi_i^\varepsilon (u) \psi_i^\varepsilon (v) w_i^\varepsilon\biggr\}.
\end{equation}
This was  first proposed by \citet{BCLR} and \citet{BLR} and may be regarded as providing an algorithm (remembering that every such decomposition gives an algorithm)
\[
\Beta^\varepsilon(u,v) = \sum_{i=1}^{r} \varphi_i^\varepsilon (u) \psi_i^\varepsilon (v) w_i^\varepsilon
\]
that approximates $\Beta(u,v)$ up to any arbitrary $\varepsilon$-accuracy. While \eqref{eq:brankbilin} relies on a limit, it may also be defined purely algebraically over any ring as in \citet[p.~522]{Knuth2} and \citet[Definition~15.19]{BCS}, which are in turn based on \citet{Bini}. Clearly $\brank(\Beta) \le \rank(\Beta)$, but one may wonder if there are indeed instances where the inequality is strict. There are explicit examples in \citet{BCLR,BLR,BCS} and \citet{Knuth2}, but they do not reveal the familiar underlying phenomenon, namely that of a difference quotient converging to a derivative \cite{DSL}:
\begin{align}\label{eq:DSL}
& \lim_{\varepsilon \to 0^+} \dfrac{(\varphi_1(u) +\varepsilon \varphi_2(u)) (\psi_1(v) +\varepsilon \psi_2(v)) (w_1 +\varepsilon w_2) -  \varphi_1(u)\psi_1(v) w}{\varepsilon} \notag \\*
&\qquad =\varphi_2(u) \psi_2(v)w_1 +\varphi_1(u) \psi_2(v)w_1+\varphi_1(u) \psi_1(v)w_2.
\end{align}
The left-hand side clearly has rank no more than two; one may  show that as long as $\varphi_1, \varphi_2$ are not collinear,
and likewise for $\psi_1,\psi_2$ and $w_1,w_2$, the right-hand side of \eqref{eq:DSL} must have rank three, \tie, it defines a bilinear operator with rank three and border rank two.

Tensor rank and border rank are purely algebraic notions defined over any vector spaces, or even modules, which are generalizations of vector spaces that we will soon introduce. However, if $\mathbb{U}$, $\mathbb{V}$, $\mathbb{W}$ are norm spaces, we may introduce various notions of norms on bilinear operators $\Beta \colon  \mathbb{U} \times \mathbb{V} \to \mathbb{W}$. We will slightly abuse notation by denoting the norms on all three spaces by $\lVert\,\cdot\,\rVert$. Recall that for a linear functional $\psi \colon  \mathbb{V} \to \mathbb{R}$, its \emph{dual norm} is just
\begin{equation}\label{eq:dualnorm}
\lVert \psi \rVert_* \coloneqq \sup_{v \ne 0} \dfrac{\lvert \psi(v) \rvert}{\lVert v \rVert}.
\end{equation}
A continuous analogue of tensor rank in \eqref{eq:trankbilin} may then be defined by
\begin{equation}\label{eq:tennuclear}
\lVert \Beta \rVert_\nu \coloneqq \inf\biggl\{  \sum_{i=1}^r \lVert\varphi_i \rVert_* \lVert\psi_i \rVert_* \lVert w\rVert \colon  \Beta(u,v) = \sum_{i=1}^{r} \varphi_i (u) \psi_i (v) w_i\biggr\}
\end{equation}
and we call this a \emph{tensor nuclear norm}. This defines a norm dual to the spectral norm in \eqref{eq:specnorm1}, which in this case becomes
\begin{equation}\label{eq:tenspectral}
\lVert \Beta \rVert_\sigma = \sup_{u,v \ne 0} \dfrac{\lVert \Beta(u,v) \rVert}{\lVert u \rVert \lVert v \rVert}.
\end{equation}
We will argue later that the tensor nuclear norm, in an appropriate sense, quantifies the optimal numerical stability of computing $\Beta$ just as tensor rank quantifies bilinear complexity.

It will be instructive to begin from some low-dimensional examples where $\mathbb{U}$, $\mathbb{V}$, $\mathbb{W}$ are of dimensions two and three.

\begin{example}[Gauss's algorithm for complex multiplication]\label{eg:Gauss}
As is fairly well known, one may multiply a pair of complex numbers with three instead of the usual four real multiplications:
\begin{align}
(a + bi)(c+ di) &= (ac - bd) + i (bc+ad) \notag\\* 
&= (ac - bd) + i[(a+b)(c+d) -ac -bd].\label{eq:cplxmult}
\end{align}
The latter algorithm is usually attributed to Gauss. Note that while Gauss's algorithm uses three multiplications, it involves five additions/subtractions where  the usual algorithm has only two, but in bilinear complexity we only care about multiplication of variables, which in this case are $a,b,c,d$. Observe that complex multiplication $\Beta_\mathbb{C} \colon  \mathbb{C} \times \mathbb{C} \to \mathbb{C}$, $(z,w) \mapsto zw$ is an $\mathbb{R}$-bilinear map when we identify $\mathbb{C}\cong \mathbb{R}^2$. So we may regard
$\Beta_\mathbb{C} \colon \mathbb{R}^2 \times \mathbb{R}^2 \to \mathbb{R}^2$ and rewrite \eqref{eq:cplxmult} as
\[
\Beta_\mathbb{C} \biggl(\begin{bmatrix} a \\ b \end{bmatrix}\!, \begin{bmatrix} c \\ d \end{bmatrix}\biggr) = 
\begin{bmatrix} ac - bd \\ bc+ad \end{bmatrix}= \begin{bmatrix} ac - bd \\(a+b)(c+d)- ac -bd \end{bmatrix}\!.
\]
The standard basis vectors in $\mathbb{R}^2$,
\[
e_1 = \begin{bmatrix} 1 \\ 0 \end{bmatrix}\!, \quad e_2 = \begin{bmatrix} 0 \\ 1 \end{bmatrix}\!,
\]
correspond to $1, i \in \mathbb{C}$ and they
have dual basis $e_1^*, e_2^* \colon  \mathbb{R}^2 \to \mathbb{R}$ given by
\[
e_1^*\biggl(\begin{bmatrix} a \\ b \end{bmatrix}\biggr) = a, \quad
e_2^*\biggl(\begin{bmatrix} a \\ b \end{bmatrix}\biggr) = b.
\]
Then the standard algorithm is given by the decomposition
\begin{align}
\Beta_\mathbb{C} (z,w)&=[e_1^*(z)e_1^*(w)-e_2^*(z) e_2^*(w)] e_1 
+  [e_1^*(z) e_2^*(w) + e_2^*(z)  e_1^*(w)] e_2 , \label{eq:standarddecomp}
\end{align}
whereas Gauss's algorithm is given by the decomposition
\begin{align}
\Beta_\mathbb{C} (z,w) &= [(e_1^*+e_2^*)(z) (e_1^*+e_2^*)(w) ]e_2 \notag  \\*
&\qquad + [e_1^*(z)e_1^*(w)] (e_1-e_2) - [e_2^*(z) e_2^*(w)] (e_1+e_2). \label{eq:Gaussdecomp}
\end{align}
One may show that $\rank(\Beta_\mathbb{C}) = 3 = \brank(\Beta_\mathbb{C})$,
\tie, Gauss's algorithm has optimal bilinear complexity whether in the exact or approximate sense. While using Gauss's algorithm for actual multiplication of complex \emph{numbers} is pointless overkill, it is actually useful in practice \cite{Higham1} as one may use it for the multiplication of complex \emph{matrices}:
\begin{equation}\label{eq:cplxmatmult}
(A + iB)(C+ iD)= (AC - BD) + i[(A+B)(C+D) -AC -BD]
\end{equation}
for any $A + iB, C + i D \in \mathbb{C}^{n \times n}$ with $A,B,C,D \in \mathbb{R}^{n \times n}$. 
\end{example}

For bilinear maps on two-dimensional vector spaces, Example~\ref{eg:Gauss} is essentially the only example. We might ask for, say, parallel evaluation of the standard inner product and the standard symplectic form on $\mathbb{R}^2$, \tie,
\begin{equation}\label{eq:gomega}
g(x,y) = x_1y_1 + x_2y_2 \quad \text{and}\quad
\omega(x,y) = x_1y_2 - x_2y_1.
\end{equation}
An algorithm similar to Gauss's algorithm would yield the result in three multiplications, and it is optimal going by either rank or border rank. For bilinear maps on three-dimensional vector spaces, a natural example is the bilinear operator $\Beta_\wedge \colon \wedge^2(\mathbb{R}^n) \times \mathbb{R}^n \to \mathbb{R}^n$ given by the $3 \times 3$ skew-symmetric matrix--vector product
\[
\begin{bmatrix}
0 & a & b\\
-a & 0 & c\\
-b & -c & 0
\end{bmatrix}
\begin{bmatrix} x \\ y \\ z \end{bmatrix} = 
\begin{bmatrix} ay+bz \\ -ax + cz \\ -bx -cy \end{bmatrix}\!.
\]
In this case $\rank(\Beta_\wedge) = 5 = \brank(\Beta_\wedge)$; see \citet[Proposition~12]{struct} and \citet[Theorem~1.3]{Makam}. For a truly interesting example, one would have to look at bilinear maps on four-dimensional vector spaces and we shall do so next.

\begin{example}[Strassen's algorithm for matrix multiplication]\label{eg:Strassen}\hspz%
  \citet{Strassen1} discovered an algorithm
   for inverting an $n \times n$ matrix in fewer than $5.64 n^{\log_2 7}$ arithmetic operations (counting both additions and multiplications). This was a huge surprise at that time as there were results proving that one may not do this with fewer multiplications than the $n^3/3$ required by Gaussian elimination \cite{KK}. The issue is that such impossibility results invariably assume that one is limited to \emph{row} and \emph{column} operations. Strassen's algorithm, on the other hand, is based on \emph{block} operations combined with a novel algorithm that multiplies two $2 \times 2$ matrices with seven multiplications:
\begin{align*}
&\begin{bmatrix}
a_1 & a_2 \\
a_3 & a_4
\end{bmatrix}
\begin{bmatrix}
b_1 & b_2 \\
b_3 & b_4
\end{bmatrix}
\\*
&\quad=
\begin{bmatrix}
a_1 b_1 + a_2 b_3 & \beta +\gamma +  (a_1+a_2 -a_3 -a_4)b_4\\
\alpha +\gamma+a_4(b_2+b_3 -b_1 -b_4) &  \alpha + \beta + \gamma
\end{bmatrix} 
\end{align*}
with
\begin{gather*}
\alpha = (a_3-a_1)(b_2 -b_4),\quad
\beta = (a_3+a_4)(b_2- b_1),\\*
\gamma =a_1 b_1 + (a_3+a_4 -a_1)(b_1 - b_2+ b_4).
\end{gather*}
As in the case of Gauss's algorithm, the saving of one multiplication comes at the cost of an increase in the number of additions/subtractions from eight to fifteen. In fact Strassen's original version had eighteen; the version presented here is the well-known but unpublished Winograd variant discussed in \citet[p.~500]{Knuth2} and \citet[equation~23.6]{HighamBook}. Recursively applying this algorithm to $2 \times 2$ block matrices produces an algorithm for multiplying $n \times n$ matrices with $O(n^{\log_2 7})$ multiplications. More generally,  the  bilinear operator defined by
\[
\Mu_{m,n,p} \colon  \mathbb{R}^{m \times n} \times \mathbb{R}^{n \times p} \to \mathbb{R}^{m \times p}, \quad (A,B)\mapsto AB,
\]
is called either the \emph{matrix multiplication tensor} or \emph{Strassen's tensor}. Every decomposition
\[
\Mu_{m,n,p}(A,B) = \sum_{i=1}^r \varphi_i(A)\psi_i(B) W_i
\]
with linear functionals $\varphi_i \colon   \mathbb{R}^{m \times n} \to \mathbb{R}$, $\psi_i \colon  \mathbb{R}^{n \times p} \to \mathbb{R}$ and $W_i \in \mathbb{R}^{m \times p}$ gives us an algorithm for multiplying an $m\times n$ matrix by an $n \times p$ matrix that requires just $r$ multiplications, with the smallest possible $r$ given by $\rank(\Mu_{m,n,p})$. The \emph{exponent of matrix multiplication} is defined to be
\[
\omega \coloneqq \inf\{ p \in\mathbb{R} \colon  \rank(\Mu_{n,n,n})=O(n^p)\}.
\]
Strassen's work showed that $\omega < \log_2 7\approx 2.807\,354\,9$, and this has been improved over the years to $\omega < 2.372\,859\,6$ at the time of writing \cite{Alman}. Recall that any linear functional $\varphi \colon  \mathbb{R}^{m \times n} \to \mathbb{R}$ must take the form $\varphi(A) = \tr(V^\tp A)$ for some matrix $V \in \mathbb{R}^{m \times n}$, a consequence of the Riesz representation theorem for an inner product space. For concreteness, when $m = n =p =2$, Winograd's variant of Strassen's algorithm is given by
\[
\Mu_{2,2,2}(A,B) = \sum_{i=1}^7 \varphi_i(A)\psi_i(B) W_i,
\]
where $\varphi_i(A) = \tr(U_i^\tp A)$ and $\psi_i(B) = \tr(V_i^\tp B)$ with
 \begin{alignat*}{7}
 U_1 & = \begin{bmatrix*}[r] -1 & 0 \\ 1 & 1 \end{bmatrix*}, &
 V_1 & = \begin{bmatrix*}[r] 1 & -1 \\ 0 & 1 \end{bmatrix*}, &
 W_1 & = \begin{bmatrix*}[r] 0 & 1 \\ 1 & 1 \end{bmatrix*}, \\*
 U_2 & = \begin{bmatrix*}[r] 1 & 0 \\ 0 & 0 \end{bmatrix*}, &
 V_2 & = \begin{bmatrix*}[r] 1 & 0 \\ 0 & 0 \end{bmatrix*}, &
 W_2 & = \begin{bmatrix*}[r] 1 & 1 \\ 1 & 1 \end{bmatrix*}, \\
 U_3 & = \begin{bmatrix*}[r] 0 & 1 \\ 0 & 0 \end{bmatrix*}, &
 V_3 & = \begin{bmatrix*}[r] 0 & 0 \\ 1 & 0 \end{bmatrix*}, &
 W_3 & = \begin{bmatrix*}[r] 1 & 0 \\ 0 & 0 \end{bmatrix*}, \\
 U_4 & = \begin{bmatrix*}[r] 1 & 0 \\ -1 & 0 \end{bmatrix*}, &
 V_4 & = \begin{bmatrix*}[r] 0& -1 \\ 0 & 1 \end{bmatrix*}, &
 W_4 & = \begin{bmatrix*}[r] 0 & 0 \\ 1 & 1 \end{bmatrix*}, \\
 U_5 & = \begin{bmatrix*}[r] 0 & 0 \\ 1 & 1 \end{bmatrix*}, &
 V_5 & = \begin{bmatrix*}[r] -1 & 1 \\ 0 & 0 \end{bmatrix*}, &
 W_5 & = \begin{bmatrix*}[r] 0 & 1 \\ 0 & 1 \end{bmatrix*}, \\*
 U_6 & = \begin{bmatrix*}[r] 1 & 1 \\ -1 & -1 \end{bmatrix*},\quad &
 V_6 & = \begin{bmatrix*}[r] 0 & 0 \\ 0 & 1 \end{bmatrix*},  &
 W_6 & = \begin{bmatrix*}[r] 0 & 1 \\ 0 & 0 \end{bmatrix*}, \\*
 U_7 & = \begin{bmatrix*}[r] 0 & 0 \\ 0 & 1 \end{bmatrix*}, &
 V_7 & = \begin{bmatrix*}[r] 1 & -1 \\ -1 & 1 \end{bmatrix*},\quad &
 W_7 & = \begin{bmatrix*}[r] 0 & 0 \\ -1 & 0 \end{bmatrix*}.
 \end{alignat*} 
\end{example}

The name `exponent of matrix multiplication' is in retrospect a misnomer: $\omega$ should rightly be called the exponent of nearly all matrix computations, as we will see next.

\begin{example}[beyond matrix multiplication]\label{eg:beyond}
This example is a summary of the highly informative book by \citet[Chapter~16]{BCS}, with a few more items drawn from \citet{Schonhage}, but stripped of the technical details to make the message clear, namely that $\omega$ pervades numerical linear algebra. Consider the following problems.\footnote{An upper Hessenberg matrix $H$ is one where $h_{ij} =0$ whenever $i > j+1$ and $\nnz(A) \coloneqq \#\{(i,j) \colon  a_{ij} \ne 0\}$ counts the number of non-zero entries.}
\begin{enumerate}[\upshape (i)]
\setlength\itemsep{3pt}
\item\label{it:inv} \emph{Inversion.} Given $A \in \GL(n)$, find $A^{-1} \in \GL(n)$.

\item \emph{Determinant.} Given $A \in \GL(n)$, find $\det(A) \in \mathbb{R}$.

\item \emph{Null basis.} Given $A \in \mathbb{R}^{n \times n}$, find a basis $v_1,\ldots, v_m \in \mathbb{R}^n$ of $\ker(A)$.

\item\label{it:LS} \emph{Linear system.} Given $A \in \GL(n)$ and $b \in \mathbb{R}^n$, find $v \in \mathbb{R}^n$ so that $Av =b$.

\item \emph{LU decomposition.} Given $A \in \mathbb{R}^{m \times n}$ of full rank, find permutation $P$, unit lower triangular  $L \in  \mathbb{R}^{m \times m}$, upper triangular $U \in \mathbb{R}^{m \times n}$ so that $PA = LU$.

\item\label{it:QR} \emph{QR decomposition.} Given $A \in \mathbb{R}^{n \times n}$, find orthogonal $Q \in  \Or(n)$, upper triangular $U \in \mathbb{R}^{n \times n}$ so that $A = QR$.

\item \emph{Eigenvalue decomposition.} Given $A \in \mathbb{S}^n$, find $Q \in \Or(n)$ and a diagonal matrix $\Lambda \in \mathbb{R}^{n \times n}$ so that $A = Q\Lambda Q^\tp$.

\item\label{it:Hess} \emph{Hessenberg decomposition.}  Given $A \in \mathbb{R}^{n \times n}$, find $Q \in \Or(n)$ and an upper Hessenberg matrix $H \in \mathbb{R}^{n \times n}$ so that $A = QH Q^\tp$.

\item \emph{Characteristic polynomial.} Given $A \in \mathbb{R}^{n \times n}$, find $(a_0,\ldots,a_{n-1}) \in \mathbb{R}^n$ so that $\det(xI - A) = x^n + a_{n-1}x^{n-1} + \dots +a_1 x + a_0$.

\item\label{it:Spar} \emph{Sparsification.} Given $A \in \mathbb{R}^{n \times n}$ and $c \in [1,\infty)$, find $X, Y \in \GL(n)$ so that $\nnz(XAY^{-1}) \le cn$.

\end{enumerate}
For each of problems~\ref{it:inv}--\ref{it:Spar}, if there is an algorithm that computes the $n \times n$ matrix product in $O(n^p)$ multiplications, then there is an algorithm that solves that problem in $O(n^p)$ or possibly  $O(n^p \log n)$ arithmetic operations, \ie\ inclusive of additions and scalar multiplications. 
By the definition of $\omega$ in Example~\ref{eg:Strassen}, there is an algorithm for the $n \times n$ matrix product in $O(n^{\omega + \varepsilon})$ multiplications for any $\varepsilon > 0$. So the  `exponents' of problems~\ref{it:inv}--\ref{it:Spar}, which may be properly defined \cite[Definition~16.1]{BCS} even though they are not bilinear operators, are all equal to or bounded by $\omega$. These results are known to be sharp, \tie, one cannot do better than  $\omega$, except for problems~\ref{it:LS}, \ref{it:QR} and  \ref{it:Hess}. In particular, it is still an open problem whether one might solve a non-singular linear system with fewer than  $O(n^\omega)$ arithmetic operations asymptotically \cite{Strassen4}. 

To get an inkling of why these results hold, a key idea, mentioned in Example~\ref{eg:Strassen} and well known to practitioners of numerical linear algebra, is to avoid working with rows and columns and work with blocks instead, making recursive use of formulas for block factorization, inversion and determinant similar to\label{pg:blockmatrixops}
\begin{gather*}
\begin{bmatrix}
A & B\\
C & D
\end{bmatrix}
=
\begin{bmatrix}
I & 0\\
CA^{-1} & I
\end{bmatrix}
\begin{bmatrix}
A & B\\
0 & S
\end{bmatrix}\!,\quad 
\det\begin{bmatrix}
A & B\\
C & D
\end{bmatrix} =\det(A)\det(S),\\
\begin{bmatrix}
A & B\\
C & D
\end{bmatrix}^{-1} = 
\begin{bmatrix}
A^{-1}+A^{-1}BS^{-1}CA^{-1} & -A^{-1}BS^{-1}\\
-S^{-1}CA^{-1} & S^{-1}
\end{bmatrix}\!,
\end{gather*}
where $S=D-CA^{-1}B$ is the Schur complement. For illustration, Strassen's algorithm for matrix inversion \cite{Strassen1} in $O(n^{\omega + \varepsilon})$ complexity computes
\[
\begin{bmatrix}
A & B\\
C & D
\end{bmatrix}^{-1} = 
\begin{bmatrix}
X_1 - X_3 X_6 X_2 & X_3 X_6\\
X_6 X_2 & - X_6
\end{bmatrix}\!,\quad
\begin{aligned}
X_1 &= A^{-1}, & X_2 &= CX_1,\\
X_3 &= X_1 B, &X_4 &= CX_3,\\
X_5 &= X_4 - D, &X_6 &= X_5^{-1},
\end{aligned}
\]
with the inversion in $X_1$ and $X_6$ performed recursively using the  same algorithm \cite[p.~449]{HighamBook}.
\end{example}

We would like to stress that the results in Examples~\ref{eg:Strassen} and \ref{eg:beyond}, while certainly easy to state, took the combined efforts of a multitude of superb scholars over many decades to establish. We refer interested readers to \citet[Sections~15.13 and 16.12]{BCS}, \citet[Section~4.6.4]{Knuth2} and \citet[Chapters~1--5]{Land2} for further information on the work involved.

The examples above are primarily about the tensor rank of a bilinear operator over vector spaces. We next consider three variations on this main theme:
\begin{itemize}
\setlength\itemsep{3pt}
\item nuclear/spectral norms instead of tensor rank,
\item modules instead of vector spaces,
\item trilinear functionals instead of bilinear operators.
\end{itemize}
An algorithm, however fast, is generally not  practical if it lacks numerical stability, as the correctness of its output may no longer be guaranteed in finite precision. Fortunately, Strassen's algorithm in Example~\ref{eg:Strassen} is only slightly less stable than the standard algorithm for matrix multiplication, an unpublished result\footnote{For $C=AB \in \mathbb{R}^{n \times n}$, we have $\lVert C - \widehat{C} \rVert \le c_n  \lVert A \rVert \lVert B \rVert \varepsilon + O(\varepsilon^2)$ with $c_n=n^2$ (standard), $cn^{\log_2(12)}$ (Strassen), $c'n^{\log_2(18)}$ (Winograd) using  $\lVert X \rVert = \max_{i,j=1,\ldots,n} \lvert x_{ij} \rvert$.}  of \citet{Brent} that has been reproduced and extended in \citet[Theorems~23.2 and 23.3]{HighamBook}. 
Numerical stability is, however, a complicated issue that depends on many factors and on hardware architecture. Designing numerically stable algorithms is as much an art as it is a science, but the six \emph{Higham guidelines} for numerical stability \cite[Section~1.18]{HighamBook} capture the most salient aspects. The second guideline, to `minimize the size of intermediate quantities relative to the final solution', together with the fifth are the two most unequivocal ones. We will discuss how such considerations naturally lead us to the definition of tensor nuclear norms in \eqref{eq:tennuclear}.

\begin{example}[nuclear norm and numerical stability]\label{eg:nuclear}
As we have discussed, an algorithm for evaluating a bilinear operator is a decomposition of the form \eqref{eq:decompbilin2},
\begin{equation}\label{eq:decompbilin2a}
\Beta(u,v) = \sum_{i=1}^{r} \varphi_i (u) \psi_i (v) w_i,
\end{equation}
where $r$ may or may not be minimal. An intermediate quantity has a natural interpretation as a rank-one term $\varphi_i (u) \psi_i (v) w_i$ among the summands and its size also has a natural interpretation as its norm. Note that such a rank-one bilinear operator has nuclear and spectral norms equal to each other and to
$\lVert \varphi_i \rVert_*\lVert \psi_i \rVert_*\lVert w_i \rVert$.
Hence the sum
\[
\sum_{i=1}^{r} \lVert \varphi_i \rVert_*\lVert \psi_i \rVert_*\lVert w_i \rVert
\]
measures the total size of the intermediate quantities in the algorithm \eqref{eq:decompbilin2a} and
its minimum value, given by the nuclear norm of $\Beta$ as defined in \eqref{eq:tennuclear}, provides a measure of optimal numerical stability in the sense of Higham's second guideline.
This was first discussed in \citet[Section~3.2]{struct}. Using complex multiplication $\Beta_\mathbb{C}$ in Example~\ref{eg:Gauss} for illustration, one may show that the nuclear norm \cite[Lemma~6.1]{nuclear} is given by
\[
\lVert \Beta_\mathbb{C}\rVert_\nu = 4.
\]
The standard algorithm \eqref{eq:standarddecomp} attains this minimum value:
\begin{align*}
&\lVert e_1^*\rVert_*\lVert e_1^*\rVert_*\lVert  e_1\rVert  + \lVert - e_2^* \rVert_*\lVert  e_2^*\rVert_*\lVert  e_1 \rVert + \lVert e_1^*\rVert_*\lVert  e_2^*\rVert_*\lVert  e_2\rVert \\*
&\qquad +\lVert e_2^* \rVert_*\lVert  e_1^*\rVert_*\lVert  e_2\rVert = 4 = \lVert \Beta_\mathbb{C}\rVert_\nu
\end{align*}
but not  Gauss's algorithm \eqref{eq:Gaussdecomp}, which has total size
\begin{align*}
&\lVert e_1^*+e_2^*\rVert_*\lVert  e_1^*+e_2^*\rVert_*\lVert  e_2 \rVert
+\lVert e_1^*\rVert_*\lVert  e_1^*\rVert_*\lVert  e_1-e_2\rVert\\*
&\qquad + \lVert - e_2^*\rVert_*\lVert  e_2^*\rVert_*\lVert  e_1+e_2\rVert =2(1+ \sqrt{2}) > \lVert \Beta_\mathbb{C}\rVert_\nu.
\end{align*}
So Gauss's algorithm is faster (by bilinear complexity) but less stable (by Higham's second guideline) than the standard algorithm. One might ask if it is possible to have an algorithm for complex multiplication that is optimal in both measures. It turns out that the decomposition\footnote{We have gone to some lengths to avoid the tensor product $\otimes$ in this section, preferring to defer it to Section~\ref{sec:tenprod}. The decompositions  \eqref{eq:standarddecomp}, \eqref{eq:Gaussdecomp}, \eqref{eq:FLdecomp} would be considerably neater if expressed in tensor product form, giving another impetus for definition~\ref{st:tensor3}.}
\begin{align}\label{eq:FLdecomp}
\Beta_\mathbb{C} (z,w) &= \dfrac{4}{3}\biggl(\biggl[\dfrac{\sqrt{3}}{2}e_1^*(z)+\dfrac{1}{2}e_2^*(z)\biggr]
\biggl[\dfrac{\sqrt{3}}{2}e_1^*(w)+\dfrac{1}{2}e_2^*(w)\biggr]\biggl[\dfrac{1}{2}e_1+\dfrac{\sqrt{3}}{2}e_2\biggr]\notag \\*
&\qquad + \biggl[\dfrac{\sqrt{3}}{2}e_1^*(z)-\dfrac{1}{2}e_2^*(z)\biggr]\biggl[\dfrac{\sqrt{3}}{2}e_1^*(w)-\dfrac{1}{2}e_2^*(w)\biggr]\biggl[\dfrac{1}{2}e_1-\dfrac{\sqrt{3}}{2}e_2\biggr] \notag \\*
&\qquad -e_2^*(z)e_2^*(w) e_1\biggr)
\end{align}
gives an algorithm that attains both $\rank(\Beta_\mathbb{C})$ and $\lVert \Beta_\mathbb{C}\rVert_\nu$. In conventional notation, this algorithm multiplies two complex numbers via
\begin{align*}
   &(a + bi)(c+ di) \\*
   &\quad = \dfrac{1}{2} \biggl[
\biggl(a + \dfrac{1}{\sqrt{3}}b\biggr)
\biggl(c + \dfrac{1}{\sqrt{3}}d\biggr)
+\biggl(a - \dfrac{1}{\sqrt{3}}b\biggr)
\biggl(c -  \dfrac{1}{\sqrt{3}}d\biggr)
-\dfrac{8}{3}bd
\biggr]\\*
&\quad\quad\, +
\dfrac{i\sqrt{3}}{2} \biggl[
\biggl(a + \dfrac{1}{\sqrt{3}}b\biggr)
\biggl(c + \dfrac{1}{\sqrt{3}}d\biggr)
-\biggl(a - \dfrac{1}{\sqrt{3}}b\biggr)
\biggl(c - \dfrac{1}{\sqrt{3}}d\biggr)
\biggr].
\end{align*}
We stress that the nuclear norm is solely a measure of stability in the sense of Higham's second guideline. Numerical stability is too complicated an issue to be adequately captured by a single number. For instance, from the perspective of cancellation errors, our algorithm above also suffers from the same issue pointed out in \citet[Section~23.2.4]{HighamBook} for Gauss's algorithm. By choosing $z=w$ and $b =\sqrt{3}/a$, our algorithm  computes
\[
\dfrac{1}{2} \biggl[
\biggl(a + \dfrac{1}{a}\biggr)^2
+\biggl(a - \dfrac{1}{a}\biggr)^2
-\dfrac{8}{a^2}
\biggr] + \dfrac{i\sqrt{3}}{2} \biggl[
\biggl(a + \dfrac{1}{a}\biggr)^2
-\biggl(a - \dfrac{1}{a}\biggr)^2
\biggr]\eqqcolon x + iy.
\]
There will be cancellation error in the computed real part $\widehat{x}$
when $\lvert a \rvert$ is small and likewise in the computed imaginary
part $\widehat{y}$ when $\lvert a \rvert$ is large. Nevertheless, as
discussed in \citet[Section~23.2.4]{HighamBook}, the algorithm is still
stable in the weaker sense of having acceptably small $\lvert x -
\widehat{x} \rvert/\lvert z \rvert$ and $\lvert y - \widehat{y}
\rvert/\lvert z \rvert$ even if $\lvert x - \widehat{x} \rvert/\lvert x
\rvert$ or $\lvert y - \widehat{y} \rvert/\lvert y \rvert$ might be large.
\end{example}

As is the case for Example~\ref{eg:Gauss}, the algorithm for complex multiplication above is useful only when applied to complex matrices. When $A,B,C,D \in \mathbb{R}^{n \times n}$, the algorithms in Examples~\ref{eg:Gauss} and \ref{eg:nuclear} provide substantial savings when used to multiply  $A + iB, C+ iD \in \mathbb{C}^{n \times n}$. This gives a good reason for extending multilinear maps and tensors to  \emph{modules} \cite{Lang}, \ie\ vector spaces whose field of scalars is replaced by other more general rings. Formally, if $R$ is a ring with multiplicative identity $1$ (which we will assume henceforth),  an $R$-module $\mathbb{V}$ is a commutative group under a group operation denoted $+ $ and has a scalar multiplication operation denoted $\cdot$ such that, 
for all $a,b\in R$ and $v,w \in \mathbb{V}$,\label{pg:module}
\begin{enumerate}[\upshape (i)]
\setlength\itemsep{3pt}
\item $a \cdot ( v + w ) = a \cdot v + a \cdot w$,
\item $( a + b ) \cdot v = a \cdot v + b \cdot v$,
\item $( a b ) \cdot v = a \cdot ( b \cdot v )$,
\item $1 \cdot v = v$.
\end{enumerate}
Clearly, when $R = \mathbb{R}$ or $\mathbb{C}$ or any other field, a module just becomes a vector space. What is useful about the notion is that it allows us to include rings of objects we would not normally consider as `scalars'. For example, in \eqref{eq:cplxmult} we regard $\mathbb{C}$ as a two-dimensional vector space over $\mathbb{R}$, but in \eqref{eq:cplxmatmult} we regard $\mathbb{C}^{n \times n}$ as a two-dimensional\footnote{Strictly speaking, the terminology over modules should be \emph{length} instead of dimension.} module over $\mathbb{R}^{n \times n}$. So in the latter the `scalars' are actually matrices, \ie\ $R=\mathbb{R}^{n \times n}$. When we consider block matrix operations on square matrices such as those on page~\pageref{pg:blockmatrixops}, we are implicitly doing linear algebra over the ring $R=\mathbb{R}^{n \times n}$, which is not even commutative.

Many standard notions in linear and multilinear algebra carry through from vector spaces to modules with little or no change. For example, the  multilinear maps and multilinear functionals of Definitions~\ref{def:tensor2} and \ref{def:tensor2a} apply verbatim to modules with the field of scalars $\mathbb{R}$ replaced by any ring $R$. In other words, the notion of tensors in the sense of definition~\ref{st:tensor2} applies equally well over modules. We will discuss three examples of multilinear maps over modules: tensor fields, fast integer multiplication algorithms and cryptographic multilinear maps.

When people speak of tensors in physics or geometry, they often really mean a tensor field. As a result one may find statements of the tensor transformation rules that bear little resemblance to our version.
The next two examples are
intended to describe the analogues of definitions~\ref{st:tensor1} and~\ref{st:tensor2} for a tensor field and show how they fit into the narrative of this article. Also, outside pure mathematics, defining a tensor field is by far the biggest reason for considering multilinear maps over modules.

\begin{example}[modules and tensor fields]\label{eg:tenfield}
A \emph{tensor field} is -- roughly speaking -- a tensor-valued function on a manifold. Let $M$ be a smooth manifold and $C^\infty(M) \coloneqq \{f \colon  M \to \mathbb{R} \colon  f$ is smooth$\}$. Note that $C^\infty(M) $ is a ring, as products and linear combinations of smooth real-valued functions are smooth (henceforth we drop the word `smooth' as all functions and fields in this example are assumed smooth). Thus $R = C^\infty(M) $ may play the role of the ring of scalars. A $0$-tensor field is just a  function in $C^\infty(M) $. A contravariant $1$-tensor field is a vector field, \ie\ a function $f$ whose value at $x \in M$ is a vector in the tangent space at $x$, denoted $\mathbb{T}_x(M)$. A covariant $1$-tensor field is a covector field, \ie\ a function $f$ whose value at $x \in M$ is a vector in the cotangent space at $x$, denoted $\mathbb{T}_x^*(M)$. Here $\mathbb{T}_x(M)$ is a vector space and $\mathbb{T}_x^*(M)$ is its dual vector space, in which case a $d$-tensor field of contravariant order $p$ and covariant order $d-p$ is simply a function $\varphi$ whose value at $x \in M$ is such a $d$-tensor, \tie, by Definition~\ref{def:tensor2a}, a multilinear functional
\[
\varphi(x)  \colon  \underbracket[0.5pt]{\mathbb{T}_x^*(M)  \times \dots \times \mathbb{T}_x^*(M)}_{p \text{ copies}}  \times \underbracket[0.5pt]{\mathbb{T}_x(M) \times \dots \times \mathbb{T}_x(M)}_{d-p \text{ copies}}  \to \mathbb{R}.
\]
This seems pretty straightforward; the catch is that $\varphi$ is not a function in the usual sense, which has a fixed domain and a fixed codomain, but the codomain of $\varphi$ depends on the point $x$ where it is evaluated. So we have only defined a tensor field at a  point $x \in M$ but we still need a way to `glue' all these pointwise definitions together. The customary way to do this is via coordinate charts and transition maps, but an alternative is to simply define the \emph{tangent bundle}
\[
\mathbb{T}(M) \coloneqq \{ (x, v) \colon  x \in M,\; v \in \mathbb{T}_x(M)\}
\]
and \emph{cotangent bundle}
\[
\mathbb{T}^*(M) \coloneqq \{ (x, \varphi) \colon  x \in M,\; \varphi \in \mathbb{T}_x^*(M)\}
\]
and observe that these are $C^\infty(M) $-modules with scalar product given by pointwise multiplication of real-valued functions with vector/covector fields. A $d$-tensor field of  contravariant order $p$ and covariant order $d-p$ is then defined to be the multilinear functional
\[
\varphi \colon  \underbracket[0.5pt]{\mathbb{T}^*(M)  \times \dots \times \mathbb{T}^*(M)}_{p \text{ copies}}  \times \underbracket[0.5pt]{\mathbb{T}(M) \times \dots \times \mathbb{T}(M)}_{d-p \text{ copies}}  \to C^\infty(M).
\]
Note that this is a multilinear functional in the sense of modules:  the `scalars' are drawn from a ring $R=C^\infty(M)$ and the `vector spaces' are the $R$-modules $\mathbb{T}(M)$ and $\mathbb{T}^*(M)$. This is another upside of defining tensors via definition~\ref{st:tensor2}; it may be easily extended to include tensor fields.

What about definition~\ref{st:tensor1}? The first thing to note is that not every result in  linear algebra over vector spaces carries over to modules. An example is the notion of basis. While some modules do have a basis -- for example,  when we speak of $\mathbb{C}^{n \times n}$ as a two-dimensional  $\mathbb{R}^{n \times n}$-module, it is with the basis $\{1,i\}$ in mind -- others such as the $C^\infty(M) $-module of $d$-tensor fields may not have a basis when $d \ge 1$. An explicit and well-known example is given by vector fields on the $2$-sphere $S^2$, which as a $C^\infty(S^2) $-module does not have a basis because of the hairy ball theorem \cite[Proposition~17.7]{Passman}.
Consequently, given a $d$-tensor field on $M$, the idea that one may just choose bases and represent it as $d$-dimensional hypermatrix with entries in $C^\infty(M) $ is false even when $d=1$. There is, however, one important exception, namely when $M = \mathbb{R}^n$, in which case, given a $d$-tensor field $\varphi$, it is possible to write down a $d$-hypermatrix
\[
A(v)= [\varphi_{ij\cdots k}(v)]_{i,j,\ldots,k=1}^n \in R^{n \times n \times \dots \times n} 
\]
with entries in $R = C^\infty(\mathbb{R}^n) $ that holds for all $v \in \mathbb{R}^n$. The analogue of the transformation rules is a bit more complicated since we are now allowed a \emph{non-linear change of coordinates} $v' =F(v)$
as opposed to merely a linear change of basis as in Section~\ref{sec:transdtensor}. Here the change-of-coordinates function $F \colon  N(v) \to \mathbb{R}^n$ is any smooth function defined on a neighbourhood $N(v) \subseteq \mathbb{R}^n$ of $v$ that is \emph{locally invertible}, \tie, the derivative $DF(v)$ as defined in Example~\ref{eg:hod} is an invertible linear map in a neighbourhood of $v$. This is sometimes called a \emph{curvilinear} change of coordinates.
The analogue of tensor transformation rule \eqref{eq:mixed2} for a $d$-tensor field on $\mathbb{R}^n$ is then
\begin{equation}\label{eq:mixed2field}
A(v') = (DF(v)^{-1},\ldots,DF(v)^{-1},DF(v)^\tp,\ldots,DF(v)^\tp) \cdot A(v).
\end{equation}
If we have $v' = Xv$ for some $X\in \GL(n)$, then $DF(v) = X$ and we recover \eqref{eq:mixed2}. More generally, if we write $v=(v_1,\ldots,v_n)$ and $v'=(v_1',\ldots,v_n')$, then $DF(v)$ and $DF(v)^{-1}$ may be written down as Jacobian matrices,
\[
DF(v) = \begin{bmatrix}
    \dfrac{\partial v'_1}{\partial v_1} & \cdots & \dfrac{\partial v'_1}{\partial v_n}\\
    \vdots                             & \ddots & \vdots\\
    \dfrac{\partial v'_n}{\partial v_1} & \cdots & \dfrac{\partial v'_n}{\partial v_n}
\end{bmatrix}\!,\quad
DF(v)^{-1} = \begin{bmatrix}
    \dfrac{\partial v_1}{\partial v'_1} & \cdots & \dfrac{\partial v_1}{\partial v'_n}\\
    \vdots                             & \ddots & \vdots\\
    \dfrac{\partial v_n}{\partial v'_1} & \cdots & \dfrac{\partial v_n}{\partial v'_n}
\end{bmatrix}\!,
\]
where the latter is a consequence of the inverse function theorem: $ DF(v)^{-1} = DF^{-1}(v')$ for $v' = F(v)$.

Take, for instance, Cartesian coordinates  $v=(x,y,z)$ and spherical coordinates $v'=(r,\theta,\varphi)$ on $\mathbb{R}^3$. Then 
\begin{alignat}{3}\label{eq:carsph}
 r &= \sqrt{x^2 + y^2 + z^2}, \quad  & x &= r \sin\theta \, \cos\varphi, \notag \\*
 \theta &= \arctan\bigl(\sqrt{x^2+y^2}/z\bigr), \quad  & y &= r \sin\theta \, \sin\varphi,\\*
 \varphi &= \arctan(y/x), \quad   & z &= r \cos\theta, \notag 
\end{alignat}
and 
\begin{gather*}
\begin{aligned}
DF(v)&= 
\dfrac{\partial(r,\theta,\varphi)}{\partial(x,y,z)}
=\begin{bmatrix}
    \dfrac{x}{r}&\dfrac{y}{r}&\dfrac{z}{r}\\[9pt]
    \dfrac{xz}{r^2\sqrt{x^2+y^2}}&\dfrac{yz}{r^2\sqrt{x^2+y^2}}&\dfrac{-(x^2+y^2)}{r^2\sqrt{x^2+y^2}}\\[14pt]
    \dfrac{-y}{x^2+y^2}&\dfrac{x}{x^2+y^2}&0
   \end{bmatrix}\!,\\[3pt]
DF(v)^{-1}  &= 
\dfrac{\partial(x,y,z)}{\partial(r,\theta,\varphi)}
  =\begin{bmatrix}
     \sin\theta\cos\varphi&r\cos\theta\cos\varphi&-r\sin\theta\sin\varphi\\
     \sin\theta \sin\varphi&r\cos\theta\sin\varphi&r\sin\theta\cos\varphi\\
     \cos\theta&-r\sin\theta&0
   \end{bmatrix}\!.
\end{aligned}
\end{gather*}
Note that either of the last two matrices may be expressed solely in terms of $x,y,z$ or solely in terms of $r,\theta,\varphi$; the forms above are chosen for convenience. The transformation rule \eqref{eq:mixed2field} then allows us to transform a hypermatrix in $x,y,z$ to one in $r,\theta,\varphi$ or {\em vice versa}.
\end{example}

\begin{example}[tensor fields over manifolds]\label{eg:tenfield2} 
  Since a smooth manifold $M$ locally resembles $\mathbb{R}^n$, we may take \eqref{eq:mixed2field} as the \emph{local} transformation rule for a tensor field $\varphi$. It describes how,  for a point $x \in M$, a hypermatrix $A(v(x))$ representing $\varphi$ in one system of local coordinates $v \colon  N(x) \to \mathbb{R}^n$
  on a neighbourhood $N(x) \subseteq M$
  is related to another hypermatrix  $A(v'(x))$  representing $\varphi$ in another system of local coordinates $v' \colon  N'(x)  \to \mathbb{R}^n$. As is the case for the tensor transformation rule \eqref{eq:mixed2}, one may use the tensor field transformation rule \eqref{eq:mixed2field} to ascertain whether or not a physical or mathematical object represented locally by a hypermatrix with entries in $C^\infty(M)$ is a tensor field; this is how one would usually show that connection forms or Christoffel symbols are not tensor fields \cite[pp.~64--65]{Simmonds}.

In geometry, tensor fields play a central role, as  geometric structures are nearly always tensor fields (exceptions such as Finsler metrics tend to be less studied). The most common ones include Riemannian metrics, which are symmetric $2$-tensor fields $g \colon  \mathbb{T}(M) \times\mathbb{T}(M) \to C^\infty(M)$, and symplectic forms,  which are alternating $2$-tensor fields $\omega \colon  \mathbb{T}(M) \times\mathbb{T}(M) \to C^\infty(M)$, toy examples of which we have seen in \eqref{eq:gomega}.  As we saw in the previous example, the change-of-coordinates maps for tensors are invertible linear maps, but for tensor fields they are locally invertible linear maps;
these are called \emph{diffeomorphisms}, \emph{Riemannian isometries}, \emph{symplectomorphisms}, {\em etc.}, depending on what geometric structures they preserve. The analogues of the finite-dimensional matrix groups in  \eqref{eq:groups}  then become infinite-dimensional Lie groups such as $\Diff(M)$, $\Iso(M,g)$, $\Symp(M,\omega)$, {\em etc.} 

Although beyond the scope of our article, there are two tensor fields that are too important to be left completely unmentioned: the Ricci curvature tensor and the Riemann curvature tensor. Without introducing additional materials, a straightforward way to define them is to use a special system of local coordinates called Riemannian \emph{normal coordinates}  \cite[Chapter~5]{Chern}. For a point $x\in M$, an  $n$-dimensional Riemannian manifold, we choose normal coordinates $(x_1,\dots,x_n)  \colon  N(x)  \to \mathbb{R}^n$ in a neighbourhood $N(x)  \subseteq M$ with $x$ at the origin $0$. Then the Riemannian metric $g$ takes values in $\mathbb{S}^n_\pp$ and we have 
\[
r_{ijkl}(x) =  \frac32 \frac{\partial^2 g_{ij}}{\partial x_k \partial x_l} (0) , \quad
\bar{r}_{ij}(x)= \frac{3}{\sqrt{\det(g)}} \frac{\partial^2 \sqrt{\det(g)}}{\partial x_k \partial x_l} (0).
\]
The Riemann and Ricci curvature tensors at $x$, respectively, are the quadrilinear and bilinear functionals
\[
R(x) \colon \mathbb{T}_x(M) \times \mathbb{T}_x(M) \times \mathbb{T}_x(M) \times \mathbb{T}_x(M) \to \mathbb{R}, \quad
\bar{R}(x) \colon \mathbb{T}_x(M) \times \mathbb{T}_x(M) \to \mathbb{R}, 
\]
given by
\[
R(x)(u,v,w,z) = \sum_{i,j,k,l=1}^n r_{ijkl}(x) \, u_i v_j w_k z_l, \quad
\bar{R}(x) (u,v) = \sum_{i=1}^n\sum_{j=1}^n \bar{r}_{ij}(x)\, u_i v_j.
\]
These are tensor fields -- note that the coefficients depend on $x$ -- so $R(x)$ will generally be a different multilinear functional at different points $x\in M$. The Ricci curvature tensor will make a brief appearance in Example~\ref{eg:sepPDE} on separation of variables for PDEs. Riemann curvature, being a $4$-tensor, is difficult to handle as is, but  when $M$ is embedded in Euclidean space, it appears implicitly in  the form of a $2$-tensor called the Weingarten map \cite[Chapter~21]{Condition} or second fundamental form \cite{Niyogi}, whose eigenvalues, called \emph{principal curvatures}, give us condition numbers.
There are other higher-order tensor fields in geometry
\cite{Dodson}
such as the torsion tensor, the  Nijenhuis tensor ($d=3$) and the Weyl curvature tensor ($d=4$), all of which are unfortunately beyond our scope.

In physics, it is probably fair to say that (i)~most physical fields are tensor fields and (ii)~most tensors are tensor fields. For (i), while there are important exceptions such as spinor fields, the most common fields such as temperature, pressure and Higgs fields are scalar fields;  electric, magnetic and flow velocity fields are vector fields; the Cauchy stress tensor, Einstein tensor and Faraday electromagnetic tensor are $2$-tensor fields; higher-order tensor fields are rarer in physics but there are also examples
such as the Cotton tensor \cite{Cotton} and the Lanczos tensor \cite{Lanczos}, both $3$-tensor fields. The last five named tensors also make the case for (ii): a `tensor' in physics almost always means a tensor field of order two or (more rarely) higher. We will describe the Cauchy stress tensor and mention a few higher-order tensors related to it in Example~\ref{eg:stress}. 
\end{example}

The above examples are analytic in nature but the next two will be algebraic. They show why it is useful to consider bilinear and more generally multilinear maps for modules over $\mathbb{Z}/m\mathbb{Z}$, the ring of integers modulo $m$.

\begin{example}[modules and integer multiplication]\label{eg:int}\hspace{-3.2pt}%
Trivially, integer multi\-pli\-cation $\Beta \colon  \mathbb{Z} \times \mathbb{Z} \to \mathbb{Z}$, $(a,b)\mapsto ab$ is a bilinear map over the $\mathbb{Z}$-module $\mathbb{Z}$, but this is not the relevant module structure that one exploits in fast integer multiplication algorithms. Instead they are based primarily on two key ideas. The first idea is that integers (assumed unsigned) may be represented as polynomials,
\[
a = \sum_{i=0}^{p-1} a_i \theta^i \eqqcolon a(\theta) \quad\text{and}\quad b = \sum_{j=0}^{p-1} b_j \theta^j \eqqcolon b(\theta) 
\]
for some number base $\theta$, and the product has coefficients given by convolutions,
\[
ab = \sum_{k=0}^{2p-2} c_k \theta^k \eqqcolon c(\theta)  \quad\text{with}\quad
c_k = \sum_{i=0}^k a_i b_{k-i}.
\]
Let  $n = 2p-1$ and pad the vectors of coefficients with enough zeros so that we may consider $(a_0,\ldots,a_{n-1})$, $(b_0,\ldots,b_{n-1})$, $(c_0,\ldots,c_{n-1}) $ on an equal footing.
The second idea is to use the discrete Fourier transform\footnote{A slight departure from \eqref{eq:dft} is that we dropped the $1/\sqrt{n}$ coefficient from our DFT and instead put a $1/n$ with our inverse DFT to avoid surds.} (DFT) for some root of unity $\omega$ to perform the convolution,
\begin{align*}
\begin{bmatrix} a_0' \\[2pt] a_1' \\[2pt] a_2'\\ \vdots \\ a_{n-1}' \end{bmatrix} 
&=
\begin{bmatrix}
1&1&1&\cdots &1 \\[2pt]
1&\omega&\omega^2&\cdots&\omega^{n-1} \\[2pt]
1&\omega^2&\omega^4&\cdots&\omega^{2(n-1)}\\ 
\vdots&\vdots&\vdots&\ddots&\vdots\\
1&\omega^{n-1}&\omega^{2(n-1)}&\cdots&\omega^{(n-1)(n-1)}
\end{bmatrix}
\begin{bmatrix} a_0 \\[2pt] a_1 \\[2pt] a_2 \\ \vdots \\ a_{n-1} \end{bmatrix}\!,
\\[4pt]
\begin{bmatrix} b_0' \\[2pt] b_1' \\[2pt] b_2'\\ \vdots \\ b_{n-1}' \end{bmatrix}
&=
\begin{bmatrix}
1&1&1&\cdots &1 \\[2pt]
1&\omega&\omega^2&\cdots&\omega^{n-1} \\[2pt]
1&\omega^2&\omega^4&\cdots&\omega^{2(n-1)}\\ 
\vdots&\vdots&\vdots&\ddots&\vdots\\
1&\omega^{n-1}&\omega^{2(n-1)}&\cdots&\omega^{(n-1)(n-1)}
\end{bmatrix}
\begin{bmatrix} b_0 \\[2pt] b_1 \\[2pt] b_2 \\ \vdots \\ b_{n-1} \end{bmatrix}\!,\\
\begin{bmatrix} c_0 \\[2pt] c_1 \\[2pt] c_2 \\ \vdots \\ c_{n-1} \end{bmatrix}
&=
\dfrac{1}{n}
\begin{bmatrix}
1&1&1&\cdots &1 \\[2pt]
1&\omega^{-1}&\omega^{-2}&\cdots&\omega^{-(n-1)} \\[2pt]
1&\omega^{-2}&\omega^{-4}&\cdots&\omega^{-2(n-1)}\\ 
\vdots&\vdots&\vdots&\ddots&\vdots\\
1&\omega^{-(n-1)}&\omega^{-2(n-1)}&\cdots&\omega^{-(n-1)(n-1)}
\end{bmatrix}
\begin{bmatrix} a_0'b_0' \\[2pt] a_1'b_1' \\[2pt] a_2'b_2'\\ \vdots \\ a_{n-1}'b_{n-1}' \end{bmatrix}\!,
\end{align*}
taking advantage of the following well-known property: a Fourier transform $\mathcal{F}$ turns convolution product $\ast$ into pointwise product $\cdot$ and the inverse Fourier transform turns it back, \tie,
\[
a \ast b = \mathcal{F}^{-1}( \mathcal{F}(a) \cdot \mathcal{F}(a) ).
\]
Practical considerations inform the way we choose  $\theta$ and $\omega$. We choose $\theta = 2^s$ so that $a_i, b_i \in \{0,1,2,\ldots,2^s-1\}$ are $s$-bit binary numbers that can be readily handled. We choose $\omega$ to be a $2p$th root of unity in a ring $\mathbb{Z}/m\mathbb{Z}$ with $m > c_k$ for all $k=0,1,\ldots,n-1$. The $2p$th root of unity ensures that the powers $1,\omega,\ldots,\omega^{n-1}$ are distinct and  $m > c_k$ prevents any `carrying' when computing $c_k$. It turns out that such a choice is always possible with, say, $m = 2^{3d}+1$. Practical considerations aside, the key ideas are to convert integer multiplication to a bilinear operator,
\begin{align*}
\Beta_1 \colon  (\mathbb{Z}/2^s\mathbb{Z})[\theta] \times  (\mathbb{Z}/2^s\mathbb{Z})[\theta] &\to  (\mathbb{Z}/2^s\mathbb{Z})[\theta], \\*
 (a(\theta), b(\theta) ) &\mapsto a(\theta) b(\theta),
\end{align*}
followed by a Fourier conversion into another bilinear operator,
\begin{align*}
\Beta_2 \colon  (\mathbb{Z}/m\mathbb{Z})^n \times (\mathbb{Z}/m\mathbb{Z})^n &\to (\mathbb{Z}/m\mathbb{Z})^n, \\
((a_0,\ldots,a_{n-1}), (b_0,\ldots,b_{n-1}) )&\mapsto (a_0b_0,\ldots,a_{n-1}b_{n-1}).
\end{align*}
In the former, the univariate polynomial ring $(\mathbb{Z}/2^s\mathbb{Z})[\theta]$ is a $\mathbb{Z}/2^s\mathbb{Z} $-module with $\theta$ regarded as the indeterminate of the polynomials. In the latter $(\mathbb{Z}/m\mathbb{Z})^n$ is a $\mathbb{Z}/m\mathbb{Z}$-module, in fact it is the direct sum of $n$ copies of $\mathbb{Z}/m\mathbb{Z}$. 
The algorithms of \citet{Karatsuba}, \citet{Cook}, \citet{Toom}, \citet{SS} and \citet{Furer} are all variations of one or both of these two ideas: incorporating a divide-and-conquer strategy, computing the DFT with fast Fourier transform, replacing the DFT over $\mathbb{Z}/m\mathbb{Z}$ with one over $\mathbb{C}$, {\em etc.}  The recent breakthrough of \citet{Harvey} that led to an $O(n \log n)$ algorithm for $n$-bit integer multiplications was achieved with a multidimensional variation: (i)~replacing the univariate polynomials $(\mathbb{Z}/2^s\mathbb{Z})[\theta]$  with an $R$-module of multivariate polynomials $R[\theta_1,\ldots,\theta_d]$ with a slightly more complicated coefficient ring $R$, and (ii)~replacing the one-dimensional DFT with a $d$-dimensional DFT,
\[
a'(\phi_1,\phi_2,\ldots,\phi_d)=\sum_{\theta_1=0}^{n_1}\cdots\sum_{\theta_d=0}^{n_d}\omega_1^{\phi_1 \theta_1}\omega_2^{\phi_2 \theta_2}\cdots \omega_d^{\phi_d \theta_d}
a(\theta_1,\theta_2,\ldots,\theta_d),
\]
and we will discuss some tensorial features of such multidimensional transforms in Example~\ref{eg:DMT}.
Their  resulting algorithm is not only the fastest ever but the fastest possible. As in the case of matrix multiplication, the implications of an  $O(n \log n)$ integer multiplication algorithm stretch far beyond integer arithmetic. Among other things, we may use the algorithm of \citet{Harvey}  to compute to $n$-bit precision $x/y$ or $\sqrt[k]{x}$ in time $O(n \log n)$ and $\rme^x$ or $\pi$ in time $O(n \log^2 n)$, for real inputs $x,y$ \cite{BrentBook}.
\end{example}

While the previous example exploits the computational efficiency of multilinear maps over modules, the next one exploits their computational intractability.

\begin{example}[cryptographic multilinear maps]\label{eg:crypt}
  Suppose Alice and Bob want to generate a (secure) common password they may both use over the (insecure) internet. One way they may do so is to pick a large prime number $p$ and pick a primitive root of unity $g \in (\mathbb{Z}/p\mathbb{Z})^\times$, \tie, $g$ generates the \emph{multiplicative} group of integers modulo $p$ in the sense that every non-zero $x\in \mathbb{Z}/p\mathbb{Z}$ may be expressed as $x =g^a$ (group theory notation) or $x \equiv g^a \pmod p$ (number theory notation) for some $a \in \mathbb{Z}$. Alice will pick a secret $a \in \mathbb{Z}$ and send $g^a$ publicly to Bob; Bob will pick a secret $b \in \mathbb{Z}$ and send $g^b$ publicly to Alice. Alice, knowing the value of $a$, may compute $g^{ab} = (g^b)^a$ from the $g^b$ she received from Bob, and Bob, knowing the value of $b$, may compute $g^{ab} = (g^a)^b$ from the $g^a$ he received from Alice. They now share the secure password $g^{ab}$. This is the renowned Diffie--Hellman key exchange. The security of the version described is based on the intractability of the discrete log problem: determining the value $a = \log_g (g^a)$ from $g^a$ and $p$ is believed to be intractable. Although the problem has a well-known polynomial-time quantum  algorithm \cite{Shor} and has recently been shown to be quasi-polynomial-time in $n$ \cite{Kleinjung} for a finite field $\mathbb{F}_{p^n}$ when $p$ is fixed (note that for $n=1$, $\mathbb{F}_p= \mathbb{Z}/p\mathbb{Z}$)
    the technology required for the former is still in its infancy, whereas the latter does not apply in our case where complexity is measured in terms of $p$ and not $n$ (for us, $n=1$ always but $p \to \infty$).

Now observe that $(\mathbb{Z}/p\mathbb{Z})^\times$ is a commutative group under the group operation of multiplication modulo $p$, and it is a $\mathbb{Z}$-module as we may check that it satisfies the properties on page~\pageref{pg:module}: for any $x,y \in \mathbb{Z}/p\mathbb{Z}$  and $a,b \in \mathbb{Z}$, 
\begin{enumerate}[\upshape (i)]
\setlength\itemsep{3pt}
\item $(xy )^a = x^a y^b$,
\item $x^{( a + b )}  = x^a x^b$,
\item $x^{( a b )}  = ( x^b)^a$,
\item $x^1 = x$.
\end{enumerate}
Furthermore, the Diffie--Hellman key exchange is a $\mathbb{Z}$-bilinear map
\[
\Beta \colon (\mathbb{Z}/p\mathbb{Z})^\times \times (\mathbb{Z}/p\mathbb{Z})^\times \to (\mathbb{Z}/p\mathbb{Z})^\times,\quad
(g^a, g^b)  \mapsto g^{ab}
\]
since, for any $\lambda, \lambda' \in \mathbb{Z}$ and $g^a,g^b \in (\mathbb{Z}/p\mathbb{Z})^\times$,
\[
\Beta(g^{\lambda a + \lambda' a'}, g^b)  = g^{(\lambda a + \lambda' a')b} 
= (g^{ab})^{\lambda} (g^{a'b})^{\lambda'}  = \Beta(g^a,g^b)^{\lambda}\Beta(g^{a'},g^b)^{\lambda'}
\]
and likewise for the second argument. That the notation is written multiplicatively with coefficients appearing in the power is immaterial; if anything it illustrates the power of abstract algebra in recognizing common structures across different scenarios. While one may express everything in additive notation by taking discrete logs whenever necessary, the notation $g^a$ serves as a useful mnemonic: anything appearing in the power is hard to extract, while using additive notation means having to constantly remind ourselves that extracting $a$ from $ag$ and $\lambda a$ is intractable for the former and trivial for the latter.

What if $d+1$ different parties need to establish a common password (say, in a $1000$-participant Zoom session)? In principle one may successively apply the two-party Diffie--Hellman key exchange $d+1$ times with the $d+1$ parties each doing $d+1$ exponentiations, which is obviously expensive. One may reduce the number of exponentiations to $\log_2 (d+1)$ with a more sophisticated protocol, but it was discovered in \citet{Joux} that with some mild assumptions a bilinear map like the one above already allows us to generate a tripartite password in just one round since
\[
\Beta(g^a, g^b)^c = \Beta(g^a, g^c)^b =\Beta(g^b, g^c)^a = \Beta(g,g)^{abc}.
\]
\citet{Boneh} generalized and formalized this idea as a \emph{cryptographic multilinear map} to allow a one-round $(d+1)$-partite password. Let $p$ be a prime and let $G_1,\ldots,G_d, G$ be cyclic groups of $p$ elements, written multiplicatively.\footnote{Of course they are all isomorphic to each other, but in actual cryptographic applications, it matters how the groups are realized: one might be an elliptic curve over a finite field and another a cyclic subgroup of the Jacobian of a hyperelliptic curve. Also, an explicit isomorphism may not be easy to identify or compute in practice.} Then, as in the $d=2$ case, these groups are $\mathbb{Z}$-modules and we may consider a $d$-linear map over $\mathbb{Z}$-modules:
\[
\Phi \colon  G_1 \times \dots \times G_d \to G.
\]
The only property that matters is the following consequence of multilinearity:
\[
\Phi(g_1^{a_1}, g_2^{a_2}, \ldots,g_d^{a_d}) = \Phi(g_1,g_2,\ldots,g_d)^{a_1 a_2 \cdots a_d}
\]
for any $g_1\in G_1,\ldots,g_d\in G_d$ and $a_1,\ldots,a_d \in \mathbb{Z}$. Slightly abusing notation, we will write $1$ for the identity element in all groups.
To exclude trivialities, we assume that $\Phi$ is non-degenerate, \tie, if $\Phi(g_1,\ldots,g_d) = 1$, then we must have $g_i =1$ for some $i$ (the converse is always true). For a non-degenerate multilinear map $\Phi$ to be \emph{cryptographic}, it needs to have two additional properties:
\begin{enumerate}[\upshape (i)]
\setlength\itemsep{3pt}
\item the discrete log problem in each of $G_1,\ldots,G_d$ is intractable,
\item  $\Phi(g_1,\ldots,g_d)$ is efficiently computable for any $g_1\in G_1,\ldots,g_d\in G_d$.
\end{enumerate}
Note that these assumptions imply that the discrete log problem in $G$ must also be intractable.
Given, say, $g_1$ and $g_1^{a_1}$, since $\Phi(g_1, 1, \ldots,1)$ and $\Phi(g_1^{a_1}, 1, \ldots,1)$ may be efficiently computed, if we could efficiently solve the discrete log problem $\Phi(g_1^{a_1}, 1, \ldots,1) = \Phi(g_1,1,\ldots,1)^{a_1}$ in $G$ to get $a_1$, we would have solved the discrete log problem in $G_1$ efficiently.

Suppose we have a cryptographic $d$-linear map, $\Phi \colon G \times \dots \times G \to G$, where we have assumed $G_1 = \dots = G_d = G$ for simplicity. For $d+1$ parties to generate a common password, the $i$th party just needs to pick a password $a_i$, perform a single exponentiation to get $g^{a_i}$, and broadcast  $g^{a_i}$ to the other parties, who will each do likewise so that every party now has $g^{a_1},\dots, g^{a_{d+1}}$. With these and the password $a_i$, the $i$th party will now compute
\[
\Phi(g^{a_1}, \ldots, g^{a_{i-1}}, g^{a_{i+1}}, \ldots,g^{a_{d+1}})^{a_i} = \Phi(g, \ldots,g)^{a_1 \cdots a_{d+1}},
\]
which will serve as the common password for the $d+1$ parties. Note that each of them would have arrived at  $\Phi(g, \ldots,g)^{a_1 \cdots a_{d+1}}$ in a different way with their own password but their results are guaranteed to be equal as a consequence of multilinearity.
There are several candidates for such cryptographic multilinear maps \cite{Graph,Lattice} and a variety of cryptographic applications that go beyond multipartite key exchange \cite{Boneh}.
\end{example}

We return to multilinear maps over familiar vector spaces. The next two examples are about trilinear functionals.

\begin{example}[trilinear functionals and self-concordance]\label{eg:self}
The definition of self-concordance is usually stated over $\mathbb{R}^n$.
Let $f \colon  \Omega \to \mathbb{R}$ be a convex $C^3$-function on an open convex subset $\Omega \subseteq \mathbb{R}^n$. Then $f$ is said to be \emph{self-concordant} at $x \in \Omega$ if
\begin{equation}\label{eq:self}
\lvert D^{3}f(x)(h,h,h)\rvert \le 2 [D^{2}f(x)(h,h)]^{3/2}
\end{equation}
for all $h \in \mathbb{R}^n$ \cite{Nesterov}. As we discussed in Example~\ref{eg:hod}, for any fixed $x \in \Omega$, the higher derivatives $D^2 f(x)$ and $D^3 f(x)$ in this case are bilinear and trilinear functionals on $\mathbb{R}^n$ given by
\begin{align*}
[D^{2}f(x)](h,h)&=\sum_{i,j=1}^{n}\dfrac{\partial^{2}f(x)}{\partial
x_i \partial x_{j}}h_{i}h_{j},\\*
[D^{3}f(x)](h,h,h)&=\sum_{i,j,k=1}
^{n}\dfrac{\partial^{3}f(x)}{\partial x_i \partial x_{j}\partial x_{k}}h_{i}h_{j}h_{k}.
\end{align*}
The affine invariance \cite[Proposition~2.1.1]{Nesterov} of self-concordance implies that self-concordance is a tensorial property in the sense of definition~\ref{st:tensor1}. For the convergence and complexity analysis of interior point methods, it goes hand in hand with the affine invariance of Newton's method in Example~\ref{eg:Newton}. Such analysis in turn allows one to establish the celebrated result that a convex optimization problem may be solved to arbitrary $\varepsilon$-accuracy  in polynomial time using interior point methods if it has a self-concordant barrier function whose first and second derivatives may be evaluated in polynomial time \cite[Chapter~6]{Nesterov}. These conditions are satisfied for many common problems including linear programming, convex quadratic programming, second-order cone programming, semidefinite programming and geometric programming \cite{Boyd}. Contrary to popular belief, polynomial-time solvability to $\varepsilon$-accuracy is not guaranteed by convexity alone: copositive and complete positive programming are convex optimization problems but both are known to be NP-hard \cite{Murty,Dickinson}.

By  Example~\ref{eg:tenfield}, we may view $D^2 f$ and $D^3 f$ as covariant tensor fields on the manifold $\Omega$ (any open subset of a manifold is itself a manifold) with 
\[
D^2 f(x) \colon  \mathbb{T}_x(\Omega) \times \mathbb{T}_x(\Omega) \to \mathbb{R}, \quad
D^3 f(x) \colon  \mathbb{T}_x(\Omega) \times \mathbb{T}_x(\Omega) \times \mathbb{T}_x(\Omega) \to \mathbb{R}
\]
for any fixed $x \in \Omega$. While this tensor field perspective is strictly speaking unnecessary, it helps us formulate \eqref{eq:self} concretely in situations when we are not working over $\mathbb{R}^n$. For instance, among the aforementioned optimization problems, semidefinite, complete positive and copositive programming require that we work over the space of symmetric matrices $\mathbb{S}^n$ with $\Omega$ given respectively by the following open convex cones:
\begin{align*}
\mathbb{S}^n_\pp &=\{A \in \mathbb{S}^n \colon  x^\tp Ax >0, \; 0\ne x\in \mathbb{R}^n \} =\{BB^\tp \in \mathbb{S}^n \colon  B \in \GL(n) \},\\*
\mathbb{S}^n_{\ppp} &=\{BB^\tp  \in \mathbb{S}^n \colon  B \in \GL(n) \cap \mathbb{R}^{n \times n}_\p \},\\*
\mathbb{S}^{n\ast}_{\ppp} &=\{A \in \mathbb{S}^n \colon  x^\tp Ax >0,\; 0\ne x\in \mathbb{R}^n_\p \}.
\end{align*}
In all three cases we have $\mathbb{T}_X(\Omega) = \mathbb{S}^n$ for all $X \in \Omega$, equipped with the usual trace inner product. To demonstrate self-concordance for the log barrier function $f(X)=-\log\det(X)$, Examples~\ref{eg:hog} and \ref{eg:logdet} give us
\begin{align*}
D^2 f(X)(H,H) &= \tr(H^\tp[\nabla^2 f(X)](H)) =\tr(HX^{-1} H X^{-1} ),\\*
D^3 f(X)(H,H,H) &= \tr(H^\tp[\nabla^3 f(X)](H,H)) =-2\tr(HX^{-1} H X^{-1} H X^{-1} ),
\end{align*}
and it follows from Cauchy--Schwarz that
\[
\lvert D^3 f(X)(H,H,H) \rvert \le 2 \lVert HX^{-1} \rVert^3 = 2 [D^2 f(X)(H,H) ]^{3/2},
\]
as required. For comparison, take the  \emph{inverse barrier} function $g(X) = \tr(X^{-1})$ which, as a convex function that blows up near the boundary of $\mathbb{S}^n_\pp $ and has easily computable gradient and Hessian  as we saw in Example~\ref{eg:hog}, appears to be a perfect barrier function for semidefinite programming. Nevertheless, using the derivatives calculated in Example~\ref{eg:hog},
\begin{align*}
D^2 g(X)(H,H) &=2 \tr( HX^{-1}HX^{-2}),\\*
D^3 g(X)(H,H,H) &= -6\tr( HX^{-1}HX^{-1}HX^{-2}),
\end{align*}
and since $6\lvert h \rvert^3/x^4 > 2(2h^2/x^3)^{3/2}$ as $x\to 0^+$, \eqref{eq:self}  will not be satisfied when  $X$ is near singular. Thus the inverse barrier function for $\mathbb{S}^n_\pp $ is not self-concordant.

What about $\mathbb{S}^n_{\ppp}$, the completely positive cone, and its dual cone $\mathbb{S}^{n\ast}_{\ppp}$, the copositive cone? It turns out that it is possible to construct convex self-concordant barrier functions for these \cite[Theorem~2.5.1]{Nesterov}, but the problem is that the gradients and Hessians of these barrier functions are not computable in polynomial time.
\end{example}

The last example, though more of a curiosity, is a personal favourite of the author \cite{norm}.

\begin{example}[spectral norm and Grothendieck's inequality]\label{eg:GI}\hspz%
One of the most fascinating inequalities in matrix and operator theory is the following. For any $A\in\mathbb{R}^{m\times n}$, there exists a constant $K_\G>0$ such
that
\begin{equation}\label{eq:GI}
\max_{\lVert x_i\rVert=\lVert y_j\rVert =1}\sum_{i=1}^{m}\sum_{j=1}^{n}a_{ij}\langle x_i ,y_j \rangle 
\leq K_\G\max_{\lvert \varepsilon_i\rvert = \lvert \delta_j\rvert =1}\sum_{i=1}^{m}\sum_{j=1}^{n}a_{ij}\varepsilon_{i}\delta_{j},
\end{equation}
with the maximum on the left taken over vectors $x_1,\ldots,x_m, y_1,\ldots,y_n \in \mathbb{R}^p$ of unit $2$-norm and that on the right taken over $\varepsilon_1,\ldots,\varepsilon_m, \delta_1,\ldots,\delta_n \in  \{-1,+1\}$. Here $p \ge m+n$ can be arbitrarily large.
What is remarkable is that the constant $K_\G$ is \emph{universal}, \ie\ independent of $m,n,p$ or the matrix $A$. It is christened the Grothendieck constant after the discoverer of the inequality \cite{Grothendieck}, which has found widespread applications in combinatorial optimization, complexity theory and quantum information \cite{Pisier}. The exact value of the constant is unknown but there are excellent bounds \cite{Davie,Krivine}:
\[
1.67696 \le K_\G \le 1.78221.
\]
Nevertheless the nature of the Grothendieck constant remains somewhat mysterious. One clue is that since the constant is independent of the matrix $A$ and  obviously also independent of the $x_i$, $y_j$, $\varepsilon_i$, $\delta_j$, as these are dummy variables that are maximized over, whatever underlying object the Grothendieck constant measures must have nothing to do with any of the quantities appearing in \eqref{eq:GI}. We will show that this object is in fact the Strassen matrix multiplication tensor $\Mu_{m,n,p}$ in  Example~\ref{eg:Strassen}. First, observe that the right-hand side of \eqref{eq:GI} is just the matrix $(\infty,1)$-norm,
\[
\lVert A \rVert_{\infty,1} = \max_{v\ne 0} \dfrac{\lVert Av \rVert_1}{\lVert v\rVert_\infty}=\max_{\varepsilon_i,\delta_j= \pm 1}\sum_{i=1}^{m}\sum_{j=1}^{n}a_{ij}\varepsilon_{i}\delta_{j},
\]
which is in fact a reason why the inequality is so useful, the $(\infty,1)$-norm being  ubiquitous  in combinatorial optimization but NP-hard to compute, and the left-hand side of \eqref{eq:GI} being readily computable via semidefinite programming. Writing $x_i$ and $y_j$, respectively,
as columns and rows of matrices $X = [x_1,\ldots,x_m] \in \mathbb{R}^{p \times m}$ and $Y =[y_1^\tp,\ldots,y_n^\tp]^\tp \in \mathbb{R}^{n \times p}$, we next observe that the constraints on the left-hand side of \eqref{eq:GI} may be expressed in terms of their $(1,2)$-norm and $(2,\infty)$-norm:
\begin{align*}
\lVert X \rVert_{1,2} &= \max_{v \ne 0} \dfrac{\lVert Xv \rVert_2}{\lVert v\rVert_1}= \max_{i=1,\ldots,m} \lVert x_i \rVert_2,\\*[4pt]
\lVert Y \rVert_{2,\infty} &= \max_{v \ne 0} \dfrac{\lVert Yv \rVert_\infty}{\lVert v\rVert_2} = \max_{j=1,\ldots,n} \lVert y_j \rVert_2,
\end{align*}
namely, as $\lVert X \rVert_{1,2} \le 1$ and $\lVert Y \rVert_{2,\infty} \le 1$. Lastly, we observe that for the standard inner product on $\mathbb{R}^p$,
\[
\sum_{i=1}^m\sum_{j=1}^n a_{ij}\langle x_i ,y_j \rangle = \tr(XAY),
\]
and thus \eqref{eq:GI} may be written as
\[
\max_{A,X,Y\ne 0}\dfrac{ \tr(XAY)}{\lVert
A\rVert_{\infty,1}\lVert X\rVert_{1,2}\lVert Y\rVert_{2,\infty}} \le K_\G.
\]
As we saw in \eqref{eq:specnorm1}, the expression on the left is precisely the spectral norm of the trilinear functional
\[
\tau_{m,n,p} \colon  \mathbb{R}^{m \times n} \times \mathbb{R}^{n \times p} \times \mathbb{R}^{p \times m} \to \mathbb{R}, \quad (A,B,C) \mapsto \tr(ABC),
\]
if we equip the three vector spaces of matrices with the matrix $(1,2)$, $(\infty,1)$ and $(2,\infty)$-norms respectively. We shall denote this norm as $\lVert \, \cdot \, \rVert_{1,2,\infty}$. The last step is to observe that as $3$-tensors, the trilinear functional $\tau_{m,n,p}$ and the  matrix multiplication tensor 
\[
\Mu_{m,n,p} \colon   \mathbb{R}^{m \times n} \times \mathbb{R}^{n \times p} \to \mathbb{R}^{m \times p}, \quad (A,B) \mapsto AB
\]
are one and the same, assuming we identify ($\mathbb{R}^{p \times n})^* = \mathbb{R}^{n \times p}$. One way to see this is that the same type of tensor transformation rule applies:
\begin{gather*}
(XAY^{-1})(YBZ^{-1}) = X(AB)Z^{-1}, \\* \tr((XAY^{-1})(YBZ^{-1})(ZCX^{-1})) =  \tr(ABC)
\end{gather*}
for any $(X,Y,Z) \in \GL(m) \times \GL(n) \times \GL(p)$. As we have {alluded,} recognizing that two different multilinear maps correspond to the same tensor without reference to the transformation rule is one benefit of definition~\ref{st:tensor3}, which we will discuss immediately after this example. Hence we have shown that
\[
K_\G = \sup_{m,n,p \in \mathbb{N}} \lVert \Mu_{m,n,p}  \rVert_{1,2,\infty}.
\]
In particular, the Grothendieck constant is an asymptotic value of the spectral norm of $\Mu_{m,n,p}$ just as the exponent of matrix multiplication $\omega$ in Example~\ref{eg:Strassen} is an asymptotic value of the tensor rank of $\Mu_{m,n,p}$.
\end{example}

\section{Tensors via tensor products}\label{sec:tenprod}

In the last two sections we have alluded to various inadequacies of definitions~\ref{st:tensor1} and~\ref{st:tensor2}. The modern definition of tensors, \ie\ definition~\ref{st:tensor3}, remedies them by regarding a tensor as an element of a tensor product of vector spaces (or, more generally, modules). We may view the development of the notion of tensor via these three definitions as a series of improvements. Definition~\ref{st:tensor1} leaves the tensor unspecified and instead just describes its change-of-basis theorems as transformation rules. Definition~\ref{st:tensor2} fixes this by supplying a multilinear map as a  candidate for the unspecified object but the problem becomes one of oversupply:  with increasing order $d$, there is an exponential multitude of different $d$-linear maps for each transformation rule. Definition~\ref{st:tensor3} fixes both issues by providing an object for definition~\ref{st:tensor1} that at the same time also neatly classifies the  multilinear maps in definition~\ref{st:tensor2} by the transformation rules they satisfy. In physics, the perspective may be somewhat different. Instead of viewing them as a series of pedagogical improvements, definition~\ref{st:tensor1} was discovered alongside (most notably) general relativity, with definition~\ref{st:tensor2} as its addendum, and definition~\ref{st:tensor3} was discovered alongside (most notably) quantum mechanics. Both remain useful in different ways for different purposes and each is used independently of the other.

One downside of definition~\ref{st:tensor3} is that it relies on the tensor product construction and/or the universal factorization property;  both have a reputation of being abstract. We will see that the reputation is largely undeserved; with the proper motivations they are far easier to appreciate than, say, definition~\ref{st:tensor1}. One reason for its `abstract' reputation is that the construction is often cast in a way that  lends itself to vast generalizations. In fact, tensor products of Hilbert spaces \cite{Reed}, modules \cite{Lang}, vector bundles \cite{Milnor}, operator algebras \cite{Takesaki}, representations \cite{Fulton}, sheaves \cite{Hartshorne}, cohomology rings \cite{Hatcher}, {\em etc.}, have become foundational materials covered in standard textbooks, with more specialized tensor product constructions, such as those of Banach spaces \cite{Ryan}, distributions \cite{Treves}, operads \cite{Operad} or more generally objects in any tensor category \cite{Etingof}, readily available in various monographs. This generality is  a feature, not a bug. In our (much more modest) context, it allows us to define tensor products of  norm spaces and inner product spaces, which in turn allows us to define norms and inner products for tensors, to view multivariate functions as tensors, an enormously useful perspective in computation, and to identify `separation of variables' as the common underlying thread in a disparate array of well-known algorithms. In physics, the importance of  tensor products cannot be over-emphasized: it is one of the fundamental postulates of quantum mechanics \cite[Section~2.2.8]{Nielsen}, the source behind many curious quantum phenomena that lie at the heart of the subject, and is indispensable in technologies such as quantum computing and quantum cryptography.\\

Considering its central role, we will discuss three common approaches to constructing tensor products, increasing in abstraction:
\begin{enumerate}[\upshape (i)]
\setlength\itemsep{3pt}
\item via tensor products of function spaces,
\item via tensor products of more general vector spaces,
\item via the universal factorization property.
\end{enumerate}
As in the case of the three definitions of tensor, each of these constructions is useful in its own way and each may be taken to be a variant of definition~\ref{st:tensor3}. We will motivate each construction with concrete examples but we defer all examples of computational relevance to Section~\ref{sec:sepvar}. The approach of defining a tensor as an element of a tensor product of vector spaces likely  first appeared in the first French edition of \citet{Bourbaki} and is now standard in graduate algebra texts \cite{DumFoot,Lang,Vinberg}. It has also caught on in physics \cite{Geroch} and in statistics \cite{McC}.
For further historical information we refer the reader to \citet{Conrad} and the last chapter of \citet{KM}.

\subsection{Tensor products of function spaces}\label{sec:tensor3a}

This is the most intuitive of our three  tensor product constructions as it reflects the way we build functions of more variables out of functions of fewer variables. In many situations, especially analytic ones, this tensor-as-a-multivariate-function view suffices and one need go no further.

For any set $X$, the set of all real-valued functions on $X$,
\[
\mathbb{R}^X \coloneqq \{ f \colon  X \to \mathbb{R} \} ,
\]
is always a vector space under pointwise addition and scalar multiplication of functions, \tie, by  defining
\[
(a_1f_1 + a_2f_2)(x) \coloneqq a_1 f_1(x) + a_2 f_2(x)\quad\text{for each } x \in X.
\]
Note that the left-hand side is a real-valued function whose value at $x$ is defined by the right-hand side, which involves only additions and multiplications of real numbers $a_1, a_2, f_1(x), f_2(x)$. For any finite set $X$,
$\dim \mathbb{R}^X = \# X$, the cardinality of~$X$.

It turns out that the converse is also true: any vector space $\mathbb{V}$ may be regarded as a function space $\mathbb{R}^X$ for some $X$. It is easiest to see this for finite-dimensional vector spaces, so we first assume that $\dim\mathbb{V} =n \in \mathbb{N}$. In this case $\mathbb{V}$ has a basis $\mathscr{B}=\{v_1,\ldots,v_n\}$ and every element $v \in \mathbb{V}$ may be written as
\begin{equation}\label{eq:uniquerep}
v = a_1 v_1 + \dots + a_n v_n\quad\text{for some unique } a_1,\ldots,a_n \in \mathbb{R}.
\end{equation}
So we simply take $X = \mathscr{B}$ and regard $v\in \mathbb{V}$ as the function $f_v \in \mathbb{R}^X$ defined by
\[
f_v \colon  X \to \mathbb{R}, \quad f_v(v_i) = a_i, \quad i =1,\ldots,n.
\]
The map $\mathbb{V} \to \mathbb{R}^{\mathscr{B}}$, $v \mapsto f_v$ is linear and invertible and so the two vector spaces are isomorphic.
The function $f_v$ depends on our choice of basis, and any two such functions $f_v$ and $f'_v$ corresponding to different bases $\mathscr{B}$ and $\mathscr{B}'$ will be related via the contravariant $1$-tensor transformation rule if we represent $f_v$ as a vector in $\mathbb{R}^n$ with coordinates given by its values
\[
 [f_v]_{\mathscr{B}} \coloneqq \begin{bmatrix} f(v_1) \\ \vdots \\ f(v_n) \end{bmatrix} =\begin{bmatrix} a_1 \\ \vdots \\ a_n \end{bmatrix} \in \mathbb{R}^n
\]
and do likewise for $f'_v$. This is of course just restating \eqref{eq:vecrep} and \eqref{eq:cob1} in slightly different language.  As long as we assume that every vector space $\mathbb{V}$ has a basis  $\mathscr{B}$, or, equivalently, the axiom of choice \cite{Blass}, the above construction for $\mathbb{V} \cong \mathbb{R}^{\mathscr{B}}$ in principle applies verbatim to infinite-dimensional vector spaces, as having a unique finite linear combination of the form \eqref{eq:uniquerep} is the very definition of a basis. When used in this sense, such a basis for a vector space is called a Hamel basis.

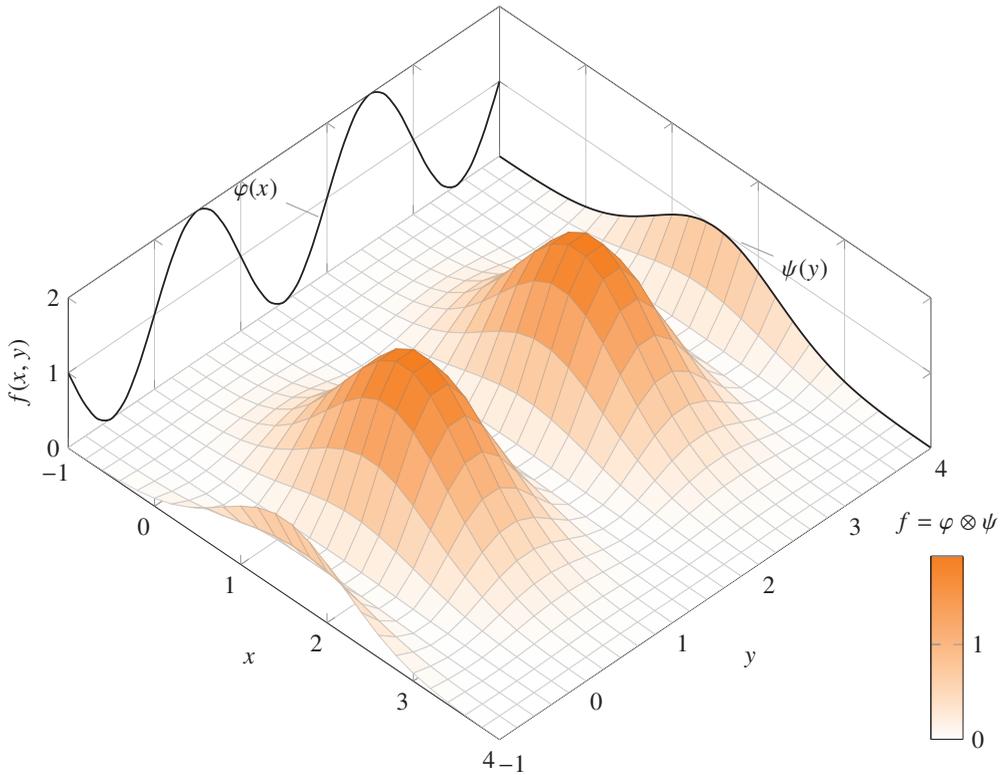
\begin{figure}[ht]
\pgfplotsset{width=7cm,compat=1.8}
\pgfplotsset{%
  colormap={whitered}{color(0cm)=(white);
    color(1cm)=(orange!75!red)}
}
\begin{tikzpicture}[scale=0.85,
  declare function={
    phi=exp(-(x-3/2)^2);
  },
  declare function={
    psi=1+sin(180*x);
  },
  declare function={
    f=exp(-(x-3/2)^2)*(1+sin(180*y));
  }]

  \begin{axis}[
    colormap name=whitered,
    width=15cm,
    view={45}{70},
    enlargelimits=false,
    grid=major,
    domain=-1:4,
    y domain=-1:4,
    samples=26,
    xlabel=$x$,
    ylabel=$y$,
    zlabel={$f(x,y)$},
    colorbar,
    colorbar style={
      at={(1,0)},
      anchor=south west,
      height=0.25*\pgfkeysvalueof{/pgfplots/parent axis height},
      title={$f= \varphi \otimes \psi$}
    }]

    \addplot3 [surf] {f};
    \addplot3 [domain=-1:4, samples=31, samples y=0, thick, smooth]
      (x, 4, {phi});
    \addplot3 [domain=-1:4, samples=31, samples y=0, thick, smooth]
      (-1, x, {psi});

    \node at (axis cs:-1,2,0.7) [pin=165:$\varphi(x)$] {};
    \node at (axis cs:1.7,4,1.0) [pin=-15:$\psi(y)$] {};
  \end{axis}
\end{tikzpicture}
\caption{Separable function $f(x,y) = \varphi(x)\psi(y)$ with $\varphi(x)=1+\sin(180x)$ and $\psi(y)=\exp\bigl(-(y-3/2)^2\bigr)$.}
\label{fig:sepfun}
\end{figure}

Since any vector space is isomorphic to a function space, by defining tensor products on function spaces we define tensor products on all vector spaces. This is the most straightforward approach to defining tensor products. In this regard, a $d$-tensor is just a real-valued function $f(x_1,\ldots,x_d)$ of $d$ variables. For two sets $X$ and $Y$, we define the \emph{tensor product} $\mathbb{R}^X \otimes \mathbb{R}^Y$ of $\mathbb{R}^X$ and $\mathbb{R}^Y$ to be  the subspace of $\mathbb{R}^{X \times Y}$ comprising all functions that can be written as a finite sum of product of two univariate functions, one of $x$ and another of $y$:
\begin{equation}\label{eq:sep1}
f(x,y) =  \sum_{i=1}^r \varphi_i(x) \psi_i(y).
\end{equation}
For univariate functions $\varphi \colon  X \to \mathbb{R}$ and $\psi \colon  Y \to \mathbb{R}$, we will write $\varphi \otimes \psi \colon  X \times Y \to \mathbb{R}$ for the bivariate function defined by
\[
(\varphi \otimes \psi) (x,y) = \varphi(x)\psi(y)\quad \text{for all } x\in X, \; y \in Y.
\]
Such a function is called a \emph{separable function}, and when $\otimes$ is used in this sense, we will refer to it a \emph{separable product}; see Figure~\ref{fig:sepfun} for a simple example. In notation,
\[
\mathbb{R}^X \otimes \mathbb{R}^Y \coloneqq \biggl\{ f \in \mathbb{R}^{X \times Y} \colon  f = \sum_{i=1}^r \varphi_i \otimes \psi_i ,\; \varphi_i \in \mathbb{R}^X, \; \psi_i \in \mathbb{R}^Y, \; r \in \mathbb{N} \biggr\}.
\]
Note that $\mathbb{R}^{X \times Y} = \{ f\colon  X \times Y \to \mathbb{R} \}$ is the set of all bivariate $\mathbb{R}$-valued functions on $X \times Y$, separable or not. In reality, when $\mathbb{V}$ is infinite-dimensional, its Hamel basis  $\mathscr{B}$ is often uncountably infinite and impossible to write down explicitly for many common infinite-dimensional vector spaces (\eg\  any complete topological vector space), and thus this approach of constructing tensor products by identifying a vector space $\mathbb{V}$  with $ \mathbb{R}^{\mathscr{B}}$ is not practically useful for most infinite-dimensional vector spaces. Furthermore, if we take this route we will inevitably have to discuss a `change of Hamel bases', a tricky undertaking that has never been done, as far as we know. As such, tensor products of infinite-dimensional vector spaces are almost always discussed in the context of \emph{topological tensor products}, which we will describe in Example~\ref{eg:infinite}.

If $X$ and $Y$ are finite sets, then in fact we have
\begin{equation}\label{eq:mv0}
\mathbb{R}^X \otimes \mathbb{R}^Y = \mathbb{R}^{X \times Y}
\end{equation}
as, clearly,
\[
\dim (\mathbb{R}^X \otimes \mathbb{R}^Y)= \# X \cdot \# Y = \# (X\times Y) =\dim  \mathbb{R}^{X \times Y},
\]
and so every function $f \colon  X \times Y \to \mathbb{R} $ can be expressed in the form \eqref{eq:sep1} for some finite~$r$. In other words, elements of $\mathbb{R}^X \otimes \mathbb{R}^Y$ give all functions of two variables $(x,y) \in X \times Y$. If $X$ and $Y$ are infinite sets, then \eqref{eq:mv0} is demonstrably false.  Nevertheless, function spaces over infinite sets almost always have additional properties imposed (\eg\  continuity, differentiability, integrability, {\em etc.})\ and/or carry additional structures (\eg\  metrics, norms, inner products, completeness in the topology induced by these, {\em etc.}). If we constrain our functions in some manner (instead of allowing them to be completely arbitrary) and limit our attention to some smaller subspace $\mathcal{F}(X) \subseteq \mathbb{R}^X$, then it turns out that the analogue of \eqref{eq:mv0} often holds, \tie,
\begin{equation}\label{eq:fun}
\mathcal{F}(X) \otimes \mathcal{F}(Y) = \mathcal{F}(X \times Y)
\end{equation}
for the right function class $\mathcal{F}$ and the right interpretation of $\otimes$, as we will see in the next two examples. For infinite-dimensional vector spaces, there are many ways to define (topological) tensor products, and the `right tensor product' is usually regarded as the one that gives  us \eqref{eq:fun}.

\begin{example}[tensor products in infinite dimensions]\label{eg:infinite}\hspz%
We write $\mathbb{R}[x_1, \dots , x_m]$ for the set of all multivariate polynomials in $m$ variables, with no restriction on degrees. This is one of the few infinite-dimensional vector spaces that has a countable Hamel basis, namely the set of all monic monomials
\[
\mathscr{B} = \{ x_1^{d_1} x_2^{d_2} \cdots x_m^{d_m} \colon  d_1,d_2,\ldots,d_m \in \mathbb{N} \cup \{0\} \}.
\]
For multivariate polynomials, it is clearly true that\footnote{One may extend \eqref{eq:mv2}  to monomials involving arbitrary real powers, \ie\ Laurent polynomials, posynomials, signomials, {\em etc.}, as they are all \emph{finite} sums of monomials.} 
\begin{equation}\label{eq:mv2}
\mathbb{R}[x_1, \dots , x_m] \otimes \mathbb{R}[y_1, \ldots, y_n]  =
\mathbb{R}[x_1, \ldots, x_m, y_1, \ldots, y_n],
\end{equation}
as polynomials are sums of finitely many monomials and monomials are always separable, for example
\[
7x_1^2 x_2 y_2^3 y_3 - 6x_2^4 y_1 y_2^5 y_3  = (7x_1^2 x_2) \cdot ( y_2^3 y_3) + (- 6x_2^4) \cdot(y_1 y_2^5 y_3).
\]
But \eqref{eq:fun} is also clearly false in general for other infinite-dimensional vector spaces. Take continuous real-valued functions  for illustration. While $\sin(x+y)=\sin x \cos y + \cos x \sin y $ and $\log(xy) = \log x + \log y$, we can never have
\[
\sin(xy) = \sum_{i=1}^r \varphi_i(x)\psi_i(y)\quad\text{or}\quad \log(x-y) = \sum_{i=1}^r \varphi_i(x)\psi_i(y)
\]
for any continuous functions $\varphi_i$, $\psi_i$ and finite $r\in \mathbb{N}$. But if we allow $r \to \infty$, then 
\[
\sin(xy) =\sum^{\infty}_{n=0} \dfrac{(-1)^n}{(2n+1)!} x^{2n+1} y^{2n+1},
\]
and by Taylor's theorem, $\sin(xy)$ can be approximated to arbitrary accuracy by sums of separable functions. Likewise, 
\begin{align*}
\log(x-y) &= \log(-y) -\dfrac{x}{y} - \dfrac{x^2}{2y^2} - \dfrac{x^3}{3y^3} -\dots -\dfrac{x^n}{ny^n} +O(x^{n+1})
\intertext{as $x \to 0$, or, more relevant for our later purposes,}
\log(x-y) &= \log(x) -\dfrac{y}{x} - \dfrac{y^2}{2x^2} - \dfrac{y^3}{3x^3} -\dots -\dfrac{y^n}{nx^n} +O\biggl(\dfrac{1}{x^{n+1}}\biggr),
\end{align*}
as $x \to \infty$. In other words, we may approximate $\log(x-y)$ by sums of separable functions when $x$ is very small or  very large;  the latter will be important when we discuss the fast multipole method in Example~\ref{eg:multpole}.

Given that \eqref{eq:mv2} holds for polynomials, the Stone--Weierstrass theorem gives an indication that we might be able to extend \eqref{eq:fun} to other infinite-dimensional spaces, say, the space of continuous functions $C(X)$ on some nice domain $X$. However, we will need to relax the definition of tensor product to allow for limits, \ie\ taking \emph{completion} with respect to some appropriate choice of norm. For instance, we may establish, for $1 \le p < \infty$,
\[
C(X) \hatotimes C(Y) = C(X \times Y), \quad L^p(X) \hatotimes L^p(Y) = L^p(X \times Y),
\]
where $X$ and $Y$ are locally compact Hausdorff spaces in the left equality and $\sigma$-finite measure spaces in the right equality \cite[Corollaries~1.14 and 1.52]{Cheney}. The \emph{topological tensor product} $C(X) \hatotimes C(Y)$ above refers to the completion~of
\begin{align*}
  C(X) \otimes C(Y) 
 = \biggl\{ f \in C(X \times Y) \colon  f = \sum_{i=1}^r \varphi_i \otimes \psi_i , \varphi_i \in C(X), \; \psi_i \in C(Y), \; r \in \mathbb{N} \biggr\}
\end{align*}
with respect to some \emph{cross-norm}, \ie\ a norm $\lVert\,\cdot\,\rVert$ that is multiplicative on rank-one elements,
\[
\lVert \varphi \otimes \psi \otimes \dots \otimes \theta \rVert = \lVert \varphi \rVert \lVert \psi \rVert  \cdots \lVert  \theta \rVert,
\]
like those in \eqref{eq:nn0}, \eqref{eq:sn0}, \eqref{eq:tennuclear}, \eqref{eq:tenspectral}; we will discuss cross-norms in greater detail in Examples~\ref{eg:norms} and \ref{eg:Hilbert} for Banach and Hilbert spaces. As we have seen, without taking completion, $C(X) \otimes C(Y)$ will be a strictly smaller subset of $C(X \times Y)$. The same goes for $L^p(X) \hatotimes L^p(Y)$, but the completion is with respect to some other cross-norm dependent on the value of $p$.

There many important infinite-dimensional vector spaces that are not Banach (and thus not Hilbert) spaces. The topology of these spaces is usually defined by families of \emph{seminorms}. Fortunately, the   construction outlined above also works with such spaces to give topological tensor products that satisfy \eqref{eq:fun}. If we limit our consideration of $X$ and $Y$ to subsets of Euclidean spaces, then for spaces of rapidly decaying {Schwartz} functions (see Example~\ref{eg:Calderon}), smooth functions, compactly supported smooth functions and holomorphic functions, we have \cite[Theorem~51.6]{Treves}
\begin{alignat*}{3}
S(X) \hatotimes S(Y) &= S(X \times Y),\quad  &
C^\infty(X) \hatotimes C^\infty(Y) &= C^\infty(X \times Y),\\*
C^\infty_c(X) \hatotimes C^\infty_c(Y) &= C^\infty_c(X \times Y),\quad &
H(X) \hatotimes H(Y) &= H(X \times Y),
\end{alignat*}
where $X$ and $Y$, respectively, are $\mathbb{R}^m$ and $\mathbb{R}^n$ (Schwartz), open subsets of $\mathbb{R}^m$ and $\mathbb{R}^n$ (smooth), compact subsets of  $\mathbb{R}^m$ and $\mathbb{R}^n$ (compactly supported) and open subsets of  $\mathbb{C}^m$ and $\mathbb{C}^n$ (holomorphic).

Similar results have been established for more complicated function spaces such as Sobolev and Besov spaces \cite{Sickel} or extended to more general domains such as having $X$ and $Y$ be smooth manifolds or Stein manifolds \cite{Akbarov}. The bottom line is that while it is not always true that all $(X,Y)$-variate functions are separable products of $X$-variate functions and $Y$-variate functions  --  $L^\infty$-functions being a notable exception  --  \eqref{eq:fun} holds true for many common function spaces and many common types of domains when $\otimes$ is interpreted as an appropriate topological tensor product $\hatotimes$.
\end{example}

To keep our promise of simplicity, our discussion in Example~\ref{eg:infinite} is necessarily terse as the technical details are quite involved. Nevertheless, the readers will get a fairly good idea of the actual construction for the special cases of Banach and Hilbert spaces in Examples~\ref{eg:norms} and \ref{eg:Hilbert}. We will briefly mention another example.

\begin{example}[tensor products of distributions]\label{eg:distributions}\hspz%
Strictly speaking, the construction of tensor products via sums of separable functions gives us \emph{contravariant} tensors, as we will see more clearly in Example~\ref{eg:hyp}. Nevertheless, we may get covariant tensors by taking the continuous dual of the function spaces in Example~\ref{eg:infinite}. For a given $\mathbb{V}$ equipped with a norm, this is the set of all linear functionals $\varphi \in \mathbb{V}^*$ whose dual norm in \eqref{eq:dualnorm} is finite. For a function space $\mathcal{F}(X)$, the continuous dual space will be denoted $\mathcal{F}'(X)$ and it is a subspace of the algebraic dual space $ \mathcal{F}(X)^*$. Elements of  $\mathcal{F}'(X)$ are called \emph{distributions} or generalized functions with the best-known example being the Dirac delta `function', intuitively defined by
\[
\int_{-\infty}^\infty f(x) \delta(x) \D x = f(0),
\]
and rigorously as the continuous linear functional $\delta \colon  C^\infty_c(\mathbb{R}) \to \mathbb{R}$, $\varphi \to \varphi(0)$, \ie\ given by evaluation at zero. The continuous duals of $C^\infty_c(X)$, $S(X)$, $C^\infty(X)$ and $H(X)$ are denoted by $D'(X)$ (distributions), $S'(X)$ (tempered distributions), $E'(X$) (compactly supported distributions) and $H'(X)$ (analytic functionals) respectively, and one may obtain the following dual analogues \cite[Theorem~51.7]{Treves}:
\begin{equation}\label{eq:sepdist}
\begin{aligned}
S'(X) \hatotimes S'(Y) &= S'(X \times Y),  &
E'(X) \hatotimes E'(Y) &= E'(X \times Y),\\
D'(X) \hatotimes D'(Y) &= D'(X \times Y),\;  &
H'(X) \hatotimes H'(Y) &= H'(X \times Y),
\end{aligned}
\end{equation}
with $X$ and $Y$ as in Example~\ref{eg:infinite}. These apparently innocuous statements conceal some interesting facts. While every linear operator $\Phi \colon  \mathbb{V} \to \mathbb{W}$ and bilinear functional $\beta \colon  \mathbb{V} \times \mathbb{W} \to \mathbb{R}$ on finite-dimensional vector spaces  may be represented by a matrix $A \in \mathbb{R}^{m \times n}$,
as  we discussed in Section~\ref{sec:multmaps}, the continuous analogue on infinite-dimensional spaces is clearly false. For instance, not every continuous linear operator $\Phi \colon  L^2(\mathbb{R}) \to L^2(\mathbb{R})$ may be expressed as
\begin{equation}\label{eq:wrong}
[\Phi(f)](x) = \int_{-\infty}^\infty K(x,y) f(y) \D y
\end{equation}
with an integral kernel $K$ in place of the matrix $A$, as the identity operator is already a counterexample. This is the impetus behind the idea of a \emph{reproducing kernel Hilbert space} that we will discuss in Example~\ref{eg:Mercer}. What \eqref{eq:sepdist} implies is that this is true for the function spaces in Example~\ref{eg:infinite}. Any continuous linear operator $\Phi \colon  S(X) \to S'(Y)$ has a unique kernel $K \in S'(X \times Y)$
and likewise with $C^\infty_c$, $C^\infty$, $H$ in place of $S$ \cite[Chapters~50 and 51]{Treves}. For bilinear functionals, \eqref{eq:sepdist} implies that for these spaces, say, $\beta \colon  C^\infty_c(X) \times C^\infty_c(Y) \to \mathbb{R}$, we may decompose it into
\begin{equation}\label{eq:bilinfundecomp}
\beta = \sum_{i=1}^\infty a_i \varphi_i \otimes \psi_i
\end{equation}
with $(a_i)_{i=1}^\infty \in l^1(\mathbb{N})$, and $\varphi_i \in D'(X)$, $\psi_i \in D'(Y)$ satisfying $\lim_{i\to \infty} \varphi_i = 0 = \lim_{i\to \infty} \psi_i$, among other properties \cite[Theorem~9.5]{Schaefer}. These results are all consequences of Schwartz's kernel theorem, which, alongside Mercer's kernel theorem in Example~\ref{eg:Mercer}, represent different approaches to resolving the falsity of \eqref{eq:wrong}. While Mercer restricted to smaller spaces (reproducing kernel Hilbert spaces), Schwartz expanded to larger ones (spaces of distributions), but both are useful. Among other things, the former allows for the $\varphi_i$ and $\psi_i$ in \eqref{eq:bilinfundecomp} to be functions instead of distributions, from which we obtain feature maps for support vector machines \cite{Scholkopf,Steinwart}; on the other hand the latter would give us distributional solutions for PDEs \cite[Theorem~52.6]{Treves}. 
\end{example}

The discussion up to this point extends to any $d$ sets $X_1, X_2,\ldots, X_d$, either by repeated application of the $d=2$ case or directly:
\begin{align}\label{eq:tp1}
  &\mathbb{R}^{X_1} \otimes\mathbb{R}^{X_2} \otimes \dots \otimes \mathbb{R}^{X_d} \notag \\*
  &\quad \coloneqq \biggl\{ f \in \mathbb{R}^{X_1 \times X_2 \times \dots \times X_d} \colon 
  f = \sum_{i=1}^r \varphi_i \otimes \psi_i \otimes \dots \otimes \theta_i, \notag \\*
  &\hspace{40pt}
  \varphi_i \in \mathbb{R}^{X_1},  \psi_i \in \mathbb{R}^{X_2},\ldots,\theta_i \in \mathbb{R}^{X_d},\; i = 1,\ldots, r\in \mathbb{N} \biggr\}.
\end{align}
Each summand is a separable function defined by
\begin{equation}\label{eq:sep}
(\varphi \otimes \psi \otimes \dots \otimes \theta)(x_1,x_2,\ldots, x_d) \coloneqq \varphi(x_1) \psi(x_2) \cdots \theta(x_d) 
\end{equation}
for $x_1 \in X_1, x_2 \in X_2,\ldots,x_d \in X_d$, and
\[
\mathbb{R}^{X_1 \times \dots \times X_d} = \{ f\colon  X_1 \times \dots \times X_d \to \mathbb{R} \}
\]
is the set of $d$-variate real-valued functions on $X_1 \times \dots \times X_d$. We will next look at some concrete examples of separable functions.

\begin{example}[multivariate\zza normal\zza and\zza quantum\zza harmonic\zza oscillator]\label{eg:Gaussian}\hspace{-3.2pt}%
The  quintessential example of an $n$-variate separable function  is a Gaussian:
\[
f(x) = \exp (x^* A x + b^* x  + c)
\]
for $A \in \mathbb{C}^{n \times n}$, $b \in \mathbb{C}^n$, $c \in \mathbb{C}$, which is separable under a change of variables $x \mapsto V^* x$ where the columns of $V \in \Un(n)$ are an orthonormal eigenbasis of $(A+A^*)/2$. Usually one requires that $(A+A^*)/2$ is negative definite and $b$ is purely imaginary to ensure that it is an $L^2$-function. This is a singularly important separable function.

In statistics, it arises as the probability density, moment generating and characteristic functions of a normal random variable $X \sim N(\mu, \Sigma)$, with the last of the three given by
\[
\varphi_X (x) = \exp\biggl( \rmi\mu^\tp x - \dfrac{1}{2} x^\tp\Sigma x\biggr).
\]
This gives the unique distribution all of whose cumulant $d$-tensors, \ie\ the $d$th derivative as in Example~\ref{eg:hod} of $\log \varphi_X$, vanish when $d \ge 3$, a fact from which the central limit theorem follows almost as an afterthought \cite[Section~2.7]{McC}. It is an undisputed fact that this is the single most important probability distribution in statistics.

In physics, it arises as Gaussian wave functions for the quantum harmonic oscillator. For a spinless particle in $\mathbb{R}^3$ of mass $\mu$ subjected to the potential
\[
V(x,y,z) = \dfrac{1}{2}\mu (\omega_x^2 x^2 + \omega_y^2 y^2+\omega_z^2 z^2),
\]
if we assume that $V$ is isotropic, \ie\ $\omega_x = \omega_y = \omega_z = \omega$, then the  states are given~by
\[
\psi_{m,n,p}(x,y,z)=\biggl(\dfrac{\beta^2}{\pi}\biggr)^{3/4} \dfrac{H_m(\beta x) H_n(\beta y) H_p(\beta z)}{\sqrt{2^{m+n+p}\, m!\, n!\, p!} }\exp\biggl[\dfrac{\beta^2}{2}(x^2 + y^2 + z^2)\biggr] ,
\]
where $\beta^2 \coloneqq \mu \omega/\hbar$ and $H_n$ denotes a Hermite polynomial of degree $n$
\cite[Chapter~VII, Complement~B]{Cohen1}. The ground state, \ie\ $m = n = p=0$, is a Gaussian but note that even excited states, \ie\ $m+n+p > 0$, are separable functions. If one subscribes to the dictum that `physics is that subset of human experience which can be reduced to coupled harmonic oscillators' (attributed to Michael Peskin), then this has a status not unlike that of the normal distribution in statistics.
\end{example}

As in the case $d=2$, when $X_1,\ldots,X_d$ are finite sets, we have
\begin{equation}\label{eq:mv1}
\mathbb{R}^{X_1} \otimes \dots \otimes \mathbb{R}^{X_d} = \mathbb{R}^{X_1 \times \dots \times X_d},
\end{equation}
\tie, $d$-tensors in $\mathbb{R}^{X_1} \otimes \dots \otimes \mathbb{R}^{X_d}$ are just functions of $d$ variables $(x_1,\ldots,x_d) \in X_1\times \dots \times X_d$. In particular, the dimension of  the tensor space is
\[
\dim (\mathbb{R}^{X_1} \otimes \dots \otimes \mathbb{R}^{X_d}) = n_1  \cdots n_d,
\]
where $n_i \coloneqq \# X_i$, $i=1,\ldots,d$. As in the $d=2$ case, constructing a tensor  product of infinite-dimensional vector spaces $\mathbb{V}_1,\ldots,\mathbb{V}_d$ by picking Hamel bases $\mathscr{B}_1,\ldots,\mathscr{B}_d$ and then defining
\begin{equation}\label{eq:hamel}
\mathbb{V}_1 \otimes \dots \otimes \mathbb{V}_d \coloneqq \mathbb{R}^{\mathscr{B}_1} \otimes \dots \otimes \mathbb{R}^{\mathscr{B}_d}
\end{equation}
is rarely practical, with multivariate polynomials just about the only exception: Clearly,
\begin{align*}
& \mathbb{R}[x_1, \dots , x_m] \otimes \mathbb{R}[y_1, \ldots, y_n]  \otimes \dots \otimes \mathbb{R}[z_1, \ldots, z_p]  \\*
&\quad = \mathbb{R}[x_1, \ldots, x_m, y_1, \ldots, y_n,\ldots,z_1,\ldots,z_p],
\end{align*}
just as in  \eqref{eq:mv2}. However, if $\hatotimes$ is interpreted to be an appropriate topological tensor product, then repeated applications of the results in Example~\ref{eg:infinite} lead us to equalities like
\begin{equation}\label{eq:mv3}
L^2(X_1) \hatotimes L^2(X_2) \hatotimes \dots \hatotimes L^2(X_d) = L^2(X_1 \times X_2 \times \dots \times X_d).
\end{equation}
We will furnish more details in Examples~\ref{eg:norms} and \ref{eg:Hilbert}. 

We now attempt to summarize the discussions up to this point with a definition that is intended to capture all examples in this section. Note that it is difficult to be more precise without restricting the scope.

\begin{definition}[tensors as multivariate functions]\label{def:tensor3a}
Let $\Fn_i(X_i) \subseteq \mathbb{R}^{X_i}$ be a subspace of real-valued functions on  $X_i$, $i =1,\ldots,d$. The topological tensor product of $\Fn_1(X_1),\ldots, \Fn_d(X_d)$ is a subspace
\[
\Fn_1(X_1) \hatotimes \dots \hatotimes \Fn_d(X_d) \subseteq \mathbb{R}^{X_1 \times \dots \times X_d}
\]
that comprises finite sums of real-valued separable functions $\varphi_1 \otimes \dots \otimes \varphi_d$ with $\varphi_i \in \Fn_i(X_i)$, $i =1,\ldots,d$, and their limits with respect to some topology. A $d$-tensor is  an element of this space.
\end{definition}

Some explanations are in order. Definition~\ref{def:tensor3a} includes the special cases of functions on finite sets and multivariate polynomials by simply taking the topology to be the discrete topology. It is not limited  to norm topology, as numerous tensor product constructions require seminorms \cite{Treves} or quasinorms \cite{Sickel}. The definition allows for tensor products of different types of function spaces so as to accommodate different regularities in different arguments. For example, a classical solution for the fluid velocity $v$ in the Navier--Stokes equation \cite{Fefferman}
\begin{equation}\label{eq:navier}
\dfrac{\partial v_i}{\partial t} +\sum_{j=1}^ 3 \dfrac{\partial v_i}{\partial x_j}v_j= -\dfrac{1}{\rho}\dfrac{\partial p}{\partial x_i} + \nu\sum_{j=1}^{3}\dfrac{\partial^2 v_i}{\partial x_j^2} +f_i, \quad i =1,2,3,
\end{equation}
may be regarded as a tensor
$v \in C^2(\mathbb{R}^3) \hatotimes C^1 [0,\infty) \otimes \mathbb{R}^3$, \tie, twice continuously differentiable in its spatial arguments, once in its temporal argument, with a discrete third argument capturing the fact that it is a three-dimensional vector field. We will have more to say about the last point in Example~\ref{eg:vecvalobj}.

Definition~\ref{def:tensor3a} allows us to view numerous types of multivariate functions as tensors. The variables involved may be continuous or discrete or a mix of both. The functions may represent anything from solutions to PDEs to target functions in machine learning. In short, Definition~\ref{def:tensor3a} is extremely versatile and is by far the most common manifestation of tensors.

We will reflect on how Definition~\ref{def:tensor3a} relates to definitions~\ref{st:tensor1} and~\ref{st:tensor2} by way of two examples.

\begin{example}[hypermatrices]\label{eg:hyp}
A hypermatrix has served as our `multi-indexed object' in definition~\ref{st:tensor1} but may now be properly defined as an element of $\mathbb{R}^{X_1 \times \dots \times X_d}$ when $X_1,\ldots,X_d$ are finite or countable  discrete sets. For any $n \in \mathbb{N}$, we will write
\[
[n] \coloneqq \{1,2,\ldots, n\}.
\]
We begin from $d =1$. While $\mathbb{R}^n$ is usually taken as a Cartesian product
$\mathbb{R}^n = \mathbb{R} \times  \dots \times \mathbb{R}$,
it is often fruitful to regard it  as the set of functions
\begin{equation}\label{eq:cn}
\mathbb{R}^{[n]} = \{f \colon  [n] \to \mathbb{R} \}.
\end{equation}
Any $n$-tuple $(a_1,\ldots,a_n) \in \mathbb{R}^n$ defines a function $f \colon  [n] \to \mathbb{R}$ with $f(i) = a_i$, $i=1,\ldots,n$, and any function $f \in \mathbb{R}^{[n]}$ defines an $n$-tuple $(f(1), \ldots, f(n)) \in \mathbb{R}^n$. So we may as well assume that  $\mathbb{R}^n \equiv \mathbb{R}^{[n]}$, given that any difference is superficial. There are two advantages. The first is that it allows us to speak of sequence spaces in the same breath, as the space of all sequences is then
\[
\mathbb{R}^{\mathbb{N}} = \{f \colon \mathbb{N} \to \mathbb{R} \} = \{(a_i)_{i=1}^\infty \colon  a_i \in \mathbb{R} \},
\]
and this in turn allows us to treat various Banach spaces of sequences $l^p(\mathbb{N})$, $c(\mathbb{N})$, $c_0(\mathbb{N})$ as special cases of their function spaces counterpart $L^p(\Omega)$, $C(\Omega)$, $C_0(\Omega)$. This perspective is of course common knowledge in real and functional analysis.

Given that there is no such thing as a `two-way Cartesian product', a second advantage is that it allows us to define matrix spaces via
\begin{equation}\label{eq:cmn}
\mathbb{R}^{m \times n} \coloneqq \mathbb{R}^{[m] \times [n]} = \{f \colon  [m] \times [n] \to \mathbb{R} \},
\end{equation}
where we regard  $(a_{ij})_{i,j=1}^{m,n}$ as shorthand for the function with $f(i,j) =  a_{ij}$, $i \in [m]$, $j \in [n]$. Why not just say that it is a two-dimensional array? The answer is that an `array' is an undefined term in mathematics; at best we may view it as a type of data structure, and as soon as we try to define it rigorously, we will end up with essentially \eqref{eq:cmn}. Again this perspective is known; careful treatments of matrices in linear algebra such as \citet[Definition~4.1.3]{Berberian} or \citet{Brualdi} would define a matrix in the same vein. An infinite-dimensional variation gives us double sequences or equivalently infinite-dimensional matrices in $l^2(\mathbb{N} \times \mathbb{N})$. These are essential in Heisenberg's matrix mechanics, the commutation relation $PQ-QP =i \hbar I$ being false for finite-dimensional matrices (trace zero on the left, non-zero on the right), requiring infinite-dimensional ones instead:
\begin{align*}
Q &= \sqrt{\dfrac{\hbar}{2\nu \omega}}
\begin{bmatrix}
0 & \sqrt{1} & 0 & 0   & \cdots \\
\sqrt{1} & 0 & \sqrt{2} & 0  & \cdots \\
0 & \sqrt{2} & 0 & \sqrt{3} &  \cdots \\
0 & 0 & \sqrt{3} & 0  & \cdots \\
\vdots & \vdots & \vdots &  \vdots & \ddots \\
\end{bmatrix}\!,\\*
P &= \sqrt{\dfrac{\hbar\nu \omega}{2}}
\begin{bmatrix}
0 & -i\sqrt{1} & 0 & 0 &  \cdots \\
i\sqrt{1} & 0 & -i\sqrt{2} & 0 &  \cdots \\
0 & i\sqrt{2} & 0 & -i\sqrt{3} &  \cdots \\
0 & 0 & i\sqrt{3} & 0 &  \cdots\\
\vdots & \vdots & \vdots &  \vdots & \ddots \\
\end{bmatrix}\!.
\end{align*}
$P$ and $Q$ are in fact the position and momentum operators for the harmonic oscillator in Example~\ref{eg:Gaussian} restricted to one dimension with all $y$- and $z$-terms dropped \cite[Chapter~18]{Jordan}.

It is straightforward to generalize. Let $n_1,\ldots, n_d \in \mathbb{N}$. Then
\begin{equation}\label{eq:cn1nd}
\mathbb{R}^{n_1 \times \dots \times n_d} \coloneqq \mathbb{R}^{[n_1] \times \dots \times [n_d]}
= \{f \colon  [n_1] \times \dots \times [n_d] \to \mathbb{R}\},
\end{equation}
\tie, a $d$-dimensional hypermatrix or more precisely an $n_1 \times \dots \times n_d$ hypermatrix is a real-valued function on $[n_1] \times \dots \times [n_d]$.  We introduce the shorthand $(a_{i_1\cdots i_d})_{i_1, \ldots, i_d=1}^{n_1,\ldots, n_d}$ for a hypermatrix with $f(i_1,\ldots, i_d) = a_{i_1  \cdots i_d}$. We say that a hypermatrix is \emph{hypercubical} when $n_1 = \dots = n_d$. These terms originated in the Russian mathematical literature, notably in \citet*{GKZpaper,GKZ}, although the notion itself has appeared implicitly as coefficients of multilinear forms in \citet{Cayley}, alongside the
  term `hyperdeterminant'. The names `multidimensional matrix' \cite{Sokolov2} or `spatial matrix' \cite{Sokolov1} have also been used for hypermatrices. We would have preferred simply `matrix' or `$d$-matrix' instead of `hypermatrix' or `$d$-hypermatrix', but the  conventional view of the term `matrix' as exclusively  two-dimensional has become too deeply entrenched to permit such a change in nomenclature.

We will allow for a further generalization of the above, namely, $X_1,\ldots,X_d$ may be any discrete {sets.}
\begin{enumerate}[\upshape (a)]
\setlength\itemsep{3pt}
\item Instead of restricting $i_1,\ldots,i_d$ to ordinal variables with values in $[n_k] = \{1,\ldots,n_k\}$, we allow for them to be nominal, \eg\  $\{\uparrow, \downarrow\}$, $\{\mathsf{A},\mathsf{C},\mathsf{G},\mathsf{T}\}$, or irreducible representations of some group (see Example~\ref{eg:neural}).
\item The discrete set may have additional structures, \eg\  they might be graphs or finite groups or combinatorial designs.
\item We permit sets of countably infinite cardinality, \eg\  $\mathbb{N}$ or $\mathbb{Z}$ (see Example~\ref{eg:Calderon}).
\end{enumerate}
Readers may have noted that we rarely invoke a $d$-hypermatrix with $d \ge 3$ in actual examples. The reason is that whatever we can do with a  $d$-hypermatrix $A \in \mathbb{R}^{n_1 \times \dots \times n_d}$
we can do without or do better with $f \colon  [n_1] \times \dots \times [n_d] \to \mathbb{R}$. If anything, hypermatrices give a misplaced sense of familiarity by furnishing a false analogy with standard two-dimensional matrices;  they are nothing alike, as we explain next.

In conjunction with specific bases, \tie, when a $d$-hypermatrix is a coordinate representation of a $d$-tensor, the values $a_{i_1 \cdots i_d}$ convey information \emph{about those basis vectors} such as  near orthogonality (Example~\ref{eg:Calderon}) or incidence (below).
But on their own, the utility of $d$-hypermatrices is far from that of the usual two-dimensional matrices, one of the most powerful and universal mathematical tools. While there is a dedicated calculus\footnote{In the sense of a system of calculation and reasoning.} for working with two-dimensional matrices  --  a rich collection of matrix operations, functions, invariants, decompositions, canonical forms, {\em etc.}, that are consistent with the tensor transformation rules and most of them readily computable   --  there is no equivalent for $d$-hypermatrices. There are no canonical forms when $d \ge 3$ \cite[Chapter~10]{Land1}, no analogue of matrix--matrix product when $d$ is odd (Example~\ref{eg:homult}) and thus no analogue of matrix inverse, and almost every calculation is computationally intractable even for small examples (see Section~\ref{sec:comp}); here `no' is in the sense of mathematically proved to be non-existent. The most useful aspect of hypermatrices is probably in supplying some analogues of basic block operations, notably in  Strassen's laser method \cite[Chapter~15]{BCS}, or of row and column operations, notably in various variations of slice rank \cite{Lovett}, although there is nothing as sophisticated as the kind of block matrix operations on page~\pageref{pg:blockmatrixops}.

Usually when we recast a problem in the form of a matrix problem, we are on our way to a solution:  matrices are the tool that gets us there. The same is not true with hypermatrices. For instance, while we may easily capture the adjacency structure of a $d$-hypergraph $G =(V, E)$, $E \subseteq \binom{V}{d}$, with a hypermatrix  $f \colon  V\times\dots \times V \to \mathbb{R}$, 
\[
f(i_1,\ldots,i_d) \coloneqq \begin{cases} 1 & \{i_1,\ldots,i_d\} \in E, \\ 0 & \{i_1,\ldots,i_d\} \notin E, \end{cases}
\]
the analogy with spectral graph theory stops here; the hypermatrix view  \cite{Friedman1,Friedman2} does not get us anywhere close to what is possible in the $d=2$ case. Almost every mathematically sensible problem involving hypermatrices is a challenge, even for $d=3$ and very small dimensions such as $4 \times 4 \times 4$. For example, we know that the $m \times n$ matrices of rank not more than $r$ are precisely the ones with vanishing $(r+1) \times (r+1)$ minors. Resolving the equivalent problem for $4 \times 4 \times 4$ hypermatrices of rank not more than $4$ requires that we first reduce it to a series of questions about matrices and then throw every tool in the matrix arsenal at them \cite{Friedland1,Friedland2}.
\end{example}

The following resolves a potential point of confusion that the observant reader might have already noticed. In Definition~\ref{def:tensor2a}, a $d$-tensor is a multilinear functional, \tie,
\begin{equation}\label{eq:mm1}
\varphi \colon  \mathbb{V}_1 \times \dots \times \mathbb{V}_d \to \mathbb{R}
\end{equation}
satisfying \eqref{eq:multilinear}, but in Definition~\ref{def:tensor3a} a $d$-tensor is merely a multivariate function
\begin{equation}\label{eq:mm2}
f \colon  X_1 \times \dots \times X_d \to \mathbb{R}
\end{equation}
that is not required to be multilinear, a requirement that in fact makes no sense as the sets $X_1,\ldots, X_d$ may not be vector spaces. One might think that \eqref{eq:mm1} is a special case of \eqref{eq:mm2} as a multilinear functional is a special case of a multivariate function, but this would be incorrect:  they are different types of tensors. By Definition~\ref{def:tensor2a}, $\varphi$ is a covariant tensor whereas, as we pointed out in Example~\ref{eg:distributions}, $f$ is a contravariant tensor. It suffices to explain this issue over finite dimensions and so we will assume that $ \mathbb{V}_1, \ldots, \mathbb{V}_d$ are finite-dimensional vector spaces and $X_1,\ldots,X_d$ are finite sets. We will use slightly different notation below to avoid having to introduce more indices.

\begin{example}[multivariate versus multilinear]\label{eg:multmult}\hspz%
The relation between \eqref{eq:mm1} and \eqref{eq:mm2} is subtler than meets the eye. Up to a choice of bases, we can construct a unique $\varphi$ from any $f$ and {\em vice versa}. One direction is easy. Given vector spaces $\mathbb{U},\mathbb{V},\ldots,\mathbb{W}$, choose any bases
\[
\mathscr{A}=\{u_1,\ldots,u_m\}, \quad \mathscr{B}=\{v_1,\ldots,v_n\},\ldots,\mathscr{C} = \{w_1,\ldots, w_p\}.
\]
Then any multilinear functional $\varphi \colon  \mathbb{U} \times \mathbb{V} \times \dots\times\mathbb{W} \to \mathbb{R}$ gives us a real-valued function
\[
f_\varphi \colon  \mathscr{A} \times \mathscr{B} \times \dots \times \mathscr{C} \to \mathbb{R}, \quad (u_i, v_j, \ldots, w_k) \mapsto  \varphi(u_i, v_j, \ldots, w_k).
\]
The other direction is a bit more involved. For any finite set $X$ and a point $x \in X$, we define the delta function $\delta_x \in \mathbb{R}^X$ given by
\begin{alignat}{3}
\delta_x \colon  X &\to \mathbb{R}, &\quad\delta_x(x')&=\begin{cases} 1 & x = x',\\ 0 & x \ne x', \end{cases}\label{eq:delta}
\intertext{for all $x' \in X$, and the \emph{point evaluation} linear functional $\varepsilon_x \in (\mathbb{R}^X)^*$ given by}
\varepsilon_x \colon  \mathbb{R}^X &\to \mathbb{R},&\varepsilon_x(f) &=f(x)\notag
\end{alignat}
for all $f \in  \mathbb{R}^X$. These are bases of their respective spaces and for any set we have the following:
\[
X = \{x_1,\ldots,x_m\},\quad
\mathbb{R}^X = \spn \{\delta_{x_1},\ldots,\delta_{x_m}\},\quad
(\mathbb{R}^X)^* = \spn \{\varepsilon_{x_1},\ldots,\varepsilon_{x_m}\}.
\]
Given any $d$ sets,
\[
X = \{x_1,\ldots,x_m\}, \quad Y=\{y_1,\ldots,y_n\},\ldots,Z=\{z_1,\ldots,z_p\},
\]
a real-valued function $f \colon  X \times Y \times \dots \times Z \to \mathbb{R}$ gives us a multilinear functional
\begin{equation*}
\varphi_f \colon  \mathbb{R}^X \times \mathbb{R}^Y \times \dots \times \mathbb{R}^Z \to \mathbb{R}
\end{equation*}
defined by
\begin{equation*}
\varphi_f = \sum_{i=1}^m \sum_{j=1}^n\cdots\sum_{k=1}^p f(x_i,y_j, \ldots, z_k) \varepsilon_{x_i} \otimes \varepsilon_{y_j} \otimes \dots \otimes \varepsilon_{z_k}.
\end{equation*}
Evaluating $\varphi_f$ on basis vectors, we get
\begin{align*}
\varphi_f(\delta_{x_i},\delta_{y_j},\ldots,\delta_{z_k}) &= f(x_i,y_j, \ldots, z_k),
\intertext{taken together with the earlier definition of $f_\varphi$,}
f_\varphi(u_i,v_j,\ldots,w_k) &= \varphi(u_i,v_j,\ldots,w_k),
\end{align*}
and we obtain
\[
\varphi_{f_\varphi} = \varphi, \quad f_{\varphi_f} = f.
\]
Roughly speaking, a real-valued multivariate function is a tensor in the sense of Definition~\ref{def:tensor2a} because one may view it as a linear combination of point evaluation linear functionals. Such point evaluation functionals are the basis behind numerical quadrature (Example~\ref{eg:quadrature}).
\end{example}

The next example is intended to serve as further elaboration for Examples~\ref{eg:Gaussian}, \ref{eg:hyp} and \ref{eg:multmult}, and as partial impetus for Sections~\ref{sec:tensor3b} and \ref{sec:tensor3c}. 

\begin{example}[quantum spin and covalent bonds]\label{eg:spin}
As we saw  in Example~\ref{eg:Gaussian}, the quantum state of a spinless particle like the Gaussian wave function $\psi_{m,n,p}$ is an element of the Hilbert space\footnote{To avoid differentiability and  integrability issues, we need $\psi_{m,n,p}$ to be in the Schwartz space $S(\mathbb{R}^3) \subseteq L^2(\mathbb{R}^3)$ if we take a rigged Hilbert space approach \cite{rigged} or the Sobolev space $H^2(\mathbb{R}^3) \subseteq L^2(\mathbb{R}^3)$ if we take an unbounded self-adjoint operators approach \cite[Chapter~7]{Teschl}, but we disregard these to keep things simple.} $L^2(\mathbb{R}^3)$. But even for a single particle, tensor products come in handy when we need to discuss particles with spin, which are crucial in chemistry. The simplest but also the most important case is when a quantum particle has spin $-\frac{1}{2}, \frac{1}{2}$, called a \emph{spin-half particle} for short, and  its quantum state is
\begin{equation}\label{eq:1electron}
\Psi \in L^2(\mathbb{R}^3) \otimes \mathbb{C}^2.
\end{equation}
Going by Definition~\ref{def:tensor3a}, this is a function $\Psi\colon  \mathbb{R}^3 \times \{-\frac{1}{2}, \frac{1}{2}\} \to \mathbb{C}$, where $\Psi(x,\sigma)$ is an $L^2(\mathbb{R}^3)$-integrable function in the argument $x$ for each  $\sigma = -\frac{1}{2}$ or $\sigma = \frac{1}{2}$.
It is also not hard to see that when we have a tensor product of an
infinite-dimensional space with a finite-dimensional one, every element must be a \emph{finite} sum of  separable terms, so $\Psi$ in \eqref{eq:1electron} may be expressed as
\begin{equation}\label{eq:1electron1}
\Psi = \sum_{i=1}^r \psi_i \otimes \chi_i\quad\text{or}\quad \Psi(x,\sigma) = \sum_{i=1}^r \psi_i (x) \chi_i(\sigma)
\end{equation}
with $\psi_i \in L^2(\mathbb{R}^3)$ and $\chi_i \colon  \{-\frac{1}{2}, \frac{1}{2}\} \to \mathbb{C}$. In physics parlance, $\Psi$ is the \emph{total wave function} of the particle and $\psi_i$ and $\chi_i$ are its spatial and spin wave functions; the variables that $\Psi$ depends on are called \emph{degrees of freedom}, with those describing position and momentum called \emph{external} degrees of freedom and others like spin called \emph{internal} degrees of freedom \cite[Chapter~II, Section~F]{Cohen1}.

This also illustrates why, as we discussed in Example~\ref{eg:hyp}, it is often desirable to view $\mathbb{C}^n$ as the set of complex-valued functions on some finite sets of $n$ elements. Here the finite set is $\{-\frac{1}{2}, \frac{1}{2}\}$ and $\mathbb{C}^2$ in \eqref{eq:1electron} is just shorthand for the \emph{spin state space}
\[
\mathbb{C}^{ \{-1/2, 1/2\} } = \biggl\{\chi \colon  \biggl\{-\frac{1}{2}, \frac{1}{2}\biggr\} \to \mathbb{C}\biggr\}.
\]
The domain  $ \mathbb{R}^3 \times \{-\frac{1}{2}, \frac{1}{2}\} $ is called the \emph{position--spin space} \cite{Pauli} and plays a role in the Pauli exclusion principle. While two classical particles cannot simultaneously occupy the exact same location in $\mathbb{R}^3$, two quantum particles can as long as they have different spins, as that means they are occupying different locations in $ \mathbb{R}^3 \times \{-\frac{1}{2}, \frac{1}{2}\} $. A consequence is that two electrons with opposite spins can occupy the same molecular orbital, and when they do, we have a covalent bond in the molecule. The Pauli exclusion principle implies the converse: if two electrons occupy the same molecular orbital, then they must have opposite spins. We will see that this is a consequence of the antisymmetry of the total wave function.

In fact, `two electrons with opposite spins occupying the same molecular orbital' is the quantum mechanical definition of a covalent bond in chemistry. When this happens, quantum mechanics mandates that we may not speak of each electron individually but need to consider both as a single entity described by a single wave function
\begin{equation}\label{eq:2electrons}
\Psi \in ( L^2(\mathbb{R}^3) \otimes \mathbb{C}^2) \hatotimes ( L^2(\mathbb{R}^3) \otimes \mathbb{C}^2),
\end{equation}
with each copy of $ L^2(\mathbb{R}^3) \otimes \mathbb{C}^2$ associated with one of the electrons. The issue with writing the tensor product in the form \eqref{eq:2electrons} is that $\Psi$ will not be a finite sum of separable terms since it is not a tensor product of two infinite-dimensional spaces. As we will see in Section~\ref{sec:tensor3c}, we may arrange the factors in a tensor product in arbitrary order and obtain isomorphic spaces
\begin{equation}\label{eq:2electrons1}
( L^2(\mathbb{R}^3) \otimes \mathbb{C}^2) \hatotimes ( L^2(\mathbb{R}^3) \otimes \mathbb{C}^2) \cong L^2(\mathbb{R}^3)  \hatotimes L^2(\mathbb{R}^3) \otimes \mathbb{C}^2 \otimes \mathbb{C}^2,
\end{equation}
but this is also obvious from Definition~\ref{def:tensor3a} by observing that for finite sums
\[
\sum_{i=1}^r \psi_i (x) \chi_i(\sigma)\varphi_i (y) \xi_i(\tau) = \sum_{i=1}^r \psi_i (x) \varphi_i (y) \chi_i(\sigma)\xi_i(\tau),
\]
and then taking completion.
Since $ L^2(\mathbb{R}^3) \hatotimes L^2(\mathbb{R}^3) = L^2(\mathbb{R}^3  \times \mathbb{R}^3)$ by \eqref{eq:mv3} and $\mathbb{C}^2 \otimes \mathbb{C}^2  = \mathbb{C}^{2 \times 2}$ by Example~\ref{eg:hyp}, again the latter is shorthand for
\[
\mathbb{C}^{ \{-1/2, 1/2\} \times \{-1/2, 1/2\} }
=  \biggl\{\chi \colon  \biggl\{\biggl(-\frac{1}{2},-\frac{1}{2}\biggr),\biggl(-\frac{1}{2},\frac{1}{2}\biggr),\biggl(\frac{1}{2},-\frac{1}{2}\biggr),\biggl(\frac{1}{2},\frac{1}{2}\biggr)\biggr\} \to \mathbb{C}\biggr\},
\]
we deduce that the state space in \eqref{eq:2electrons1} is $L^2(\mathbb{R}^6) \otimes \mathbb{C}^{2 \times 2}$. The total wave function for two electrons  may now be expressed as a finite sum of separable terms
\[
\Psi(x,y,\sigma,\tau) = \sum_{i=1}^r \varphi_i(x,y) \chi_i(\sigma,\tau),
\]
nicely separated into external degrees of freedom $(x,y) \in \mathbb{R}^3 \times \mathbb{R}^3$ and internal ones $(\sigma,\tau)$ as before. Henceforth we will write $\Psi(x,\sigma; y,\tau)$ if we view $\Psi$ as an element of the space on the left of \eqref{eq:2electrons1} and $\Psi(x,y,\sigma,\tau)$ if we view it as an element of the space on the right. The former emphasizes that $\Psi$ is a wave function of two electrons, one with degrees of freedom $(x,\sigma)$ and the other  $(y,\tau)$; the latter emphasizes the separation of the wave function into external and internal degrees of freedom.

Electrons are \emph{fermions}, which are quantum particles with antisymmetric wave functions, \tie,
\begin{equation}\label{eq:antisymm1}
\Psi(x,\sigma; y,\tau) = -\Psi(y,\tau; x, \sigma).
\end{equation}
For simplicity, suppose we have a separable wave function
\[
\Psi(x,\sigma; y,\tau) = \Psi(x,y,\sigma,\tau) = \varphi(x,y) \chi(\sigma,\tau).
\]
Aside from their spins, identical electrons in the same orbital are indistinguishable and we must have $\varphi(x,y) = \varphi(y,x)$, so the antisymmetry \eqref{eq:antisymm1} must be a consequence of  $\chi(\sigma,\tau) = -\chi(\tau,\sigma)$.  As in Example~\ref{eg:multmult}, a basis for $\mathbb{C}^2$ is given by the delta functions $\delta_\mh $ and $\delta_\ph $ and thus a basis for $\mathbb{C}^{2\times 2}$ is given by
\[
\delta_\mh \otimes \delta_\mh , \quad
\delta_\mh \otimes \delta_\ph , \quad 
\delta_\ph \otimes \delta_\mh , \quad 
\delta_\ph \otimes \delta_\ph .
\] 
We may thus express
\[
\chi = a \delta_\mh \otimes \delta_\mh  + b \delta_\mh \otimes \delta_\ph  + c\delta_\ph \otimes \delta_\mh  +d \delta_\ph \otimes \delta_\ph  ,
\]
and if $\chi(\sigma,\tau) = -\chi(\tau,\sigma)$ for all $\sigma, \tau$, then it implies that $a = d = 0$ and $b = -c$. Choosing $b = 1/\sqrt{2} $ to get unit norms, we have
\[
\chi = \dfrac{1}{\sqrt{2}} ( \delta_\mh \otimes \delta_\ph  -\delta_\ph \otimes \delta_\mh  ).
\]
In other words, the two electrons in the same molecular orbital must have opposite spins  --  either $(-\frac{1}{2},\frac{1}{2})$ or $(\frac{1}{2}, -\frac{1}{2})$.

The discussion above may be extended to more particles ($d > 2$) and more spin values ($s > \frac{1}{2}$). When we have $d$ quantum particles with spins $-\frac{1}{2}, \frac{1}{2}$, the total wave function satisfies
\begin{equation}\label{eq:delectrons}
\Psi \in (L^2(\mathbb{R}^3) \otimes \mathbb{C}^2) \hatotimes \dots  \hatotimes (L^2(\mathbb{R}^3) \otimes \mathbb{C}^2) \cong  L^2(\mathbb{R}^{3d}) \otimes  \mathbb{C}^{2 \times \dots \times 2}
\end{equation}
via the same argument. Again we may group the degrees of freedom either by particles or by external/internal, \tie,
\[
\Psi(x_1, \sigma_1; x_2, \sigma_2; \dots; x_d, \sigma_d) = \Psi(x_1, x_2, \ldots, x_d, \sigma_1, \sigma_2,\ldots, \sigma_d),
\]
depending on whether we want to view $\Psi$ as an element of the left or right space in \eqref{eq:delectrons}, with antisymmetry usually expressed in the former,
\begin{equation}\label{eq:fermion}
\Psi(x_1, \sigma_1; \dots; x_d, \sigma_d) = (-1)^{\sgn(\pi)} \Psi(x_{\pi(1)}, \sigma_{\pi(1)}; \dots; x_{\pi(d)}, \sigma_{\pi(d)})
\end{equation}
for all $\pi \in \mathfrak{S}_d$, and separability usually expressed in the latter,
\[
 \Psi(x_1, \ldots, x_d, \sigma_1, \ldots, \sigma_d) = \sum_{i=1}^r \psi_i(x_1,\ldots,x_d) \chi_i(\sigma_1,\ldots,\sigma_d).
\]
As before, $\mathbb{C}^{2 \times \dots \times 2}$ is to be regarded as
\[
\mathbb{C}^{ \{-1/2, 1/2\} \times \dots \times \{-1/2, 1/2\}  }
=  \biggl\{\chi \colon  \biggl\{-\frac{1}{2},\frac{1}{2}\biggr\}^n \to \mathbb{C}\biggr\}.
\]

We may easily extend our discussions to a spin-$s$ quantum particle for any
\[
s \in \biggl\{0, \frac{1}{2}, 1, \frac{3}{2}, 2,\frac{5}{2},\ldots\biggr\},
\]
 in which case its state space would be
\[
L^2(\mathbb{R}^3) \otimes \mathbb{C}^{2s+1},
\]
where $\mathbb{C}^{2s+1}$ is shorthand for the space of functions $\chi \colon  \{-s, -s+1,\ldots, s-1, s\} \to \mathbb{C}$. The particle is called a fermion or a \emph{boson} depending on whether $s$ is a half-integer or an integer. The key difference is that the total wave function of $d$ identical fermions is antisymmetric as in \eqref{eq:fermion} whereas that of $d$ identical bosons is \emph{symmetric}, \tie,~
\begin{equation}\label{eq:boson}
\Psi(x_1, \sigma_1; \dots; x_d, \sigma_d) = \Psi(x_{\pi(1)}, \sigma_{\pi(1)}; \dots; x_{\pi(d)}, \sigma_{\pi(d)})
\end{equation}
for all $\pi \in \mathfrak{S}_d$. Collectively, \eqref{eq:fermion} and \eqref{eq:boson} are the modern form of the Pauli exclusion principle. Even more generally, we may consider $d$ quantum particles with different spins $s_1,\ldots,s_d$ and state space
\[
(L^2(\mathbb{R}^3) \otimes \mathbb{C}^{2s_1 +1}) \hatotimes \dots  \hatotimes (L^2(\mathbb{R}^3) \otimes \mathbb{C}^{2s_d + 1}) \cong  L^2(\mathbb{R}^{3d}) \otimes  \mathbb{C}^{(2s_1 + 1) \times \dots \times (2s_d +1)}.
\]
We recommend the very lucid exposition in \citet[Sections~46--48]{Faddeev} and \citet[Chapter~4]{Takhtajan} for further information.

For most purposes in applied and computational mathematics, $s=1/2$, \tie, the spin-half case  is literally the only one that matters, as particles that constitute ordinary matter (electrons, protons, neutrons, quarks, {\em etc.}) are all spin-half fermions. Spin-zero (\eg\  Higgs) and spin-one bosons (\eg\  W and Z, photons, gluons) are already the realm of particle physics, and higher-spin particles, whether fermions or bosons, tend to be hypothetical  (\eg\  gravitons, gravitinos).
\end{example}

In Example~\ref{eg:spin} we have implicitly relied on the following three background assumptions for a physical system in quantum mechanics.\label{pg:postulates}
\begin{enumerate}[\upshape (a)]
\setlength\itemsep{3pt}
\item\label{it:state} The state space of a  system is a Hilbert space.
\item\label{it:composite} The state space of a composite system is the tensor product of the state spaces of the component systems.
\item\label{it:iden} The state space of a composite system of identical particles lies either in the symmetric or alternating tensor product of the state spaces of the component systems.
\end{enumerate}
These assumptions are standard; see \citet[Chapter~III, Section~B1]{Cohen1}, \citet[Section~2.2.8]{Nielsen} and \citet[Chapter~XIV, Section~C]{Cohen2}.

Example~\ref{eg:spin} shows the simplicity of Definition~\ref{def:tensor3a}. If our quantum system has yet other degrees of freedom, say, flavour or colour, we just include the corresponding wave functions as factors in the separable product. For instance, the total wave function of a proton or neutron  is a finite sum of separable terms of the form \cite[Section~5.6.1]{Griffiths}
\[
\psi_{\text{spatial}} \otimes \psi_{\text{spin}} \otimes \psi_{\text{flavour}} \otimes \psi_{\text{colour}}.
\]

Nevertheless, we would like more flexibility than Definition~\ref{def:tensor3a} provides. For instance, instead of writing the state of a spin-half particle as a sum of real-valued functions in \eqref{eq:1electron1}, we might prefer to view it as a $\mathbb{C}^2$-valued $L^2$-vector field on $\mathbb{R}^3$ indexed by spin:
\begin{equation}\label{eq:l2vecfield}
  \Psi(x) = \begin{bmatrix} \psi_\mh(x) \\ \psi_\ph(x) \end{bmatrix}
\end{equation}
with
\[
\lVert \Psi \rVert^2 = \int_{\mathbb{R}^3} \lvert \psi_\mh(x) \rvert^2 \D x +  \int_{\mathbb{R}^3} \lvert \psi_\ph(x) \rvert^2 \D x < \infty. 
\]
In other words, we want  to allow  $L^2(\mathbb{R}^3; \mathbb{C}^2)$ to be a possible interpretation for $L^2(\mathbb{R}^3) \otimes \mathbb{C}^2$. The next two sections will show how to accomplish this.

\subsection{Tensor products of abstract vector spaces}\label{sec:tensor3b}

Most readers' first encounter with scalars and vectors will likely take the following forms:
\begin{itemize}
\setlength\itemsep{3pt}
\item a scalar is an object that only has magnitude,
\item a vector is an object with both magnitude and direction.
\end{itemize}
In this section we will construct tensors of arbitrary order starting from these. For a vector $v$ in this sense, we will denote its magnitude by $\lVert v \rVert$ and its direction by $\widehat{v}$. Note that here $\lVert \,\cdot\, \rVert$ does not need to be a norm, just some notion of magnitude. Likewise we will write $\lvert a \rvert \in [0,\infty)$  for the magnitude of a scalar $a$ and $\sgn(a)$ for its sign if it has one.

If we go along this line of reasoning, what should the next object be? One possible answer is that it ought to be an object with a magnitude and \emph{two} directions;  this is called a \emph{dyad} and we have in fact encountered it in Section~\ref{sec:multmaps}. Henceforth it is straightforward to generalize to objects with a magnitude and an arbitrary number of directions:
\begin{alignat*}{3}
&\textsc{object} &\ & \textsc{property}  \\*
&\text{scalar $a$} && \text{has magnitude $\lvert a \rvert$}\\
&\text{vector $v$} && \text{has magnitude $\lVert v \rVert$ and a direction $\widehat{v}$}\\
&\text{dyad $v \otimes w$} &&  \text{has magnitude $\lVert v \otimes w \rVert$ and two directions $\widehat{v}$, $\widehat{w}$}\\
&\text{triad $u \otimes v \otimes w$} &&  \text{has magnitude $\lVert u \otimes v \otimes w \rVert$ and three directions $\widehat{u}$, $\widehat{v}$, $\widehat{w}$}\\
&\quad\vdots & & \quad \vdots \\*
&\text{$d$-ad $u \otimes v \otimes \!\cdots\! \otimes w$} &&   \text{has magnitude $\lVert u \otimes v \otimes\!\cdots\!\otimes w \rVert$ and $d$ directions $\widehat{u},\widehat{v}, \ldots, \widehat{w}.$}
\end{alignat*}
Such objects are called \emph{polyads} and have been around for over a century;  the first definition of tensor rank was in fact in terms of polyads \cite{Hitch1}. We have slightly modified the notation by adding the tensor product symbol $\otimes$ between vectors to bring it in line with modern treatments in algebra, which we will describe later. Note that the $\otimes$ in  $v \otimes w$ is solely used as a delimiter; we could have written it as $vw$ as in \citet{Hitch1} and \citet{Morse} or $\lvert v \rangle\lvert w \rangle$ or $\lvert v, w \rangle$ in Dirac notation \cite[Chapter~II, Section~F2c]{Cohen1}. Indeed, $v$ and $w$ may not have coordinates and there is no `formula for $v\otimes w$'.  \label{pg:noformula}
Also, the order matters as $v$ and $w$ may denote vectors of different nature, so in general
\begin{equation}\label{eq:noncomm}
v \otimes w \ne w \otimes v.
\end{equation}
The reason for this would be clarified with Example~\ref{eg:stress}. Used in this sense, $\otimes$ is called an \emph{abstract tensor product} or \emph{dyadic product} in the older literature \cite{Chou,Morse,Tai}.

With hindsight, we see that the extension of a scalar in the sense of an object with magnitude, and a vector in the sense of one with magnitude and direction, is not a tensor but a rank-one tensor. To get tensors that are not rank-one we will have to look at the arithmetic of dyads and polyads. As we have discussed, each dyad is essentially a placeholder for three pieces of information:
\[
(\text{magnitude},\; \text{first direction},\; \text{second direction})
\]
or, in our notation $(\lVert v \otimes w \rVert,\widehat{v}, \widehat{w})$. We are trying to devise a consistent system of arithmetic for such objects that (a) preserves their information content, (b)~merges information whenever possible, and (c) is consistent with scalar and vector arithmetic.

Scalar products are easy. The role of scalars is that they scale vectors, as reflected in its name. We expect the same for dyads. If any of the vectors in a dyad is scaled by $a$, we expect its magnitude to be scaled by $a$ but the two directions to remain unchanged. So we require
\begin{equation}\label{eq:distrib1}
(av) \otimes w  = v\otimes (aw) \eqqcolon a v \otimes w,
\end{equation}
where the last term is defined to be the common value of the first two. This seemingly innocuous property is the main reason tensor products, not direct sums, are used to combine quantum state spaces. In quantum mechanics, a quantum state is described not so much by a vector $v$ but the entire one-dimensional subspace spanned by $v$; thus \eqref{eq:distrib1} ensures that when combining two quantum states in the form of two one-dimensional subspaces, it does not matter which (non-zero) vector in the subspace we pick to represent the state. On the other hand, for direct sums,
\begin{equation}\label{eq:nondistrib}
(av) \oplus w \ne v \oplus (aw) \ne a (v \oplus w)
\end{equation}
in general; for example,  with $v = (1,0) \in \mathbb{R}^2$, $w = 1 \in \mathbb{R}^1$ and $a = 2 \in \mathbb{R}$, we get $(2,0,1) \ne (1,0,2) \ne (2,0,2)$.

To ensure consistency with the usual scalar and vector arithmetic, one also needs the assumption that scalar multiplication is always distributive and associative in the following sense:
\[
( a + b) v \otimes w =  a v \otimes w +  b v \otimes w, \quad (ab) v \otimes w = a(bv \otimes w).
\]

Addition of dyads is trickier. If 
$a_1,a_2$ are scalars, then $a_1+ a_2$ is a scalar; if $v_1,v_2$ are vectors, then $v_1+ v_2$ is a vector according to either the parallelogram law of addition (for geometrical or physical vectors) or the axioms of vector spaces (for abstract vectors). But for dyads $v_1 \otimes w_1$ and  $v_2 \otimes w_2$, their sum $v_1 \otimes w_1+v_2 \otimes w_2$ is no longer a dyad;  in general there is no way to simplify the sum further and it should  be regarded as an object with two magnitudes and four
directions:\footnote{We could of course add corresponding vectors $v_1 + v_2$ and $w_1+w_2$, but that leads to direct sums of vector spaces, which is inappropriate because of \eqref{eq:nondistrib}.}
\begin{align*}
&(\text{magnitude 1},\; \text{first direction 1},\; \text{second direction 1}) \\*
&\quad \text{ \& } (\text{magnitude 2},\; \text{first direction 2},\; \text{second direction 2}).
\end{align*}
Just as the $\otimes$ is used as a delimiter where order matters, the $+$ in the sum of dyads is used as a different type of delimiter where order does not matter:
\[
v_1 \otimes w_1+v_2 \otimes w_2 = v_2 \otimes w_2  + v_1 \otimes w_1,
\]
\tie, $+$ is commutative. We also want $+$ to be associative so that we may unambiguously define, for any finite $r \in \mathbb{N}$, a sum of $r$ dyads
\[
v_1 \otimes w_1 + v_2 \otimes w_2 + \dots + v_r \otimes w_r
\]
that we will call a \emph{dyadic}. To ensure consistency with the usual scalar and vector arithmetic, we assume throughout this article that addition $+$ is always commutative and associative. This assumption is also perfectly reasonable as a dyadic is just a placeholder for the information contained in its dyad summands:
\begin{align*}
& (\text{magnitude 1},\; \text{first direction 1},\; \text{second direction 1}) \\*
&\quad \text{ \& } (\text{magnitude 2},\; \text{first direction 2},\; \text{second direction 2})
\text{ \& } \cdots\\*
&\quad \quad\text{ \& } (\text{magnitude }r,\; \text{first direction }r,\; \text{second direction }r).
\end{align*}
There is one scenario where a sum can be simplified and the directional information in two dyads merged. Whenever $v_1 = v_2$ or $w_1 = w_2$, we want
\begin{equation}\label{eq:distrib2}
\begin{aligned}
v \otimes w_1 + v \otimes w_2 &= v \otimes (w_1 + w_2), \\
v_1  \otimes w + v_2 \otimes w &= (v_1 + v_2) \otimes w,
\end{aligned}
\end{equation}
and these may be applied to combine any pair of summands in a dyadic that share a common factor. The following example provides further motivation for why we would like these arithmetical properties to hold.

\begin{example}[stress tensor]\label{eg:stress}
In physics, a dyadic is essentially a vector whose components are themselves vectors. The standard example is the \emph{stress} tensor, sometimes called the Cauchy stress tensor to distinguish it from other related models. The stress $\sigma$ at a point in space\footnote{We assume `space' here means $\mathbb{R}^3$ with  coordinate axes labelled $x,y,z$. The older literature would use $\mathrm{i}$, $\mathrm{j}$, $\mathrm{k}$ instead of $e_x$, $e_y$, $e_z$.} has three components in the directions given by unit vectors $e_x$, $e_y$, $e_z$:
\begin{equation}\label{eq:stress1}
\sigma = \sigma_x e_x + \sigma_y e_y + \sigma_z e_z.
\end{equation}
Now the deviation from a usual linear combination of vectors is that each of the coefficients $\sigma_x$, $\sigma_y$, $\sigma_z$ is not a scalar but a vector. This is the nature of stress:  in every direction, the stress in that direction has a normal component in that direction and two shear components in the plane perpendicular to it. For example,
\begin{equation}\label{eq:ex}
\sigma_x = \sigma_{xx} e_x + \sigma_{yx} e_y + \sigma_{zx} e_z,
\end{equation}
the component $\sigma_{xx}$ in the direction of $e_x$ is called \emph{normal stress} whereas the two components  $\sigma_{yx}$ and $\sigma_{zx}$ in the plane perpendicular to $e_x$, \ie\ $\spn \{e_y,e_z\}$, is called \emph{shear stress}. Here $\sigma_{xx}$, $\sigma_{yx}$, $\sigma_{zx}$ are scalars and \eqref{eq:ex} is an honest linear combination. It is easiest to represent $\sigma$ pictorially, as in \citet[Figure~4.10]{Borg} or \cite[Figure~1.5]{Chou},
adapted here as Figure~\ref{fig:stress}, where the point is blown up into a cube to show the fine details.
\begin{figure}[ht]
\centering
 \begin{tikzpicture}[%
        3d view = {28}{28}
    ]
    
        \draw[-{Latex[scale = 0.7]}] (-4, -4, 0) -- (-3, -4, 0)
            node[right] {\footnotesize$e_x$};
        \draw[-{Latex[scale = 0.7]}] (-4, -4, 0) -- (-4, -5.25, 0)
            node[left] {\footnotesize$e_y$};
        \draw[-{Latex[scale = 0.7]}] (-4, -4, 0) -- (-4, -4, 1)
            node[above] {\footnotesize$e_z$};

        \begin{scope}[canvas is xz plane at y = -2]
        
            \draw[
                fill = lightgray
            ] (-2, -2) rectangle (2, 2);
            
            \draw[-Latex] (-1.75, 0) -- (1.75, 0)
                node[below left] {\footnotesize$\sigma_{xy}$};
            \draw[-Latex] (0, -1.75) -- (0, 1.75)
                node[below left] {\footnotesize$\sigma_{zy}$};
        
        \end{scope}     
    
        \draw[-Latex] (0, -2, 0) -- (0, -4, 0)
            node[left] {\footnotesize$\sigma_{yy}$};
    
        \begin{scope}[canvas is yz plane at x = 2]
    
            \draw[
                fill = lightgray
            ] (-2, -2) rectangle (2, 2);
            
            \draw[-Latex] (1.75, 0) -- (-1.75, 0)
                node[below right] {\footnotesize$\sigma_{yx}$};
            \draw[-Latex] (0, -1.75) -- (0, 1.75)
                node[below right] {\footnotesize$\sigma_{zx}$};
    
        \end{scope} 

        \draw[-Latex] (2, 0, 0) -- (3.5, 0, 0)
            node[right] {\footnotesize$\sigma_{xx}$};

        \begin{scope}[canvas is xy plane at z = 2]

            \draw[
                fill = lightgray
            ] (-2, -2) rectangle (2, 2);

            \draw[-Latex] (-1.75, 0) -- (1.75, 0)
                node[above]{\footnotesize$\sigma_{xz}$};
            \draw[-Latex] (0, 1.75) -- (0, -1.75)
                node[above left]{\footnotesize$\sigma_{yz}$};

        \end{scope} 

        \draw[-Latex] (0, 0, 2) -- (0, 0, 3.5)
            node[above] {\footnotesize$\sigma_{zz}$};

    \end{tikzpicture}
\caption{Depiction of the stress tensor. The cube represents a single point in $\mathbb{R}^3$ and should be regarded as infinitesimally small. In particular, the origins of the three coordinate frames on the surface of the cube are the same point.}
\label{fig:stress}
\end{figure}
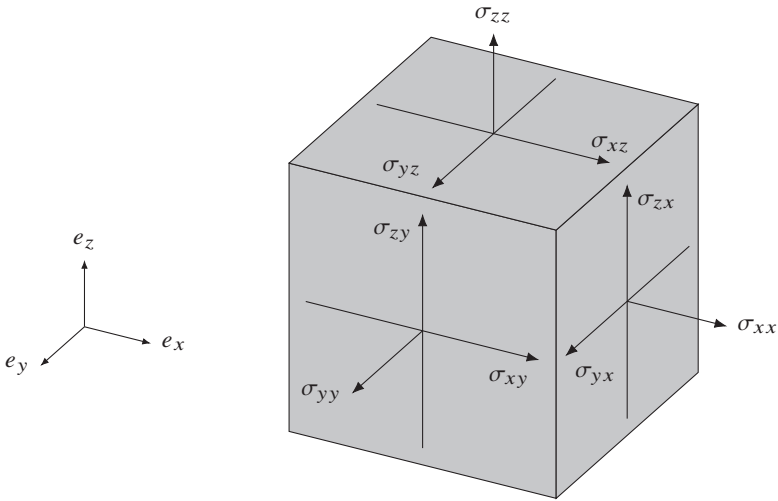
The coefficients $\sigma_x$,  $\sigma_y$, $\sigma_z$ are called \emph{stress vectors} in the directions $e_x$, $e_y$, $e_z$ respectively. As Figure~\ref{fig:stress} indicates, each stress vector has a normal component pointing in the same direction as itself and two shear components in the plane perpendicular to it:
\begin{align*}
\sigma_y &= \sigma_{xy} e_x + \sigma_{yy} e_y + \sigma_{zy} e_z,\\*
\sigma_z &= \sigma_{xz} e_x + \sigma_{yz} e_y + \sigma_{zz} e_z.
\end{align*}
Since the `coefficients' $\sigma_x$,  $\sigma_y$, $\sigma_z$ in the `linear combination' in \eqref{eq:stress1} are vectors, to conform to modern notation we just insert $\otimes$:
\[
\sigma = \sigma_x \otimes e_x + \sigma_y \otimes e_y + \sigma_z \otimes e_z.
\]
If we then plug the expressions for $\sigma_x$, $\sigma_y$, $\sigma_z$ into that of $\sigma$ and use  \eqref{eq:distrib2}, we~get
\begin{align}\label{eq:stress2}
\sigma & = \sigma_{xx} e_x \otimes e_x + \sigma_{yx} e_y \otimes  e_x + \sigma_{zx} e_z \otimes  e_x \notag \\*
&\quad +  \sigma_{xy} e_x \otimes  e_y + \sigma_{yy} e_y\otimes  e_y + \sigma_{zy} e_z \otimes  e_y \notag \\*
&\quad\quad + \sigma_{xz} e_x \otimes  e_z + \sigma_{yz} e_y \otimes  e_ z+ \sigma_{zz} e_z \otimes  e_z.
\end{align}
Without the $\otimes$, this is the expression of stress as a dyadic in \citet[pp.~70--71]{Morse}. This also explains why we want \eqref{eq:noncomm}, \tie, $\otimes$ should be non-commutative. For instance, $e_z \otimes  e_y$, the $z$-shear component for the $y$-normal direction, and $e_z \otimes  e_y$, the $y$-shear component for the $z$-normal direction mean completely different things.

If we set a basis $\mathscr{B} =\{e_x,e_y,e_z\}$, we may represent $\sigma$ as a $3\times 3$ matrix, \tie,
\[
\Sigma \coloneqq [\sigma]_{\mathscr{B}} = \begin{bmatrix}
\sigma_{xx} & \sigma_{xy} & \sigma_{xz} \\
\sigma_{yx} & \sigma_{yy} & \sigma_{yz} \\
\sigma_{zx} & \sigma_{zy} & \sigma_{zz} \\
\end{bmatrix}
\]
with normal stresses on the diagonal and shear stresses on the off-diagonal.
If we have a different basis $\mathscr{B}' =\{e'_x,e'_y,e'_z\}$ with a different representation
\[
\Sigma' \coloneqq [\sigma]_{\mathscr{B}'} =
\begin{bmatrix}
\sigma'_{xx} & \sigma'_{xy} & \sigma'_{xz} \\
\sigma'_{yx} & \sigma'_{yy} & \sigma'_{yz} \\
\sigma'_{zx} & \sigma'_{zy} & \sigma'_{zz} \\
\end{bmatrix}\!,
\]
then by plugging the change-of-basis relations
\[
\left\{
\begin{aligned}
e'_x &= c_{xx} e_x + c_{yx} e_y+c_{zx} e_z,\\
e'_y &= c_{xy} e_x + c_{yy} e_y+c_{zy} e_z,\\
e'_z &= c_{xz} e_x + c_{yz} e_y+c_{zz} e_z,
\end{aligned}\right. \quad
C \coloneqq \begin{bmatrix}
c_{xx} & c_{xy} & c_{xz} \\
c_{yx} & c_{yy} & c_{yz} \\
c_{zx} & c_{zy} & c_{zz} \\
\end{bmatrix}
\]
into
\begin{align*}
\sigma & = \sigma_{xx}' e_x' \otimes e_x' + \sigma_{yx}' e_y' \otimes  e_x' + \sigma_{zx}' e_z' \otimes e_x' \\*
&\quad +  \sigma_{xy}' e_x' \otimes  e_y' + \sigma_{yy}' e_y\otimes  e_y' + \sigma_{zy}' e_z' \otimes  e_y' \\*
&\quad\quad + \sigma_{xz}' e_x' \otimes  e_z' + \sigma_{yz}' e_y' \otimes  e_ z' + \sigma_{zz}' e_z' \otimes  e_z',
\end{align*}
applying \eqref{eq:distrib2} and comparing with \eqref{eq:stress2}, we get
\[
\begin{bmatrix}
\sigma_{xx} & \sigma_{xy} & \sigma_{xz} \\
\sigma_{yx} & \sigma_{yy} & \sigma_{yz} \\
\sigma_{zx} & \sigma_{zy} & \sigma_{zz} \\
\end{bmatrix}=
\begin{bmatrix}
c_{xx} & c_{xy} & c_{xz} \\
c_{yx} & c_{yy} & c_{yz} \\
c_{zx} & c_{zy} & c_{zz} \\
\end{bmatrix} \begin{bmatrix}
\sigma_{xx}' & \sigma_{xy}' & \sigma_{xz}' \\
\sigma_{yx}' & \sigma_{yy}' & \sigma_{yz}' \\
\sigma_{zx}' & \sigma_{zy}' & \sigma_{zz}' \\
\end{bmatrix} \begin{bmatrix}
c_{xx} & c_{yx} & c_{zx} \\
c_{xy} & c_{yy} & c_{zy} \\
c_{xz} & c_{yz} & c_{zz} \\
\end{bmatrix}\!. 
\]
By the table on page~\pageref{transrulesummary2},
the fact that any two coordinate representations of Cauchy stress $\sigma$  satisfy the transformation rule $\Sigma' = C^{-1}\Sigma C^{-\tp}$ says that it is a contravariant $2$-tensor. This provides yet another example why we do not want to identify a $2$-tensor with the matrix that represents it:  stripped of the basis vectors, the physical connotations of \eqref{eq:stress1} and \eqref{eq:ex} that we discussed earlier are irretrievably lost.

What we have described above is the behaviour of stress at a single point in space, assumed to be $\mathbb{R}^3$, and in Cartesian coordinates. As we mentioned in Example~\ref{eg:tenfield2}, stress is a tensor field and  $e_x, e_y, e_z$ really form a basis of the tangent space $\mathbb{T}_v(\mathbb{R}^3)$ at the point $v= (x,y,z) \in \mathbb{R}^3$, and we could have used any coordinates. Note that our discussion above merely used stress in a purely nominal way:  our focus is on the dyadic that represents stress. The same discussion  applies {\em mutatis mutandis} to any contravariant $2$-tensors describing inertia, polarization, strain, tidal force, viscosity, {\em etc.}\ \cite{Borg}.  We refer the readers to \citet[Chapter~4]{Borg}, \citet[Chapter~1]{Chou} and \citet[Chapter~2]{Irgens} for the actual physical details regarding stress, including how to derive \eqref{eq:stress1} from first principles.

Stress is a particularly vital notion, as many tensors describing physical properties in various constitutive equations \cite[Table~1]{Hartmann} are defined as derivatives with respect to it. In fact, this is one way in which tensors of higher order arise in physics -- as \emph{second}-order derivatives of real-valued functions with respect to first- and second-order tensors. For instance, in Cartesian coordinates, the piezo-electric tensor, the piezo-magnetic tensor and the elastic tensor are represented as hypermatrices $D ,Q\in \mathbb{R}^{3 \times 3 \times 3}$ and $S \in \mathbb{R}^{3 \times 3 \times 3 \times 3}$, where
\[
d_{ijk} = - \frac{\partial^2 G}{\partial \sigma_{ij} \partial e_k},\quad
q_{ijk} = - \frac{\partial^2 G}{\partial \sigma_{ij} \partial h_k},\quad
s_{ijkl} = - \frac{\partial^2 G}{\partial \sigma_{ij} \partial \sigma_{kl}},
\]
and $i,j,k,l = 1,2,3$. Here $G=G(\Sigma, E,H, T)$ is the \emph{Gibbs potential}, a real-valued function of the second-order stress tensor $\Sigma$, the first-order tensors  $E$ and $H$ representing electric and magnetic field respectively, and the zeroth-order tensor $T$ representing temperature \cite[Section~3]{Hartmann}.
\end{example}

It is straightforward to extend the above  construction of a dyadic to arbitrary dimensions and arbitrary vectors in abstract vector spaces $\mathbb{U}$ and $\mathbb{V}$. A dyadic is a `linear combination' of vectors in $\mathbb{V}$,
\[
\sigma = \sigma_1 v_1 + \dots + \sigma_n v_n,
\]
whose coefficients $\sigma_j$ are vectors in $\mathbb{U}$:
\[
\sigma_j = \sigma_{1j} u_1 + \dots + \sigma_{mj} u_j, \quad j =1,\ldots,n.
\]
Thus by \eqref{eq:distrib1} and \eqref{eq:distrib2} we have
\[
\sigma = \sum_{i=1}^m\sum_{j=1}^n \sigma_{ij} \, u_i v_j = \sum_{i=1}^m\sum_{j=1}^n \sigma_{ij}\, u_i \otimes v_j
\]
in old (pre-$\otimes$) and modern notation respectively. Note that the coefficients $\sigma_{ij}$ are now scalars and a dyadic is an honest linear combination of dyads with scalar coefficients. We denote the set of all such dyadics as $\mathbb{U} \otimes \mathbb{V}$.
 
We may recursively extend the construction to arbitrary order. With a third vector space $\mathbb{W}$, a \emph{triadic} is a `linear combination' of vectors in $\mathbb{W}$,
\[
\tau = \tau_1 w_1 + \dots + \tau_p w_p,
\]
whose coefficients $\tau_k$ are dyadics in  $\mathbb{U} \otimes \mathbb{V}$:
\[
\tau_k = \sum_{i=1}^m\sum_{j=1}^n \tau_{ijk} \, u_i v_j = \sum_{i=1}^m\sum_{j=1}^n \tau_{ijk} \, u_i \otimes v_j, \quad k=1,\ldots,p.
\]
Thus we have
\[
\tau = \sum_{i=1}^m\sum_{j=1}^n\sum_{j=1}^p \tau_{ijk} \, u_i v_j w_k = \sum_{i=1}^m\sum_{j=1}^n\sum_{j=1}^p \tau_{ijk} \, u_i \otimes v_j \otimes w_k
\]
in old and modern notation respectively. Henceforth we will ditch the old notation. Again a triadic is an honest linear combination of triads and the set of all triadics will be denoted $\mathbb{U} \otimes \mathbb{V} \otimes \mathbb{W}$. Strictly speaking, we have constructed  $(\mathbb{U} \otimes \mathbb{V}) \otimes \mathbb{W}$ and an analogous construction considering `linear combinations' of dyadics in $\mathbb{V} \otimes \mathbb{W}$ with coefficients in $\mathbb{U}$ would give us $\mathbb{U} \otimes (\mathbb{V} \otimes \mathbb{W})$, but it does not matter for us as we have implicitly imposed that
\begin{equation}\label{eq:assoc1}
(u \otimes v) \otimes w = u \otimes (v \otimes w) \eqqcolon u \otimes v \otimes w
\end{equation}
for all $u \in \mathbb{U}$, $v \in \mathbb{V}$, $w \in \mathbb{W}$.

The observant reader may have noticed that aside from \eqref{eq:assoc1} we have also implicitly imposed the third-order analogues of \eqref{eq:distrib1} and \eqref{eq:distrib2} in the construction above. For completeness we will state them formally but in a combined form:\label{pg:construction}
\begin{align}\label{eq:distrib3}
(\lambda u + \lambda' u') \otimes v \otimes w &= \lambda u \otimes v \otimes w + \lambda' u' \otimes v \otimes w, \notag \\*
u \otimes (\lambda v + \lambda' v') \otimes w &= \lambda u \otimes v \otimes w + \lambda' u \otimes v' \otimes w,\\*
u \otimes v \otimes (\lambda w + \lambda' w')  &= \lambda u \otimes v \otimes w + \lambda' u \otimes v \otimes w' \notag
\end{align}
for all vectors $u, u' \in \mathbb{U}$, $v, v' \in \mathbb{V}$, $w, w' \in \mathbb{W}$ and scalars $\lambda, \lambda'$.

The construction described above makes the set  $\mathbb{U} \otimes \mathbb{V} \otimes \mathbb{W}$ into an algebraic object with scalar product, tensor addition $+$ and tensor product $\otimes$ interacting in a consistent manner according to the algebraic rules \eqref{eq:assoc1} and \eqref{eq:distrib3}. We will call it a \emph{tensor product} of $\mathbb{U}$, $\mathbb{V}$, $\mathbb{W}$, or an \emph{abstract tensor product} if we need to emphasize that it refers to this particular construction. Observe that this is an extremely general construction.
\begin{enumerate}[\upshape (i)]
\setlength\itemsep{3pt}
\item It does not depend what we call `scalars' so long as we may add and multiply them, \tie, this construction works for any arbitrary modules $\mathbb{U}$, $\mathbb{V}$, $\mathbb{W}$ over a ring $R$.

\item It does not require $\mathbb{U}$, $\mathbb{V}$, $\mathbb{W}$ to be finite-dimensional (vector spaces) or finitely generated (modules);  the whole construction is about specifying how linear combinations behave under $\otimes$, and {there} is no need for elements to be linear combinations of some basis or generating set.

\item It does not call for separate treatments of covariant and mixed tensors; the definition is agnostic to having some or all of $\mathbb{U}$, $\mathbb{V}$, $\mathbb{W}$  replaced by their dual spaces $\mathbb{U}^*$, $\mathbb{V}^*$, $\mathbb{W}^*$.
\end{enumerate}
It almost goes without saying that the construction can be readily extended to arbitrary order $d$. We may consider `linear combinations' of dyadics with dyadic coefficients to get $d=4$,  `linear combinations' of triadics with dyadic coefficients to get $d =5$, {\em etc.} But looking at the end results of the $d=2$ and $d=3$ constructions, if we have $d$ vector spaces $\mathbb{U}, \mathbb{V},\ldots,\mathbb{W}$, a more direct way is to simply consider the set of all $d$-adics, \ie\ finite sums of $d$-ads,
\begin{align}\label{eq:tenprodvecsp}
  &\mathbb{U} \otimes \mathbb{V} \otimes \dots \otimes \mathbb{W} 
 \coloneqq \biggl\{
\sum_{i=1}^r u_i \otimes v_i \otimes \dots \otimes w_i \colon  
u_i \in \mathbb{U}, v_i\in \mathbb{V},\ldots,w_i \in \mathbb{W}, \; r \in \mathbb{N} \biggr\},
\end{align}
and decree that $\otimes$ is associative, $+$ is associative and commutative, and $\otimes$ is distributive over $+$ in the sense of
\begin{equation}\label{eq:distrib4}
\begin{aligned}
(\lambda u + \lambda' u') \otimes v \otimes \dots \otimes w &= \lambda u \otimes v \otimes  \dots \otimes w + \lambda' u' \otimes v \otimes  \dots \otimes w,\\
u \otimes (\lambda v + \lambda' v') \otimes  \dots \otimes w &= \lambda u \otimes v \otimes  \dots \otimes w + \lambda' u \otimes v' \otimes  \dots \otimes w,\\
\vdots \quad\quad\quad & \quad\quad\quad\quad\quad \vdots\\
u \otimes v \otimes  \dots \otimes (\lambda w + \lambda' w')  &= \lambda u \otimes v \otimes  \dots \otimes w + \lambda' u \otimes v \otimes  \dots \otimes w'
\end{aligned}
\end{equation}
for all vectors $u, u' \in \mathbb{U}, v, v' \in \mathbb{V}, \ldots, w, w' \in \mathbb{W}$, and scalars $\lambda, \lambda'$.

 We have used the old polyadic terminology to connect our discussions to the older  engineering and physics literature, and to show that the modern abstract construction\footnote{This construction is equivalent to the one in \citet[Chapter~XVI, Section~1]{Lang}, likely adapted from \citet[Chapter~II, Section~3]{Bourbaki}, and widely used in pure mathematics. Take the free module generated by all rank-one elements $u \otimes v \otimes\dots\otimes w$ and quotient it by the submodule generated by elements of the forms $(\lambda u + \lambda' u') \otimes v \otimes \dots \otimes w - \lambda u \otimes v \otimes \dots  \otimes w - \lambda' u' \otimes v \otimes \dots  \otimes w$, {\em etc.}, to ensure that \eqref{eq:distrib4} holds.} is firmly rooted in  concrete physical considerations. We will henceforth ditch the obsolete names $d$-ads and $d$-adics and use their modern names, rank-one $d$-tensors and $d$-tensors. Note that these arithmetic rules also allow us to define, for any $T \in \mathbb{U} \otimes \mathbb{V} \otimes \dots \otimes \mathbb{W}$ and any $T' \in  \mathbb{U}' \otimes \mathbb{V}' \otimes \dots \otimes \mathbb{W}'$, an element
\begin{equation}\label{eq:outprod3}
T \otimes T' \in  \mathbb{U} \otimes \mathbb{V} \otimes \dots \otimes \mathbb{W} \otimes \mathbb{U}' \otimes \mathbb{V}' \otimes \dots \otimes \mathbb{W}'.
\end{equation}
We emphasize that $\mathbb{U} \otimes \mathbb{V} \otimes \dots \otimes \mathbb{W}$ does not mean the set of rank-one tensors but the set of finite sums of them as in \eqref{eq:tenprodvecsp}. This mistake is common enough to warrant a displayed equation to highlight the pitfall: 
\begin{equation}\label{eq:pitfall}
\mathbb{U} \otimes \mathbb{V} \otimes \dots \otimes \mathbb{W} \ne \{
u \otimes v \otimes \dots \otimes w \colon  u \in \mathbb{U}, v\in \mathbb{V},\ldots,w \in \mathbb{W} \},
\end{equation}
unless all but one of $\mathbb{U}, \mathbb{V}, \ldots, \mathbb{W}$ are one-dimensional.

For easy reference, we state a special case of the above abstract tensor product construction, but where we have both vector spaces and dual spaces, as another definition for tensors.

\begin{definition}[tensors as elements of tensor spaces]\label{def:tensor3b}\hspz%
Let $p\le d$ be non-negative integers and let $\mathbb{V}_1,\ldots,\mathbb{V}_d$ be vector spaces. A \emph{tensor} of \emph{contravariant order} $p$ and \emph{covariant order} $d-p$ is an element
\[
T \in \mathbb{V}_1 \otimes \dots \otimes \mathbb{V}_p \otimes \mathbb{V}_{p+1}^* \otimes \dots \otimes \mathbb{V}_d^*.
\]
The set above is as defined in \eqref{eq:tenprodvecsp} and is called a \emph{tensor product} of vector spaces or a \emph{tensor space} for short. The tensor $T$ is said to be of \emph{type} $(p,d-p)$ and of \emph{order}  $d$. A tensor of type $(d,0)$, \ie\ $T \in \mathbb{V}_1 \otimes \dots \otimes \mathbb{V}_d$, is called a contravariant $d$-tensor, and a tensor of type $(0,d)$, \ie\ $T \in \mathbb{V}_1^* \otimes \dots \otimes \mathbb{V}_d^*$ is called a covariant $d$-tensor.
\end{definition}

Vector spaces have bases and a tensor space $\mathbb{V}_1  \otimes \dots \otimes \mathbb{V}_d$ has a \emph{tensor product basis} given by
\begin{align}\label{eq:tenprodbasis2}
\mathscr{B}_1 \otimes \mathscr{B}_2 \otimes \dots \otimes \mathscr{B}_d
\coloneqq \{u \otimes v \otimes \dots \otimes w \colon  
u \in \mathscr{B}_1, v \in \mathscr{B}_2,\ldots, w \in \mathscr{B}_d\},
\end{align}
where $\mathscr{B}_i$ is any basis of  $\mathbb{V}_i$, $i=1,\ldots,d$. This is the \emph{only} way to get a basis of rank-one tensors on a tensor space. If $\mathscr{B}_1 \otimes \dots \otimes \mathscr{B}_d$ is a basis for $\mathbb{V}_1  \otimes \dots \otimes \mathbb{V}_d$, then $\mathscr{B}_1,\ldots,\mathscr{B}_d$ must be bases of  $\mathbb{V}_1,\ldots, \mathbb{V}_d$, respectively.  More generally, for mixed tensor spaces, $\mathscr{B}_1 \otimes \dots \otimes \mathscr{B}_p  \otimes \mathscr{B}_{p+1}^* \otimes \dots \otimes \mathscr{B}_d^*$ is a basis for $\mathbb{V}_1 \otimes \dots \otimes \mathbb{V}_p \otimes \mathbb{V}_{p+1}^* \otimes \dots \otimes \mathbb{V}_d^*$, where $\mathscr{B}^*$ denotes the dual basis of $\mathscr{B}$. There is no need to assume finite dimension;  these bases may well be uncountable. If the vector spaces are finite-dimensional, then \eqref{eq:tenprodbasis2} implies that
\begin{equation}\label{eq:dim}
\dim (\mathbb{V}_1 \otimes \dots  \otimes \mathbb{V}_d )=  \dim (\mathbb{V}_1) \cdots \dim (\mathbb{V}_d).
\end{equation}
Aside from \eqref{eq:distrib1},  relation \eqref{eq:dim} is another reason the tensor product is the proper operation for combining quantum states of different systems: if a quantum system is made up of superpositions (\ie\  linear combinations) of $m$ distinct (\ie\  linearly independent) states and another is made up of superpositions of $n$ distinct states, then it is physically reasonable for the combined system to be made up of superpositions of $mn$ distinct states \cite[p.~94]{Nielsen}. As a result, we expect the combination of an $m$-dimensional quantum state space and an $n$-dimensional quantum state space to be an $mn$-dimensional quantum state space, and by \eqref{eq:dim}, the tensor product fits the bill perfectly.
This of course is hand-waving and assumes finite-dimensionality. For a fully rigorous justification, we refer readers to \citet{Daubechies1,Daubechies2,Daubechies3}.

Note that there is no contradiction between \eqref{eq:pitfall} and \eqref{eq:tenprodbasis2}:  in the former the objects are vector spaces; in the latter they are bases of vector spaces. The tensor product symbol $\otimes$ is used in at least a dozen different ways, but fortunately there is little cause for confusion as its meaning is almost always unambiguous from the context.

One downside of Definition~\ref{def:tensor3a} in Section~\ref{sec:tensor3a} is that it requires us to convert everything into real-valued functions before we can discuss tensor products. What we gain from the abstraction in Definition~\ref{def:tensor3b} is generality: the construction above allows us to form tensor products of any objects \emph{as is}, so long as they belong to some vector space (or module).

\begin{example}[concrete tensor products]\label{eg:concretetenprod}
As we mentioned on page~\pageref{pg:noformula}, when $\mathbb{U}$, $\mathbb{V}$, $\mathbb{W}$ are abstract vector spaces, there is no further `meaning' one may attach to $\otimes$, just that it satisfies arithmetical properties such as \eqref{eq:distrib3}. But with abstraction comes generality:  it also means that any associative product operation that satisfies \eqref{eq:distrib3} may play the role of $\otimes$. As soon as we pick specific $\mathbb{U}$, $\mathbb{V}$, $\mathbb{W}$  --  Euclidean spaces, function spaces, distribution spaces, measure spaces, {\em etc.}  --  we may start taking tensor products of these spaces with concrete formulas for $\otimes$. We will go over the most common examples.
\begin{enumerate}[\upshape (i)]
\setlength\itemsep{3pt}
\item\label{it:outprod} Outer products of vectors $a = (a_1,\ldots,a_m) \in \mathbb{R}^m$, $b = (b_1,\ldots,b_n) \in \mathbb{R}^n$:
\[
a \otimes b \coloneqq a b^\tp = 
\begin{bmatrix}
a_1 b_1 &\cdots &a_1 b_n \\
\vdots & \ddots & \vdots \\
a_m b_1 &\cdots & a_m b_n
\end{bmatrix} \in \mathbb{R}^{m \times n}.
\]

\item\label{it:kronprod} Kronecker products of matrices $A \in \mathbb{R}^{m \times n}$, $B \in \mathbb{R}^{p \times q}$:
\[
A \otimes B \coloneqq 
\begin{bmatrix}
a_{11} B &\cdots &a_{1n} B \\
\vdots & \ddots & \vdots \\
a_{m1} B &\cdots & a_{mn} B
\end{bmatrix} \in \mathbb{R}^{mp \times nq}.
\]

\item\label{it:sepprod} Separable products of functions $f \colon  X \to \mathbb{R}$, $g  \colon  Y \to \mathbb{R}$:
\[
f \otimes g \colon  X \times Y \to \mathbb{R},\quad (x,y) \mapsto f(x) g(y).
\]

\item\label{it:tenprod} Separable products of kernels $K \colon  X \times X' \to \mathbb{R}$, $H  \colon  Y \times Y' \to \mathbb{R}$:
\begin{align*}
K \otimes H \colon  (X \times X') \times (Y \times Y') &\to \mathbb{R},\\*
((x,x'), (y,y')) &\mapsto K(x,x') H(y,y').
\end{align*}
\end{enumerate}
Hence these apparently different products of different objects may all be regarded as tensor products. In particular, associativity \eqref{eq:assoc1} is satisfied and we may extend these products to an arbitrary number of factors.
Readers who recall Example~\ref{eg:hyp} would see that \ref{it:outprod} is a special case of \ref{it:sepprod} with $X = [m]$, $Y=[n]$, and \ref{it:kronprod} is a special case of \ref{it:tenprod} with $X = [m]$, $X'=[n]$, $Y = [p]$, $Y'=[q]$. Nevertheless, \ref{it:outprod}  and  \ref{it:kronprod}  may also be independently derived from Definition~\ref{def:tensor3b}. More generally \ref{it:sepprod} recovers Definition~\ref{def:tensor3a}.

The outer product in \ref{it:outprod} comes from Definition~\ref{def:tensor3b}. Take two vectors represented in, say, Cartesian and spherical coordinates:
\begin{equation}\label{eq:uv}
u = a_1 e_x + a_2 e_y + a_3 e_z, \quad v = b_1 e_r + b_2 e_\theta + b_3 e_\phi
\end{equation}
respectively. Consider the contravariant $2$-tensor $u \otimes v$. It follows from \eqref{eq:distrib3} that
\begin{align*}
u \otimes v & = a_1 b_1\, e_x \otimes e_r + a_1b_2 \, e_x \otimes e_\theta + a_1  b_3 \, e_x \otimes e_\phi \\*
&\quad +a_2 b_1\, e_y \otimes e_r + a_2b_2 \, e_y \otimes e_\theta + a_2  b_3 \, e_y \otimes e_\phi \\*
&\quad +a_3 b_1\, e_z \otimes e_r + a_3b_2 \, e_z \otimes e_\theta + a_3  b_3 \, e_z \otimes e_\phi.
\end{align*}
The coefficients are captured by the outer product of their coefficients, \tie,
\[
a\otimes b = \begin{bmatrix} a_1 b_1 & a_1 b_2 & a_1 b_3 \\ a_2 b_1 & a_2 b_2 & a_2 b_3 \\ a_3 b_1 & a_3 b_2 & a_3 b_3 \end{bmatrix} = [u \otimes v]_{\mathscr{B}_1\otimes \mathscr{B}_2} ,
\]
which is a (hyper)matrix representation of $u \otimes v $ in the tensor product basis
\begin{align}\label{eq:tensorprodbasis1}
\mathscr{B}_1\otimes \mathscr{B}_2 & \coloneqq \{e_x \otimes e_r,\;  e_x \otimes e_\theta,\;  e_x \otimes e_\phi,\; e_y \otimes e_r, \notag \\*
&\quad\ \ \; e_y \otimes e_\theta,\;   e_y \otimes e_\phi,\;   e_z \otimes e_r,\;  e_z \otimes e_\theta,\; e_z \otimes e_\phi\}
\end{align}
of the two bases $\mathscr{B}_1 = \{e_x, e_y, e_z\}$ and $\mathscr{B}_2 = \{e_r, e_\theta, e_\phi\}$. If we have another pair of vectors
\[
u' = a'_1 e_x + a'_2 e_y + a'_3 e_z, \quad v' = b'_1 e_r + b'_2 e_\theta + b'_3 e_\phi,
\]
then $u \otimes v + u' \otimes v'$, or more generally $\lambda u \otimes v + \lambda' u' \otimes v'$ for any scalars $\lambda, \lambda'$, obeys
\[
[\lambda u \otimes v + \lambda' u' \otimes v' ]_{\mathscr{B}_1\otimes \mathscr{B}_2}
= \lambda a \otimes b + \lambda' a' \otimes b'.
\]
In other words, the arithmetic of $2$-tensors in $\mathbb{V} \otimes \mathbb{W}$ mirrors the arithmetic of rank-one matrices in $\mathbb{R}^{m \times n}$, $m = \dim \mathbb{W}$ and $n = \dim \mathbb{V}$. It is easy to see that this statement extends to higher order, \tie, the arithmetic of $d$-tensors in $\mathbb{U} \otimes \mathbb{V} \otimes \dots \otimes \mathbb{W}$ mirrors the arithmetic of rank-one $d$-hypermatrices in $\mathbb{R}^{m \times n \times \dots \times p}$:
\[
\biggl[\sum_{i=1}^r \lambda_ i u_i \otimes v_i \otimes \dots \otimes w_i \biggr]_{\mathscr{B}_1\otimes \mathscr{B}_2 \otimes \dots \otimes \mathscr{B}_d}
=\sum_{i=1}^r \lambda_ i [u_i]_{\mathscr{B}_1} \otimes [v_i]_{\mathscr{B}_2} \otimes \dots\otimes [w_i ]_{\mathscr{B}_d}.
\]
Note that the tensor products  on the left of the equation are abstract tensor products and those on the right are outer products.

A few words of caution are in order whenever hypermatrices are involved. Firstly, hypermatrices are dependent on the  bases: if we change $u$ in \eqref{eq:uv} to spherical coordinates with the same numerical values,
\[
u = a_1 e_r + a_2 e_\theta + a_3 e_\phi,
\]
then $\lambda a \otimes b + \lambda' a' \otimes b'$, perfectly well-defined as a (hyper)matrix, is a meaningless thing to compute as we are adding coordinate representations in different bases. Secondly, hypermatrices are oblivious to covariance and contravariance: if we change $u$ in \eqref{eq:uv} to a dual vector with exactly the same coordinates,
\[
u = a_1 e_x^* + a_2 e_y^* + a_3 e_z^*,
\]
then again $\lambda a \otimes b + \lambda' a' \otimes b'$ is meaningless as we are adding coordinate representations of incompatible objects in different vector spaces. Both are of course obvious points but they are often obfuscated when higher-order tensors are brought into the picture.
\end{example}

The Kronecker product above may also be viewed as another manifestation of Definition~\ref{def:tensor3b}, but a more fruitful way is to deduce it as a tensor product of linear operators that naturally follows from Definition~\ref{def:tensor3b}. This is important and sufficiently interesting to warrant separate treatment.

\begin{example}[Kronecker product]\label{eg:kron}
Given linear operators
\[
\Phi_1 \colon  \mathbb{V}_1 \to \mathbb{W}_1 \quad\text{and}\quad \Phi_2 \colon  \mathbb{V}_2 \to \mathbb{W}_2,
\]
forming the tensor products $\mathbb{V}_1  \otimes \mathbb{V}_2$ and $\mathbb{W}_1  \otimes \mathbb{W}_2$ automatically\footnote{This is called \emph{functoriality} in category theory. The fact that the tensor product of vector spaces in Definition~\ref{def:tensor3b} gives a tensor product of linear operators on these spaces says that the tensor product in Definition~\ref{def:tensor3b} is functorial.} gives us a linear operator
\[
\Phi_1 \otimes \Phi_2 \colon  \mathbb{V}_1  \otimes \mathbb{V}_2 \to \mathbb{W}_1  \otimes \mathbb{W}_2
\]
defined on rank-one elements by
\[
\Phi_1 \otimes \Phi_2( v_1 \otimes v_2) \coloneqq \Phi_1 (v_1)\otimes \Phi_2( v_2)
\]
and extended linearly to all elements of $ \mathbb{V}_1  \otimes \mathbb{V}_2$, which are all finite linear combinations of rank-one elements. The way it is defined, $\Phi_1 \otimes \Phi_2$ is clearly unique and we will call it the Kronecker product of linear operators $\Phi_1$ and $\Phi_2$. It extends easily to an arbitrary number of linear operators via
\begin{equation}\label{eq:kronprod2}
\Phi_1 \otimes\dots\otimes \Phi_d( v_1 \otimes\dots\otimes v_d) \coloneqq \Phi_1 (v_1)\otimes\dots\otimes \Phi_d( v_d),
\end{equation}
and the result obeys \eqref{eq:distrib4}:
\[
\begin{aligned}
(\lambda \Phi_1 + \lambda' \Phi_1') \otimes \Phi_2 \otimes \dots \otimes \Phi_d &= \lambda \Phi_1 \otimes \Phi_2 \otimes  \dots \otimes \Phi_d + \lambda' \Phi_1' \otimes \Phi_2 \otimes  \dots \otimes \Phi_d,\\
\Phi_1 \otimes (\lambda \Phi_2 + \lambda' \Phi_2') \otimes  \dots \otimes \Phi_d &= \lambda \Phi_1 \otimes \Phi_2 \otimes  \dots \otimes \Phi_d + \lambda' \Phi_1 \otimes \Phi_2' \otimes  \dots \otimes \Phi_d,\\
\vdots \quad\quad\quad & \quad\quad\quad\quad\quad \vdots\\
\Phi_1 \otimes \Phi_2 \otimes  \dots \otimes (\lambda \Phi_d + \lambda' \Phi_d')  &= \lambda \Phi_1 \otimes \Phi_2 \otimes  \dots \otimes \Phi_d + \lambda' \Phi_1 \otimes \Phi_2 \otimes  \dots \otimes \Phi_d'.
\end{aligned}
\]
In other words, the Kronecker product defines a tensor product in the sense of Definition~\ref{def:tensor3b} on the space of linear operators,
\[
\Lin(\mathbb{V}_1;\mathbb{W}_1) \otimes \dots \otimes \Lin(\mathbb{V}_d;\mathbb{W}_d),
\]
that for finite-dimensional vector spaces equals
\[
\Lin(\mathbb{V}_1\otimes \dots \otimes \mathbb{V}_d;\mathbb{W}_1 \otimes \dots \otimes \mathbb{W}_d).
\]
In addition, as linear operators they may be composed and Moore--Penrose-inverted, and have adjoints and ranks, images and null spaces, all of which work in tandem with the Kronecker product:
\begin{align}
(\Phi_1 \otimes\dots\otimes \Phi_d)  (\Psi_1 \otimes\dots\otimes \Psi_d) &= 
\Phi_1  \Psi_1 \otimes\dots\otimes \Phi_d \Psi_d, \label{eq:kronprodprop1}\\*
(\Phi_1 \otimes\dots\otimes \Phi_d)^\dag &= \Phi_1^\dag \otimes\dots\otimes \Phi_d^\dag,\label{eq:kronprodprop2}\\
(\Phi_1 \otimes\dots\otimes \Phi_d)^* &= \Phi_1^* \otimes\dots\otimes \Phi_d^*,\label{eq:kronprodprop3}\\
\rank(\Phi_1 \otimes\dots\otimes \Phi_d) &= \rank(\Phi_1)\cdots \rank(\Phi_d), \label{eq:kronprodprop4}\\*
\im(\Phi_1 \otimes\dots\otimes \Phi_d) &= \im(\Phi_1) \otimes \cdots \otimes \im(\Phi_d). \label{eq:kronprodprop4a}
\end{align}
Observe that the $\otimes$ on the left of \eqref{eq:kronprodprop4a} is the Kronecker product of operators whereas that on the right is the tensor product of vector spaces. For null spaces, we have
\[
\ker(\Phi_1 \otimes \Phi_2) = \ker(\Phi_1) \otimes \mathbb{V}_2 + \mathbb{V}_1 \otimes  \ker(\Phi_2),
\]
when $d=2$ and more generally
\begin{align*}
\ker(\Phi_1 \otimes \dots\otimes \Phi_d) & = \ker(\Phi_1) \otimes \mathbb{V}_2 \otimes \cdots \otimes \mathbb{V}_d  
 + \mathbb{V}_1 \otimes  \ker(\Phi_2) \otimes \cdots \otimes \mathbb{V}_d\\*
&\quad +\cdots + \mathbb{V}_1 \otimes \mathbb{V}_2 \otimes \cdots \otimes \ker(\Phi_d).
\end{align*}
Therefore injectivity and surjectivity are preserved by taking Kronecker products.
If $\mathbb{V}_i =\mathbb{W}_i$, $i =1,\ldots,d$, then the eigenpairs of $\Phi_1 \otimes\dots\otimes \Phi_d$ are exactly those given by
\begin{equation}\label{eq:kronprodprop5}
(\lambda_1\lambda_2\cdots \lambda_d, v_1 \otimes v_2 \otimes \dots \otimes v_d),
\end{equation}
where $(\lambda_i,v_i)$ is an eigenpair of $\Phi_i$, $i=1,\ldots,d$. All of these are straightforward consequences of \eqref{eq:kronprod2} \cite[Section~13.2]{Berberian}.

The Kronecker product of matrices simply expresses how the matrix representing $\Phi_1 \otimes \dots \otimes \Phi_d$ is related to those representing $\Phi_1, \ldots, \Phi_d$. If we pick bases $\mathscr{B}_i$ for $\mathbb{V}_i$ and $\mathscr{C}_i$ for $\mathbb{W}_i$, then each linear operator $\Phi_i$ has a matrix representation
\[
[\Phi_i]_{\mathscr{B}_i, \mathscr{C}_i}  = X_i \in \mathbb{R}^{m_i \times n_i}
\]
as in \eqref{eq:matreplinop}, $i =1,\ldots,d$, and
\[
[\Phi_1 \otimes \dots\otimes \Phi_d]_{\mathscr{B}_1 \otimes \dots \otimes \mathscr{B}_d, \;\mathscr{C}_1 \otimes \dots \otimes \mathscr{C}_d} = X_1 \otimes \dots \otimes X_d \in \mathbb{R}^{m_1m_2 \cdots m_d \times n_1 n_2 \cdots n_d}.
\]
Note that the tensor products on the left of the equation are Kronecker products as defined in \eqref{eq:kronprod2}  and those on the right are matrix Kronecker products as defined in Example~\ref{eg:concretetenprod}\ref{it:kronprod}, applied $d$ times in any order (order does not matter as $\otimes$ is associative). The tensor product of $d$ bases is as defined in \eqref{eq:tenprodbasis2}.

Clearly the matrix Kronecker product inherits the properties \eqref{eq:kronprodprop1}--\eqref{eq:kronprodprop5}, and when $m_i = n_i$, the eigenvalue property in particular gives
\begin{align}
\tr(X_1 \otimes \dots \otimes X_d) &= \tr(X_1) \cdots \tr(X_d), \label{eq:kronprodtr} \\*
\det(X_1 \otimes \dots \otimes X_d) &= \det(X_1)^{p_1} \cdots \det(X_d)^{p_d}, \notag
\end{align}
where $p_i = (n_1 n_2 \cdots n_d)/n_i$, $i=1,\ldots,d$. Incidentally these last two properties also apply to Kronecker products of operators with  coordinate independent definitions of trace (see Example~\ref{eg:cftr}) and determinant \cite[Section~7.2]{Berberian}.
\end{example}

Next we will see that the multilinear matrix multiplication in \eqref{eq:mmm} and the matrix Kronecker product in Example~\ref{eg:concretetenprod}\ref{it:kronprod} are intimately related.

\begin{example}[multilinear matrix multiplication revisited]\label{eg:mmm}\hspace{-8pt}%
  There is one caveat to obtaining the neat formula for the matrix Kronecker product in Example~\ref{eg:concretetenprod}\ref{it:kronprod} from the operator Kronecker product in Example~\ref{eg:kron}, namely, we must choose the \emph{lexicographic order} for the tensor product bases $\mathscr{B}_1 \otimes \dots \otimes \mathscr{B}_d$ and $\mathscr{C}_1 \otimes \dots \otimes \mathscr{C}_d$, recalling that `basis' in this article always means ordered basis. This is the most common and intuitive way  to totally order a Cartesian product of totally ordered sets. For our purposes, we just need to define it for $[n_1] \times [n_2] \times \dots \times [n_d]$ since our tensor product bases and hypermatrices are all indexed by this. Let
  \[
(i_1, \ldots, i_d), (j_1, \ldots, j_d) \in [n_1]  \times \dots \times [n_d].
  \]
  The lexicographic order $\prec$ is defined as
\[
(i_1, i_2, \ldots, i_d) \prec (j_1, j_2, \ldots, j_d)\quad\text{if and only if}\quad
\begin{cases}
i_1 = j_1, \ldots, i_{k-1} = j_{k-1}, \\
i_k < j_k \text{ for some } k \in [d].
\end{cases}
\]
Take $d =2$, for instance. If $\mathscr{B}_1 = \{u_1, u_2\}$, $\mathscr{B}_2 = \{v_1, v_2, v_3\}$, we have to order
\[
\mathscr{B}_1  \otimes \mathscr{B}_2 = \{ u_1 \otimes v_1, \; u_1 \otimes v_2, \; u_1 \otimes v_3,\;
u_2 \otimes v_1, \; u_2 \otimes v_2, \; u_2 \otimes v_3\}
\]
and not, say, $\{ u_1 \otimes v_1, u_2 \otimes v_1, u_1 \otimes v_2, u_2 \otimes v_2, u_1 \otimes v_3, u_2 \otimes v_3\} \eqqcolon \mathscr{A}$. While this is a trifling point, the latter ordering is unfortunately what is used in the $\vect$ operation common in mathematical software where one concatenates the $n$ columns of an $m \times n$ matrix into a vector of dimension $mn$. The coordinate representation of the $2$-tensor
\begin{align*}
T & = a_{11} u_1 \otimes v_1 + a_{12}  u_1 \otimes v_2 + a_{13}  u_1 \otimes v_3 \\*
&\quad + a_{21} u_2 \otimes v_1 + a_{22}  u_2 \otimes v_2 + a_{23}  u_2 \otimes v_3
\end{align*}
in $\mathscr{B}_1  \otimes \mathscr{B}_2$ (lexicographic) and $\mathscr{A}$ (non-lexicographic) order are
\[
[T]_{\mathscr{B}_1  \otimes \mathscr{B}_2} = \begin{bmatrix}
a_{11} \\ a_{12} \\ a_{13} \\ a_{21} \\ a_{22} \\ a_{23} \end{bmatrix} \in \mathbb{R}^6,\quad
[T]_{\mathscr{A}} = \begin{bmatrix}
a_{11} \\ a_{21} \\ a_{12} \\ a_{22} \\ a_{13} \\ a_{23} \end{bmatrix}
= \vect\biggl( \begin{bmatrix}  a_{11} & a_{12} & a_{13} \\ a_{21} & a_{22} & a_{23} \end{bmatrix} \biggr)\in \mathbb{R}^6
\]
respectively, where $\vect$ denotes the usual `column-concatenated' operation stacking the columns of an $m \times n$ matrix to obtain a vector of dimension $mn$. 
To see why we favour lexicographic order, define the vector space isomorphism
\[
\vect_\ell \colon  \mathbb{R}^{n_1 \times n_2 \times \dots \times n_d} \to \mathbb{R}^{n_1  n_2 \cdots n_d}
\]
that takes a $d$-hypermatrix $A = [a_{i_1i_2\cdots i_d}]$ to a vector
\[
\vect_\ell(A) = (a_{11\cdots 1}, a_{11\cdots 2},\ldots, a_{n_1n_2\cdots n_d}) \in \mathbb{R}^{n_1  n_2 \cdots n_d},
\]
whose coordinates are those in $A$ ordered lexicographically. For any matrices $X_1 \in \mathbb{R}^{m_1 \times n_1},\ldots, X_d \in \mathbb{R}^{m_d \times n_d}$, it is routine to check that
\begin{equation}\label{eq:mmmkron}
\vect_\ell ( (X_1,\ldots,X_d) \cdot A ) = (X_1\otimes\dots\otimes X_d)  \vect_\ell(A).
\end{equation}
In other words, the lexicographic $\vect_\ell$ operator takes multilinear matrix multiplication to the matrix Kronecker product. They are the same operation represented on different but isomorphic vector spaces $\mathbb{R}^{n_1 \times \dots \times n_d}$ and $\mathbb{R}^{n_1\cdots n_d}$, respectively. Incidentally,
$\vect_\ell$ is an example of a `forgetful map':  it forgets the hypermatrix structure. This loss of information is manifested in that, given $(X_1,\ldots,X_d)$, we can compute $X_1\otimes \dots\otimes X_d$ but we cannot uniquely recover $(X_1,\ldots,X_d)$ from their matrix Kronecker product.

For $d=2$, $\vect_\ell$  is just the `row-concatenated $\vect$ operation', \ie\ $\vect_\ell(A) = \vect(A^\tp)$ for any $A \in \mathbb{R}^{m \times n}$, and we have
\[
\vect_\ell(XAY^\tp) = (X \otimes Y)\vect_\ell(A), \quad \vect(XAY^\tp) = (Y^\tp \otimes X)\vect(A).
\]
The neater form for $\vect_\ell$, especially when extended to higher order, is why we prefer lexicographic ordering.
\end{example}

As we discussed in Example~\ref{eg:homult}, there is no formula that extends the usual matrix--matrix multiplication to $3$-hypermatrices or indeed any $d$-hypermatrices for odd $d \in \mathbb{N}$. Such a formula, however, exists for any even $d$, thanks to the Kronecker product.

\begin{example}[`multiplying higher-order tensors' revisited]\label{eg:homult2}\hspz%
Consider the following $2d$-hypermatrices:
\[
A \in \mathbb{R}^{m_1 \times n_1 \times m_2 \times n_2 \times \dots \times m_d \times n_d},\quad
B \in \mathbb{R}^{n_1 \times p_1 \times n_2 \times p_2 \times \dots \times n_d \times p_d}.
\]
Their product, denoted $AB$ to be consistent with the usual matrix--matrix product, is  the $2d$-hypermatrix $C \in \mathbb{R}^{m_1 \times p_1 \times m_2 \times p_2 \times \dots \times m_d \times p_d}$, where
\[
c_{i_1 k_1 \cdots i_d k_d} \coloneqq \sum_{j_1=1}^{n_1}\cdots\sum_{j_d=1}^{n_d} a_{i_1 j_1 \cdots i_d j_d}
b_{j_1 k_1 \cdots j_d k_d}.
\]
This clearly extends the matrix--matrix product, which is just the case when $d=1$.
We have in fact already encountered Example~\ref{eg:kron} in a different form: this is \eqref{eq:kronprodprop1}, the composition of two Kronecker products of $d$ linear operators. To view it in hypermatrix form, we just use the fact that $\Lin(\mathbb{V}; \mathbb{W}) = \mathbb{V}^* \otimes \mathbb{W}$. Indeed, we may also view it as a tensor contraction of
\[
T \in (\mathbb{U}^*_1 \otimes \mathbb{V}_1 ) \otimes \dots \otimes (\mathbb{U}^*_d \otimes \mathbb{V}_d ),\quad
T' \in (\mathbb{V}^*_1 \otimes \mathbb{W}_1 ) \otimes \dots \otimes (\mathbb{V}^*_d \otimes \mathbb{W}_d ),
\]
from which we obtain
\[
\langle T, T' \rangle \in (\mathbb{U}^*_1 \otimes \mathbb{W}_1 ) \otimes \dots \otimes (\mathbb{U}^*_d \otimes \mathbb{W}_d ).
\]
This is entirely within expectation as the first fundamental theorem of invariant theory, the result that prevents the existence of a $d$-hypermatrix--hypermatrix product for odd $d$, also tells us that essentially the only way to multiply tensors without increasing order is via contractions. This product is well-defined for $2d$-tensors and does not depend on bases: in terms of Kronecker products of matrices,
\begin{align*}
& ( (X_1A_1 Y_1^{-1}) \otimes \dots \otimes (X_dA_d Y_d^{-1}) )
( (Y_1B_1 Z_1^{-1}) \otimes \dots \otimes (Y_d B_d Z_d^{-1}) ) \\*
&\quad = (X_1A_1B_1 Z_1^{-1}) \otimes \dots \otimes (X_d A_dB_d Z_d^{-1}),
\end{align*}
and in terms of hypermatrix--hypermatrix products,
\begin{align*}
& ((X_1,Y_1^{-1},\ldots, X_d, Y_d^{-1}) \cdot A) ((Y_1,Z_1^{-1},\ldots, Y_d, Z_d^{-1}) \cdot B)\\*
&\quad = (X_1,Z_1^{-1},\ldots, X_d, Z_d^{-1}) \cdot (AB),
\end{align*}
\tie, they satisfy the higher-order analogue of \eqref{eq:matprod}.
Frankly, we do not see any advantage in formulating such a product as $2d$-hypermatrices, but an abundance of disadvantages.
\end{example}

Just as it is possible to deduce the tensor transformation rules in definition~\ref{st:tensor1} from definition~\ref{st:tensor2}, we can do likewise with definition~\ref{st:tensor3} in the form of Definition~\ref{def:tensor3b}.

\begin{example}[tensor transformation rules revisited]\label{eg:tentransf3}
  As stated on page~\pageref{pg:construction}, the tensor product construction in this section does not require that  $\mathbb{U}, \mathbb{V}, \ldots, \mathbb{W}$ have bases, \tie, they could be modules, which do not have bases in general
  (those with bases are called \emph{free modules}). Nevertheless, in the event where they do have bases, their change-of-basis theorems would lead us directly to the transformation rules discussed in Section~\ref{sec:trans}.

We first remind the reader of a simple notion, discussed in standard linear algebra textbooks such as \citet[Section~3.9]{Berberian} and \citet[Section~2.6]{FIS} but often overlooked. Any linear operator $\Phi \colon  \mathbb{V} \to \mathbb{W}$ induces a \emph{transpose} linear operator on the dual spaces defined as
\[
\Phi^\tp \colon  \mathbb{W}^* \to \mathbb{V}^*, \quad \Phi(\varphi) \coloneqq \varphi \circ \Phi
\]
for any linear functional $\varphi \in  \mathbb{W}^* $. Note that the composition $\varphi \circ \Phi \colon  \mathbb{V} \to \mathbb{R}$ is indeed a linear functional in $\mathbb{V}^*$. The reason for its name is that
\[
[\Phi]_{\mathscr{B},\mathscr{C}} = A \in \mathbb{R}^{m \times n}\quad \text{if and only if} \quad
[\Phi^\tp]_{\mathscr{C}^*,\mathscr{B}^*} = A^\tp \in \mathbb{R}^{n \times m}.
\]
One may show that $\Phi$ is injective if and only if $\Phi^\tp$ is subjective, and $\Phi$ is surjective if and only if $\Phi^\tp$ is injective.\footnote{Not a typo. Injectivity and surjectivity are in fact dual notions in this sense.} So if $\Phi \colon  \mathbb{V} \to \mathbb{W}$ is invertible, then so is $\Phi^\tp  \colon  \mathbb{W}^* \to \mathbb{V}^*$ and its inverse is a linear operator,
\[
\Phi^{-\tp} \colon  \mathbb{V}^* \to \mathbb{W}^*.
\]
Another name for `invertible linear operator' is vector space isomorphism, especially when used in the following context.
Any basis $\mathscr{B} = \{v_1,\ldots,v_n\}$ of $\mathbb{V}$ gives us two vector space isomorphisms, namely (i)~$\Phi_{\mathscr{B}} \colon  \mathbb{V}  \to \mathbb{R}^n$  that takes {$v = a_1 v_1 + \dots +a_n v_n \in \mathbb{V}$} to its coordinate representation in $\mathscr{B}$, and
(ii)~$\Phi_{\mathscr{B}^*} \colon  \mathbb{V}^*  \to \mathbb{R}^n$  that takes $\varphi = b_1 v_1^* + \dots +b_n v_n^* \in \mathbb{V}^*$ to its coordinate representation in $\mathscr{B}^*$:
\begin{equation}
\Phi_{\mathscr{B}}(v) = \begin{bmatrix} a_1\\ \vdots \\ a_n\end{bmatrix} \in \mathbb{R}^n,\quad
  \Phi_{\mathscr{B}^*}(\varphi) = \begin{bmatrix} b_1\\ \vdots \\ b_n\end{bmatrix} \in \mathbb{R}^n, \label{eq:Phibasis}
    \end{equation}
  and as expected
  \[
\Phi_{\mathscr{B}^*} = \Phi_{\mathscr{B}}^{-\tp}. \notag
\]
This is of course just paraphrasing what we have already discussed in \eqref{eq:vecrep} and thereabouts. Note that we do not distinguish between $\mathbb{R}^n$ and its dual space $(\mathbb{R}^n)^*$.

Given vector spaces $\mathbb{V}_1,\ldots,\mathbb{V}_d$ and $\mathbb{W}_1,\ldots,\mathbb{W}_d$ and invertible linear operators $\Phi_i \colon  \mathbb{V}_i \to \mathbb{W}_i$, $i =1,\ldots,d$, take the Kronecker product of $\Phi_1,\ldots,\Phi_p$, $\Phi_{p+1}^{-\tp},\ldots,\Phi_d^{-\tp}$ as defined in \eqref{eq:kronprod2} of the previous example. Then the linear operator
\begin{align*}
& \Phi_1 \otimes \dots \otimes \Phi_p \otimes \Phi_{p+1}^{-\tp}\otimes \dots\otimes \Phi_d^{-\tp} \colon  \mathbb{V}_1 \otimes \dots \otimes \mathbb{V}_p \otimes \mathbb{V}_{p+1}^* \otimes \dots \otimes \mathbb{V}_d^*\\*
&\quad  \to \mathbb{W}_1 \otimes \dots \otimes \mathbb{W}_p \otimes \mathbb{W}_{p+1}^* \otimes \dots \otimes \mathbb{W}_d^*
\end{align*}
must also be invertible by \eqref{eq:kronprodprop2}.

For each $i=1,\ldots,d$, let $\mathscr{B}_i$ and $\mathscr{B}_i'$ be two different  bases of $\mathbb{V}_i$; let $\mathbb{W}_i = \mathbb{R}^{n_i}$ with $n_i = \dim \mathbb{V}_i$ and
\[
\Phi_i \coloneqq \Phi_{\mathscr{B}_i}, \quad \Psi_i \coloneqq \Phi_{\mathscr{B}_i'}
\]
as in \eqref{eq:Phibasis}. Since we assume $\mathbb{R}^n = (\mathbb{R}^n)^*$, we have
\begin{equation}\label{eq:fudge1}
\mathbb{W}_1 \otimes \dots \otimes \mathbb{W}_p \otimes \mathbb{W}_{p+1}^* \otimes \dots \otimes \mathbb{W}_d^* \cong \mathbb{R}^{n_1} \otimes \dots \otimes \mathbb{R}^{n_d}
\cong \mathbb{R}^{n_1 \times \dots \times n_d},
\end{equation}
where the last $\cong$ is a consequence of using Example~\ref{eg:concretetenprod}\ref{it:outprod}. With these choices, the vector space isomorphisms 
\[
\begin{aligned}
&  \Phi_1 \otimes \dots \otimes \Phi_p \otimes \Phi_{p+1}^{-\tp}\otimes \dots\otimes \Phi_d^{-\tp} \colon \\
&\hspace{48pt}  \mathbb{V}_1 \otimes \dots \otimes \mathbb{V}_p \otimes \mathbb{V}_{p+1}^* \otimes \dots \otimes \mathbb{V}_d^* &\to \mathbb{R}^{n_1 \times \dots \times n_d},\\
&  \Psi_1 \otimes \dots \otimes \Psi_p \otimes \Psi_{p+1}^{-\tp}\otimes \dots\otimes \Psi_d^{-\tp} \colon \\
&\hspace{48pt}  \mathbb{V}_1 \otimes \dots \otimes \mathbb{V}_p \otimes \mathbb{V}_{p+1}^* \otimes \dots \otimes \mathbb{V}_d^* &\to \mathbb{R}^{n_1 \times \dots \times n_d}
\end{aligned}
\]
give us hypermatrix representations of $d$-tensors with respect to $\mathscr{B}_1,\ldots,\mathscr{B}_d$ and $\mathscr{B}_1',\ldots,\mathscr{B}_d'$ respectively. 
Let $T \in \mathbb{V}_1 \otimes \dots \otimes \mathbb{V}_p \otimes \mathbb{V}_{p+1}^* \otimes \dots \otimes \mathbb{V}_d^* $ and
\begin{align*}
\Phi_1 \otimes \dots \otimes \Phi_p \otimes \Phi_{p+1}^{-\tp}\otimes \dots\otimes \Phi_d^{-\tp} (T) &= A \in \mathbb{R}^{n_1 \times \dots \times n_d},\\*
\Psi_1 \otimes \dots \otimes \Psi_p \otimes \Psi_{p+1}^{-\tp}\otimes \dots\otimes \Psi_d^{-\tp} (T) &= A' \in \mathbb{R}^{n_1 \times \dots \times n_d}.
\end{align*}
Then the relation between the hypermatrices $A$ and $A'$ is given by
\[
A = \Phi_1\Psi_1^{-1} \otimes \dots \otimes \Phi_p\Psi_p^{-1} \otimes (\Phi_{p+1} \Psi_{p+1}^{-1})^{-\tp} \otimes \dots\otimes (\Phi_d \Psi_d^{-1})^{-\tp} (A'),
\]
where we have used \eqref{eq:kronprodprop1} and \eqref{eq:kronprodprop2}. Note that each $\Phi_i\Psi_i^{-1} \colon  \mathbb{R}^{n_i} \to \mathbb{R}^{n_i} $ is an invertible linear operator and so it must be given by $\Phi_i\Psi_i^{-1}(v)= X_i v$ for some $X_i \in \GL(n_i)$. Using \eqref{eq:mmmkron}, the above relation between $A$ and $A'$ in terms of multilinear matrix multiplication is just
\[
A = (X_1,\ldots,X_p, X_{p+1}^{-\tp},\ldots,X_d^{-\tp}) \cdot A',
\]
which gives us the tensor transformation rule \eqref{eq:mixed1}.
To get the isomorphism in \eqref{eq:fudge1}, we had to identify $\mathbb{R}^n$ with $(\mathbb{R}^n)^*$, and this is where we lost all information pertaining to covariance and contravariance; a hypermatrix cannot perfectly represent a tensor, which is one reason why we need to look at the transformation rules to ascertain the tensor.
\end{example}

The example above essentially discusses change of basis without any mention of bases. In pure mathematics, sweeping such details under the rug is generally regarded as a good thing; in applied and computational mathematics, such details are usually unavoidable and we will work them out below.

\begin{example}[representing tensors as hypermatrices]\label{eg:hypmatrep}\hspace{-7.7pt}%
Let $\mathscr{B}_1 = \{u_1,\ldots,u_m\}$, $\mathscr{B}_2=\{v_1,\ldots,v_n\},\ldots,\mathscr{B}_d = \{w_1,\ldots,w_p\}$ be any bases of finite-dimensional vector spaces  $\mathbb{U} ,\mathbb{V},\ldots, \mathbb{W}$. Then any $d$-tensor  $T \in \mathbb{U} \otimes \mathbb{V} \otimes \dots \otimes \mathbb{W}$ may be expressed as
\begin{equation}\label{eq:hypmatrep}
T = \sum_{i=1}^m \sum_{j=1}^n \cdots \sum_{k=1}^p  a_{ij\cdots k}\, u_i \otimes v_j \otimes \cdots \otimes w_k,
\end{equation}
and the $d$-hypermatrix of coefficients
\[
[T]_{\mathscr{B}_1 ,\mathscr{B}_2,\ldots,\mathscr{B}_d} \coloneqq  [a_{ij\cdots k}]_{i,j,\ldots,k=1}^{m,n,\ldots,p} = A \in \mathbb{R}^{m \times n \times \dots \times p}
\]
is said to \emph{represent} or to be a \emph{coordinate representation} of $T$ with respect to  $\mathscr{B}_1 ,\mathscr{B}_2,\ldots,\mathscr{B}_d$. That \eqref{eq:hypmatrep} always holds is easily seen as follows. First it holds for a rank-one tensor $T =u \otimes v \otimes \dots \otimes w$ as we may express $u,v,\ldots,w$ as a linear combination of basis vectors and apply  \eqref{eq:distrib4} to get
\begin{align}\label{eq:outprod2}
&  (a_1 u_1 + \dots + a_m u_m) \otimes (b_1 v_1 + \dots + b_n v_n) \otimes \dots \otimes (c_1 w_1 + \dots + c_p w_p)\notag  \\*
  &\quad =  \sum_{i=1}^m \sum_{j=1}^n \cdots \sum_{k=1}^p  a_i b_j\cdots c_k\, u_i \otimes v_j \otimes \cdots \otimes w_k,
\end{align}
which has the form in \eqref{eq:hypmatrep}. Since an arbitrary $T$ is simply a finite sum of rank-one tensors as in \eqref{eq:tenprodvecsp}, the distributivity in \eqref{eq:distrib4} again ensures that
\begin{align}\label{eq:lincomb2}
&\lambda \sum_{i=1}^m \sum_{j=1}^n \cdots \sum_{k=1}^p  a_{ij\cdots k}\, u_i \otimes v_j \otimes \cdots \otimes w_k \notag \\*
&\qquad + \lambda' \sum_{i=1}^m \sum_{j=1}^n \cdots \sum_{k=1}^p  b_{ij\cdots k}\, u_i \otimes v_j \otimes \cdots \otimes w_k \notag \\*
& \quad =\sum_{i=1}^m \sum_{j=1}^n \cdots \sum_{k=1}^p (\lambda  a_{ij\cdots k}  + \lambda' b_{ij\cdots k} )\, u_i \otimes v_j \otimes \cdots \otimes w_k.
\end{align}
By \eqref{eq:outprod2} and \eqref{eq:lincomb2}, any $T$ must be expressible in the form \eqref{eq:hypmatrep}. For a different choice of bases $\mathscr{B}_1' ,\mathscr{B}_2',\ldots,\mathscr{B}_d'$ we obtain a different hypermatrix representation $A'  \in \mathbb{R}^{m \times n \times \dots \times p}$ for $T$ that is related to $A$ by a multilinear matrix multiplication  as in Example~\ref{eg:tentransf3}.

Note that \eqref{eq:outprod2} and \eqref{eq:lincomb2}, respectively, give the formulas for outer product and linear combination:
\begin{align*}
(a_1,\ldots,a_m) \otimes (b_1,\ldots,b_n) \otimes \dots \otimes (c_1,\ldots,c_p) &= [ a_i b_j\cdots c_k],\\*
\lambda [ a_{ij\cdots k} ] + \lambda'  [ b_{ij\cdots k} ] &=  [ \lambda  a_{ij\cdots k}  + \lambda' b_{ij\cdots k} ].
\end{align*}
These also follow from our definition of hypermatrices as real-valued functions in Example~\ref{eg:hyp}. 

If $T' \in \mathbb{U}' \otimes \mathbb{V}' \otimes \dots \otimes \mathbb{W}'$ is a $d'$-tensor, expressed as in \eqref{eq:hypmatrep} as
\[
T' =  \sum_{i'=1}^{m'} \sum_{j'=1}^{n'} \cdots \sum_{k'=1}^{p'}   b_{i'j'\cdots k'}\, u'_{i'} \otimes v'_{j'} \otimes \dots \otimes w'_{k'},
\]
with bases of $\mathbb{U}',\mathbb{V}',\ldots, \mathbb{W}'$
\[
\mathscr{B}_1' = \{u'_1,\ldots,u'_{m'}\},\quad \mathscr{B}_2'=\{v_1',\ldots,v'_{n'}\},\ldots,\mathscr{B}_{d'}' = \{w_1',\ldots,w'_{p'}\},
\]
then the outer product of $T$ and $T'$ in \eqref{eq:outprod3} may be expressed as
\[
T \otimes T' = \sum_{i,j,\ldots,k,i',j',\ldots,k'=1}^{m,n\ldots,p,m',n'\ldots,p'}  a_{ij\cdots k} b_{i'j'\cdots k'}\,  u_i \otimes v_j \otimes \cdots \otimes w_k \otimes u'_{i'} \otimes v'_{j'} \otimes \cdots \otimes w'_{k'},
\]
which is a $(d+d')$-tensor. In terms of their hypermatrix representations,
\[
A \otimes B = [a_{ij\cdots k} b_{i'j'\cdots k'}]_{i,j,\ldots,k,i',j',\ldots,k'=1}^{m,n\ldots,p,m',n'\ldots,p'}  \in \mathbb{R}^{m \times n \times \dots \times p \times m' \times n' \times \dots \times p'}.
\]
It is straightforward to show that with respect to the Kronecker product or multilinear matrix multiplication,
  \[
\begin{aligned}
&  [\Phi_1 \otimes \dots \otimes \Phi_d(T)] \otimes [\Psi_1 \otimes \dots \otimes \Psi_{d'}(T')] \\
  &\hspace{48pt} =
\Phi_1 \otimes \dots \otimes \Phi_d \otimes \Psi_1 \otimes \dots \otimes \Psi_{d'}(T\otimes T'),\\
&[(X_1,\ldots,X_d) \cdot A] \otimes [(Y_1,\ldots,Y_{d'}) \cdot B]\\*
&\hspace{48pt} = (X_1,\ldots,X_d,Y_1,\ldots,Y_{d'}) \cdot A \otimes B,
\end{aligned} 
\]
with notation as in Examples~\ref{eg:kron}, \ref{eg:mmm} and \ref{eg:tentransf3}.

The above discussion extends verbatim to covariant and mixed tensors by replacing some or all of these vector spaces with their duals and bases with dual bases. It also extends to separable Hilbert spaces with orthonormal bases $\mathscr{B}_1 ,\ldots,\mathscr{B}_d$ and $A \in l^2(\mathbb{N}^d)$  an infinite-dimensional hypermatrix as discussed in Example~\ref{eg:hyp}. Nevertheless we do not recommend working exclusively with such hypermatrices. The reason is simple:  there is a better representation of the tensor $T$ sitting in plain sight, namely \eqref{eq:hypmatrep}, which contains not just all information in the hypermatrix $A$ but also the bases $\mathscr{B}_1 ,\ldots,\mathscr{B}_d$. Indeed, the representation of a tensor in \eqref{eq:hypmatrep}, coupled with the arithmetic rules for manipulating $\otimes$ like \eqref{eq:distrib4}, is more or less the standard method for working with tensors when bases are involved. Take $d=2$ for illustration. We can tell at a glance from
\[
\beta = \sum_{i=1}^m \sum_{j=1}^n a_{ij} \,u_i^* \otimes v_j^*,\quad
\Phi = \sum_{i=1}^m \sum_{j=1}^n a_{ij} \, u_i \otimes v_j^*
\]
that $\beta$ is a bilinear functional and $\Phi$ a linear operator (the dyad case is discussed in Example~\ref{eg:stress}), information that is lost if one just looks at the matrix $A = [a_{ij}]$. Furthermore, the $u_i$ and $v_j$ are sometimes more important than the $a_{ij}$. We might have
\begin{equation}\label{eq:quadrature1}
\beta = \varepsilon_\mrt \otimes \varepsilon_\mrt + \varepsilon_\mrt \otimes \varepsilon_\prt + \varepsilon_\prt \otimes \varepsilon_\mrt+ \varepsilon_\prt \otimes \varepsilon_\prt,
\end{equation}
where the $a_{ij}$ are all ones, and all the important information is in the basis vectors, which are the point evaluation functionals we encountered in Example~\ref{eg:multmult},
\tie, $\varepsilon_x (f) = f(x)$ for any $f \colon  [-1,1]\to \mathbb{R}$. The decomposition \eqref{eq:quadrature1} is in fact a
four-point Gauss quadrature formula on the domain $[-1,1] \times [-1,1]$. We will revisit this in Example~\ref{eg:quadrature}.
\end{example}

We next discuss tensor rank in light of Definition~\ref{def:tensor3b}.

\begin{example}[multilinear rank and tensor rank]\label{eg:mrank}\hspace{-6.5pt}%
For any $d$ vector spaces $\mathbb{U}, \mathbb{V}$, $\ldots, \mathbb{W}$ and any subspaces $\mathbb{U}' \subseteq \mathbb{U}, \mathbb{V}' \subseteq \mathbb{V},\ldots,\mathbb{W}' \subseteq \mathbb{W}$, it follows from \eqref{eq:tenprodvecsp}  that
\[
\mathbb{U}' \otimes \mathbb{V}' \otimes \dots \otimes \mathbb{W}' \subseteq \mathbb{U} \otimes \mathbb{V} \otimes \dots \otimes \mathbb{W}.
\]
Given a non-zero $d$-tensor $T \in \mathbb{U} \otimes \mathbb{V} \otimes \dots \otimes \mathbb{W}$, one might ask for a smallest tensor subspace that contains $T$, \tie,
\begin{align}\label{eq:mrank}
\mrank(T) &\coloneqq \min \{(\dim \mathbb{U}',\dim \mathbb{V}',\ldots, \dim \mathbb{W}') \in \mathbb{N}^d \colon  \\*
&\hspace{36pt} T\in \mathbb{U}' \otimes \mathbb{V}' \otimes \dots \otimes \mathbb{W}',\; \mathbb{U}' \subseteq \mathbb{U}, \mathbb{V}' \subseteq \mathbb{V},\ldots,\mathbb{W}' \subseteq \mathbb{W}\}.\notag
\end{align}
This gives the \emph{multilinear rank} of $T$, first studied in \citet[Theorem~II]{Hitch1}. Let $\mathbb{U}', \mathbb{V}',\ldots,\mathbb{W}'$ be subspaces that attain the minimum in \eqref{eq:mrank} and let $ \{u_1,\ldots,u_p\},\{v_1,\ldots,v_q\},\ldots,\{w_1,\ldots,w_r\}$ be any respective bases. Then, by \eqref{eq:hypmatrep}, we get a multilinear rank decomposition \cite[Corollary~I]{Hitch1}:
\begin{equation}\label{eq:multdecomp}
T = \sum_{i=1}^p \sum_{j=1}^q \cdots \sum_{k=1}^r c_{ij\cdots k} \, u_i \otimes v_j \otimes \dots \otimes w_k.
\end{equation}
For a given $T$, the vectors on the right of \eqref{eq:multdecomp} are not unique but the subspaces $\mathbb{U}', \mathbb{V}',\ldots,\mathbb{W}'$ that they span are unique, and thus $\mrank(T) $ is well-defined.\footnote{While $\mathbb{N}^d$ is only partially ordered, this shows that the minimum in \eqref{eq:mrank} is unique.} 
One may view these subspaces as generalizations of the row and column spaces of a matrix and multilinear rank as a generalization of the row  and column ranks.
Note that the coefficient hypermatrix $C = [c_{ij\cdots k}] \in \mathbb{R}^{p \times q \times \dots \times r}$ is not unique either, but two such coefficient hypermatrices $C$ and $C'$ must be related by the tensor transformation rule \eqref{eq:mixed1}, as we saw in Example~\ref{eg:hypmatrep}.

The decomposition in \eqref{eq:multdecomp} is a sum of rank-one tensors, and if we simply use this as impetus and define
\begin{align*}
 \rank(T)& \coloneqq \min\biggl\{ r \in \mathbb{N} \colon  T = \sum_{i=1}^r u_i \otimes v_i \otimes \dots \otimes w_i, \; u_1,\ldots,u_r \in \mathbb{U},\\*
 &\hspace{38pt} v_1,\ldots,v_r \in \mathbb{V}, \ldots, w_1,\ldots,w_r \in \mathbb{W}\biggr\},
\end{align*}
we obtain \emph{tensor rank} \cite[equations~2 and $2_a$]{Hitch1}, which we have encountered in a more restricted form in \eqref{eq:trankbilin}. For completeness, we define
\[
\mrank(T) \coloneqq (0,0,\ldots,0)\quad \Leftrightarrow \quad \rank(T) \coloneqq 0\quad \Leftrightarrow \quad T = 0.
\]
Note that
\[
\mrank(T) = (1,1,\ldots,1)\quad \Leftrightarrow \quad \rank(T) = 1\quad \Leftrightarrow \quad
T = u \otimes v \otimes \dots \otimes w,
\]
\tie, a rank-one tensor has rank one irrespective of the rank we use. 
Also, we did not exclude infinite-dimensional vector spaces from our discussion.
Note that by construction $\mathbb{U} \otimes \mathbb{V} \otimes \dots \otimes \mathbb{W}$ comprises
\emph{finite} linear combinations of rank-one tensors, and so by its definition, $\mathbb{U} \otimes \mathbb{V} \otimes \dots \otimes \mathbb{W}$ is the set of all tensors of \emph{finite rank}. Furthermore,
whether or not a tensor has finite rank  is independent of the rank we use.
\end{example}

Vector spaces become enormously more interesting when equipped with inner products and norms; it is the same with tensor spaces. The tensor product construction leading to Definition~\ref{def:tensor3b} allows us to incorporate them easily.

\begin{example}[tensor\, products\, of\, inner\, product\, and\, norm\, spaces]\label{eg:norms}\hspace{-1.7pt}%
 We let $\mathbb{U}, \mathbb{V}, \ldots, \mathbb{W}$ be $d$ vector spaces equipped with either norms
\[
\lVert\,\cdot\,\rVert_1\colon  \mathbb{U} \to [0,\infty),\; \lVert\,\cdot\,\rVert_2\colon  \mathbb{V} \to [0,\infty), \ldots, \lVert\,\cdot\,\rVert_d \colon  \mathbb{W} \to [0,\infty)
\]
or inner products
\[
\langle\, \cdot\,, \cdot \,\rangle_1 \colon  \mathbb{U} \times \mathbb{U} \to \mathbb{R},\; \langle\, \cdot\,, \cdot \,\rangle_2\colon  \mathbb{V} \times \mathbb{V} \to \mathbb{R}, \ldots, \langle\, \cdot\,, \cdot \,\rangle_d\colon  \mathbb{W} \times \mathbb{W} \to \mathbb{R}.
\]
We would like to define a norm or an inner product on the tensor space $\mathbb{U} \otimes \mathbb{V} \otimes \dots \otimes \mathbb{W}$ as defined in \eqref{eq:tenprodvecsp}. Note that we do not assume that these vector spaces are finite-dimensional.

The inner product is easy. First define it on rank-one tensors,
\begin{equation}\label{eq:innprod1}
\langle u \otimes v \otimes \dots \otimes w, u' \otimes v' \otimes \dots \otimes w' \rangle\coloneqq \langle u, u' \rangle_1 \langle v, v' \rangle_2 \cdots \langle w, w' \rangle_d
\end{equation}
for all $u,u' \in \mathbb{U}$, $v,v' \in \mathbb{V},\ldots, w,w'\in \mathbb{W}$. Then extend it bilinearly, \tie, decree that for rank-one tensors $S,T,S',T'$ and scalars $\lambda,\lambda'$,
\begin{align*}
\langle \lambda S + \lambda' S', T \rangle &\coloneqq \lambda  \langle S, T \rangle + \lambda' \langle S', T \rangle,\\*
\langle S, \lambda T +  \lambda' T' \rangle &\coloneqq  \lambda \langle S, T \rangle + \lambda' \langle S, T' \rangle.
\end{align*}
As $\mathbb{U} \otimes \mathbb{V} \otimes \dots \otimes \mathbb{W}$ comprises finite linear combinations of rank-one tensors by definition, this defines $\langle\, \cdot\,, \cdot \,\rangle$ on the whole tensor space. It is then routine to check that the axioms of inner product are satisfied by  $\langle\, \cdot\,, \cdot \,\rangle$. This construction applies verbatim to other bilinear functionals such as the Lorentzian scalar product in Example~\ref{eg:cartesian}.

Norms are slightly trickier because there are multiple ways to define them on tensor spaces. We will restrict ourselves to the three most common constructions. If we have an inner product as above, then the norm induced by the inner product,
\[
\lVert T \rVert_\F \coloneqq \sqrt{\langle S, T \rangle},
\]
is called the Hilbert--Schmidt norm. Note that by virtue of \eqref{eq:innprod1}, the Hilbert--Schmidt norm would satisfy
\begin{equation}\label{eq:normprod1}
\lVert u \otimes v \otimes \dots \otimes w \rVert = \lVert u \rVert_1 \lVert v \rVert_2 \cdots \lVert  w \rVert_d,
\end{equation}
a property that we will require of all tensor norms. As we mentioned in Example~\ref{eg:infinite}, such norms are called cross-norms. We have to extend \eqref{eq:normprod1} so that $\lVert\,\cdot\,\rVert$  is defined on not just rank-one tensors but on all finite sums of these. Naively we might just define it to be the sum of the norms of its rank-one summands, but that does not work as we get different values with different sums of rank-one tensors. Fortunately, taking the infimum over all possible finite sums
\begin{equation}\label{eq:nuclear2}
\lVert T \rVert_\nu \coloneqq \inf\biggl\{\sum_{i=1}^r \lVert u_i \rVert_1 \lVert v_i \rVert_2 \cdots \lVert  w_i \rVert_d \colon 
T =  \sum_{i=1}^r u_i \otimes v_i \otimes \dots \otimes w_i \biggr\}
\end{equation}
fixes the issue and gives us a norm. This is the \emph{nuclear norm}, special cases of which we have encountered in \eqref{eq:tennuclear} and in Examples~\ref{eg:nuclear} and \ref{eg:infinite}.

There is an alternative. We may also take the supremum over all rank-one tensors in the dual tensor space:
\begin{equation}\label{eq:spectral2}
\lVert T\rVert_\sigma \coloneqq \sup \biggl\{ \dfrac{\lvert \varphi \otimes \psi \otimes \dots \otimes \theta (T) \rvert}{\lVert \varphi \rVert_{1\ast} \lVert \psi \rVert_{2\ast} \cdots \lVert \theta \rVert_{d\ast} } \colon   \varphi \in \mathbb{U}^*,\psi \in \mathbb{V}^*,\ldots,\theta \in \mathbb{W}^* \biggr\},
\end{equation}
where  $\lVert\,\cdot\,\rVert_{i\ast}$ is the dual norm of $\lVert\,\cdot\,\rVert_i$ as defined by \eqref{eq:dualnorm}, $i=1,\ldots,d$. This is  the \emph{spectral norm}, special cases of which we have encountered in \eqref{eq:specnorm1} and \eqref{eq:tenspectral} and in Example~\ref{eg:GI}. There are many more cross-norms \cite[Chapter~4]{Diestel} but these three are the best known. The nuclear and spectral norms are also special in that any cross-norm $\lVert \,\cdot \,\rVert$ must satisfy
\begin{equation}\label{eq:smallestlargest}
\lVert T \rVert_\sigma \le \lVert T \rVert \le \lVert T \rVert_\nu \quad\text{for all } T \in \mathbb{U} \otimes \mathbb{V} \otimes \dots \otimes \mathbb{W};
\end{equation}
and conversely any norm $\lVert \,\cdot \,\rVert$  that satisfies \eqref{eq:smallestlargest} must be a cross-norm \cite[Proposition~6.1]{Ryan}. In this sense, the nuclear and spectral norms, respectively, are the largest and smallest cross-norms. If $\mathbb{V}$ is a Hilbert space, then by the Riesz representation theorem, a linear functional $\varphi \colon  \mathbb{V} \to \mathbb{R}$ takes the form $\varphi = \langle v, \cdot \rangle$ for some $v \in \mathbb{V}$ and the Hilbert space norm on $\mathbb{V}$ and its dual norm may be identified. Thus, if $\mathbb{U}, \mathbb{V},\ldots, \mathbb{W}$ are Hilbert spaces, then the spectral norm in \eqref{eq:spectral2} takes the form
\begin{equation}\label{eq:spectral3}
\lVert T\rVert_\sigma = \sup \biggl\{ \dfrac{\lvert \langle T, u \otimes v \otimes \dots \otimes w \rangle \rvert}{\lVert u \rVert_1 \lVert v \rVert_2 \cdots \lVert w \rVert_d }\colon  u \in \mathbb{U}, v \in \mathbb{V},\ldots, w \in \mathbb{W} \biggr\}.
\end{equation}

The definition of inner products on tensor spaces fits perfectly with other tensorial notions such as the tensor product basis and the Kronecker product discussed in Example~\ref{eg:kron}. If $\mathscr{B}_1,\mathscr{B}_2,\ldots,\mathscr{B}_d$ are orthonormal bases for $\mathbb{U},\mathbb{V},\ldots,\mathbb{W}$,  then the tensor product basis $\mathscr{B}_1 \otimes \mathscr{B}_2 \otimes \dots \otimes \mathscr{B}_d$ as defined in \eqref{eq:tenprodbasis2} is an orthonormal basis for $\mathbb{U}  \otimes \mathbb{V} \otimes \dots \otimes \mathbb{W}$. Here we have no need to assume finite dimension or separability:  these orthonormal bases may be uncountable. If $\Phi_1 \colon  \mathbb{U} \to \mathbb{U}, \ldots,\Phi_d \colon  \mathbb{W} \to \mathbb{W}$ are orthogonal linear operators in the sense that
\[
\langle\Phi_1(u),\Phi_1(u')\rangle_1 = \langle u,u'\rangle_1,
 \ldots,
 \langle\Phi_d(w),\Phi_d(w')\rangle_d= \langle w,w'\rangle_d
\]
for all $u,u' \in \mathbb{U},\ldots, w,w'\in \mathbb{W}$, then the Kronecker product $\Phi_1  \otimes \dots \otimes \Phi_d \colon  \mathbb{U}\otimes \dots \otimes \mathbb{W} \to  \mathbb{U}\otimes \dots \otimes \mathbb{W}$ is an orthogonal linear operator, \tie,
\[
 \langle \Phi_1 \otimes \dots \otimes \Phi_d(S), \Phi_1 \otimes \dots \otimes \Phi_d(T) \rangle = \langle S, T\rangle
\]
for all $S,T \in \mathbb{U} \otimes \dots \otimes \mathbb{W}$. For Kronecker products of operators defined on norm spaces, we will defer the discussion to the next example.

Let $\mathbb{U} ,\mathbb{V},\ldots, \mathbb{W}$ be finite-dimensional inner product spaces. Let
$\mathscr{B}_1 ,\mathscr{B}_2,\ldots,\mathscr{B}_d $ be orthonormal bases. Then, as we saw in Example~\ref{eg:hypmatrep}, a $d$-tensor in  $T \in \mathbb{U} \otimes \mathbb{V} \otimes \dots \otimes \mathbb{W}$ may be represented by a $d$-hypermatrix $A \in \mathbb{R}^{m \times n \times \dots \times p}$, in which case the expressions for norms and inner product may be given in terms of $A$:
\[
\begin{aligned}
\lVert A \rVert _\F &=\biggl(\sum_{i=1}^m\sum_{j=1}^n\cdots\sum_{k=1}^p\, \lvert a_{ij \cdots k} \rvert^2 \biggr)^{1/2},\quad \langle A, B\rangle=\sum_{i=1}^m\sum_{j=1}^n\cdots\sum_{k=1}^p a_{ij \cdots k} b_{ij \cdots k},\\
\lVert A \rVert_\nu &= \inf\biggl\{\sum_{i=1}^r \lvert \lambda_i\rvert \colon  
A =  \sum_{i=1}^r \lambda_i u_i \otimes v_i \otimes \dots \otimes w_i, \;\lVert u \rVert = \lVert v \rVert =\dots = \lVert w \rVert = 1\biggr\},\\
\lVert A \rVert _\sigma &= \sup \{\lvert\langle A, u \otimes v \otimes \dots \otimes w \rangle\rvert \colon   \lVert u \rVert = \lVert v \rVert =\dots = \lVert w \rVert = 1 \}. 
\end{aligned} 
\]
Here the nuclear and spectral norms are in fact dual norms (not true in general)
\begin{equation}\label{eq:dualnorms}
\lVert A \rVert_{\nu\ast} =   \sup_{\lVert B \rVert_\nu \le 1} \lvert \langle A, B \rangle \rvert  =\lVert A \rVert_\sigma
\end{equation}
and thus satisfy
\[
\lvert \langle A, B \rangle \rvert \le \lVert A \rVert_\sigma \lVert B \rVert_\nu.
\]
For $d = 2$, these inner product and norms become the usual trace inner product $\langle A, B\rangle=\tr(A^\tp B)$ and spectral, Frobenius, nuclear norms of matrices:
\[
\lVert A \rVert_\nu = \sum_{i=1}^r \sigma_i(A), \quad
\lVert A \rVert _\F = \biggl(\sum_{i=1}^r \sigma_i(A)^2\biggr)^{1/2},\quad
\lVert A \rVert _\sigma = \max_{i=1,\ldots,r} \sigma_i(A),
\]
where $\sigma_i(A)$ denotes the $i$th singular value of $A \in \mathbb{R}^{m \times n}$ and $r = \rank(A)$.

We emphasize that the discussion in the above paragraph requires the bases in Example~\ref{eg:hypmatrep} to be orthonormal. Nevertheless, the values of the inner product and norms do not depend on which orthonormal bases we choose. In the terminology of Section~\ref{sec:trans} they are invariant under multilinear matrix multiplication by $(X_1,X_2,\ldots,X_d) \in \Or(m) \times \Or(n) \times \dots \times \Or(p)$, or equivalently, they are defined on Cartesian tensors.
More generally, if $\langle X_i v, X_i v \rangle_i = \langle  v, v \rangle_i$, $i=1,\ldots,d$, then
\[
\langle (X_1,X_2,\ldots,X_d) \cdot A, (X_1,X_2,\ldots,X_d) \cdot B\rangle = \langle A, B\rangle,
\]
and thus $\lVert (X_1,X_2,\ldots,X_d) \cdot A \rVert _\F = \lVert A \rVert _\F$; if $\lVert X_i v \rVert_i = \lVert v \rVert_i$, $i=1,\ldots,d$, then
\[
\lVert (X_1,X_2,\ldots,X_d) \cdot A \rVert_\nu = \lVert A \rVert_\nu, \quad \lVert (X_1,X_2,\ldots,X_d) \cdot A \rVert _\sigma = \lVert A \rVert _\sigma.
\]
This explains why, when discussing definition~\ref{st:tensor1} in conjunction with inner products or norms, we expect the change-of-basis matrices in the transformation rules to preserve these inner products or norms.
\end{example}

Inner product and norm spaces become enormously more interesting when they are completed into Hilbert and Banach spaces. The study of cross-norms was in fact started by \citet{Grothendieck2} and \citet{Schatten} in order to define tensor products of Banach spaces, and this has grown into a vast subject \cite{Defant,Diestel,Cheney,Ryan,Treves} that we are unable to survey at any reasonable level of detail. The following example is intended to convey an idea of how cross-norms allow one to complete the tensor spaces  in a suitable manner.

\begin{example}[tensor\zzb product\zzb of\zzb Hilbert\zzb and\zzb Banach\zzb spaces]\label{eg:Hilbert}\hspz%
  We revisit our\linebreak
  discussion in Example~\ref{eg:infinite},
properly defining topological tensor products $\hatotimes$
  using continuous and integrable functions as illustrations.
\begin{enumerate}[\upshape (i)]
\setlength\itemsep{3pt}
\item If $X$ and $Y$ are compact Hausdorff topological spaces, then by the Stone--Weierstrass theorem, $C(X) \otimes C(Y)$ is a dense subset of $C(X \times Y)$ with respect to the uniform norm
\begin{equation}\label{eq:unifnorm}
\lVert f \rVert_\infty = \sup_{(x,y) \in X \times Y} \, \lvert f(x,y) \rvert.
\end{equation}

\item\label{it:int} If $X$ and $Y$ are $\sigma$-finite measure spaces, then by the  Fubini--Tonelli theorem, $L^1(X) \otimes L^1(Y)$ is a dense subset of $L^1(X \times Y)$ with respect to the $L^1$-norm
\begin{equation}\label{eq:L1norm}
\lVert f \rVert_1 = \int_{X \times Y} \lvert f(x,y) \rvert  \D x \D y.
\end{equation}
\end{enumerate}
Here a tensor product of vector spaces is as defined in \eqref{eq:tenprodvecsp}, and as we just saw in \eqref{eg:norms}, there are several ways to equip it with a norm and with respect to any norm we may  complete it (\ie\  by  adding the limits of all Cauchy sequences) to obtain Banach and Hilbert spaces out of norm and inner product spaces. The completed space depends on the choice of norms;  with a judicious choice, we get
\begin{equation}\label{eq:ctsint}
C(X) \hatotimes_\sigma C(Y) = C(X \times Y),\quad L^1(X) \hatotimes_\nu L^1(Y) = L^1(X \times Y),
\end{equation}
as was first discovered in \citet{Grothendieck2}. Here $\hatotimes_\sigma$ and $\hatotimes_\nu $ denote completion in the spectral and nuclear norm respectively and are called the \emph{injective tensor product} and \emph{projective tensor product} respectively. To be clear, the first equality in \eqref{eq:ctsint} says that if we equip  $C(X) \otimes C(Y)$ with the spectral norm \eqref{eq:spectral2} and complete it to obtain $C(X) \hatotimes_\sigma C(Y)$, then the resulting space is $ C(X \times Y)$ equipped with the uniform norm \eqref{eq:unifnorm}, and likewise for the second equality. In particular, \eqref{eq:ctsint} also tells us that the uniform and spectral norms are equal on $C(X\times Y)$, and likewise for the $L^1$-norm in \eqref{eq:L1norm} and nuclear norm in \eqref{eq:nuclear2}. For $f \in L^1(X \times Y)$,
\begin{align*}
  &\int_{X \times Y} \lvert f(x,y) \rvert  \D x \D y \\*
  & \quad = \inf \biggl\{ \sum_{i=1}^r \int_X \lvert \varphi_i(x) \rvert  \D x\int_Y \lvert \psi_i(y) \rvert  \D y\colon  
 {\varphi_i \in L^1(X), \; \psi_i \in L^1(Y), \; r \in \mathbb{N} \biggr\},}
\end{align*}
a fact that can be independently verified by considering simple functions\footnote{Finite linear combinations of characteristic functions of measurable sets.} and taking limits.

These are examples of topological tensor products that involve completing the (algebraic) tensor product in \eqref{eq:tenprodvecsp} with respect to a choice of cross-norm to obtain a complete topological vector space. These also suggest why it is desirable to have a variety of different cross-norms, and with each a different topological tensor product, as the `right' cross-norm to choose for a class of functions $\Fn(X)$ is usually the one that gives us $\Fn(X) \hatotimes \Fn(Y) = \Fn(X \times Y)$ as we discussed after \eqref{eq:fun}. For example, to get the corresponding result for $L^2$-functions, we have to use the Hilbert--Schmidt norm
\begin{equation}\label{eq:L2tenprod}
L^2(X) \hatotimes_\F L^2(Y) = L^2(X \times Y).
\end{equation}
Essentially the proof relies on the fact that the completion of $L^2(X) \otimes L^2(Y)$ with respect to the  Hilbert--Schmidt norm is the closure of $L^2(X) \otimes L^2(Y)$ as a subspace of $L^2(X \times Y)$, and as the orthogonal complement of $L^2(X) \otimes L^2(Y)$ is zero, its closure is the whole space \cite[Theorem~1.39]{Cheney}. Nevertheless, such results are not always possible: there is no cross-norm that will complete $L^\infty(X) \otimes L^\infty(Y)$ into $L^\infty(X \times Y)$ for all $\sigma$-finite $X$ and $Y$ \cite[Theorem~1.53]{Cheney}. Also, we should add that a `right' cross-norm that guarantees \eqref{eq:ctsint} may be less interesting than a `wrong' cross-norm that gives us a new tensor space. For instance, had we used $\hatotimes_\nu$ to form a tensor product of $C(X)$ and $C(Y)$, we would have obtained a smaller subset:
\[
C(X) \hatotimes_\nu C(Y) \subsetneq C(X \times Y).
\]
This smaller tensor space of continuous functions on $X \times Y$, more generally the tensor space $C(X_1) \hatotimes_\nu \dots \hatotimes_\nu  C(X_n)$, is called the \emph{Varopoulos algebra} and it turns out to be very interesting and useful in harmonic analysis \cite{Varopoulos1,Varopoulos2}.

A point worth highlighting is the difference between $L^2(X) \otimes L^2(Y)$ and $L^2(X) \hatotimes_\F L^2(Y) $. While the former contains only \emph{finite} sums of separable functions
\[
\sum_{i=1}^r f_i \otimes g_i
\]
with $f_i \in L^2(X)$, $g_i \in L^2(Y)$, the latter includes all convergent infinite series of separable functions
\[
\sum_{i=1}^\infty f_i \otimes g_i \coloneqq \lim_{r\to\infty} \sum_{i=1}^r f_i \otimes g_i,
\]
\tie, the right-hand side converges to some limit in $L^2(X \times Y)$ in the Hilbert--Schmidt norm $\lVert\,\cdot\,\rVert_\F$.  The equality in \eqref{eq:L2tenprod} says that every function in $L^2(X \times Y)$ is given by such a limit, but if we had taken completion with respect to some other norms such as $\lVert\,\cdot\,\rVert_\nu$ or $\lVert\,\cdot\,\rVert_\sigma$, their topological tensor products $L^2(X) \hatotimes_\nu L^2(Y)$ or $L^2(X) \hatotimes_\sigma L^2(Y)$ would in general be smaller or larger than $L^2(X \times Y)$ respectively.
Whatever the choice of $\hatotimes$, the (algebraic) tensor product $\mathbb{U} \otimes \mathbb{V} \otimes \dots \otimes \mathbb{W}$ should be regarded as the subset of all finite-rank tensors in the topological tensor product $\mathbb{U} \hatotimes \mathbb{V} \hatotimes \dots \hatotimes \mathbb{W}$.
\end{example}

The tensor product of Hilbert spaces invariably refers to the topological tensor product with respect to the Hilbert--Schmidt norm because the result is always a Hilbert space. In particular, the $\hatotimes$ in \eqref{eq:mv3} should be interpreted as $\hatotimes_\F$. However, as we pointed out above, other topological tensor products with respect to other cross-norms may also be very interesting.

\begin{example}[trace-class, Hilbert--Schmidt,  compact operators]\label{eg:compact}
For a separable Hilbert space $\mathbb{H}$ with inner product $\langle \,\cdot\,, \cdot\,\rangle$ and induced norm $\lVert \, \cdot \, \rVert$, completing $\mathbb{H} \otimes \mathbb{H}^*$ with respect to the nuclear, Hilbert--Schmidt and spectral norms, we obtain different types of bounded linear operators on $\mathbb{H}$ as follows:
\begin{alignat*}{3}
&\text{trace-class} &  \mathbb{H} \hatotimes_\nu \mathbb{H}^*&= \biggl\{ \Phi \in \Bd (\mathbb{H}) \colon  \sum_{i \in I}\sum_{j \in I} \, \lvert \langle \Phi  (e_i), f_j \rangle \rvert < \infty \biggr\},\\*
&\text{Hilbert--Schmidt}\quad &  \mathbb{H} \hatotimes_\F \mathbb{H}^*&= \biggl\{ \Phi \in \Bd (\mathbb{H})  \colon  \sum_{i \in I} \, \lVert \Phi(e_i) \rVert^2 < \infty \biggr\},\\
  &\text{compact} &  \mathbb{H} \hatotimes_\sigma \mathbb{H}^*&=\biggl\{ \Phi \in \Bd (\mathbb{H}) \colon  \begin{aligned}
    & X \subseteq\mathbb{H} \text{ bounded}\\*
&\Rightarrow \overline{\Phi(X)} \subseteq\mathbb{H} \text{ compact} \end{aligned} \biggr\}.
\end{alignat*}
The series convergence is understood to mean `for some orthonormal bases $\{e_i \colon  i \in I\}$ and $\{f_i \colon  i \in I\}$ of $\mathbb{H}$', although for trace-class operators the condition could be simplified to $\sum_{i \in I}\, \langle \Phi  (e_i), e_j \rangle < \infty $ provided that `for some' is replaced by `for all'. See \citet[p.~278]{Schaefer} for the trace-class result  and \citet[Theorems~48.3]{Treves} for the compact result.

If $\mathbb{H}$ is separable, then such operators are  characterized by their having a \emph{Schmidt decomposition}:
\begin{equation}\label{eq:schmidt}
\Phi = \sum_{i=1}^\infty \sigma_i  u_i \otimes v_i^*
\quad\text{or}\quad
\Phi(x) = \sum_{i=1}^\infty \sigma_i \langle v_i, x \rangle u_i \quad \text{for all } x \in \mathbb{H},
\end{equation}
where  $\{ u_i \colon  i \in \mathbb{N}\}$ and $\{ v_i \colon  i \in \mathbb{N}\}$ are orthonormal sets, $\sigma_i \ge 0$,  and
\[
\sum_{i=1}^\infty \sigma_i < \infty,\quad 
\sum_{i=1}^\infty \sigma_i^2 < \infty, \quad 
\lim_{i \to \infty} \, \sigma_i = 0
\]
depending on whether they are trace-class, Hilbert--Schmidt or compact, respectively \cite[Chapter~VI, Sections 5--6]{Reed}. Note that $\mathbb{H}$ and $\mathbb{H}^*$ are naturally isomorphic by the Riesz representation theorem and we have identified $v^*$ with $\langle v, \cdot\, \rangle$. As expected, we also have
\[
\lVert\Phi \rVert_\nu = \sum_{i=1}^\infty \sigma_i , \quad 
\lVert\Phi \rVert_\F =\biggl(\sum_{i=1}^\infty \sigma_i^2\biggr)^{1/2}, \quad 
\lVert\Phi \rVert_\sigma = \sup_{i \in \mathbb{N}} \, \sigma_i.
\]
In other words, this is the infinite-dimensional version of the relation between the various norms and matrix singular values in Example~\ref{eg:norms}, and the Schmidt decomposition is an infinite-dimensional generalization of singular value decomposition. Unlike the finite-dimensional case, where one may freely speak of the nuclear, Frobenius and spectral norms of any matrix, for infinite-dimensional $\mathbb{H}$ we have
\[
\mathbb{H} \otimes \mathbb{H}^* \subsetneq 
\mathbb{H} \hatotimes_\nu \mathbb{H}^* \subsetneq  \mathbb{H} \hatotimes_\F \mathbb{H}^* \subsetneq  \mathbb{H} \hatotimes_\sigma \mathbb{H}^*, \quad \lVert\Phi \rVert_\sigma \le \lVert\Phi \rVert_\F \le \lVert\Phi \rVert_\nu,
\]
and the inclusions are strict, for example, a compact operator $\Phi$ may have $\lVert\Phi \rVert_\nu =\infty$.
By our discussion at the end of Example~\ref{eg:Hilbert}, $\mathbb{H} \otimes \mathbb{H}^*$ is the subset of finite-rank operators in any of these larger spaces.
The inequality relating the three norms is a special case of \eqref{eq:smallestlargest} as $\lVert\,\cdot\, \rVert_\F$ is a cross-norm. On $\mathbb{H} \hatotimes_\nu \mathbb{H}^*$, nuclear and spectral norms are dual norms as in the finite-dimensional case \eqref{eq:dualnorms}:
\[
\lVert \Phi \rVert_{\nu\ast} =   \sup_{\lVert \Psi \rVert_\nu \le 1} \lvert \langle \Phi, \Psi \rangle_\F \rvert  =\lVert \Phi \rVert_\sigma,
\]
where the inner product will be defined below in \eqref{eq:HSip}.

The reason for calling $\Phi \in \mathbb{H} \hatotimes_\nu \mathbb{H}^*$ trace-class is that its \emph{trace} is well-defined and always finite:
\[
\tr(\Phi) \coloneqq \sum_{i=1}^\infty  \langle \Phi(e_i), e_i \rangle
\]
for any orthonormal basis $\{e_i \colon  i \in \mathbb{N} \}$ of $\mathbb{H}$. As trace is independent of the choice of orthonormal basis, it is truly a tensorial notion. Moreover $\tr(\lvert\Phi\rvert) = \lVert \Phi \rVert_\nu$, where $\lvert \Phi \rvert \coloneqq \sqrt{\Phi^* \Phi}$. To see why the expression above is considered a trace, we may, as discussed in Example~\ref{eg:hypmatrep}, represent $\Phi$ as an infinite-dimensional matrix with respect to the chosen orthonormal basis,
\begin{equation}\label{eq:infmatrep}
\Phi = \sum_{i=1}^\infty \sum_{j=1}^\infty a_{ij}\, e_i \otimes e_j^*,
\end{equation}
where $A = (a_{ij})_{i,j=1}^\infty \in l^2(\mathbb{N}^2)$, and now observe that
\begin{equation}\label{eq:trace}
\tr(\Phi)=
\sum_{k=1}^\infty  \langle \Phi(e_k), e_k \rangle =
\sum_{k=1}^\infty  \biggl\langle \sum_{i=1}^\infty  a_{ik} e_i  , e_k \biggr\rangle =
\sum_{k=1}^\infty  a_{kk}.
\end{equation}

One may show that an operator is trace-class if and only if it is a product of two Hilbert--Schmidt operators. A consequence is that
\begin{equation}\label{eq:HSip}
\langle \Phi, \Psi \rangle_\F \coloneqq \tr(\Phi^*\Psi)
\end{equation}
is always finite and defines an inner product on $\mathbb{H} \hatotimes_\F \mathbb{H}^*$ that gives  $\lVert\,\cdot\, \rVert_\F$ as its induced  norm:
\[
\langle \Phi, \Phi \rangle_\F = \tr(\Phi^*\Phi) = \lVert\Phi \rVert_\F^2.
\]

While $(\mathbb{H} \hatotimes_\nu \mathbb{H}^*, \lVert\,\cdot\, \rVert_\nu )$, $(\mathbb{H} \hatotimes_\F \mathbb{H}^*, \lVert\,\cdot\, \rVert_\F )$, $(\mathbb{H} \hatotimes_\sigma \mathbb{H}^*, \lVert\,\cdot\, \rVert_\sigma )$ are all Banach spaces, only $\mathbb{H} \hatotimes_\F \mathbb{H}^* $ is a Hilbert space with the inner product in \eqref{eq:HSip}. So a topological tensor product of Hilbert spaces may or may not be a Hilbert space: it depends on the choice of cross-norms. Also, the choice of norm is critical in determining whether we have a Banach space: for instance, $\mathbb{H} \hatotimes_\nu \mathbb{H}^*$ is complete with respect to $\lVert\,\cdot\, \rVert_\nu$ but not $\lVert\,\cdot\, \rVert_\sigma$.

Taking topological tensor products of infinite-dimensional spaces has some unexpected properties. When $ \mathbb{H} $ is infinite-dimensional, the identity operator $I_\mathbb{H} \in \Bd(\mathbb{H})$ is not in $\mathbb{H} \hatotimes \mathbb{H}^*$ regardless of whether we use $\hatotimes_\nu$, $\hatotimes_\F$ or $\hatotimes_\sigma$, as $I_\mathbb{H}$ is non-compact (and thus neither Hilbert--Schmidt nor trace-class). Hence $\mathbb{H} \hatotimes \mathbb{H}^* \ne \Bd(\mathbb{H})$, no matter which notion of topological tensor product $\hatotimes$ we use.
Incidentally, the strong and weak duals of $\Bd(\mathbb{H})$, \ie\ the set of linear functionals on $\Bd(\mathbb{H})$ continuous under the strong and weak operator topologies, are equal to each other and to $\mathbb{H} \otimes \mathbb{H}^*$, \ie\ the set of finite-rank operators \cite[Chapter~IV, Section~1]{Dunford}.
Taking the (continuous) dual of the compact operators gives us the trace-class operators
\begin{equation}\label{eq:dualcpttr}
(\mathbb{H} \hatotimes_\sigma \mathbb{H}^*)^* \cong \mathbb{H} \hatotimes_\nu \mathbb{H}^*,
\end{equation}
so taking dual can change the topological tensor product.

We have framed our discussion in terms of $\mathbb{H} \hatotimes \mathbb{H}^*$ as this is the most common scenario, but it may be extended to other situations.

\mypara{Different Hilbert spaces} It applies verbatim to $\mathbb{H}_1 \hatotimes \mathbb{H}_2^*$ as operators in $\Bd(\mathbb{H}_1; \mathbb{H}_2)$.

\mypara{Covariant and contravariant} We covered the mixed case $\mathbb{H}_1 \hatotimes \mathbb{H}_2^*$, but the covariant case $\mathbb{H}_1^* \hatotimes \mathbb{H}_2^*$ and contravariant case $\mathbb{H}_1 \hatotimes \mathbb{H}_2$ follow from replacing $\mathbb{H}_1$ with $\mathbb{H}_1^*$ or $\mathbb{H}_2$ with $\mathbb{H}_2^*$.

\mypara{Banach spaces} Note that \eqref{eq:nuclear2} and \eqref{eq:spectral2} are defined without any reference to inner products, so for $\mathbb{H}_1 \hatotimes_\sigma \mathbb{H}_2^*$ and $\mathbb{H}_1 \hatotimes_\nu \mathbb{H}_2^*$ we do not need $\mathbb{H}_1$ and  $\mathbb{H}_2$ to be Hilbert spaces;  compact and trace-class operators may be defined for any pair of Banach spaces, although the latter are usually called nuclear operators in this context.

\mypara{Higher order} One may define order-$d \ge 3$ analogues of bounded, compact, Hilbert--Schmidt, trace-class operators in a straightforward manner, but corresponding results are more difficult;  one reason is that Schmidt decomposition \eqref{eq:schmidt} no longer holds \cite{Benyi,CKP1,CKP2}.
\end{example}

We remind the reader that many of the topological vector spaces considered in Examples~\ref{eg:infinite} and \ref{eg:distributions}, such as $S(X)$, $C^\infty(X)$, $C^\infty_c(X)$,
$H(X)$, $S'(X)$, $E'(X)$, $D'(X)$ and $H'(X)$, are not Banach spaces (and thus not Hilbert spaces) but so-called \emph{nuclear spaces}. Nevertheless, similar ideas apply to yield topological tensor products. In fact nuclear spaces have the nice property that topological tensor products with respect to the nuclear and spectral norms, \ie\ $\hatotimes_\nu$ and $\hatotimes_\sigma$, are always equal and thus independent of the choice of cross-norms.

We mention two applications of Example~\ref{eg:compact} in machine learning (Mercer kernels) and quantum mechanics (density matrices).

\begin{example}[Mercer kernels]\label{eg:Mercer}
Hilbert--Schmidt operators, in the form of Mercer kernels \cite{Mercer}, are of special significance to machine learning. Here $\mathbb{H} = L^2(X)$ for a set $X$ of machine learning significance (where a training set is sampled from) that has some appropriate measure (usually Borel) and topology (usually locally compact). For any $L^2$-integral kernel $K \colon  X \times X \to \mathbb{R}$, \ie~$K \in L^2(X \times X) = L^2(X) \hatotimes_\F L^2(X)$, the integral operator $\Phi \colon  L^2(X) \to L^2(X)$, 
\[
\Phi(f)(x)  = \int_X K(x,y) f(y) \D y,
\]
is Hilbert--Schmidt; in fact every Hilbert--Schmidt operator on $L^2(X)$ arises in this manner and we always have $\lVert \Phi \rVert_\F = \lVert K \rVert_{L^2}$ \cite[p.~267]{Conway}. By our discussion above, $\Phi$ has a Schmidt decomposition as in \eqref{eq:schmidt}:
\begin{equation}\label{eq:HS1}
\Phi = \sum_{i=1}^\infty \sigma_i \varphi_i \otimes \psi_i^*,
\end{equation}
where $\psi^*(f) = \langle \psi, f\rangle = \int_X \psi(x) f(x)  \D x$ by the Riesz representation. Hence
\begin{equation}\label{eq:HS2}
K = \sum_{i=1}^\infty \sigma_i \varphi_i \otimes \psi_i ,
\end{equation}
and we see that the difference between \eqref{eq:HS1} and \eqref{eq:HS2} is just one of covariance and contravariance:  every $K \in  L^2(X) \hatotimes_\F L^2(X)$ gives us a $\Phi \in L^2(X) \hatotimes_\F L^2(X)^*$ and {\em vice versa}. Mercer's kernel theorem is essentially the statement that if $K$ is continuous and symmetric positive semidefinite, \ie\ $[K(x_i,x_j)]_{i,j=1}^n \in \mathbb{S}^n_\p$ for any $x_1,\ldots,x_n \in X$ and any $n \in \mathbb{N}$, then $\varphi_i $ and $\psi_i$ may be chosen to be equal in \eqref{eq:HS1} and \eqref{eq:HS2}. In this case we obtain a \emph{feature map}:
\[
F \colon  X \to l^2(\mathbb{N}), \quad x \mapsto ( \sqrt{\sigma_i} \varphi_i(x) )_{i=1}^\infty,
\]
and in this context $l^2(\mathbb{N})$ is called the \emph{feature space}.  It follows from \eqref{eq:HS2} that $K(x,y) = \langle F(x), F(y) \rangle_{l^2(\mathbb{N})}$.  Assuming that $\sigma_i > 0$ for all $i \in \mathbb{N}$, the subspace of $L^2(X)$ defined by
\[
L^2_K(X) \coloneqq\biggl\{ \sum_{i=1}^\infty a_i \varphi_i \in L^2(X) \colon  ( a_i/\sqrt{\sigma_i} )_{i=1}^\infty \in l^2(\mathbb{N}) \biggr\},
\]
equipped with the inner product
\[
\biggl\langle \sum_{i=1}^\infty a_i \varphi_i , \sum_{i=1}^\infty b_i \varphi_i \biggr\rangle_K \coloneqq \sum_{i=1}^\infty \dfrac{a_i b_i }{ \sigma_i} = a^\tp \Sigma^{-1} b,
\]
is called the \emph{reproducing kernel Hilbert space} associated with the kernel $K$ \cite{Smale}. The last expression is given in terms of  infinite-dimensional matrices, as described in Example~\ref{eg:hyp}, \ie\ $a,b \in l^2(\mathbb{N})$ and $\Sigma \in l^2(\mathbb{N} \times \mathbb{N})$, to show its connection with the finite-dimensional case.

These notions are used to great effect in machine learning \cite{Hofmann,Scholkopf,Steinwart}. For instance, as depicted in Figure~\ref{fig:svm}, the feature map allows one to unfurl a complicated set of points $x_1,\ldots,x_n$ in a lower-dimensional space $X$ into a simpler set of points $F(x_1),\ldots,F(x_n)$, say, one that may be partitioned with a hyperplane, in an infinite-dimensional feature space $l^2(\mathbb{N})$.
\begin{figure}[ht]
\centering
\begin{tikzpicture}[>=stealth',x=1cm,y=1cm,scale=0.9]
  
%draw[color=gray] (0,0) grid (6,6);
\draw (0,0) rectangle (6,6) node at (5.7,5.65) {\normalsize$X$};
% \draw line
\draw[color=blue,line width=2pt]
  (2,6) .. controls (3,5.5) and (3,5) .. 
  (3,5) .. controls (3,4) and (2,2.5) .. 
  (2,2) .. controls (2,1) and (2.8,1) .. 
  (3,1) .. controls (3.5,1) and (3.5,2) .. 
  (4,2) .. controls (4.5,2) and (6,0) .. 
  (6,0);
% \draw left dashed line
\draw[dashed] 
  (1.5,6) .. controls (2.5,5.5) and (2.5,5) .. 
  (2.5,5) .. controls (2.5,4) and (1.5,2.5) .. 
  (1.5,2) .. controls (1.5,.5) and (2.8,.5) .. 
  (3,.5) .. controls (3.75,.5) and (3.5,1.5) .. 
  (4,1.5) .. controls (4.5,1.5) and (5.5,0) .. 
  (5.5,0);
% \draw right dashed line
\draw[dashed] 
  (2.5,6) .. controls (3.5,5.5) and (3.5,5) .. 
  (3.5,5) .. controls (3.5,4) and (2.5,2.5) .. 
  (2.5,2) .. controls (2.5,1.5) and (2.8,1.5) .. 
  (3,1.5) .. controls (3.25,1.5) and (3.5,2.5) .. 
  (4,2.5) .. controls (4.5,2.5) and (6,0.5) .. 
  (6,0.5);
%\draw[color=gray] (2,6) -- (3,5) -- (2,2) -- (3,1) -- (4,2) -- (6,0);
%\draw[color=gray] (1.5,6) -- (2.5,5) -- (1.5,2) -- (3,.5)-- (4,1.5)-- (5.5,0);
%\draw[color=gray] (2.5,6) -- (3.5,5) -- (2.5,2) -- (3,1.5)-- (4,2.5)-- (6,0.5);
  
%\draw[color=gray] (7,0) grid (13,6);
\draw (7,0) rectangle (13,6) node at (12.4,5.6) {\normalsize$l^2(\mathbb{N})$};
% \draw line
\draw[color=blue,line width=2pt] (8.5,6) -- (12,0);
% \draw dashed line
\draw[dashed]  (8,6) -- (11.5,0);
\draw[dashed]  (9,6) -- (12.5,0);
  
\draw[->,thick] (5,3) -- (8,3) node [above,pos=.5] {\normalsize$F$};
  
\def\positive{{%
{2.3,5.3},
{3.5,.7},
{1.5,2},
{1.2,2.1},
{1.8,.8},
{1,5.5},
{1.2,5.8},
{.75,.2},
{2,4},
{5, 0.5},
{1.5,3},
{2.3,.5},
{9.3,3.3},
{11,.8},
{8.5,2},
{7.2,4.1},
{8.8,.8},
{8,5.5},
{8.2,5},
{7.75,.2},
{9,4.2},
{12, 0.5},
{8.5,3},
{9.3,.5},
}}
  
% \draw positive dots
\foreach \i in {0,...,20} {
  \pgfmathparse{\positive[\i][0]}\let \x \pgfmathresult;
  \pgfmathparse{\positive[\i][1]}\let \y \pgfmathresult;
  \fill[black] (\x,\y) circle (2pt);
}
  
\def\negative{{%
{4,2.5},
{3.5,5},
{2.6,1.6},
{4.5,5.2},
{5.5,3.7},
{3.9,4.7},
{5,2.7},
{3.5,4.2},
{5.8,.9},
{10.75,3},
{10.5,5},
{11.6,1.6},
{11.5,5.2},
{12.5,3.7},
{10.9,4.7},
{12,2.7},
{10.5,4.2},
{12.8,.9},
}}
  
% \draw negative dots
\foreach \i in {0,...,16} {
  \pgfmathparse{\negative[\i][0]}\let \x \pgfmathresult;
  \pgfmathparse{\negative[\i][1]}\let \y \pgfmathresult;
  \draw[black] (\x,\y) circle (2pt);
}
  
\end{tikzpicture}
\caption{Depiction of the feature map $F \colon  X \to l^2(\mathbb{N})$ in the context of support-vector machines.}
\label{fig:svm}
\end{figure}
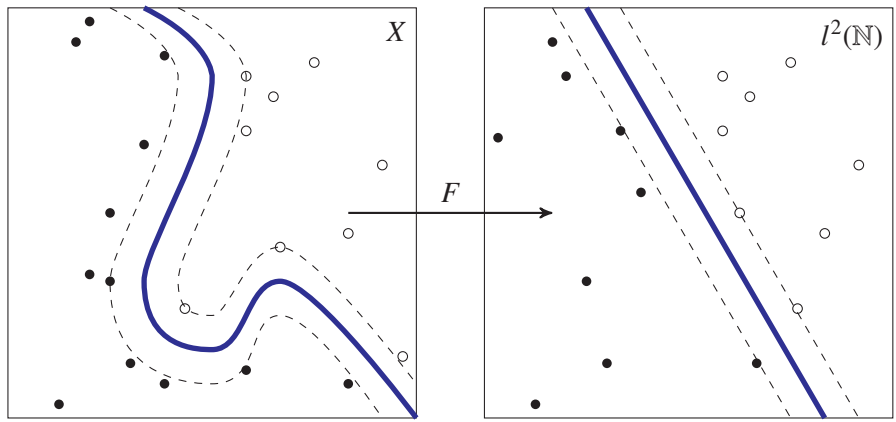
\end{example}

The previous example is about symmetric positive semidefinite Hilbert--Schmidt operators; the next one will be about Hermitian positive semidefinite trace-class operators.

\begin{example}[density operators]\label{eg:density}\hspace{-1.9pt}%
These are sometimes called density matrices, and are important alike in quantum physics \cite{Blum}, quantum chemistry \cite{Davidson} and quantum computing \cite{Nielsen}. Let $\mathbb{H}$ be a Hilbert space, assumed separable for notational convenience, and let $\langle\,\cdot\,,\cdot\,\rangle$ be its inner product. A trace-class operator $\rho \in \mathbb{H} \hatotimes_\nu \mathbb{H}^*$ is called a density operator if it is Hermitian $\rho^* = \rho$, positive semidefinite $\lvert \rho \rvert = \rho$, and of unit trace $\tr(\rho) = 1$ (or, equivalently, unit nuclear norm as $\tr(\rho) = \lVert \rho \rVert_\nu$ in this case). By the same considerations as in Example~\ref{eg:Mercer}, its Schmidt decomposition may be chosen to have $u_i = v_i$, \tie,
\begin{equation}\label{eq:mixedstate}
\rho = \sum_{i=1}^\infty \sigma_i v_i \otimes v_i^*, \quad \tr(\rho) = \sum_{i=1}^\infty  \sigma_i = 1.
\end{equation}
A key reason density operators are important is that the postulates of quantum mechanics may be reformulated in terms of them \cite[Section~2.4.2]{Nielsen}.  Instead of representing a quantum state as a one-dimensional subspace of $\mathbb{H}$ spanned by a non-zero vector $v \in \mathbb{H}$, we let it be represented by the projection operator onto that subspace, \tie,
\[
v \otimes v^* \colon  \mathbb{H} \to \mathbb{H}, \quad w \mapsto v^*(w)v = \langle v, w\rangle v,
\]
where we have assumed that $v$ has unit norm (if not,  normalize).
This is evidently a density operator and such rank-one density operators are called \emph{pure states}. An advantage is that one may now speak of \emph{mixed states}: they are infinite convex linear combination of pure states as in \eqref{eq:mixedstate}, \ie\ the density operators. Note that this is not possible with states represented as vectors in $\mathbb{H}$:  a linear combination of two vectors would just be another vector in $\mathbb{H}$. There is a well-known condition for telling pure states apart from mixed states without involving a Schmidt decomposition, namely, $\rho$ is pure if and only if $\tr(\rho^2) = 1$ and $\rho$ is mixed if and only if $\tr(\rho^2) < 1$. Consequently, the value of $\tr(\rho^2)$ is called the purity.

In this density operator formulation, a quantum system is completely described by a single density operator, and for $d$ quantum systems described by density operators $\rho_1,\ldots,\rho_d$, the composite system is described by their  Kronecker product $\rho_1 \otimes \dots \otimes \rho_d$ as defined in Example~\ref{eg:kron}. The reader should compare these with \ref{it:state} and \ref{it:composite} on  page~\pageref{pg:postulates}. A Kronecker product of density operators remains a density operator because of \eqref{eq:kronprodprop3}, \eqref{eq:kronprodprop5} and the fact that \eqref{eq:kronprodtr} holds for trace-class operators:
\[
\tr(\rho_1 \otimes \dots \otimes \rho_d) = \tr(\rho_1) \cdots \tr(\rho_d).
\]

We end  with a word about the  role of compact operators $\mathbb{H} \hatotimes_\sigma \mathbb{H}^*$ in this discussion. For any density operator $\rho$, the linear functional $\omega_\rho \colon \Bd(\mathbb{H})  \to \mathbb{C}$, $\Phi \mapsto \tr(\rho \Phi)$ satisfies $\omega_\rho(I_\mathbb{H}) = 1$ and $\omega_\rho(\lvert \Phi \rvert) \ge 0$ for all $\Phi \in  \Bd(\mathbb{H})$. A linear functional $\omega$ with the last two properties is called a \emph{positive linear functional}, abbreviated $\omega \succeq 0$. It turns out that every positive linear functional on the space of compact operators is of the form $\omega_\rho$ for some density operator $\rho$  \cite[Corollaries~4.14 and 4.15]{Landsman}:
\[
\{ \omega \in  (\mathbb{H} \hatotimes_\sigma \mathbb{H}^*)^* \colon  \omega \succeq 0 \}
\cong \{ \rho \in \mathbb{H} \hatotimes_\nu \mathbb{H}^* \colon  \rho \text{ density operator} \},
\]
a very useful addendum to \eqref{eq:dualcpttr}. With this observation, we may now use positive linear functionals to represent quantum states. The advantage is that it permits ready generalization: $\mathbb{H} \hatotimes_\sigma \mathbb{H}^*$ has the structure of a $C^*$-algebra and one may replace it with other $C^*$-algebra  $\mathcal{A}$ and represent quantum states as positive linear functionals $\omega \colon  \mathcal{A} \to \mathbb{C}$.
\end{example}

The evolution of the definition of quantum states from vectors in Hilbert spaces to density operators to positive linear functionals on $C^*$-algebras is not unlike the three increasingly sophisticated definitions of tensors or the three increasingly abstract definitions of tensor products. As in the case of tensors and tensor products, each definition of quantum states is useful in its own way, and all three remain in use.

Our discussion of tensor product will not be complete without mentioning the tensor algebra. We will make this the last example of this section.

\begin{example}[tensor algebra]\label{eg:tenalg}
Let $\mathbb{V}$ be a vector space and $v \in \mathbb{V}$. We introduce the shorthand
\[
\mathbb{V}^{\otimes d} \coloneqq \overbracket[0.5pt]{\mathbb{V}  \otimes\dots\otimes \mathbb{V}}^{d \text{ copies}}, \quad v^{\otimes d} \coloneqq \overbracket[0.5pt]{v \otimes\dots\otimes v}^{d \text{ copies}}
\]
for any $d \in \mathbb{N}$. We also define $\mathbb{V}^{\otimes 0} \coloneqq \mathbb{R}$ and $v^{\otimes 0} \coloneqq 1$. There is a risk of construing erroneous relations  from notation like this. We caution that $\mathbb{V}^{\otimes d}$ is \emph{not} the set of tensors of the form $v^{\otimes d}$, as we pointed out in \eqref{eq:pitfall}, but neither is it the set of linear combinations of such tensors;  $\mathbb{V}^{\otimes d}$ contains all finite linear combinations of any $v_1 \otimes \dots \otimes v_d$, including but not limited to those of the form $v^{\otimes d}$. As $\mathbb{V}^{\otimes d}$, $d = 0, 1,2,\ldots,$ are all vector spaces, we may form their direct sum to obtain an infinite-dimensional vector space:
\begin{equation}\label{eq:tenalg0}
\Ten(\mathbb{V}) = \bigoplus_{k=0}^\infty \mathbb{V}^{\otimes k} = \mathbb{R} \oplus \mathbb{V} \oplus \mathbb{V}^{\otimes 2} \oplus \mathbb{V}^{\otimes 3} \oplus \cdots ,
\end{equation}
called the \emph{tensor algebra} of $\mathbb{V}$. Those unfamiliar with direct sums may simply regard  $\Ten(\mathbb{V}) $ as the set of all \emph{finite} sums of the tensors of any order:
\begin{equation}\label{eq:tenalg1}
\Ten(\mathbb{V})  = \biggl\{\sum_{k=0}^d T_k \colon  T_k \in \mathbb{V}^{\otimes k}, \; d \in \mathbb{N} \biggr\}.
\end{equation}
What does it mean to add a tensor $T \in \mathbb{V}^{\otimes 2}$ to another tensor $T' \in \mathbb{V}^{\otimes 5}$, \ie\ of different orders? This is exactly like the addition of polyads we discussed at the beginning of this section.  Namely there is no formula for actually adding $T$ to $T'$, a sum like $T + T'$ cannot be further simplified, it is just intended as a place holder for the two tensors $T$ and $T'$, we could have written instead $(0,0,T,0,0,T',0,\ldots)$, and indeed another possible way to view $\Ten(\mathbb{V})$ is that it is the set of all sequences with finitely many non-zero terms:
\[
\{ (T_k)_{k=0}^\infty \colon   T_k \in \mathbb{V}^{\otimes k}, \; T_k = 0 \text{ for all but finitely many }k\} .
\]
However, we prefer the form in \eqref{eq:tenalg1}. 
The reason $\Ten(\mathbb{V})$ is called an algebra is that it inherits a product operation given by the product of tensors in \eqref{eq:outprod3}:
\[
\mathbb{V}^{\otimes d} \times \mathbb{V}^{\otimes d'} \ni (T,T') \mapsto T \otimes T' \in \mathbb{V}^{\otimes (d + d')},
\]
which can be extended to $\Ten(\mathbb{V})$ by
\[
\biggl(\sum_{j=0}^d T_j \biggr) \otimes \biggl(\sum_{k=0}^{d'} T'_k \biggr) \coloneqq \sum_{j=0}^d \sum_{k=0}^{d'} T_j \otimes T'_k.
\]
Nothing we have discussed up to this point involves a basis of $\mathbb{V}$, so it all applies to modules as well. If we now take advantage of the fact that $\mathbb{V}$ is a vector space and further assume that it is finite-dimensional, then $\mathbb{V}$ has a basis $\mathscr{B} = \{e_1,\ldots,e_n\}$ and we may represent any element in $\Ten(\mathbb{V})$ as
\[
\sum_{k = 0}^d\; \sum_{i_1,i_2, \dots ,i_k=1}^n a_{i_1i_2 \cdots i_k} e_{i_1}\otimes e_{i_2} \otimes \dots \otimes e_{i_k}
\]
with $d\in \mathbb{N}$. Another way to represent this is as a \emph{non-commutative polynomial}
\[
f(X_1,X_2,\ldots,X_n) =
\sum_{k = 0}^d\; \sum_{i_1,i_2, \dots ,i_k=1}^n a_{i_1i_2 \cdots i_k} X_{i_1} X_{i_2} \cdots X_{i_k},
\]
in $n$ non-commutative variables $X_1, X_2,\ldots,X_n$. Here non-commutative means $X_i X_j \ne X_j X_i$ if $i \ne j$, \ie\ just like \eqref{eq:noncomm}, and we require other arithmetic rules for manipulating $X_1,X_2,\ldots,X_n$ to mirror those in \eqref{eq:distrib4}. The set of all such polynomials is usually denoted $\mathbb{R}\langle X_1,\ldots,X_n\rangle$. Hence, as algebras,
\begin{equation}\label{eq:tenalg2}
\Ten(\mathbb{V}) \cong \mathbb{R}\langle X_1,\ldots,X_n\rangle ,
\end{equation}
where the isomorphism is given by the map that takes $e_{i_1}\otimes e_{i_2} \otimes \dots \otimes e_{i_k} \mapsto X_{i_1} X_{i_2} \cdots X_{i_k}$. Such a representation also applies to any $d$-tensor in $\mathbb{V}^{\otimes d}$ as a special case, giving us a \emph{homogeneous} non-commutative polynomial
\[
\sum_{i_1,i_2, \dots ,i_d=1}^n a_{i_1i_2 \cdots i_d} X_{i_1} X_{i_2} \cdots X_{i_d} ,
\]
where all terms have degree $d$. This is a superior representation compared to merely representing the tensor as a hypermatrix $(a_{i_1i_2 \cdots i_d} )_{i_1,i_2, \dots ,i_d=1}^n \in \mathbb{R}^{n \times n \times \dots \times n}$ because, like usual polynomials, we can take derivatives and integrals of non-commutative polynomials and evaluate them on matrices, \tie, we can take $A_1,A_2,\ldots,A_n \in \mathbb{R}^{m \times m}$ and plug them into $f(X_1,X_2,\ldots,X_n)$ to get $f(A_1,A_2,\ldots,A_n) \in \mathbb{R}^{m \times m}$. This last observation alone leads to a rich subject, often called `non-commutative sums-of-squares', that has many engineering applications \cite{Helton}.

What if we need to speak of an infinite sum of tensors of different orders? This is not just of theoretical interest; as we will see in Example~\ref{eg:multpole}, a multipole expansion is such an infinite sum. A straightforward way to do this would be to replace the direct sums in \eqref{eq:tenalg0} with direct products; the difference, if the reader recalls, is that while the elements of a direct sum are zero for all but a finite number of summands, those of a direct product may contain an infinite number of non-zero summands. If there are just a finite number of vector spaces $\mathbb{V}_1,\dots,\mathbb{V}_d$, the direct sum $\mathbb{V}_1 \oplus \dots \oplus \mathbb{V}_d$ and the direct product $\mathbb{V}_1 \times \dots \times \mathbb{V}_d$ are identical, but if we have infinitely many vector spaces, a direct product is much larger than a direct sum. For instance, a direct sum of countably many \emph{two}-dimensional vector spaces $\bigoplus_{k \in \mathbb{N}}\mathbb{V}_k$ has countable dimension whereas a direct product $\prod_{k \in \mathbb{N}}\mathbb{V}_k$ has uncountable dimension.

Nevertheless, instead of introducing new notation for the direct product of $\mathbb{V}^{\otimes k}$, $k =0,1,2,\ldots,$ we may simply consider the dual vector space of $\Ten(\mathbb{V})$,
\begin{equation}\label{eq:tenalgdual0}
\Ten(\mathbb{V})^*  = \prod_{k=0}^\infty \mathbb{V}^{\ast \otimes k}
= \biggl\{\sum_{k=0}^\infty \varphi_k \colon  \varphi_k \in \mathbb{V}^{\ast \otimes k} \biggr\},
\end{equation}
noting that taking duals turns a direct sum into a direct product. We call this the \emph{dual tensor algebra}.  While it is true that $(\mathbb{V}^{\otimes k } )^*= \mathbb{V}^{\ast \otimes k}$ -- this follows from the `calculus of tensor products' that we will present in Example~\ref{eg:calc} -- we emphasize that the dual tensor algebra of any non-zero vector space $\mathbb{V} $ is not the same as the tensor algebra of its dual space $\mathbb{V}^*$, \tie,
\[
\Ten(\mathbb{V})^* \ne \Ten(\mathbb{V}^*).
\]
Indeed, the dual vector space of any infinite-dimensional vector space must have dimension strictly larger than that of the vector space \cite[Chapter~IX, Section 5]{Jacobson2}. In this case, if $\mathbb{V}$ is finite-dimensional, then by what we discussed earlier, $\Ten(\mathbb{V}^*)$ has countable dimension whereas $\Ten(\mathbb{V})^*$ has uncountable dimension. Although what we have defined in \eqref{eq:tenalgdual0} is the direct product of $\mathbb{V}^{\ast \otimes k}$, $k =0,1,2,\ldots,$ the direct product of $\mathbb{V}^{\otimes k}$, $k =0,1,2,\ldots,$ is simply given by the dual tensor algebra of the dual space:
\[
\Ten(\mathbb{V}^*)^*  = \prod_{k=0}^\infty \mathbb{V}^{\otimes k}
= \biggl\{\sum_{k=0}^\infty T_k \colon  T_k \in \mathbb{V}^{\otimes k} \biggr\}.
\]
Just as $\Ten(\mathbb{V})$ may be regarded as the ring of non-commutative polynomials, the same arguments that led to \eqref{eq:tenalg2} also give
\[
\Ten(\mathbb{V}^*)^* \cong \mathbb{R}\llangle X_1,\ldots,X_n\rrangle,
\]
\tie, $\Ten(\mathbb{V}^*)^*$ may be regarded as the ring of non-commutative power series.

One issue with these constructions is that taking direct sums results in a space $\Ten(\mathbb{V})$ that is sometimes too small whereas taking direct products leads to a space $\Ten(\mathbb{V}^*)^*$ that is often too big. However, these constructions are purely algebraic; if $\mathbb{V}$ is equipped with a norm or inner product, we can often use that to obtain a tensor algebra of the desired size. Suppose $\mathbb{V}$ has an inner product $\langle \, \cdot\,, \cdot\, \rangle$; then as in Example~\ref{eg:norms}, defining
\[
\langle u \otimes v \otimes \dots \otimes w, u' \otimes v' \otimes \dots \otimes w' \rangle\coloneqq \langle u, u' \rangle \langle v, v' \rangle \cdots \langle w, w' \rangle
\]
for any $ u \otimes v \otimes \dots \otimes w, u' \otimes v' \otimes \dots \otimes w' \in \mathbb{V}^{\otimes d}$, any $d \in \mathbb{N}$, and extending bilinearly to all of $\Ten(\mathbb{V})$ makes it into an inner product space. We will denote the inner products on $\mathbb{V}^{\otimes d}$ and $\Ten(\mathbb{V})$ by $\langle \, \cdot\,, \cdot\, \rangle$ and the corresponding induced norms by $\lVert\,\cdot\,\rVert$; there is no ambiguity since they are all consistently defined.  Taking completion gives us a Hilbert space
\[
\widehat{\Ten}(\mathbb{V}) \coloneqq  \biggl\{\sum_{k=0}^\infty T_k \colon  T_k \in \mathbb{V}^{\otimes k}, \;  \sum_{k=0}^\infty \lVert T_k \rVert^2 < \infty \biggr\}
\]
with inner product
\[
\biggl\langle \sum_{k=0}^\infty T_k , \sum_{k=0}^\infty T_k' \biggr\rangle \coloneqq \sum_{k=0}^\infty \langle T_k, T_k' \rangle.
\]
Alternatively, one may also describe $\widehat{\Ten}(\mathbb{V})$ as the \emph{Hilbert space direct sum} of the Hilbert spaces $\mathbb{V}^{\otimes k}$, $k =0,1,2,\ldots.$ 
If $\mathscr{B}$ is an orthonormal basis of $\mathbb{V}$, then the tensor product basis \eqref{eq:tenprodbasis2}, denoted $\mathscr{B}^{\otimes k}$, is an orthonormal basis on $\mathbb{V}^{\otimes k}$, and  $\bigcup_{k=0}^\infty \mathscr{B}^{\otimes k}$ is a countable orthonormal basis on $\widehat{\Ten}(\mathbb{V})$, \ie\ it is a separable Hilbert space. This can be used in a `tensor trick' for linearizing a non-linear function with convergent Taylor series. For example, for any $x,y \in \mathbb{R}^n$,
\begin{align*}
\exp(\langle x, y\rangle) &= \sum_{k=0}^\infty \dfrac{1}{k!} \langle x,y\rangle^k = \sum_{k=0}^\infty \dfrac{1}{k!} \langle x^{\otimes k},y^{\otimes k}\rangle\\*
&= \biggl\langle \sum_{k=0}^\infty \dfrac{1}{\sqrt{k!}}  x^{\otimes k}, \sum_{k=0}^\infty \dfrac{1}{\sqrt{k!}} y^{\otimes k}\biggr\rangle = \langle S(x), S(y) \rangle,
\end{align*}
where the map
\[
S \colon  \mathbb{R}^n \to \widehat{\Ten}(\mathbb{R}^n), \quad S(x) =  \sum_{k=0}^\infty \dfrac{1}{\sqrt{k!}}  x^{\otimes k}
\]
is well-defined since
\[
\lVert S(x) \rVert^2 = \sum_{k=0}^\infty \biggl\lVert \dfrac{1}{\sqrt{k!}}  x^{\otimes k} \biggr\rVert^2
= \sum_{k=0}^\infty \dfrac{1}{k!} \lVert x \rVert^{2k} = \exp(\lVert x \rVert^2) < \infty.
\]
This was  used, with $\sin$ in place of $\exp$,  in an elegant proof \cite{Krivine} of Grothendieck's inequality that also yields
\[
K_\G \le \dfrac{\pi}{2 \log (1+ \sqrt{2})} \approx 1.78221,
\]
the upper bound we saw in Example~\ref{eg:GI}. This is in fact the best explicit upper bound for the Grothendieck constant over $\mathbb{R}$, although it is now known that it is not sharp.
\end{example}

\subsection{Tensors via the universal factorization property}\label{sec:tensor3c}

We now come to the most abstract and general way to define tensor products and thus tensors: the universal factorization property. It is both a potent tool and an illuminating concept. Among other things, it unifies almost every notion of tensors we have discussed thus far. For instance, an immediate consequence is that multilinear maps (Definitions~\ref{def:tensor2} and~\ref{def:tensor2a}) and sums of separable functions (Definitions~\ref{def:tensor3a} and~\ref{def:tensor3b}) are one and the same notion, unifying definitions~\ref{st:tensor2} and~\ref{st:tensor3}.

The straightforward way to regard the universal factorization property is that it is a property of multilinear maps. It says that the multilinearity in any multilinear map $\Phi$ is the result of a special multilinear map $\sigma_\otimes$ that is an omnipresent component of all multilinear {maps: $\Phi = F \circ \sigma_\otimes$  --  once $\sigma_\otimes$ is factored out of $\Phi$, what remains, \ie\ the component $F$, is just a linear map.} Any multilinearity
  contained in $\Phi$ is all due to $\sigma_\otimes$.

However, because of its universal nature, we may turn the property on its head and use it to define the tensor product:  a tensor product is whatever satisfies the universal factorization property. More precisely, whatever vector space one can put in the box below that makes the diagram
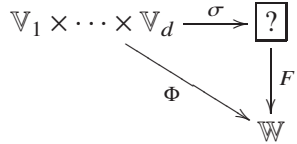
\[
\xymatrix{
\mathbb{V}_1\times\dots\times \mathbb{V}_d\ar[r]^-{\sigma}\ar[rd]_{\Phi} & \fbox{?}\ar[d]^{F}\\
 & \mathbb{W}}
\]
commute for all vector spaces $\mathbb{W}$ with $\Phi$ multilinear and $F$ linear is defined to be a tensor product of $\mathbb{V}_1,\ldots, \mathbb{V}_d$. One may show that this defines the tensor product (and the map $\sigma$) uniquely up to vector space isomorphism.

We now supply the details. The special multilinear map in question is $\seg \colon  \mathbb{V}_1\times \cdots \times \mathbb{V}_d\to \mathbb{V}_1\otimes \dots\otimes \mathbb{V}_d$,
\begin{equation}\label{eq:seg}
\seg(v_1,\ldots,v_d) = v_1\otimes\dots\otimes v_d,
\end{equation}
called the \emph{Segre map} or Segre embedding.\footnote{There is no standard notation. It has been variously denoted as $\sigma$ \cite{Harris}, as $\varphi$ \cite{Bourbaki,Lang}, as $\otimes$ \cite{Conrad} and as $\operatorname{Seg}$ \cite{Land1}.} Here the tensor space $\mathbb{V}_1\otimes \dots\otimes \mathbb{V}_d$ is as defined in Section~\ref{sec:tensor3b} and \eqref{eq:seg} defines a multilinear map by virtue of \eqref{eq:distrib4}.

\mypara{Universal factorization property} Every multilinear map is the composition of a linear operator with the Segre map.
\medskip

This is also called the universal mapping property, or the universal property of the tensor product. Formally, it says that if $\Phi$ is a multilinear
map from
$\mathbb{V}_1\times \cdots \times \mathbb{V}_d$ into a vector space $\mathbb{W}$, then $\Phi$ determines a unique linear
map $F_\Phi$ from $\mathbb{V}_1\otimes\dots\otimes \mathbb{V}_d$ into $\mathbb{W}$,
that makes the following diagram commute:
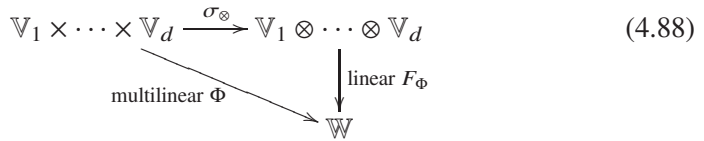
\begin{equation}\label{eq:commdiag}
\xymatrix{
\mathbb{V}_1\times\dots\times \mathbb{V}_d\ar[r]^{\seg}\ar[rd]_{\text{multilinear }\Phi \quad} & \mathbb{V}_1\otimes\dots\otimes \mathbb{V}_d\ar[d]^{\text{linear }F_\Phi}\\
 & \mathbb{W}}
\end{equation}
or, equivalently, $\Phi$ can be factored as
\begin{equation}\label{eq:factor0}
\Phi = F_\Phi \circ \seg,
\end{equation}
where $F_\Phi$ is a `linear component' dependent on $\Phi$ and $\seg$ is a `multilinear component'  independent of $\Phi$. The component $\seg$ is universal, \ie\ common to all multilinear maps.
This can be stated in less precise but more intuitive terms as follows. 
\begin{enumerate}[\upshape (a)]
\setlength\itemsep{3pt}
\item\label{it:factor} All the `multilinearness' in a multilinear map $\Phi$ may be factored out of $\Phi$, leaving behind only the `linearness' that is encapsulated in $F_\Phi$.
\item\label{it:same} The `multilinearness' extracted from any multilinear map is identical, \tie, all multilinear maps are multilinear because they contain a copy of the Segre map $\seg$.
\end{enumerate}
Actually, $\seg$ does depend on $d$ and so we should have said it is universal for $d$-linear maps. An immediate consequence of the universal factorization property is that
\[
\Mult^d(\mathbb{V}_1,\ldots,\mathbb{V}_d; \mathbb{W}) = \Lin(\mathbb{V}_1 \otimes \dots \otimes \mathbb{V}_d; \mathbb{W}).
\]
In principle one could use the universal factorization property to avoid multilinear maps entirely and discuss only linear maps, but of course there is no reason to do that: as we saw in Section~\ref{sec:mult}, multilinearity is a very useful notion in its own right. The universal factorization property should be viewed as a way to move back and forth between the multilinear realm and the linear one.  While we have assumed contravariant tensors for simplicity, the discussions above apply verbatim to mixed tensor spaces $\mathbb{V}_1\otimes \dots\otimes \mathbb{V}_p\otimes \mathbb{V}_{p+1}^*\otimes \dots\otimes \mathbb{V}_d^*$ and with modules in place of vector spaces. 

Pure mathematicians accustomed to commutative diagrams would swear by \eqref{eq:commdiag} but it expresses the same thing as \eqref{eq:factor0}, which is more palatable to applied and computational mathematicians as it reminds us of matrix factorizations like Example~\ref{eg:qr}, where we factor a matrix $A \in \mathbb{R}^{m \times n}$ into  $A=QR$ with an orthogonal component $Q \in \Or(m)$ that captures the `orthogonalness' in $A$, leaving behind a triangular component $R$. There is one big difference, though:  both $Q$ and $R$ depend on $A$, but in \eqref{eq:factor0}, whatever our choice of $\Phi$, the multilinear component will always be $\seg$. The next three examples will bear witness to this intriguing property.

\begin{example}[trilinear functionals]\label{eg:trilin}
The universal factorization property explains why  \eqref{eq:tl} and \eqref{eq:distrib3} look alike. A trilinear functional $\tau \colon \mathbb{U} \times \mathbb{V} \times \mathbb{W} \to \mathbb{R}$ can be factored as
\begin{equation}\label{eq:factor1}
\tau = F_\tau \circ \seg
\end{equation}
and the multilinearity of $\seg$ in \eqref{eq:distrib3} accounts for that of $\tau$ in \eqref{eq:tl}. For example, take the first equation in \eqref{eq:tl}. We may view it as arising from
\begin{align*}
\tau(\lambda u +\lambda' u',v, w) &= F_\tau ( \seg (\lambda u +\lambda' u',v, w)) &&\text{by \eqref{eq:factor1},}\\*
&=  F_\tau ((\lambda u +\lambda' u')  \otimes v  \otimes  w ) && \text{by \eqref{eq:seg},}\\
&=  F_\tau (\lambda u\otimes v  \otimes  w +  \lambda' u'  \otimes v  \otimes  w) && \text{by \eqref{eq:distrib3},}\\
&=  \lambda F_\tau(u\otimes v  \otimes  w )+  \lambda' F_\tau(u'  \otimes v  \otimes  w) && F_\tau\text{ is linear,}\\
&=  \lambda  F_\tau ( \seg ( u, v,  w ) )+  \lambda' F_\tau ( \seg ( u', v,  w ) ) && \text{by \eqref{eq:seg},}\\*
&= \lambda\tau(u, v,w) +\lambda' \tau(u', v, w) &&\text{by \eqref{eq:factor1},}
\end{align*}
and similarly for the second and third equations in \eqref{eq:tl}. 
More generally, if we reread Section~\ref{sec:structmultmaps} with the hindsight of this section, we will realize that it is simply a discussion of the universal factorization property without invoking the $\otimes$ symbol. 
\end{example}

We next see how the universal factorization property ties together the three matrix products that have made an appearance in this article.

\begin{example}[matrix, Hadamard and Kronecker products]\label{eg:variousprod}
Let us consider the standard matrix product in \eqref{eq:mm} and  Hadamard product in \eqref{eq:hp} on $2\times 2$ matrices. Let $\mu, \eta \colon  \mathbb{R}^{2 \times 2} \times \mathbb{R}^{2 \times 2}  \to \mathbb{R}^{2 \times 2} $ be, respectively,
\begin{align*}
\mu\biggl(
\begin{bmatrix}
a_{11} & a_{12} \\
a_{21} & a_{22}
\end{bmatrix}\!,
\begin{bmatrix}
b_{11} & b_{12} \\
b_{21} & b_{22}
\end{bmatrix}
\biggr)
&=
\begin{bmatrix}
a_{11} b_{11} + a_{12} b_{21} & a_{11} b_{12} + a_{12} b_{22}\\
a_{21} b_{11} + a_{22} b_{21} & a_{21} b_{12} + a_{22} b_{22}
\end{bmatrix}\!,\\
\eta\biggl(
\begin{bmatrix}
a_{11} & a_{12} \\
a_{21} & a_{22}
\end{bmatrix}\!,
\begin{bmatrix}
b_{11} & b_{12} \\
b_{21} & b_{22}
\end{bmatrix}
\biggr)
&=
\begin{bmatrix}
a_{11} b_{11}  &  a_{12} b_{12}\\
a_{21} b_{21} &  a_{22} b_{22}
\end{bmatrix}\!.
\end{align*}
They look nothing alike and, as discussed in Section~\ref{sec:trans}, are of entirely different natures. But the universal factorization property tells us that, by virtue of the fact that both are bilinear operators from  $\mathbb{R}^{2 \times 2} \times \mathbb{R}^{2 \times 2}$  to  $\mathbb{R}^{2 \times 2} $, the `bilinearness' in $\mu$ and  $\eta$ is the same, encapsulated in the Segre map $\seg \colon  \mathbb{R}^{2 \times 2} \times \mathbb{R}^{2 \times 2} \to \mathbb{R}^{2 \times 2} \otimes \mathbb{R}^{2 \times 2} \cong \mathbb{R}^{4 \times 4}$,
\begin{align*}
\seg\biggl(
\begin{bmatrix}
a_{11} & a_{12} \\
a_{21} & a_{22}
\end{bmatrix}\!,
\begin{bmatrix}
b_{11} & b_{12} \\
b_{21} & b_{22}
\end{bmatrix}
\biggr)
&=
\begin{bmatrix}
a_{11} & a_{12} \\
a_{21} & a_{22}
\end{bmatrix} \otimes
\begin{bmatrix}
b_{11} & b_{12} \\
b_{21} & b_{22}
\end{bmatrix}\\
&=
\begin{bmatrix}
a_{11} b_{11} & a_{11} b_{12} & a_{12}b_{11} & a_{12}b_{12} \\
a_{11} b_{21} & a_{11} b_{22} & a_{12}b_{21} & a_{12}b_{22}\\
a_{21}b_{11} & a_{21}b_{12} & a_{22}b_{11} & a_{22}b_{12} \\
a_{21}b_{21} & a_{21}b_{22} &a_{22}b_{21} & a_{22}b_{22}
\end{bmatrix}\!,
\end{align*}
\ie\ the Kronecker product  in Example~\ref{eg:concretetenprod}\ref{it:kronprod}.  The difference between $\mu$ and $\eta$ is due entirely to  their linear components $F_\mu, F_\eta \colon  \mathbb{R}^{4 \times 4} \to \mathbb{R}^{2 \times 2}$,
\begin{align*}
F_\mu
\left( 
\begin{bmatrix}
c_{11} & c_{12} & c_{13} & c_{14} \\
c_{21} & c_{22} & c_{23} & c_{24} \\
c_{31} & c_{32} & c_{33} & c_{34} \\
c_{41} & c_{42} & c_{43} & c_{44}
\end{bmatrix}
\right)
&=
\begin{bmatrix}
c_{11} + c_{23} & c_{12} + c_{24} \\
c_{31} + c_{43} & c_{32} + c_{44}
\end{bmatrix}\!,\\
F_\eta
\left( 
\begin{bmatrix}
c_{11} & c_{12} & c_{13} & c_{14} \\
c_{21} & c_{22} & c_{23} & c_{24} \\
c_{31} & c_{32} & c_{33} & c_{34} \\
c_{41} & c_{42} & c_{43} & c_{44}
\end{bmatrix}
\right)
&=
\begin{bmatrix}
c_{11} & c_{14} \\
c_{41} & c_{44}
\end{bmatrix}\!,
\end{align*}
which are indeed very different. As a sanity check, the reader may like to verify from these formulas that
\[
\mu = F_\mu \circ \seg, \quad \eta = F_\eta \circ \seg\!.
\]
More generally, the same argument shows that any products we define on $2\times 2$ matrices must share something in common, namely the Kronecker product that takes a pair of $2\times 2$ matrices to a $4 \times 4$ matrix. The Kronecker product accounts for the `bilinearness' of all conceivable matrix products.
\end{example}

We will include an infinite-dimensional example for variety.

\begin{example}[convolution I]\label{eg:L1}
We recall the convolution bilinear operator 
\[
\ast \colon  L^1(\mathbb{R}^n)\times L^1(\mathbb{R}^n)  \to L^1(\mathbb{R}^n), \quad (f,g) \mapsto f \ast g
\]
from Example~\ref{eg:Calderon}.
As we saw in Example~\ref{eg:Hilbert}, $L^1(\mathbb{R}^n) \otimes L^1(\mathbb{R}^n) = L^1(\mathbb{R}^n \times \mathbb{R}^n)$ when we interpret $\otimes$ as $\hatotimes_\nu$, the topological tensor product with respect to the nuclear norm, and so the Segre map is
\[
\seg \colon  L^1(\mathbb{R}^n)\times L^1(\mathbb{R}^n)  \to L^1(\mathbb{R}^n \times \mathbb{R}^n), \quad (f,g) \mapsto f \otimes g,
\]
where, as a reminder, $f \otimes g(x,y) = f(x)g(y)$, where $\otimes$ is interpreted as the separable product. Commutative diagrams like \eqref{eq:commdiag} are useful when we have three or more maps between multiple spaces, as they allow us to assemble all maps into a single picture, which helps to keep track of what gets mapped where. Here we have
\[
\xymatrix{
L^1(\mathbb{R}^n)\times L^1(\mathbb{R}^n) \ar[r]^-{\seg}\ar[rd]_{\ast} & L^1(\mathbb{R}^n \times \mathbb{R}^n) \ar[d]^{F_\ast}\\
 & L^1(\mathbb{R}^n)}
\]
and we would like to find the unique linear map $F_\ast$ that makes the diagram commute, \ie\ $f \ast g = F_*( f \otimes g)$. Writing out this expression, we see that $F_\ast = G \circ H$ is a composition of the two linear maps defined by
\begin{alignat*}{3}
G \colon  L^1(\mathbb{R}^n \times \mathbb{R}^n) &\to L^1(\mathbb{R}^n), & G(K)(x) &\coloneqq \int_{\mathbb{R}^n} K(x,y)  \D y,\\*
H \colon  L^1(\mathbb{R}^n \times \mathbb{R}^n) &\to L^1(\mathbb{R}^n \times \mathbb{R}^n),\quad & H(K)(x,y) &\coloneqq K(x-y,y)
\end{alignat*}
for all $K \in L^1(\mathbb{R}^n \times \mathbb{R}^n) $.
\end{example}

Strictly speaking, the space $\mathbb{V}_1\otimes \dots\otimes \mathbb{V}_d$ has already been fixed by Definition~\ref{def:tensor3b}. In order to allow some flexibility in \eqref{eq:commdiag}, we will have to allow for $\mathbb{V}_1\otimes \dots\otimes \mathbb{V}_d$ to be replaced by any vector space $\mathbb{T}$ and $\seg$ by any multilinear map
\[
\segt \colon  \mathbb{V}_1\times \cdots \times \mathbb{V}_d\to \mathbb{T}
\]
as long as the universal factorization property holds for this choice of $\mathbb{T}$ and $\segt$, \tie,
\begin{equation}\label{eq:commdiag1}
\xymatrix{
\mathbb{V}_1\times\dots\times \mathbb{V}_d\ar[r]^-{\segt}\ar[rd]_{\text{multilinear }\Phi \quad} & \mathbb{T}\ar[d]^{\text{linear }F_\Phi}\\
 & \mathbb{W}}
\end{equation}

\begin{definition}[tensors via universal factorization]\label{def:tensor3c}\hspz%
Let $\mathbb{T}$ be a vector space and let $\segt \colon  \mathbb{V}_1\times \cdots \times \mathbb{V}_d\to \mathbb{T}$ be a multilinear map with the universal factorization property, \tie, for any vector space $\mathbb{W}$ and any multilinear map $\Phi \colon  \mathbb{V}_1\times\dots\times \mathbb{V}_d \to \mathbb{W}$, there is a unique linear map $F_\Phi \colon  \mathbb{T} \to \mathbb{W}$, called the \emph{linearization} of $\Phi$, such that the factorization $\Phi = F_\Phi \circ \segt$ holds. Then $\mathbb{T}$ is unique up to vector space isomorphisms and is called a \emph{tensor space} or more precisely a \emph{tensor product} of the vector spaces $\mathbb{V}_1,\ldots,\mathbb{V}_d$. An element of $\mathbb{T}$ is called a $d$-\emph{tensor}.
\end{definition}

As usual, one may easily incorporate covariance and contravariance by replacing some of the vector spaces with their dual spaces.

In Example~\ref{eg:trilin}, $\mathbb{T} = \mathbb{U} \otimes \mathbb{V} \otimes \mathbb{W}$, as defined in   Definition~\ref{def:tensor3b}, $\segt = \seg$, as defined in \eqref{eq:seg}, and we recover \eqref{eq:commdiag}. But for the subsequent two examples we have in fact already used Definition~\ref{def:tensor3c} implicitly. In Example~\ref{eg:variousprod}, $\mathbb{T} = \mathbb{R}^{4 \times 4}$ and  $\segt \colon  \mathbb{R}^{2\times 2} \times \mathbb{R}^{2 \times 2} \to \mathbb{T}$  takes a pair of $2 \times 2$ matrices to their Kronecker product.  In Example~\ref{eg:L1}, $\mathbb{T} = L^1(\mathbb{R}^n \times \mathbb{R}^n)$ and $\segt \colon  L^1(\mathbb{R}^n) \times  L^1(\mathbb{R}^n) \to \mathbb{T}$ takes a pair of $L^1$-functions to their separable product.

Definition~\ref{def:tensor3c} is intended to capture as many notions of tensor products as possible. The Segre map $\sigma_\mathbb{T}$ may be an abstract tensor product of vectors, an outer product of vectors in $\mathbb{R}^n$, a Kronecker product of operators or matrices, a separable product of functions, {\em etc.}; the $\otimes$ in \eqref{eq:commdiag1} could be abstract tensor products of vector spaces, Kronecker products of operator spaces, topological tensor products of topological vector spaces, {\em etc.} In the last case we require the maps $\Phi$, $\sigma_\mathbb{T}$, $F_\Phi$ to be also \emph{continuous};  without this requirement the topological tensor product just becomes the algebraic tensor product. Recall that continuity takes a simple form for a multilinear map, being equivalent to the statement that its spectral norm as defined in \eqref{eq:specnorm1} is finite, \tie, continuous multilinear maps are precisely bounded multilinear maps.
Let us verify that the maps in Example~\ref{eg:L1} are continuous.

\begin{example}[convolution II]\label{eg:L1a}
We follow the same notation in Example~\ref{eg:L1} except that we will write $\Beta_* \colon   L^1(\mathbb{R}^n) \times L^1(\mathbb{R}^n)  \to L^1(\mathbb{R}^n)$, $(f,g) \mapsto f \ast g$ for convolution as a bilinear map. We show that $\lVert \Beta_* \rVert_\sigma = \lVert \sigma_\otimes \rVert_\sigma = \lVert F_* \rVert_\sigma = 1$ and thus these maps are continuous. First recall from Example~\ref{eg:Calderon} that $\lVert \Beta_* \rVert_\sigma =1$. Next, since  $\lVert f \otimes g \rVert_1 = \lVert f  \rVert_1 \lVert g \rVert_1$, we have
\[
\lVert \sigma_\otimes \rVert_\sigma = \sup_{f,g \ne 0} \dfrac{\lVert f\otimes g\rVert_1}{\lVert f \rVert_1 \lVert g\rVert_1} = 1.
\]
Lastly, for $K \in  L^1(\mathbb{R}^n \times \mathbb{R}^n)$,
\[
\lVert F_* (K)\rVert_1 = \int_{\mathbb{R}^n} \int_{\mathbb{R}^n} \lvert K(x-y,y) \rvert  \D y  \D x = \lVert K \rVert_1
\]
by Fubini and a change of variables. Thus we get
\[
\lVert F_* \rVert_\sigma = \sup_{K \ne 0} \dfrac{\lVert F_* (K) \rVert_1}{\lVert K\rVert_1} = 1.
\]
Note that even though we use the \emph{spectral} norm to ascertain continuity of the linear and multilinear maps between spaces, the tensor product space above and in Example~\ref{eg:L1} is formed with $\hatotimes_\nu$, the topological tensor product with respect to the \emph{nuclear} norm.
Indeed, $\hatotimes_\nu$
is an amazing tensor product in that for any norm spaces $\mathbb{U}, \mathbb{V}, \mathbb{W}$, if we have the commutative diagram
\[
\xymatrix{
\mathbb{U} \times \mathbb{V} \ar[r]^-{\seg}\ar[rd]_{\Phi} & \mathbb{U} \hatotimes_\nu \mathbb{V}\ar[d]^{F_\Phi}\\
 & \mathbb{W}}
\]
then the bilinear map $\Phi \colon  \mathbb{U} \times \mathbb{V} \to \mathbb{W}$ is continuous if and only if its linearization $F_\Phi \colon  \mathbb{U} \hatotimes_\nu \mathbb{V} \to \mathbb{W}$ is continuous \cite[Chapter~3]{Defant}.
\end{example}

We will now examine a non-example to see that Definition~\ref{def:tensor3c} is not omnipotent. The tensor product of Hilbert spaces is an important exception that does not satisfy the universal factorization property. The following example is adapted from \citet{Garrett}.

\begin{example}[no universal factorization for Hilbert spaces]\label{eg:cftr}
We begin with a simple observation. For a finite-dimensional vector space $\mathbb{V}$, the map defined by $\beta \colon \mathbb{V} \times \mathbb{V}^* \to \mathbb{R}$, $(v,\varphi) \mapsto \varphi(v)$ is clearly bilinear. So we may apply universal factorization to get
\[
\xymatrix{
\mathbb{V}\times \mathbb{V}^* \ar[r]^-{\seg}\ar[rd]_{\beta} & \mathbb{V} \otimes \mathbb{V}^* \ar[d]^{\tr}\\
 & \mathbb{R}}
\]
The linear map that we obtained, $\tr \colon \mathbb{V}\otimes \mathbb{V}^* \to \mathbb{R}$, is in fact a coordinate-free way to define trace. To see this,  note that if we choose $\mathbb{V} = \mathbb{R}^n$, then for any $a \in \mathbb{R}^n$ and $b^\tp \in \mathbb{R}^{1 \times n} = \mathbb{R}^{n\ast}$, 
\[
\beta(a,b^\tp ) = b^\tp a = \tr(b^\tp a) = \tr(ab^\tp) = \tr(a \otimes b).
\]
This observation is key to our subsequent discussion. We will see that when $\mathbb{V}$ is replaced by an infinite-dimensional Hilbert space and we require all maps to be continuous, then there is no continuous linear map that can take the place of trace.

Let $\mathbb{H}$ be a separable infinite-dimensional Hilbert space. As discussed in Examples~\ref{eg:Hilbert} and \ref{eg:compact},
the tensor product of $\mathbb{H}$ and $\mathbb{H}^*$
is interpreted to be the topological tensor product $\mathbb{H} \hatotimes_\F \mathbb{H}^*$  with respect to the Hilbert--Schmidt norm $\lVert \, \cdot \,\rVert_\F$. So we must have $\seg \colon  \mathbb{H}\times \mathbb{H}^* \to  \mathbb{H} \hatotimes_\F \mathbb{H}^*$, $(v, \varphi) \mapsto v \otimes \varphi$. Take $\mathbb{W} = \mathbb{C}$ and consider the continuous bilinear functional $\beta \colon   \mathbb{H}\times \mathbb{H}^* \to \mathbb{C}$, $(v,\varphi) \mapsto \varphi(v)$. We claim that there is no continuous linear map $F_\beta \colon   \mathbb{H} \hatotimes_\F \mathbb{H}^* \to \mathbb{C}$ such that $\beta = F_\beta \circ \seg$, \tie, we cannot make the diagram
\[
\xymatrix{
\mathbb{H}\times \mathbb{H}^* \ar[r]^-{\seg}\ar[rd]_{\beta} & \mathbb{H} \hatotimes_\F \mathbb{H}^* \ar[d]^{F_{\beta}?}\\
 & \mathbb{C}}
\]
commute.
To this end, suppose $F_\beta$ exists. Take any orthonormal basis $\{e_i \colon  i \in \mathbb{N} \}$. Then
\[
\beta( e_i, e_j^* ) = \delta_{ij}, \quad \seg (e_i, e_j^*) = e_i \otimes e_j^* ,
\]
and since $\beta = F_\beta \circ \seg$, we must have
\[
F_\beta (e_i \otimes e_j^*) = \delta_{ij}.
\]
Now take any $\Phi \in  \mathbb{H} \hatotimes_\F \mathbb{H}^*$ and express it as in \eqref{eq:infmatrep}. Then since $F_\beta$ is continuous and linear, we have
\[
F_\beta(\Phi) = \sum_{k=1}^\infty a_{kk} = \tr(\Phi),
\]
\tie, $F_\beta$ must be the trace as defined in \eqref{eq:trace}. This is a contradiction as trace is unbounded on $ \mathbb{H} \hatotimes_\F \mathbb{H}^*$;  just take
\[
\Phi = \sum_{k=1}^\infty k^{-1} e_k \otimes e_k^*,
\]
 which is Hilbert--Schmidt as $\lVert \Phi \rVert_\F^2 = \sum_{k=1}^\infty k^{-2} < \infty$ but  $\tr(\Phi) = \sum_{k=1}^\infty k^{-1} = \infty$. Hence no such map exists.

A plausible follow-up question is that instead of $\hatotimes_\F$, why not take topological tensor products with respect to the nuclear norm $\hatotimes_\nu$? As we mentioned at the end of Example~\ref{eg:L1a}, this automatically guarantees continuity of all maps. The problem, as we pointed out in Example~\ref{eg:compact}, is that the resulting tensor space, \ie\ the trace-class operators $\mathbb{H} \hatotimes_\nu \mathbb{H}^*$, is not a Hilbert space. In fact, one may show that no such Hilbert space exists. Suppose there is a Hilbert space $\mathbb{T}$ and a Segre map $\segt \colon  \mathbb{H}\times \mathbb{H}^* \to \mathbb{T}$ so that the universal factorization property in \eqref{eq:commdiag1} holds. Then we may choose $\mathbb{W} = \mathbb{H} \hatotimes_\F \mathbb{H}^*$ and $\Phi = \seg$ to get
\[
\xymatrix{
\mathbb{H}\times \mathbb{H}^* \ar[r]^-{\segt}\ar[rd]_{\seg} & \mathbb{T} \ar[d]^{F_\otimes}\\
 &  \mathbb{H} \hatotimes_\F \mathbb{H}^*}
\]
from which one may deduce  (see \citealt{Garrett} for details)  that if $F_\otimes$ is a continuous linear map, then
\[
\langle \, \cdot\,, \cdot \, \rangle_\mathbb{T} = c \langle \, \cdot\,, \cdot \, \rangle_\F
\]
for some $c > 0$, \tie, $F_\otimes$ is up to a constant multiple an isometry (Hilbert space isomorphism). So $\mathbb{T} $ and $ \mathbb{H} \hatotimes_\F \mathbb{H}^*$ are essentially the same Hilbert space up to scaling and thus the same trace argument above shows that $\mathbb{T} $  does not satisfy the universal factorization property.

While the goal of this example is to demonstrate that Hilbert spaces generally do not satisfy the universal factorization  property, there is a point worth highlighting in this construction. As we saw at the beginning, if we do not care about continuity, then everything goes through without a glitch: the \emph{evaluation bilinear functional} $\beta \colon  \mathbb{V}\times \mathbb{V}^*  \to \mathbb{R}$, $(v,\varphi) \mapsto \varphi(v)$ satisfies the universal factorization property
\[
\xymatrix{
\mathbb{V}\times \mathbb{V}^* \ar[r]^-{\seg}\ar[rd]_{\beta} & \mathbb{V} \otimes \mathbb{V}^* \ar[d]^{\tr}\\
 & \mathbb{R}}
\]
with the linear functional $\tr \colon \mathbb{V} \otimes \mathbb{V}^*  \to \mathbb{R}$ given by
\[
\tr\biggl( \sum_{i=1}^r  v_i \otimes \varphi_i \biggr) = \sum_{i=1}^r  \varphi_i( v_i )
\]
for any $r\in \mathbb{N}$.  In the language of Definition~\ref{def:tensor3c}, trace is the linearization of the evaluation bilinear functional, a definition adopted in both algebra \cite[Proposition~5.7]{Lang} and analysis \cite[Section~1.3]{Ryan}. Note that here $\mathbb{V}$ can be any vector space, not required to be finite-dimensional or complete or normed, and the maps involved are not required to be continuous. The trace as defined above is always a finite sum and thus finite-valued.
\end{example}

The next example clarifies an occasional ambiguity in the definition of higher-order derivatives of vector-valued functions.

\begin{example}[linearizing higher-order derivatives]\label{eg:linhod}\hspz%
Recall from Example~\ref{eg:hod} that the $d$th derivative of $f \colon  \Omega \to \mathbb{W}$ at a point $v \in \Omega \subseteq \mathbb{V} $ is a multilinear map denoted $D^d f(v)$ as in \eqref{eq:hod}. Applying the universal factorization property
\[
\xymatrix{
\overbracket[0.5pt]{\mathbb{V}  \times\dots\times \mathbb{V}}^{d \text{ copies}}
\ar[r]^-{\seg}\ar[rd]_{D^d f(v)} & \overbracket[0.5pt]{\mathbb{V}  \otimes\dots\otimes \mathbb{V}}^{d \text{ copies}} \ar[d]^{\partial^d f(v)}\\
 &  \mathbb{W}}
\]
we obtain a linear map on a tensor space,
\[
\partial^d f(v) \colon  \mathbb{V}^{\otimes d} \to \mathbb{W},
\]
that may be regarded as an alternative candidate for the $d$th derivative of $f$. In the language of Definition~\ref{def:tensor3c}, $\partial^d f(v)$ is the linearization of $D^d f(v)$. While we have denoted them  differently for easy distinction, it is not uncommon to see $D^d f(v)$ and $\partial^d f(v)$ used interchangeably, sometimes in the same sentence. 

Take the log barrier $f \colon  \mathbb{S}^n_\pp \to \mathbb{R}$, $f(X) = -\log \det (X)$ discussed in Examples~\ref{eg:logdet} and \ref{eg:self}. We have
\begin{alignat*}{3}
D^2f(X) \colon  \mathbb{S}^n \times \mathbb{S}^n &\to \mathbb{R}, &(H_1,H_2) &\mapsto \tr(X^{-1} H_1 X^{-1}H_2),\\*
\seg \colon  \mathbb{S}^n  \times \mathbb{S}^n &\to \mathbb{S}^n  \otimes \mathbb{S}^n, \quad &(H_1, H_2) &\mapsto H_1 \otimes H_2,
\end{alignat*}
where the latter denotes the Kronecker product as discussed in Example~\ref{eg:kron}. The universal factorization property
\[
\xymatrix{
\mathbb{S}^n  \times \mathbb{S}^n
\ar[r]^-{\seg}\ar[rd]_{D^d f(X)} & \mathbb{S}^n  \otimes \mathbb{S}^n \ar[d]^{\partial^d f(X)}\\
 &  \mathbb{R}}
\]
gives us
\[
\partial^2 f(X) \colon  \mathbb{S}^n \otimes \mathbb{S}^n \to \mathbb{R}, \quad (H_1,H_2) = \vect(X^{-1})^\tp (H_1 \otimes H_2) \vect(X^{-1}).
\]
Note that here it does not matter whether we use $\vect$ or $\vect_\ell$, since the matrices involved are all symmetric.

Taylor's theorem for a vector-valued function $f \colon  \Omega \subseteq \mathbb{V} \to \mathbb{W}$ in Example~\ref{eg:hod} may alternatively be expressed in the form
\begin{align*}
f(v)& =f(v_0)+ [\partial f(v_0)](v-v_0)+  \dfrac{1}{2}[\partial^2f(v_0)] (v-v_0)^{\otimes 2} + \cdots \\*
 &\quad\quad \cdots +\dfrac{1}{d!}[\partial^d f(v_0)](v-v_0)^{\otimes d}+ R(v-v_0),
\end{align*}
where  $v_0 \in \Omega$  and $\lVert R(v - v_0) \rVert/\lVert v - v_0 \rVert^d \to 0$  as $v \to v_0$. Here the `Taylor coefficients' are linear maps $\partial^d f(v) \colon  \mathbb{V}^{\otimes d} \to \mathbb{W}$. For the important special case $\mathbb{W} = \mathbb{R}$, \ie\ a real-valued function $f \colon  \mathbb{V} \to \mathbb{R}$, we have that $\partial^d f(v) \colon  \mathbb{V}^{\otimes d} \to \mathbb{R}$ is a linear functional. So the Taylor coefficients are covariant $d$-tensors
\[
\partial^d f(v) \in \mathbb{V}^{\ast \otimes d}.
\]
By our discussions in Example~\ref{eg:tenalg} and assuming that $f \in C^\infty(\Omega)$, we obtain an element of the dual tensor algebra
\[
\partial f (v)  \coloneqq \sum_{k=0}^\infty \dfrac{1}{k!}  \partial^k f(v) \in \Ten(\mathbb{V})^*
\]
that we will call the tensor Taylor series.

A slight variation is when $\mathbb{V}$ is equipped with an inner product $\langle \, \cdot\,, \cdot \,\rangle$. We may use the Riesz representation theorem to write, for $f \in C^{d+1}(\Omega)$,
\begin{align}\label{eq:taylor0}
f(v) & =f(v_0)+ \langle \partial f(v_0), v-v_0 \rangle +  \dfrac{1}{2}\langle \partial^2f(v_0), (v-v_0)^{\otimes 2}\rangle + \cdots  \notag \\*
 &\quad  \cdots +\dfrac{1}{d!} \langle \partial^d f(v_0), (v-v_0)^{\otimes d} \rangle + R(v-v_0),
\end{align}
where the inner products on $\mathbb{V}^{\otimes d}$ are defined as in \eqref{eq:innprod1}. Here the Taylor coefficients are contravariant $d$-tensors $\partial^d f(v) \in \mathbb{V}^{\otimes d}$. If $f$ is furthermore analytic, we may rewrite this as an inner product between maps into the tensor algebra:
\begin{equation}\label{eq:taylor}
f(v) = \biggl\langle \sum_{k=0}^\infty \dfrac{1}{k!}  \partial^k f(v_0), \sum_{k=0}^\infty (v-v_0)^{\otimes k}\biggr\rangle = \langle \partial f (v_0), S(v-v_0) \rangle,
\end{equation}
where the maps $\partial f  \colon  \Omega \to \widehat{\Ten}(\mathbb{V})$ and $S \colon  \Omega \to \widehat{\Ten}(\mathbb{V})$ are given by
\[
\partial f (v)  = \sum_{k=0}^\infty \dfrac{1}{k!}  \partial^k f(v),\quad
S(v)=  \sum_{k=0}^\infty  v^{\otimes k}.
\]
The tensor geometric series $S(v)$ is clearly a well-defined element of $\widehat{\Ten}(\mathbb{V})$ whenever $\lVert v \rVert < 1$ and the tensor Taylor series $\partial f (v) $ is well-defined as long as $\lVert \partial^k f (v)  \rVert \le B$ for some uniform bound $B > 0$. We will see this in action when we discuss multipole expansion in Example~\ref{eg:multpole}, which is essentially \eqref{eq:taylor} applied to $f(v) = 1/\lVert v \rVert$.
\end{example}

Another utility of the universal factorization property is in defining various vector-valued objects, \ie\ taking values in a vector space $\mathbb{V}$. The  generality here is deliberate, intended to capture a myriad of possibilities for $\mathbb{V}$, which could refer to anything from $\mathbb{C}$ as a two-dimensional vector space over $\mathbb{R}$ to an infinite-dimensional function space, a space of matrices or hypermatrices, a Hilbert or Banach or Fr\'echet space, a space of linear or multilinear operators, a tensor space itself, {\em etc.}

\begin{example}[vector-valued objects]\label{eg:vecvalobj}
For any real vector space $\mathbb{V}$ and any set $X$, we denote the set of all functions taking values in $\mathbb{V}$ by
\[
\mathbb{V}^X \coloneqq \{ f \colon  X \to \mathbb{V} \}.
\]
This is a real vector space in the obvious way: $\lambda f + \lambda' f'$ is defined to be the function whose value at $x \in X$ is given by the vector $\lambda f(x) + \lambda' f'(x) \in \mathbb{V}$. It is also consistent with our notation $\mathbb{R}^X$ for the vector space of all real-valued functions on $X$ in Section~\ref{sec:tensor3a}.  Consider the bilinear map
\[
\Phi \colon  \mathbb{R}^X \times \mathbb{V} \to \mathbb{V}^X, \quad (f,v) \mapsto f\cdot v,
\]
where, for any fixed $f \in \mathbb{R}^X $ and  $v \in \mathbb{V}$, the function $f \cdot v $ is defined to be $f \cdot v \colon  X \to  \mathbb{V}$, $x \mapsto f(x) v$. Applying the universal factorization property, we have
\begin{equation}\label{eq:RXV}
\xymatrix{
\mathbb{R}^X  \times \mathbb{V}
\ar[r]^-{\seg}\ar[rd]_{\Phi} & \mathbb{R}^X  \otimes \mathbb{V} \ar[d]^{F_\Phi}\\
 &  \mathbb{V}^X}
\end{equation}
with the linearization of $\Phi$ given by
\[
F_\Phi \biggl( \sum_{i=1}^r f_i \otimes v_i \biggr) =\sum_{i=1}^r f_i \cdot v_i.
\]
It is straightforward \cite[Section~1.5]{Ryan} to show that $F_\Phi$ gives an isomorphism $\mathbb{R}^X  \otimes \mathbb{V}  \cong  \mathbb{V}^X$. The standard practice is to identify the two spaces, \tie, replace $\cong$ with $=$ and take $f \otimes v$ to mean $f \cdot v$ as defined above.

In reality we are often more interested in various subspaces $\Fn(X) \subseteq \mathbb{R}^X$ whose elements satisfy properties such as continuity, differentiability, integrability, summability, {\em etc.}, and we would like to have a corresponding space of vector-valued functions $\Fn(X;\mathbb{V}) \subseteq \mathbb{V}^X$ that also have these properties. If $\mathbb{V}$ is finite-dimensional, applying the universal factorization property \eqref{eq:RXV} with the subspaces $\Fn(X) $ and $\Fn(X; \mathbb{V}) $ in place of $\mathbb{R}^X$ and $\mathbb{V}^X$ gives us
\[
\Fn(X; \mathbb{V}) = \Fn(X) \otimes \mathbb{V},
\]
and this already leads to some useful constructions.

\mypara{Extension of scalars} Any real vector space $\mathbb{U}$ has a \emph{complexification}  $\mathbb{U}_\mathbb{C}$ with elements given by
\[
u_1 \otimes 1 + u_2 \otimes i, \quad u_1, u_2 \in \mathbb{U},
\]
and scalar multiplication by $a + bi \in \mathbb{C}$ given by
\[
(a+ib)(u_1 \otimes 1 + u_2 \otimes i ) = (au_1 - bu_2) \otimes 1 + (bu_1 + au_2) \otimes i.
\]
This construction is just
\[
\mathbb{W} \otimes \mathbb{C} = \mathbb{W}_\mathbb{C} .
\]

\mypara{Polynomial matrices}  In mechanical engineering, one of the simplest and most fundamental system is the mass--spring--damper model described by
\[
q(\lambda) = \lambda^2 M + \lambda C + K,
\]
where $M, C, K \in \mathbb{R}^{n \times n}$ represent mass, damping and stiffness respectively. This a polynomial matrix, which may be viewed as either a matrix of polynomials or a polynomial with matrix coefficients. This construction is just
\[
\mathbb{R}[\lambda] \otimes \mathbb{R}^{n \times n}  = (\mathbb{R}[\lambda])^{n \times n} = (\mathbb{R}^{n \times n})[\lambda].
\]

\mypara{$C^k$- and $L^p$-vector fields} In \eqref{eq:navier} we viewed the solution to a system of partial differential equations as an element of $C^2(\mathbb{R}^3) \hatotimes C^1 [0,\infty) \otimes \mathbb{R}^3$, and in \eqref{eq:l2vecfield} we viewed an $L^2$-vector field as an element of $L^2(\mathbb{R}^3) \otimes \mathbb{C}^2$. These are consequences~of
\[
C^k(\Omega) \otimes \mathbb{R}^n = C^k(\Omega; \mathbb{R}^n) \quad \text{and} \quad
L^p(\Omega) \otimes \mathbb{C}^n = L^p(\Omega; \mathbb{C}^n) 
\]
where the latter spaces refer to the sets of vector-valued functions
\[
f = (f_1,\ldots,f_n) \colon  \Omega \to \mathbb{R}^n\ \text{or}\ \mathbb{C}^n
\]whose components $f_i$ are $k$-times continuously differentiable or $L^p$-integrable on $\Omega$ respectively.
\medskip

On the other hand, if $\mathbb{V}$ is an infinite-dimensional topological vector space, then the results depend on the choice of topological tensor product and the algebraic construction above will need to be adapted to address issues of convergence. We present a few examples involving Banach-space-valued functions, mostly drawn from \citet[Chapters~2 and 3]{Ryan}. In the following, $\mathbb{B}$ denotes a Banach space with norm $\lVert\,\cdot\,\rVert$.

\mypara{Absolutely summable sequences} This is a special case of the $L^1$-functions below but worth a separate statement:
\[
l^1(\mathbb{N}) \hatotimes_\nu \mathbb{B} = l^1(\mathbb{N};\mathbb{B}) = \biggl\{ (x_i)_{i=1}^\infty \colon  x_i \in \mathbb{B}, \sum_{i=1}^\infty \lVert x_i \rVert < \infty \biggr\}.
\]
If we use $\hatotimes_\sigma$ in place of  $\hatotimes_\nu$, we obtain a larger subspace, as we will see next.

\mypara{Unconditionally summable sequences} Absolute summability implies unconditional summability but the converse is generally false when $\mathbb{B} $ is infinite-dimen\-sional. We have
\[
l^1(\mathbb{N}) \hatotimes_\sigma \mathbb{B} = \biggl\{ (x_i)_{i=1}^\infty \colon  x_i \in \mathbb{B}, \sum_{i=1}^\infty \varepsilon_i x_i < \infty \text{ for any } \varepsilon_i =\pm 1 \biggr\}.
\]
The last condition is also equivalent to $\sum_{i=1}^\infty x_{\sigma(i)} < \infty$ for any bijection $\sigma \colon  \mathbb{N} \to \mathbb{N}$.

\mypara{Integrable functions} Let $\Omega$ be $\sigma$-finite. Then
\[
L^1(\Omega) \hatotimes_\nu \mathbb{B}  = L^1(\Omega; \mathbb{B}) = \{ f \colon  \Omega \to \mathbb{B} \colon  \lVert f \rVert_1 < \infty \}.
\]
This is called the Lebesgue--Bochner space. 
Here the $L^1$-norm is defined as $\lVert f \rVert_1 = \int_\Omega \lVert f(x) \rVert  \D x$. In fact the proof of the second half of \eqref{eq:ctsint} is via
\[
L^1(X) \hatotimes_\nu L^1(Y) = L^1(X; L^1(Y)) = L^1(X \times Y).
\]

\mypara{Continuous functions}  Let $\Omega$ be compact Hausdorff. Then
\[
C(\Omega) \hatotimes_\sigma \mathbb{B}  = C(\Omega; \mathbb{B}) = \{ f \colon  \Omega \to \mathbb{B} \colon  f \text{ continuous}\}.
\]
Here continuity is with respect to  $\lVert f \rVert_\infty = \sup_{x \in \Omega} \lVert f(x) \rVert$. In fact the proof of the first half of \eqref{eq:ctsint} is via
\[
C(X) \hatotimes_\sigma C(Y) = C(X; C(Y)) = C(X \times Y).
\]

\mypara{Infinite-dimensional diagonal matrices} For finite-dimensional square matrices, the subspace of diagonal matrices $\diag(\mathbb{R}^{n \times n}) \coloneqq \spn \{e_i \otimes e_i \colon  i =1,\ldots,n\} \cong \mathbb{R}^n$ is independent of any other structure on $\mathbb{R}^{n \times n}$. For infinite-dimensional matrices, the
closed linear span $\overline{\spn}\{e_i \otimes e_i \colon i \in \mathbb{N} \}$ 
depends on both the choice of norms and the choice of tensor products:
\[
\diag(l^p(\mathbb{N}) \hatotimes_\nu l^p(\mathbb{N})) \cong
\begin{cases}
l^1(\mathbb{N} ) &1 \le p \le 2,\\
l^{p/2}(\mathbb{N} ) & 2< p < \infty,\\
c_0(\mathbb{N} ) & p = \infty.
\end{cases}
\]
We refer readers to \citet{Holub} for other  cases with $\diag(l^p(\mathbb{N}) \hatotimes_\nu l^q(\mathbb{N}))$, or $\hatotimes_\sigma$ in place of $\hatotimes_\nu$, and to \citet{Arias} for the higher-order case $l^{p_1}(\mathbb{N} ) \otimes \dots \otimes l^{p_d}(\mathbb{N})$ with $d > 2$.

\mypara{Partial traces} For  separable Hilbert spaces $\mathbb{H}_1$ and $\mathbb{H}_2$, the partial traces are the continuous linear maps
\begin{align*}
\tr_1 \colon  ( \mathbb{H}_1 \hatotimes_\F \mathbb{H}_2 ) \hatotimes_\nu ( \mathbb{H}_1 \hatotimes_\F \mathbb{H}_2 )^* &\to \mathbb{H}_2 \hatotimes_\F \mathbb{H}_2^*,\\*
\tr_2 \colon  ( \mathbb{H}_1 \hatotimes_\F \mathbb{H}_2 ) \hatotimes_\nu ( \mathbb{H}_1 \hatotimes_\F \mathbb{H}_2 )^* &\to \mathbb{H}_1 \hatotimes_\F \mathbb{H}_1^*.
\end{align*}
Note that the domain of these maps are trace-class operators on the Hilbert space $\mathbb{H}_1 \hatotimes_\F \mathbb{H}_2$ by Example~\ref{eg:compact}.  Choose any orthonormal basis on $\mathbb{H}_1$ to obtain a Hilbert space isomorphism $\mathbb{H}_1 \hatotimes_\F\mathbb{H}_1^* \cong l^2(\mathbb{N} \times \mathbb{N})$. Then
\begin{align*}
( \mathbb{H}_1 \hatotimes_\F \mathbb{H}_2 ) \hatotimes_\nu ( \mathbb{H}_1 \hatotimes_\F \mathbb{H}_2 )^* &\cong ( \mathbb{H}_1 \hatotimes_\F \mathbb{H}_1^* ) \hatotimes_\nu ( \mathbb{H}_2 \hatotimes_\F \mathbb{H}_2 )^*\\*
&\cong  l^2(\mathbb{N} \times \mathbb{N})  \hatotimes_\nu ( \mathbb{H}_2 \hatotimes_\F \mathbb{H}_2^* )
\end{align*}
and the last space is the set of infinite matrices of the form
\begin{equation}\label{eq:dblinfmat}
\begin{bmatrix}
\Phi_{11} & \Phi_{12} & \cdots & \Phi_{1j} & \cdots \\
\Phi_{21} & \Phi_{22} & \cdots & \Phi_{2j} & \cdots \\
\vdots & \vdots & \ddots & \vdots \\
\Phi_{i1}& \Phi_{i2} & \cdots & \Phi_{ij} & \cdots \\
\vdots  & \vdots & & \vdots  & \ddots
\end{bmatrix}\!,
\end{equation}
where each $\Phi_{ij}  \colon  \mathbb{H}_2 \to \mathbb{H}_2^*$ is a Hilbert--Schmidt operator for any $i,j \in \mathbb{N}$. Given a trace-class operator $\Phi \colon  \mathbb{H}_1 \hatotimes_\F \mathbb{H}_2 \to \mathbb{H}_1 \hatotimes_\F \mathbb{H}_2$, we may express it in the form \eqref{eq:dblinfmat} and then define the first partial trace as
\[
\tr_1(\Phi) \coloneqq \sum_{i=1}^\infty \Phi_{ii}.
\]
The other partial trace $\tr_2$
may be similarly defined by reversing the roles of $\mathbb{H}_1$ and $\mathbb{H}_2$. Partial traces are an indispensable tool for working with the density operators discussed in Example~\ref{eg:density}; we refer readers to \citet[Chapter~III, Complement~E, Section~5b]{Cohen1} and \citet[Sections~2.4.3]{Nielsen} for more information.
\end{example}

Our discussion about partial traces contains quite a bit of hand-waving; we will justify some of it with the next example, where we discuss some properties of tensor products that we have used liberally above.

\begin{example}[calculus of tensor products]\label{eg:calc}
The universal factorization property is particularly useful for establishing other properties \cite[Chapter~I]{Greub} of the tensor product  operation $\otimes$ such as how it interacts with itself,
\[
\mathbb{U} \otimes \mathbb{V} \cong \mathbb{V} \otimes \mathbb{U},\quad
\mathbb{U} \otimes (\mathbb{V} \otimes \mathbb{W}) \cong ( \mathbb{U} \otimes \mathbb{V} ) \otimes \mathbb{W} \cong \mathbb{U} \otimes \mathbb{V} \otimes \mathbb{W},
\]
with direct sum $\oplus$,
\[
\mathbb{U} \otimes (\mathbb{V} \oplus \mathbb{W}) \cong  (\mathbb{U} \otimes \mathbb{V} ) \oplus (\mathbb{U} \otimes \mathbb{W}),\quad
(\mathbb{U} \oplus\mathbb{V}) \otimes \mathbb{W} \cong  (\mathbb{U} \otimes \mathbb{W} ) \oplus (\mathbb{V} \otimes \mathbb{W}),
\]
and with intersection $\cap$,
\[
(\mathbb{V} \otimes \mathbb{W}) \cap (\mathbb{V}' \otimes \mathbb{W}') =
(\mathbb{V}  \cap \mathbb{V}') \otimes ( \mathbb{W} \cap \mathbb{W}'),
\]
as well as how it interacts with linear and multilinear maps,
\[
\Lin(\mathbb{U} \otimes \mathbb{V}; \mathbb{W}) \cong
\Lin(\mathbb{U}; \Lin(\mathbb{V}; \mathbb{W})) \cong
\Mult^2(\mathbb{U},\mathbb{V}; \mathbb{W}),
\]
with duality $\mathbb{V}^* = \Lin(\mathbb{V}; \mathbb{R})$ a special case,
\[
(\mathbb{V} \otimes \mathbb{W})^* \cong \mathbb{V}^* \otimes \mathbb{W}^*, \quad \mathbb{V}^* \otimes \mathbb{W} \cong \Lin(\mathbb{V}; \mathbb{W}),
\]
and with the Kronecker product,
\[
\Lin(\mathbb{V}; \mathbb{W}) \otimes \Lin(\mathbb{V}'; \mathbb{W}') \cong \Lin(\mathbb{V} \otimes \mathbb{V}'; \mathbb{W} \otimes \mathbb{W}').
\]
Collectively, these properties form a system of calculus for manipulating tensor products of vector spaces.  When pure mathematicians speak of \emph{multilinear algebra}, this is often what they have in mind, \tie, the subject is less about manipulating individual tensors and more about manipulating whole spaces of tensors. 
Many, but not all, of these properties may be extended to other vector-space-like objects such as modules and vector bundles, or to vector spaces with additional structures such as metrics, products and topologies. One needs to exercise some caution in making such extensions. For instance, if $\mathbb{V}$ and $\mathbb{W}$ are infinite-dimensional Hilbert spaces, then
\[
\mathbb{V}^* \otimes \mathbb{W} \not\cong \Bd(\mathbb{V}; \mathbb{W}),
\]
no matter what notion of tensor product or topological tensor product one uses for $\otimes$, \tie, in this context $\Bd(\mathbb{V}; \mathbb{W})$ is not the infinite-dimensional analogue of $\Lin(\mathbb{V}; \mathbb{W})$. As we saw in
Examples~\ref{eg:compact} and \ref{eg:cftr},
the proper interpretation of $\otimes$ for Hilbert spaces is more subtle.

To establish these properties, the `uniqueness up to vector space isomorphisms' in Definition~\ref{def:tensor3c} is the tool of choice. This says that if we have two tensor spaces $\mathbb{T}$ and $\mathbb{T}'$, both satisfying the universal factorization property for a multilinear map~$\Phi$,
\[
\xymatrix{
\mathbb{V}_1\times\dots\times \mathbb{V}_d\ar[r]^-{\segt}\ar[rd]_{\Phi} & \mathbb{T}\ar[d]^{F_\Phi}\\
 & \mathbb{W}}\quad
\xymatrix{
\mathbb{V}_1\times\dots\times \mathbb{V}_d\ar[r]^-{\sigma_{\mathbb{T}'}}\ar[rd]_{\Phi} & \mathbb{T}'\ar[d]^{F'_\Phi}\\
 & \mathbb{W}}
\]
then we must have $\mathbb{T} \cong \mathbb{T}'$ as vector spaces. For each of the above properties, say, $\mathbb{U} \otimes \mathbb{V} \cong \mathbb{V} \otimes \mathbb{U}$, we plug in the vector spaces on both sides of the purported isomorphism $\cong$ as our $\mathbb{T}$ and $\mathbb{W}$ and make a judicious choice for $\Phi$, say, $\mathbb{U} \times \mathbb{V} \to \mathbb{V} \otimes \mathbb{U}$, $(u, v) \mapsto v \otimes u$, to obtain $F_\Phi$ as the required isomorphism.

While we have said early on in \eqref{eq:noncomm}  that $u \otimes v \ne v \otimes u$, note that here we did not say $\mathbb{U} \otimes \mathbb{V} = \mathbb{V} \otimes \mathbb{U}$ but $\mathbb{U} \otimes \mathbb{V} \cong \mathbb{V} \otimes \mathbb{U}$. We sometimes \emph{identify} isomorphic spaces, \ie\ replace $\cong$ with $=$, but this depends on the context. In Navier--Stokes \eqref{eq:navier}, we do not want to say $C^2(\mathbb{R}^3) \hatotimes C^1 [0,\infty) = C^1 [0,\infty) \hatotimes C^2(\mathbb{R}^3)$, since the order of the variables in a solution $v(x,y,z,t)$ carries a specific meaning. On the other hand, it does not matter whether we treat the quantum states of a composite system of two particles as  $\mathbb{H}_1 \hatotimes \mathbb{H}_2$ or  $\mathbb{H}_2 \hatotimes \mathbb{H}_1$:  either the particles are distinguishable with $\mathbb{H}_1$ and $\mathbb{H}_2$ distinct and the difference becomes merely a matter of which particle we label as first and which as second (completely arbitrary); or, if they are indistinguishable, then $\mathbb{H}_1 =\mathbb{H}_2$ and the question does not arise.

The reader may remember from our discussions in Section~\ref{sec:multmaps} about a major shortcoming of definition~\ref{st:tensor2}, that for each type of tensor there are many different types of multilinear maps. Definition~\ref{def:tensor3c} neatly classifies them in that we may use the calculus of $\otimes$ described above to reduce any type of multilinear map into a tensor of type $(p,d-p)$ for some $p$ and $d$. A noteworthy special case is the multilinear functionals in Definition~\ref{def:tensor2a} where we have
\[
\Mult(\mathbb{V}_1^*, \ldots, \mathbb{V}_p^*, \mathbb{V}_{p+1}, \ldots, \mathbb{V}_d; \mathbb{R}) \cong
\mathbb{V}_1\otimes \dots\otimes \mathbb{V}_p\otimes \mathbb{V}_{p+1}^*\otimes \dots\otimes \mathbb{V}_d^*.
\]
For instance, recall the matrix--matrix product in Example~\ref{eg:Strassen} and triple product trace in Example~\ref{eg:GI}, respectively,
\[
\Mu_{m,n,p}(A,B) = AB, \quad \tau_{m,n,p}(A,B,C) =\tr(ABC),
\]
for $A \in \mathbb{R}^{m \times n}$, $B \in \mathbb{R}^{n \times p}$, $C \in \mathbb{R}^{p \times m}$. Since
\begin{align*}
\Mu_{m,n,p} \in \Mult(\mathbb{R}^{m \times n}, \mathbb{R}^{n \times p}; \mathbb{R}^{m \times p}) &\cong (\mathbb{R}^{m \times n})^* \otimes (\mathbb{R}^{n \times p})^* \otimes \mathbb{R}^{m \times p},\\*
\tau_{m,n,p} \in \Mult(\mathbb{R}^{m \times n}, \mathbb{R}^{n \times p}, \mathbb{R}^{p \times m}; \mathbb{R}) &\cong (\mathbb{R}^{m \times n})^* \otimes (\mathbb{R}^{n \times p})^* \otimes (\mathbb{R}^{p \times m})^*,
\end{align*}
and $(\mathbb{R}^{p \times m})^* \cong  \mathbb{R}^{m \times p}$, the bilinear operator $\Mu_{m,n,p}$ and the trilinear functional $\tau_{m,n,p}$ may be regarded as elements of the same tensor space. Checking their values on the standard bases shows that they are in fact the same tensor. More generally, the $(d-1)$-linear map and $d$-linear functional given by
\begin{align*}
\Mu_{n_1,n_2,\ldots,n_d}(A_1,A_2,\ldots,A_{d-1}) &= A_1 A_2 \cdots A_{d-1}, \\*
\tau_{n_1,n_2,\ldots,n_d}(A_1,A_2,\ldots,A_d) &= \tr(A_1 A_2 \cdots A_d) ,
\end{align*}
for matrices
\[
A_1 \in \mathbb{R}^{n_1 \times n_2}, A_2 \in \mathbb{R}^{n_2 \times n_3},\ldots, A_{d-1} \in \mathbb{R}^{n_{d-1} \times n_d}, A_d \in \mathbb{R}^{n_d \times n_1},
\]
are one and the same tensor.  The tensor $\tau_{n_1,n_2,\ldots,n_d}$ will appear again in Example~\ref{eg:tennet} in the form of the \emph{matrix product states} tensor network.
\end{example}

Definition~\ref{def:tensor3c} was originally due to \citet{Bourbaki} and, as is typical of their style, provides the most general and abstract definition of a tensor, stated in language equally abstruse. As a result it is sometimes viewed with trepidation by students attempting to learn tensor products from standard algebra texts, although it has also become the prevailing definition of a tensor in modern mathematics \cite{DumFoot,KM,Lang,Vinberg}. As we saw in this section, however, it is a simple and practical idea.

\subsection{Tensors in computations III: separability}\label{sec:sepvar}

The notion of separability is central to the utility of definition~\ref{st:tensor3} in its various forms. We begin with its best-known manifestation: the separation-of-variables technique.

We will first present an abstract version to show the tensor perspective and then apply it to PDEs, finite difference schemes and integro-differential equations.

\begin{example}[separation of variables: abstract]\label{eg:sepabs}
Let $\Phi \colon  \mathbb{V}_1 \otimes  \dots \otimes \mathbb{V}_d \to \mathbb{V}_1 \otimes  \dots \otimes \mathbb{V}_d$ be a linear operator of the form
\begin{equation}\label{eq:sepop} 
\Phi = \Phi_1 \otimes I_2 \otimes \dots \otimes I_d + I_1 \otimes \Phi_2 \otimes \dots \otimes I_d + \dots + I_1 \otimes I_2 \otimes \dots \otimes \Phi_d,
\end{equation}
where $I_d$ is the identity operator on $\mathbb{V}_d$ and $\otimes$ is the Kronecker product in Example~\ref{eg:kron}. The separation-of-variables technique essentially transforms a homogeneous linear system into a collection of eigenproblems:
\begin{equation}\label{eq:septransf}
\Phi( v_1 \otimes v_2 \otimes \dots \otimes v_d) = 0 \quad \longrightarrow \quad
\left\{
\begin{aligned}
\Phi_1 (v_1) &= \lambda_1 v_1,\\
\Phi_2 (v_2) &= \lambda_2 v_2,\\
&\vdots \\
\Phi_d (v_d) &= -(\lambda_1 + \dots +\lambda_{d-1}) v_d,
\end{aligned}\right.
\end{equation}
and $\Phi$ being linear, any sum, linear combination or, under the right conditions, even integral of $ v_1 \otimes v_2 \otimes \dots \otimes v_d$ is also a solution. The constants $\lambda_1,\ldots,\lambda_{d-1}$ are called separation constants. The technique relies only on one easy fact about tensor products: for any non-zero $v \in \mathbb{V}$ and $w\in \mathbb{W}$,
\begin{equation}\label{eq:easy}
v \otimes w = v' \otimes w' \quad \Rightarrow \quad v = \lambda v',\; w = \lambda^{-1} w'
\end{equation}
for some non-zero $\lambda \in \mathbb{R}$. More generally, for any $d$ non-zero vectors $v_1 \in \mathbb{V}_1, v_2\in \mathbb{V}_2,\ldots, v_d \in \mathbb{V}_d$, we have
\begin{align*} 
& v_1 \otimes v_2 \otimes \dots \otimes v_d = v'_1 \otimes v'_2 \otimes \dots \otimes v'_d \\*
&\quad \Rightarrow \quad
v_1 = \lambda_1 v'_1,  v_2 = \lambda_2 v'_2, \ldots, v_d = \lambda_d d'_d, \quad
\lambda_1\lambda_2 \cdots \lambda_d = 1,
\end{align*}
although this more general version is not as useful for separation of variables, which just requires repeated applications of \eqref{eq:easy}.

Take $d =3$ for illustration. Dropping subscripts to avoid clutter, we have
\begin{equation}\label{eq:sepproof}
( \Phi \otimes I \otimes I + I \otimes \Psi \otimes I + I \otimes I \otimes \Theta)(u \otimes v \otimes w) = 0
\end{equation}
or equivalently $\Phi(u) \otimes v \otimes w + u \otimes \Psi(v) \otimes w + u \otimes v \otimes \Theta(w) = 0$. 
Since $\Phi(u) \otimes (v \otimes w) = u \otimes [ -\Psi(v) \otimes w -  v \otimes \Theta(w) ]$,
applying \eqref{eq:easy} we get
\[
\Phi(u) = \lambda u, \quad v \otimes w =  - \lambda^{-1} [ \Psi(v) \otimes w +  v \otimes \Theta(w) ].
\]
Rearranging the second equation, $\Psi(v) \otimes w =  v \otimes [-\Theta(w) - \lambda  w ]$, and
applying \eqref{eq:easy} again, we get
\[
\Psi(v) = \mu v, \quad \Theta(w) = -(\mu +\lambda ) w.
\]
Thus we have transformed \eqref{eq:sepproof} into three eigenproblems:
\[
\left\{
\begin{aligned}
\Phi (u) &= \lambda u,\\
\Psi (v) &= \mu v,\\
\Theta (w) &= (-\mu -\lambda)w.
\end{aligned}\right.
\]

The technique applies widely. The operators $\Phi_i \colon  \mathbb{V}_i \to \mathbb{V}_i$ involved can be any operator in any coordinates. While they are most commonly differential operators such~as
\[
\Phi_\theta (f) =  \dfrac{1}{r^2 \sin^2 \phi} \dfrac{\partial^2 f}{\partial \theta^2}, \quad \Phi_\sigma (f) = \dfrac{1}{\sigma^2 +\tau^2} \dfrac{\partial^2 f}{\partial \sigma^2}
\]
in Example~\ref{eg:sepPDE}, they may well be finite difference or integro-differential operators,
\[
\Phi_k (u_{k,n}) = u_{k-1,n} + (1-2r) u_{k,n} + r u_{k+1,n}, \quad \Phi_x (f) = a \dfrac{\partial^2 f}{\partial x^2} + b \int_0^x f,
\]
as in Examples~\ref{eg:findiff} and \ref{eg:ide} respectively. The important thing is that they must take the form\footnote{One sometimes hears of `Kronecker sums' in the $d=2$ case: $\Phi_1 \oplus \Phi_2 \coloneqq \Phi_1 \otimes I_2 + I_1 \otimes \Phi_2$. This has none of the usual characteristics of a sum:  it is not associative, not commutative, and $\Phi_1 \oplus \Phi_2 \oplus \Phi_3 \ne \Phi_1 \otimes I_2 \otimes I_3 + I_1 \otimes \Phi_2 \otimes I_3 + I_1 \otimes I_2 \otimes \Phi_3$ regardless of the order $\oplus$ is performed. As such we avoid using $\oplus$ in this sense.} in \eqref{eq:sepop}:  there should not be any `cross-terms' like $\Phi_1 \otimes I_2 \otimes \Phi_3 $, \eg\  $\partial^2 f/\partial x_1 \partial x_3$, because if so then the first variable and third variable will not be separable, as we will see.
\end{example}

The abstract generality of Example~\ref{eg:sepabs} makes the separation-of-variables technique appear trivial, although as we will see in Examples~\ref{eg:findiff} and \ref{eg:ide}, it also makes the technique widely applicable. Once we restrict to specific classes of problems like PDEs, there are much deeper results about the technique. Nevertheless, elementary treatments tend to focus on auxiliary issues such as Sturm--Liouville theory and leave out what we think are the most important questions: Why and when does the technique work?  Fortunately these questions have been thoroughly investigated and addressed in \citet{Koornwinder,Miller2,Miller} and \citet{Schobel}. The explanation is intimately related to tensors but in the sense of  tensor fields in Example~\ref{eg:tenfield2}:  the metric tensor, Ricci tensor, Nijenhuis tensor, Weyl tensor all play a role. Unfortunately, the gap between these results and our brief treatment is too wide for us to convey anything beyond a glimpse of what is involved.

\begin{example}[separation of variables: PDEs]\label{eg:sepPDE}
The technique of separation of variables for solving PDEs evidently depends on the choice of coordinates. It works perfectly for the one-dimensional wave equation in the form
\begin{equation}\label{eq:wave}
\dfrac{\partial^2 f}{\partial x^2}
- \dfrac{\partial^2 f}{\partial t^2}
=0.
\end{equation}
The standard recipe is to use the ansatz $f(x,t) = \varphi(x) \psi(t) $, obtain  $\varphi''(x)/\varphi(x) = \psi''(t)/\psi(t)$, and argue that since the left-hand side is a function of $x$ and the right-hand side of $t$, they must both be equal to some separation constant $-\omega^2$, and we obtain two ODEs $\varphi'' + \omega^2 \varphi = 0$, $\psi'' + \omega^2 \psi = 0$. We may deduce the same thing in one step from \eqref{eq:septransf}:
\[
[\partial^2_x \otimes I + I \otimes (-\partial_t^2)](\varphi \otimes \psi) = 0 \quad \longrightarrow \quad 
\left\{
\begin{aligned}
\partial^2_x \varphi &= -\omega^2 \varphi,\\
-\partial^2_t \psi &= \omega^2 \psi.
\end{aligned}\right.
\]
The ODEs have solutions $\varphi(x) = a_1 \,\rme^{\omega x} + a_2 \,\rme^{-\omega x}$, $\psi(x) = a_3 \,\rme^{\omega t} + a_4 \,\rme^{-\omega t}$ if $\omega  \ne 0$, or $\varphi(x) = a_1  + a_2 x$, $\psi(t) = a_3  + a_4 t$ if $\omega  = 0$, and any finite linear combinations of  $\varphi \otimes \psi$ give us solutions for \eqref{eq:wave}. Nevertheless, a change of coordinates $\xi  = x - t$, $\eta = x + t$ transforms \eqref{eq:wave} into
\begin{equation}\label{eq:wave2}
\dfrac{\partial^2 f}{\partial \xi \partial \eta} = 0.
\end{equation}
Note that the operator here is $\partial_x \otimes \partial_y$ and does not have the required form in \eqref{eq:sepop}. 
Indeed, the solution of \eqref{eq:wave2} is easily seen to take the form $f (\xi, \eta) = \varphi(\xi) + \psi(\eta)$,  so the usual multiplicatively separable ansatz $f (\xi, \eta) = \varphi(\xi)\psi(\eta)$ will not work.

More generally, the same argument works in any dimension $n$ to separate the spatial coordinates $x=(x_1,\ldots,x_n)$ (which do not need to be Cartesian) from the temporal ones. Applying \eqref{eq:septransf} to the $n$-dimensional wave equation,
\begin{equation}\label{eq:wave3}
\Delta f - \dfrac{\partial^2 f}{\partial t^2} = 0,
\end{equation}
we get
\[
[\Delta \otimes I + I \otimes (-\partial_t^2)](\varphi \otimes \psi) = 0  \quad \longrightarrow \quad 
\left\{
\begin{aligned}
\Delta \varphi &= -\omega^2 \varphi,\\
-\partial^2_t \psi &= \omega^2 \psi,
\end{aligned}\right.
\]
with separation constant $-\omega^2$. In the remainder of this example, we will focus on the first equation, called the $n$-dimensional  Helmholtz equation,
\begin{equation}\label{eq:helm1}
\Delta f + \omega^2 f = 0,
\end{equation}
which may also be obtained by taking the Fourier transform of \eqref{eq:wave3} in time. For $n = 2$, \eqref{eq:helm1} in Cartesian and polar coordinates are given by
\begin{equation}\label{eq:helm2}
\dfrac{\partial^2 f}{\partial x^2}
+ \dfrac{\partial^2 f}{\partial y^2}
+ \omega^2 f
=0, \quad \dfrac{\partial^2 f}{\partial r^2}
+ \dfrac{1}{r} \dfrac{\partial f}{\partial r} + \dfrac{1}{r^2}\dfrac{\partial^2 f}{\partial \theta^2}
+ \omega^2 f
=0
\end{equation}
respectively. Separation of variables works in both cases but gives entirely different solutions. Applying \eqref{eq:septransf} to $\partial^2_x \otimes I + I \otimes (\partial_y^2 + \omega^2 I)$  gives us
\[
\dfrac{\dd^2 \varphi}{\dd x^2} + k^2 \varphi = 0, \quad \dfrac{\dd^2 \psi}{\dd y^2} + (\omega^2 - k^2) \psi = 0,
\]
with separation constant $k^2$ and therefore the solution
\begin{align*}
f_k(x,y) &\coloneqq a_1 \,\rme^{\rmi [kx + (\omega^2 - k^2)^{1/2} y]} + a_2 \,\rme^{\rmi [-kx + (\omega^2 - k^2)^{1/2} y]} \\*
&\quad\ + a_3 \,\rme^{\rmi [kx - (\omega^2 - k^2)^{1/2} y]} + a_4  \,\rme^{\rmi [-kx - (\omega^2 - k^2)^{1/2} y]}.
\end{align*}
Applying \eqref{eq:septransf} to $[(r^2 \partial^2_r + \omega^2 r^2 I) \otimes I + I \otimes \partial_\theta^2](\varphi \otimes \psi) = 0 $  gives us
\[
r^2\dfrac{\dd^2 \varphi}{\dd r^2} + r\dfrac{\dd \varphi}{\dd r} +(\omega^2r^2 - k^2) \varphi = 0, \quad \dfrac{\dd^2 \psi}{\dd \theta^2} +  k^2\psi = 0,
\]
with separation constant $k^2$ and therefore the solution
\[
f_k(r,\theta) \coloneqq
a_1 \,\rme^{\rmi k\theta} J_k(\omega r) + a_2 \,\rme^{-\rmi k\theta} J_k(\omega r) + a_3 \,\rme^{\rmi k\theta} J_{-k}(\omega r) + a_4 \,\rme^{-\rmi k\theta} J_{-k}(\omega r)
\]
where $J_k$ is a Bessel function.  Any solution of the two-dimensional Helmholtz equation in Cartesian coordinates is a sum or integral of $f_k(x,y)$  over $k$ and any solution in polar coordinates is one of $f_k(r,\theta)$ over $k$. Analytic solutions in different coordinate systems provide different insights. There are exactly two more such coordinate systems where separation of variables works;  we call these \emph{separable coordinates}. For $n=2$, \eqref{eq:helm1} has exactly four systems of separable coordinates: Cartesian, polar, parabolic and elliptic. For $n = 3$, there are exactly eleven \cite{Eisenhart}.

The fundamental result that allows us to deduce these numbers is the  St\"ackel condition: the $n$-dimensional Helmholtz equation in coordinates $x_1,\ldots,x_n$ can be solved using the separation-of-variables technique if and only if (a) the Euclidean metric tensor $g$ is a diagonal matrices in this coordinate system, and (b) if $g =\diag(g_{11},\ldots,g_{nn})$, then there exists an invertible matrix of the form
\[
S = \begin{bmatrix}
s_{11}(x_1) & s_{12} (x_1) & \cdots & s_{1n}(x_1) \\
s_{21}(x_2) & s_{22} (x_2) & \cdots & s_{2n}(x_2) \\
\vdots & \vdots & & \vdots \\
s_{n1}(x_n) & s_{n2} (x_n) & \cdots & s_{nn}(x_n)
\end{bmatrix}
\]
with
\begin{equation}\label{eq:stackel}
g_{jj}^{-1} = (S^{-1})_{1j}, \quad j =1,\ldots,n.
\end{equation}
Here $S$ is called a St\"ackel matrix for $g$ in coordinates $x_1,\ldots,x_n$; note that the $i$th row of $S$ depends only on the $i$th coordinate. The Euclidean metric tensor $g$ is a covariant $2$-tensor field that parametrizes the Euclidean inner product at different points in a coordinate system. As we saw in Example~\ref{eg:tenfield}, any $2$-tensor field over $\mathbb{R}^n$  may be represented as an $n \times n$ matrix of functions. 
Take $n =3$; the Euclidean metric tensor on $\mathbb{R}^3$ in Cartesian, cylindrical,  spherical and parabolic\footnote{Given by $x = \sigma \tau \cos \varphi$, $y = \sigma \tau \sin \varphi$, $z =(\sigma^2 - \tau^2)/2$.} coordinates are
\begin{alignat*}{3}
g(x,y,z) &= \begin{bmatrix} 1 & 0 & 0 \\ 0 & 1 & 0 \\ 0 & 0 & 1 \end{bmatrix}\!, 
&g(r,\theta,z) &= \begin{bmatrix} 1 & 0 & 0 \\ 0 & r^2 & 0 \\ 0 & 0 & 1 \end{bmatrix}\!, \\*
g(r,\theta,\phi) &= \begin{bmatrix} 1 & 0 & 0 \\ 0 & r^2 & 0 \\ 0 & 0 & r^2 \sin^2 \theta \end{bmatrix}\!,\quad
&g(\sigma,\tau,\phi) &= \begin{bmatrix} \sigma^2 + \tau^2 & 0 & 0 \\ 0 & \sigma^2 + \tau^2 & 0 \\ 0 & 0 & \sigma^2\tau^2 \end{bmatrix}\!.
\end{alignat*}
The St\"ackel condition is satisfied for these four coordinate systems using the following matrices and their inverses:
\begin{alignat*}{3} 
&\text{Cartesian} & S& =
\begin{bmatrix}
0 & -1 & -1\\
0 & 1 & 0\\
1 & 0 & 1
\end{bmatrix}\!, 
&S^{-1} &=
\begin{bmatrix}
1 & 1 & 1\\
0 & 1 & 0\\
-1 & -1 & 0
\end{bmatrix}\!, \\*
&\text{cylindrical}\quad & S& =
\begin{bmatrix}
0 & -\frac{1}{r^2} & -1\\[3pt]
0 & 1 & 0\\[3pt]
1 & 0 & 1
\end{bmatrix}\!, 
&S^{-1} &=
\begin{bmatrix}
1 & \frac{1}{r^2} & 1\\[3pt]
0 & 1 & 0\\[3pt]
-1 & -\frac{1}{r^2} & 0
\end{bmatrix}\!, \\
&\text{spherical} & S& =
\begin{bmatrix}
1 & -\frac{1}{r^2} & 0\\[3pt]
0 & 1 & -\frac{1}{\sin^2 \theta} \\[3pt]
0 & 0 & 1
\end{bmatrix}\!, 
\quad &S^{-1} &=
\begin{bmatrix}
1 & \frac{1}{r^2} & \frac{1}{(r^2 \sin^2 \theta)}\\[3pt]
0 & 1 & \frac{1}{\sin^2 \theta}\\[3pt]
0 & 0 & 1
\end{bmatrix}\!, \\*
&\text{parabolic} & S& =
\begin{bmatrix}
\sigma^2 & -1 & -\frac{1}{\sigma^2}\\[3pt]
\tau^2 & 1 & -\frac{1}{\tau^2}\\[3pt]
0 & 0 & 1
\end{bmatrix}\!, 
\quad&S^{-1} &=
\begin{bmatrix}
\frac{1}{(\sigma^2 +\tau^2)}& \frac{1}{(\sigma^2 +\tau^2)} & \frac{1}{(\sigma^2\tau^2)} \\[3pt]
-\frac{\tau^2}{(\sigma^2 +\tau^2)} & \frac{\sigma^2}{(\sigma^2 +\tau^2)} & \frac{1}{\tau^2 - 1/\sigma^2} \\[3pt]
0 & 0 & 1
\end{bmatrix}
\end{alignat*}
The matrices on the left are St\"ackel matrices for $g$ in the respective coordinate system, \tie, the entries in the first row of $S^{-1}$ are exactly the reciprocal of the entries on the diagonal of $g$.

What we have ascertained is that the three-dimensional Helmholtz equation can be solved by separation of variables in these four coordinate systems,   without writing down a single differential equation. In fact, with  more effort, one can show that there are exactly eleven such separable coordinate systems:
\begin{enumerate}[\upshape (i),itemindent=4pt]
\setlength\itemsep{3pt}
\item Cartesian,
\item cylindrical,
\item spherical,
\item parabolic,
\item paraboloidal,
\item ellipsoidal,
\item conical,
\item\label{it:ps} prolate spheroidal,
\item oblate spheroidal,
\item\label{it:ec} elliptic cylindrical,
\item parabolic cylindrical.
\end{enumerate}
These eleven coordinate systems have been thoroughly studied in engineering \cite{Moon,Moon2}, where they are used for different tasks:\ \ref{it:ps}~for modelling radiation from a slender spheroidal antenna, \ref{it:ec}~for heat flow in a bar of elliptic cross section, {\em etc.} More generally, the result can be extended to a Riemannian manifold $M$ of arbitrary dimension: a system of local coordinates is separable if and only if both the Riemannian metric tensor $g$ and the Ricci curvature tensor $\bar{R}$ are diagonal matrices in those coordinates,\footnote{We did not need to worry about the Ricci tensor because
$\mathbb{R}^n$  
  is a so-called Einstein manifold where $g$ and $\bar{R}$ differ by a scalar multiple.  See Example~\ref{eg:tenfield2} for a cursory discussion of $g$ and $\bar{R}$.} and $g$ satisfies the St\"ackel condition \cite{Eisenhart}.  In particular, we now see why separation of variables does not apply to PDEs with mixed partials like \eqref{eq:wave2}: when the metric tensor is diagonal, $g =\diag(g_{11},\ldots,g_{nn})$, the Laplacian takes the form
\begin{equation}\label{eq:laplace}
\Delta = 
\sum_{i=1}^n \dfrac{1}{\sqrt{\det(g)}} \dfrac{\partial}{\partial x_i} \dfrac{\sqrt{\det(g)}}{g_{ii}} \dfrac{\partial f}{\partial x_i} 
\end{equation}
and thus does not contain any mixed derivatives. We may use \eqref{eq:laplace} to write down the Helmholtz equation on $\mathbb{R}^3$ in cylindrical, spherical and parabolic coordinates:
\begin{align*}
\dfrac{\partial^2 f}{\partial r^2} + \dfrac{1}{r} \dfrac{\partial f}{\partial r} + \dfrac{1}{r^2} \dfrac{\partial^2 f}{\partial \theta^2} + \dfrac{\partial^2 f}{\partial z^2} + \omega^2 f
&=0,\\*
\dfrac{\partial^2 f}{\partial r^2} + \dfrac{2}{r} \dfrac{\partial f}{\partial r} + \dfrac{1}{r^2 \sin^2 \phi} \dfrac{\partial^2 f}{\partial \theta^2} + \dfrac{\cos \phi}{r^2 \sin^2 \phi} \dfrac{\partial f}{\partial \phi} +\dfrac{1}{r^2} \dfrac{\partial^2 f}{\partial \phi^2} +\omega^2 f
&=0,\\*
\dfrac{1}{\sigma^2 +\tau^2} \biggl[\dfrac{\partial^2 f}{\partial \sigma^2} + \dfrac{1}{\sigma} \dfrac{\partial f}{\partial \sigma} + \dfrac{\partial^2 f}{\partial \tau^2} + \dfrac{1}{\tau} \dfrac{\partial f}{\partial \tau}\biggr] + 
\dfrac{1}{\sigma^2 \tau^2} \dfrac{\partial^2 f}{\partial \phi^2}  + \omega^2 f
&=0.
\end{align*}
It is not so obvious that these are amenable to separation of variables, speaking to the power of the St\"ackel condition. While we have restricted to the Helmholtz equation for concreteness, the result has been generalized to higher-order semilinear PDEs \cite[Theorem~3.8]{Koornwinder}.
\end{example}

Separation of variables is manifestly about tensors in the form of definition~\ref{st:tensor3} and particularly Definition~\ref{def:tensor3a}, as it involves, in one way or another, a rank-one tensor $\varphi \otimes \psi \otimes \dots \otimes \theta$ or a sum of these. It may have come as a surprise that tensors in the sense of definition~\ref{st:tensor1} also play a major role in Example~\ref{eg:sepPDE}, but this should be expected as the tensor transformation rules originally came from a study of the relations between different coordinate systems.

By our discussion in Example~\ref{eg:hyp}, whether a function is of continuous variables or discrete ones, whether the variables are written as arguments or indices, should not make a substantive difference. Indeed, as we saw in Example~\ref{eg:sepabs}, separation of variables ought to apply verbatim in cases where the variables are finite or discrete. Our next example, adapted from \citet[Example~3.2.2]{Thomas}, is one where the rank-one tensors are of the form $u = a \otimes b$ with $a = (a_0,\ldots,a_m) \in \mathbb{R}^{m+1}$ and $b =(b_n)_{n=0}^\infty \in c_0(\mathbb{N} \cup \{0\})$.

\begin{example}[separation of variables: finite difference]\label{eg:findiff}
Consider the recurrence relation
\begin{equation}\label{eq:recur}
\left\{
\begin{aligned}
u_{k,n+1} &= r u_{k-1,n} + (1-2r) u_{k,n} + r u_{k+1,n}  &k &=1,\ldots,m-1,\\
u_{0,n+1} &= 0 = u_{m,n+1}&&\\
u_{k,0} &= f\biggl(\dfrac{k}{m}\biggr) &k &=0,1,\ldots,m,
\end{aligned}
\right.
\end{equation}
with $n =0,1,2,\ldots,$ and $r >0$ is some fixed constant. This comes from a forward-time centred space discretization scheme applied to an initial--boundary value problem for a one-dimensional heat equation on $[0,1]$, although readers lose nothing by simply treating this example as one of solving recurrence relations.  To be consistent with our notation, we have written $u_{k,n} = u(x_k,t_n)$ instead of the more common $u_k^n$ with the time index in superscript. Applying \eqref{eq:septransf} to \eqref{eq:recur} with
\[
\Phi_k(u_{k,n}) \coloneqq r u_{k-1,n} + (1-2r) u_{k,n} + r u_{k+1,n}, \quad \Psi_n (u_{k,n}) \coloneqq u_{k,n+1},
\]
we get
\[
[\Phi_k \otimes I + I \otimes (-\Psi_n)](a \otimes b) = 0\quad \longrightarrow \quad 
\left\{
\begin{aligned}
\Phi_k (a_k) &= \lambda a_k,\\
-\Psi_n (b_n) &= -\lambda b_n,
\end{aligned}\right.
\]
with separation constant $\lambda$. We write these out in full:
\begin{alignat*}{2}
r a_{k-1} + (1-2r) a_k + r a_{k+1} &= \lambda a_k, \quad & k &=1,\ldots,m-1,\\*
b_{n+1} &= \lambda b_n, & n & =0,1,2,\dots.
\end{alignat*}
The second equation is trivial to solve:  $b_n = \lambda^n b_0$. 
Noting that the boundary conditions $u_{0,n+1} = 0 = u_{m,n+1}$ give $a_0 = 0 =a_m$, we see that the first equation is a tridiagonal  eigenproblem:
\[
\begin{bmatrix}
1 - 2r &       r &            &            &\\
        r & 1-2r & r          &            &\\
          &       r & 1 - 2r  &          r &\\
          &          & \ddots & \ddots & \\
          &          &            & r          & 1-2r
\end{bmatrix}
\begin{bmatrix}
a_1\\
a_2 \\
a_3 \\
\vdots \\
a_{m-1}
\end{bmatrix}
=
\lambda
\begin{bmatrix}
a_1\\
a_2 \\
a_3 \\
\vdots \\
a_{m-1}
\end{bmatrix}\!.
\]
The eigenvalues and eigenvectors of a tridiagonal Toeplitz matrices have well-known closed-form expressions:
\[
\lambda_j = 1 - 4r \sin^2 \biggl( \dfrac{j \pi}{2m} \biggr), \quad a_{jk} = \sin \biggl( \dfrac{j k\pi}{m} \biggr), \quad j,k =1,\ldots, m-1,
\]
where $a_{jk}$ is the $k$th coordinate of the $j$th eigenvector. Hence we get
\[
u_{k,n} = \sum_{j=1}^{m-1} c_j b_0 \biggl[  1 - 4r \sin^2 \biggl( \dfrac{j \pi}{2m} \biggr) \biggr]^n \sin \biggl( \dfrac{j k\pi}{m} \biggr).
\]
The initial condition $u_{k,0} = f(k/m)$, $k =0,1,\ldots,m$, may be used to determine the coefficients $c_1,\ldots,c_{m-1}$ and $b_0$.

The von Neumann stability analysis is a special case of such a discrete separation of variables with $a_k = \rme^{\rmi jk \pi/m}$ and $b_n = \lambda^n$, \tie,
\[
u_{k,n} = \lambda^n \,\rme^{\rmi jk \pi/m}.
\]
Substituting this into the recursion in \eqref{eq:recur} and simplifying gives us
\[
\lambda = 1 - 4r \sin^2 \biggl( \dfrac{j \pi}{2m} \biggr),
\]
and since we require $\lvert \lambda \rvert \le 1$ for stability, we get that $r \le 1/2$, an analysis familiar to readers who have studied finite difference methods.
\end{example}

The next example is adapted from \citet{Kostoglou}, just to show the range of applicability of Example~\ref{eg:sepabs}.

\begin{example}[separation of variables: integro-differential equations]\label{eg:ide}\hspz%
Let us consider the following integro-differential equation arising from the study of heterogeneous heat transfer:
\begin{equation}\label{eq:ide}
\dfrac{\partial f}{\partial t} = a \dfrac{\partial^2 f}{\partial x^2} + b \int_0^x f(y,t)  \D y - f,
\end{equation}
with $a,b \ge 0$ and $f  \colon  [0,1] \times \mathbb{R} \to \mathbb{R}$. Note that at this juncture, if we simply differentiate both sides to eliminate the integral, we will introduce mixed derivatives and thus prevent ourselves from using separation of variables. Nevertheless, our interpretation of separation of variables in Example~\ref{eg:sepabs} allows for integrals. If we let
\[
\Phi_x(f) \coloneqq  \dfrac{\partial^2 f}{\partial x^2} -f + b \int_0^x f(y,t)  \D y , \quad \Psi_t(f)  \coloneqq -\dfrac{\partial f}{\partial t},
\]
then \eqref{eq:septransf} gives us
\[
[\Phi_x \otimes I + I \otimes \Psi_t](\varphi \otimes \psi) = 0 \quad \longrightarrow \quad 
\left\{
\begin{aligned}
\Phi_x (\varphi) &= \lambda \varphi,\\
\Psi_t (\psi) &= -\lambda \psi,
\end{aligned}\right.
\]
with separation constant $\lambda$. Writing these out in full, we have
\[
a \dfrac{\dd^2 \varphi}{\dd x^2} + (\lambda -1) \varphi + b \int_0^x \varphi(y)  \D y = 0, \quad
\dfrac{\dd\psi}{\dd t} + \lambda \psi  = 0.
\]
The second equation is easy: $\psi(t) = c \,\rme^{-\lambda t}$ for an arbitrary constant $c$ that could be determined with an initial condition. \citet{Kostoglou} solved the first equation  in a convoluted manner involving Laplace transforms and partial fractions, but this is unnecessary; at this point it is harmless to simply differentiate and eliminate the integral. With this, we obtain a third-order homogeneous ODE with constant coefficients, $a \varphi''' + (\lambda -1) \varphi' + b \varphi =0$, whose solution is standard. Physical considerations show that its characteristic polynomial
\[
r^3 +\frac{\lambda -1}{a r + \frac{b}{a} =0}
\]
must have one real and two complex roots {$r_1, r_2 \pm i r_3$}, and thus the solution is given by $c_1\,\rme^{r_1x} + \rme^{r_2 x} (c_2 \cos r_3x + c_3 \sin r_3 x)$ with arbitrary constants $c_1,c_2,c_3$ that could be determined with appropriate boundary conditions.
\end{example}

The last four examples are about exploiting separability in the structure of the solutions; the next few are about exploiting separability in the structure of the problems.

\begin{example}[separable ODEs]
The easiest ODEs to solve are probably the separable ones,
\begin{equation}\label{eq:sepODE}
\dfrac{\dd y}{\dd x} = f(x)g(y),
\end{equation}
with special cases $\dd y/\dd x = f(x)$ and $\dd y/\dd x = g(y)$ when one of the functions is constant. Solutions are, at least in principle, given by direct integration,
\[
\int \dfrac{\dd y}{g(y)} = \int f(x) \D x + c,
\]
even though closed-form expressions still depend on having closed-form integrals. The effort is not so much in solving \eqref{eq:sepODE} but in seeking a transformation  into \eqref{eq:sepODE}, \tie, given an inseparable ODE
\[
\dfrac{\dd y}{\dd x} = K(x,y),
\]
we would like to find a differentiable $\Phi \colon  \mathbb{R}^2 \to \mathbb{R}^2$, $(x,y) \mapsto (u,v)$, so that
\begin{equation}\label{eq:sepODEtransf}
\dfrac{\dd v}{\dd u} = f(u) g(v).
\end{equation}
These transformations may be regarded as a special case of the transformation rules for tensor fields (two-dimensional vector fields in this case) discussed in Example~\ref{eg:tenfield}. For instance,
\begin{equation}\label{eq:polar}
\dfrac{\dd y}{\dd x} = \dfrac{y-x}{y+x}\quad \xrightarrow{\begin{bsmallmatrix} x \rule[-.3\baselineskip]{0pt}{\baselineskip} \\ y \end{bsmallmatrix} = \begin{bsmallmatrix} r \cos \theta \\ r \sin \theta \end{bsmallmatrix}}  \quad \dfrac{\dd r}{\dd \theta} =  -r
\end{equation}
is the two-dimensional version of the tensor transformation rule in \eqref{eq:carsph}. There is a wide variety of scenarios where this is possible:
\begin{alignat*}{5}
\dfrac{\dd y}{\dd x} &= f( ax + by +c), & v &= ax + by + c,& \dfrac{\dd v}{\dd x} &= a + b f(v), \\*
\dfrac{\dd y}{\dd x} &= f\biggl(\dfrac{y}{x}\biggr), & v &= \dfrac{y}{x},  & \dfrac{\dd v}{\dd x} &=  \dfrac{f(v) - v}{x},\\
\dfrac{\dd y}{\dd x} &= f(x)y + g(x), & v &= y\,\rme^{-\int f(x)  \D x},  & \dfrac{\dd v}{\dd x} &= g(x)\,\rme^{-\int f(x)  \D x} ,\\
\dfrac{\dd^n y}{\dd x^n} &= f\biggl(\dfrac{\dd^{n-1} y}{\dd x^{n-1}}\biggr)g(x), & v &= \dfrac{\dd^{n-1} y}{\dd x^{n-1}},& \dfrac{\dd v}{\dd x} &= f(v)g(x),\\*
\dfrac{\dd y}{\dd x} &= f\biggl(\dfrac{a x + b y + c}{a' x + b' y + c'} \biggr),\quad & \begin{bmatrix} u\\ v \end{bmatrix} &= \begin{bmatrix} x - \dfrac{b'c + bc'}{a'b - ab'} \\[8pt] y- \dfrac{a'c + ac'}{a'b - ab'} \end{bmatrix}\!,  \quad & \dfrac{\dd v}{\dd u} &=  f\biggl( \dfrac{a + bv/u}{a' + b'v/u}\biggr),
\end{alignat*}
noting that the last case reduces to the second, and if $u = x$, we do not introduce a new variable \cite[Chapter~1]{Walter}. Nevertheless, unlike the case of PDEs discussed in Example~\ref{eg:sepPDE}, there does not appear to be a systematic study of when transformation to \eqref{eq:sepODEtransf} is possible. Something along the lines of the St\"ackel condition, but for ODEs, would be interesting and useful.
\end{example}

There is a close analogue of separable ODEs for integral equations with kernels that satisfy a finite version of  \eqref{eq:HS2}. As we will see below, they have neat and simple solutions, but unfortunately such kernels are not too common in practice; in fact they are regarded as a degenerate case in the study of integral equations. The following discussion is adapted from \citet[Chapter~2]{Kanwal}.

\begin{example}[separable integral equations]\label{eg:sepIE}
Let us consider Fredholm integral equations of the first and second kind:
\begin{equation}\label{eq:Fred}
g(x) = \int_a^b K(x,y) f(y)  \D y, \quad f(x) = g(x) + \lambda \int_a^b K(x,y) f(y)  \D y ,
\end{equation}
with given constants $a < b$ and functions $g \in C[a,b]$ and $K \in C([a,b] \times [a,b])$. The goal is to solve for  $f \in C[a,b]$ and, in the second case, also $\lambda \in \mathbb{R}$. The kernel $K$ is said to be degenerate or separable if
\begin{equation}\label{eq:degenker}
  K(x,y) = \sum_{i=1}^n \varphi_i(x) \psi_i(y)
\end{equation}
for some $\varphi_1,\ldots,\varphi_n, \psi_1,\ldots,\psi_n \in C[a,b]$, assumed known and linearly independent. We set
\[
v_i \coloneqq \int_a^b \psi_i(y) f(y)  \D y, \quad b_i \coloneqq \int_a^b \psi_i(x) g(x)  \D x, \quad
a_{ij} \coloneqq \int_a^b \varphi_j(x) \psi_i(x)  \D x
\]
for $i,j=1,\ldots,n$. Plugging \eqref{eq:degenker} into \eqref{eq:Fred} and integrating with respect to $y$ gives~us
\begin{equation}\label{eq:Fred2}
g(x) = \sum_{j=1}^n v_j \varphi_j(x), \quad   f(x) = g(x) + \lambda \sum_{i=1}^n v_j \varphi_j(x).
\end{equation}
Multiplying by $\psi_i(x)$ and integrating with respect to $x$ gives us
\[
b_i =  \sum_{j=1}^n a_{ij}v_j, \quad v_i  = b_i +  \lambda \sum_{j=1}^n a_{ij} v_j
\]
for $i=1,\ldots,n$, or, in matrix forms,
\[
Av = b, \quad (I-\lambda A)v = b
\]
respectively. For the latter, any $\lambda^{-1}$ not an eigenvalue of $A$ gives a unique solution for $v$ and thus for $f$ by \eqref{eq:Fred2}; for the former, considerations similar to those for obtaining minimum-norm solutions in numerical linear algebra allow us to obtain~$f$. 

It is straightforward to extend the above to Fredholm integro-differential equations  \cite[Chapter~16]{Wazwaz} such as
\[
\dfrac{\dd^kf}{\dd x^k} = g(x) + \lambda \int_a^b K(x,y) f(y)  \D y.
\]
Note that if $k =0$, this reduces to the second equation in \eqref{eq:Fred}, but it is of a different nature from the integro-differential considered in Example~\ref{eg:ide}, which is of Volterra type and involves partial derivatives. 

In principle, the discussion extends to the $L^2$-kernels in Examples~\ref{eg:compact} and \ref{eg:Mercer} with $n = \infty$, although in this case we would just be  reformulating Hilbert--Schmidt operators in terms of infinite-dimensional matrices in $l^2(\mathbb{N}\times \mathbb{N})$;  the resulting equations are no easier to analyse or solve \cite[Chapter~7]{Kanwal}.
\end{example}

While separable kernels are uncommon in integral equations, they are quite common in multidimensional integral transforms.

\begin{example}[separable integral transforms]\label{eg:septransf}
  Whenever one iteratively applies a univariate integral transform
to each variable of a multivariate function  separately, one obtains an iterated integral
\[
\int_{a_1}^{b_1} \cdots \int_{a_d}^{b_d} f(x_1,\dots,x_d) K_1(x_1, y_1)\cdots K_d(x_d, y_d) \, \dd x_1 \cdots\dd x_d
\]
with respect to the separable kernel $K_1 \otimes \dots \otimes K_d$ in Example~\ref{eg:concretetenprod}\ref{it:tenprod}. We will just mention two of the most common.
The multidimensional Fourier transform is given by
\[
\mathcal{F}(f)(\xi_1,\ldots,\xi_d) = \int_{-\infty}^\infty \cdots \int_{-\infty}^\infty
f(t_1,\ldots,t_d)\, \rme^{-\rmi \xi_1 t_1- \cdots -\rmi \xi_d t_d} \, \dd t_1 \cdots\,\dd t_d
\]
and the multidimensional Laplace transform is given by
\[
\mathcal{L}(f)(s_1,\ldots,s_d) = \int_0^\infty \cdots \int_0^\infty
f(t_1,\ldots,t_d)\, \rme^{-s_1t_1 - \cdots - s_d t_d} \, \dd t_1 \cdots \,\dd t_d.
\]
We point readers to \citet[Chapter~9]{Osgood} and \citet[Chapter~13]{Bracewell}  for properties and applications of the former, and \citet{Debnath} and \citet{Jaeger} for the latter.
\end{example} 

We briefly discuss another closely related notion of separability before our next examples, which will require it.

\begin{example}[additive separability]
Our interpretation of separability thus far is multiplicative separability of the form
\[
  f(x_1,x_2,\ldots,x_d) =\varphi_1(x_1)\varphi_2(x_2) \cdots \varphi_d(x_d)
  \]
  and not additive separability of the form
\begin{equation}\label{eq:addsep}
g(x_1,x_2,\ldots,x_d) = \varphi_1(x_1) + \varphi_2(x_2) + \dots + \varphi_d(x_d).
\end{equation}
However, as we saw in the solution of the one-dimensional wave equation in Example~\ref{eg:sepPDE}, what is multiplicatively separable in one coordinate system is additively separable in another. The two notions are intimately related although the exact relation depends on the context. In PDEs, the relation is simple \cite[Section~3.1.1]{Miller2}, at least in principle. If
\[
F( D^d f(x),\  D^{d-1} f(x), \ldots, D f(x), f(x), x) = 0
\]
has a multiplicatively separable solution, then the transformation $g(x) \coloneqq \rme^{f(x)}$ gives us a PDE
\[
G( D^d g(x),\  D^{d-1} g(x), \ldots, D g(x), g(x), x) = 0
\]
with an additively separable solution of the form \eqref{eq:addsep}. The St\"ackel condition in Example~\ref{eg:sepPDE} was in fact originally proved for an additively separable solution of the form \eqref{eq:addsep} for the Hamilton--Jacobi equation.

We describe a readily checkable condition due to  \citet{Scheffe} for deciding additive separability, which is important in statistics, particularly in generalized additive models \cite{Hastie} and in the analysis of variance \cite[Section~4.1]{ScheffeBook}.
Given any $C^2$-function $f \colon  \Omega \subseteq \mathbb{R}^d \to \mathbb{R}$, there exists a transformation $\varphi_0 \colon  \mathbb{R} \to \mathbb{R}$ so that $g = \varphi_0 \circ f$ has a decomposition of the form \eqref{eq:addsep} if and only if its Hessian takes the form
\[
\nabla^2f(x)= h(f(x)) \nabla f(x) \nabla f(x)^\tp
\]
for some function $h \colon  \mathbb{R} \to \mathbb{R}$. In this case
\[
\varphi_0(t) = c \int \rme^{-\int h(t)  \D t} \D t + c_0, \quad \varphi_i(t) = c \int  \rme^{-\int h(t)  \D t} \dfrac{\partial f}{\partial x_i}(t)   \D t + c_i,
\]
for $i=1,\ldots,d$, with constants $c >0$ and $c_0,c_1,\ldots,c_d \in \mathbb{R}$ arising from the indefinite integrals. For instance, the condition may be used to verify that
\[
\begin{aligned}
  f(x,y,z) &=  x\sqrt{(1-\smash{y^2})(1-z^2)} +  y\sqrt{(1-x^2)(1-z^2)} \\*
  &\quad \quad+ z\sqrt{(1-x^2)(1-\smash{y^2})} - xyz,\\
f(a,b,c,d) &= \dfrac{a + b + c + d + abc + abd + acd + bcd}{1 + ab + ac + ad + bc + bd + cd + abcd},
\end{aligned}
\]
are both additively separable with $\varphi_i(t) = \sin^{-1}(t)$, $i=0,1,2,3$, in the first case and $\varphi_i(t) = \log(1+t) - \log(1-t)$, $i=0,1,2,3,4$, in the second.
\end{example}

The next two examples show that it is sometimes more natural to formulate separability in the additive form.

\begin{example}[separable convex optimization]
Let $A \in \mathbb{R}^{m \times n}$, $b \in \mathbb{R}^m$ and $c \in \mathbb{R}^n$.
We consider the linear programming problem
\begin{alignat*}{3}
&\text{LP} &&
\begin{tabular}{lll}
&minimize & $c^\tp x$\\
&subject to & $Ax \le b$, $x \ge 0$
\end{tabular}&
\intertext{and integer linear programming problem}
&\text{ILP} &&
\begin{tabular}{lll}
&minimize & $c^\tp x$\\
&subject to & $Ax \le b$, $x \ge 0$, $x \in \mathbb{Z}^n$
\end{tabular}&
\end{alignat*}
where, as is customary in optimization, inequalities between vectors are interpreted in a coordinate-wise sense.  As we will be discussing complexity and will assume finite bit length inputs with $A$, $b$, $c$ having rationals entries, by clearing denominators, there is no loss of generality in assuming that they have integer entries, \tie, $A \in \mathbb{Z}^{m \times n}$, $b \in \mathbb{Z}^m$, $c \in \mathbb{Z}^n$.

LP is famously solvable in polynomial time to arbitrary $\varepsilon$-relative accuracy \cite{Khachiyan}, with its time complexity improved over the years, notably in \citet{Karmarkar} and \citet{Vaidya}, to its current bound in \citet{CLS} that is essentially in terms of $\omega$, the exponent of matrix multiplication we saw in Example~\ref{eg:Strassen}. The natural question is: Polynomial in what? The aforementioned complexity bounds invariably involve the bit length of the input,
\[
L \coloneqq \sum_{i=1}^m \sum_{j=1}^n \log_2 (\lvert a_{ij} \rvert + 1) + 
\sum_{i=1}^m \log_2 (\lvert b_i \rvert + 1) + \log_2(mn) + 1,
\]
but ideally we want an algorithm that runs in time polynomial in the size of the structure, rather than the size of the actual numbers involved. This is called a strongly polynomial-time algorithm. Whether it exists for LP is still famously unresolved \cite[Problem~9]{SmaleProb}, although a result of \citet{Tardos} shows that it is plausible: there is a polynomial-time algorithm for LP whose time complexity is independent of the vectors $b$ and $c$ and depends only on $m,n$ and the largest subdeterminant\footnote{Recall from Example~\ref{eg:hyp} that a matrix is a function $A \colon  [m] \times [n] \to \mathbb{R}$ and here $A_{\sigma,\tau}$ is the function restricted to the subset $\sigma \times \tau \subseteq [m] \times [n]$. Recall from page~\pageref{eg:rank} that a non-square determinant is identically zero.}
of $A$,
\[
\Delta \coloneqq \max \{ \lvert \det A_{\sigma \times \tau} \rvert \colon  \sigma \subseteq [m],\, \tau \subseteq [n] \}.
\]
Note that $L$ depends on $A,b,c$ and $\Delta$ only on $A$.

We will now assume that we have a box constraint $0 \le x \le 1$ in our LP and ILP. In this case ILP becomes a zero-one ILP with $x \in \{0,1\}^n$; such an ILP is also known to have time complexity that does not depend on $b$ and $c$ although it generally still depends on the entries of $A$ and not just $\Delta$ \cite{Frank}. We will let the time complexity of the box-constrained LP be $\operatorname{LP}(m, n, \Delta)$ and let that of the zero-one ILP be $\operatorname{ILP}(m, n, A)$. The former has time complexity
polynomial in $m,n,\Delta$ by \citet{Tardos}; the latter is  polynomial-time if $A$ is \emph{totally unimodular}, \ie\ when $\Delta = 1$ \cite[Chapter~19]{Schrijver}.

Observe that the linear objective
\[
c^\tp x = c_1 x_1 + \dots + c_n x_n
\]
is  additively separable. It turns out that this separability, not so much that it is linear, is the key to polynomial-time solvability. Consider an objective function $f \colon  \mathbb{R}^n \to \mathbb{R}$ of the form
\[
f(x) = f_1(x_1) + \dots + f_n(x_n),
\]
where $f_1,\ldots, f_n \colon  \mathbb{R} \to \mathbb{R}$ are all convex. Note that if we set $f_i(x_i) = c_i x_i$, which is linear and therefore convex, we recover the linear objective. A surprising result of \citet{Hochbaum} shows that the separable convex programming problem
\begin{alignat*}{3}
&\text{SP} &&
\begin{tabular}{rll}
&minimize & $ f_1(x_1) + \dots +  f_n(x_n)$\\
&subject to & $Ax \le b$, $0 \le x \le 1$
\end{tabular}&
\intertext{and its zero-one variant}
&\text{ISP} &&
\begin{tabular}{rll}
&minimize & $ f_1(x_1) + \dots +  f_n(x_n)$\\
&subject to & $Ax \le b$, $x \in \{0,1\}^n$
\end{tabular}&
\end{alignat*}
are solvable to $\varepsilon$-accuracy in time complexities
\begin{align*}
\operatorname{SP}(m,n,\Delta) &= \log_2 \biggl(\dfrac{\beta}{2\varepsilon}\biggr) \operatorname{LP}(m, 8n^2 \Delta, \Delta),\\*
\operatorname{ISP}(m,n,A) &= 
\log_2 \biggl(\dfrac{\beta}{2n \Delta}\biggr) \operatorname{LP}(m, 8n^2 \Delta,\Delta) + \operatorname{ILP}(m, 4n^2 \Delta, A \otimes \bbone_{4n\Delta}^\tp)
\end{align*}
respectively. Here $\bbone_n \in \mathbb{R}^n$ denotes the vector of all ones, $\otimes$ the Kronecker product, and $\beta > 0$ is some constant. A consequence is that SP is always polynomial-time solvable and ISP is polynomial-time solvable for totally unimodular $A$. The latter is a particularly delicate result, and even a slight deviation yields an NP-hard problem \cite{Baldick}. While the objective $f = f_1 + \dots + f_n$ is convex, we remind readers of Example~\ref{eg:self}:  it is a mistake to think that all convex optimization problems have polynomial-time algorithms.

Note that we could have stated the results above in terms of multiplicatively separable functions: if $g_1,\ldots,g_n \colon  \mathbb{R} \to \mathbb{R}_\pp$ are log convex, and $g \colon \mathbb{R}^n \to \mathbb{R}_\pp$ is defined by $g(x_1,\ldots,x_n) = g_1(x_1)\cdots g_n(x_n)$, then SP and ISP are equivalent to
\[
\begin{tabular}{lll}
&minimize & $\log g(x_1,\ldots,x_n)$\\
&subject to & $Ax \le b$, $0 \le x \le 1$,
\end{tabular}
\begin{tabular}{lll}
&minimize & $\log g(x_1,\ldots,x_n)$\\
&subject to & $Ax \le b$, $x \in \{0,1\}^n$,
\end{tabular}
\]
although it would be unnatural to discuss LP and ILP in these forms.
\end{example}

We will discuss another situation where additive separability arises.

\begin{example}[separable Hamiltonians]
The Hamilton equations
\begin{equation}\label{eq:Ham}
\dfrac{\dd p_i}{\dd t} = - \dfrac{\partial H}{\partial x_i} (x,p), \quad \dfrac{\dd x_i}{\dd t} = - \dfrac{\partial H}{\partial p_i} (x,p), \quad i=1,\ldots,n,
\end{equation}
are said to be separable if the Hamiltonian $H \colon  \mathbb{R}^n \times \mathbb{R}^n \to \mathbb{R}$ is additively separable:~
\begin{equation}\label{eq:sepHam}
H(x,p) = V(x) + T(p),
\end{equation}
\tie, the kinetic energy $T \colon  \mathbb{R}^n\to \mathbb{R}$ depends only on momentum $p=(p_1,\ldots,p_n)$ and the potential energy $V \colon  \mathbb{R}^n\to \mathbb{R}$ depends only on position $x=(x_1,\ldots,x_n)$, a common scenario. In this case the system \eqref{eq:Ham} famously admits a finite difference scheme that is both \emph{explicit}, \tie, iteration depends only on quantities already computed in a previous step, and \emph{symplectic}, \tie, it conforms to the tensor transformation rules with change-of-coordinates matrices from $\Sp(2n,\mathbb{R})$ on page~\pageref{eq:groups}. Taking $n =1$ for illustration, the equations in \eqref{eq:Ham} are
\[
\dfrac{\dd p}{\dd t} = - V'(x), \quad \dfrac{\dd x}{\dd t} = - T'(p),
\]
and with backward Euler on the first, forward Euler on the second, we obtain the finite difference scheme
\[
p^{(k+1)} = p^{(k)} +  V'(x^{(k+1)})\Delta t, \quad x^{(k+1)} = x^{(k)} - T'(x^{(k)})\Delta t,
\]
which is easily seen to be explicit as $x^{(k+1)}$ may be computed before $p^{(k+1)}$; to show that it is a symplectic integrator requires a bit more work and the reader may consult \citet[pp.~580--582]{Stuart}. Note that we could have written \eqref{eq:sepHam} in an entirely equivalent multiplicatively separable form with
$\rme^{H(x,p)} = \rme^{V(x)} \,\rme^{T(p)}$, although our subsequent discussion would be somewhat awkward.
\end{example}

We now present an example with separability structures in both problem and solution and where both additive and multiplicative separability play a role.

\begin{example}[one-electron\zzc and\zzc Hartree--Fock\zzc approximations]\label{eg:HF}
The time-dependent Schr\"odinger equation for a system of $d$ particles in $\mathbb{R}^3$ is
\[
i\hbar\dfrac{\partial}{\partial t} f(x,t) = \biggl[ - \dfrac{\hbar^2}{2m}\Delta + V(x)\biggr] f(x,t).
\]
Here $x = (x_1,\ldots,x_d) \in \mathbb{R}^{3n}$ represents the positions of the $d$ particles, $V$ a real-valued function representing potential, and
\[
\Delta = \Delta_1 + \Delta_2 + \dots + \Delta_d
\]
with each $\Delta_i \colon  L^2(\mathbb{R}^3) \to L^2(\mathbb{R}^3)$ a copy of the Laplacian on $\mathbb{R}^3$ corresponding to the $i$th particle. Note that these do not need to be in Cartesian coordinates. For instance, we might have
\[
\Delta_i =  \dfrac{1}{r_i^2} \dfrac{\partial}{\partial r_i} \biggl(r_i^2 \dfrac{\partial }{\partial r_i} \biggr) + \dfrac{1}{r_i^2 \sin \theta_i} \dfrac{\partial}{\partial \theta_i} \biggl(\sin \theta_i \dfrac{\partial }{\partial \theta_i} \biggr) + \dfrac{1}{r_i^2 \sin^2 \theta_i} \dfrac{\partial^2 }{\partial \phi_i^2}
\]
with $x_i = (r_i,\theta_i,\phi_i)$, and $i=1,\ldots,d$.

We will drop the constants, which are unimportant to our discussions, and just keep their signs:
\begin{equation}\label{eq:Schro}
(-\Delta + V) f -i \partial_t f = 0.
\end{equation}
Separation of variables \eqref{eq:septransf} applies to give~us
\[
(-\Delta + V) \otimes I + I \otimes (-i\partial_t) \quad \longrightarrow \quad 
\left\{
\begin{aligned}
(-\Delta + V )\varphi &= E \varphi,\\
-i\partial_t \psi &= -E \psi,
\end{aligned}\right.
\]
where we have written our separation constant as $-E$. The second equation is trivial to solve, $\psi(t) = \rme^{-\rmi E t}$, and as \eqref{eq:Schro} is linear, the solution $f$ is given by a linear combination of $\psi \otimes \varphi$ over all possible values of $E$ and our main task is to determine $\varphi$ and $E$ from the first equation, called the time-independent Schr\"odinger equation for $d$ particles. Everything up to this stage is similar to our discussion for the wave equation in Example~\ref{eg:sepPDE}. The difference is that we now have an extra $V$ term: $(E,\varphi)$ are an eigenpair of $-\Delta + V$.

The motivation behind the one-electron, Hartree--Fock and other approximations is as follows. If the potential $V$ is additively separable,
\begin{equation}\label{eq:seppot}
V(x) = V_1(x_1) + V_2(x_2) + \dots + V_d(x_d),
\end{equation}
then the eigenfunction $\varphi$ is multiplicatively separable,
\[
\varphi(x) = \varphi_1(x_1)\varphi_2(x_2)\cdots \varphi_d(x_d),
\]
and the eigenvalue $E$ is additively separable,
\[
E = E_1 + E_2 + \dots + E_d.
\]
The reason is that with \eqref{eq:seppot} we have
\[
\sum_{i=1}^d (-\Delta_i + V_i) \varphi  - E \varphi= 0,
\]
and we may apply \eqref{eq:septransf} to get
\begin{align}\label{eq:npart}
&[(-\Delta_1 + V_1) \otimes I \otimes \dots  \otimes I 
+ I \otimes (-\Delta_2 + V_2)  \otimes \dots \otimes I 
\notag \\*
&\quad  + I \otimes \dots \otimes I \otimes (-\Delta_d + V_d - E) ](\varphi_1 \otimes \varphi_2 \otimes \dots \otimes \varphi_d) = 0 \notag \\*
&\quad \longrightarrow \quad 
\left\{
\begin{aligned}
(-\Delta_1 + V_1 )\varphi_1 &= E_1 \varphi_1,\\
(-\Delta_2 + V_2 )\varphi_2 &= E_2 \varphi_2,\\
&\vdots\\
(-\Delta_d + V_d )\varphi_d &= (E-E_1-\dots-E_{d-1}) \varphi_d.
\end{aligned}\right. 
\end{align}
Note that we may put the $-E \varphi$ term with any $-\Delta_i + V_i$; we just chose $i =d$ for convenience. If we write $E_d \coloneqq E-E_1-\dots-E_{d-1}$ and $\varphi = \varphi_1 \otimes \dots \otimes \varphi_d$, we obtain the required expressions. So \eqref{eq:npart} transforms a $d$-particle Schr\"odinger equation into $d$ one-particle Schr\"odinger equations.

If one could solve the $d$-particle Schr\"odinger equation, then one could in principle determine all chemical and physical properties of atoms and molecules, so the above discussion seems way too simple, and indeed it is:  the additivity in \eqref{eq:seppot} can only happen if the particles do not interact. In this context an operator of the form in \eqref{eq:sepop} is called a non-interacting Hamiltonian, conceptually useful but unrealistic. It is, however, the starting point from which various approximation schemes are developed. To see the necessity of approximation, consider what appears to be a slight deviation from \eqref{eq:seppot}. Let $V$ include only pairwise interactions:
\[
V(x) = \sum_{i=1}^d V_i(x_i) + \sum_{i<j} V_{ij}(x_i, x_j),
\]
\tie, no higher-order terms of the form $V_{ijk}(x_i,x_j,x_k)$, and we may even fix $V_{ij}(x_i, x_j) = 1/\lVert x_i - x_j \rVert$. But the equation $(-\Delta + V )\varphi = E \varphi$ then becomes computationally intractable in multiple ways \cite{Whitfield}.

The one-electron approximation and Hartree--Fock approximation are  based on the belief that
\[
V(x) \approx V_1(x_1) +  \dots + V_d(x_d) \quad \Rightarrow \quad
\left\{
\begin{aligned}
\varphi(x) &\approx \varphi_1(x_1)\cdots \varphi_d(x_d),\\
E &\approx E_1 +\dots + E_d,
\end{aligned}\right.
\]
with `$\approx$' interpreted differently and with different tools:  the former uses perturbation theory and the latter calculus of variations. While these approximation methods are only tangentially related to the topic of this section, they do tell us how additive and multiplicative separability can be approximated and thus we will briefly describe the key ideas. Our discussion below is based on \citet{Fischer}, \citet[Chapters~33, 34, 50, 51]{Faddeev} and \citet[Chapter~12]{Hannabuss}.  

In most scenarios, $V(x)$ will have an additively separable component $V_1(x_1)+\dots + V_d(x_d)$ and an additively inseparable component comprising the higher-order interactions $V_{ij}(x_i,x_j),V_{ijk}(x_i,x_j,x_k),\dots.$ We will write
\[
H_0 \coloneqq \sum_{i=1}^d (-\Delta_i + V_i), \quad H_1 \coloneqq -\Delta + V = H_0 + W,
\]
with $W$ accounting for the inseparable terms. The one-electron approximation method uses perturbation theory to recursively add correction terms to the eigenpairs of $H_0$ to obtain those of $H_1$. Ignoring rigour, the method is easy to describe. Set $H_\varepsilon \coloneqq (1-\varepsilon) H_0 + \varepsilon W$ and plug the power series
\[
\lambda_\varepsilon = \lambda^{(0)} + \varepsilon \lambda^{(1)} + \varepsilon^2 \lambda^{(2)} + \cdots,\quad
\psi_\varepsilon = \psi^{(1)} + \varepsilon \psi^{(1)} + \varepsilon^2 \psi^{(2)} + \cdots
\]
into $H_\varepsilon \psi_\varepsilon = \lambda_\varepsilon \psi_\varepsilon$ to see that
\begin{align}
H_0 \psi^{(0)} &= \lambda^{(0)} \psi^{(0)}, \label{eq:perturb0}\\*
H_0 \psi^{(1)} + W \psi^{(0)} &= \lambda^{(0)} \psi^{(1)} + \lambda^{(1)} \psi^{(0)},\notag\\
&\ \ \vdots\notag \\*
H_0 \psi^{(k)} + W \psi^{(k-1)} &= \lambda^{(0)} \psi^{(k)} + \lambda^{(1)} \psi^{(k-1)} + \dots + \lambda^{(k)} \psi^{(0)}. \label{eq:perturbk}
\end{align}
By \eqref{eq:perturb0}, $\lambda^{(0)}$ and $\psi^{(0)}$ form an eigenpair of $H_0$. The subsequent terms $\lambda^{(k)}$ and $\psi^{(k)}$ are interpreted as the $k$th-order corrections to the eigenpair of $H_0$, and as $\lambda^{(k)}$ and $\psi^{(k)}$ may be recursively calculated from \eqref{eq:perturbk}, the hope is that
\[
\lambda^{(0)} + \lambda^{(1)} + \lambda^{(2)} + \cdots,\quad
\psi^{(0)} + \psi^{(1)} + \psi^{(2)} + \cdots
\]
would converge to the desired eigenpair of $H_1$. For instance, if $H_0$ has an orthonormal basis of eigenfunctions $\{\psi_i \colon  i \in \mathbb{N} \}$ with  eigenvalues $\{\lambda_i \colon  i \in \mathbb{N}\}$ and $\lambda_j$ is a simple eigenvalue, then \cite[Theorem~12.4.3]{Hannabuss}
\[
\lambda_j^{(1)} = \langle W \psi_j, \psi_j \rangle, \quad
\psi_j^{(1)} = \sum_{i \ne j} \dfrac{\langle W\psi_j, \psi_i \rangle}{\lambda_j - \lambda_i} \psi_i,\quad
\lambda_j^{(2)} = \sum_{i \ne j} \dfrac{\lvert\langle W\psi_j, \psi_i \rangle\rvert^2}{\lambda_j - \lambda_i} .
\]

Whereas one-electron approximation considers perturbation of additive separability in the problem (the Hamiltonian), Hartree--Fock approximation considers variation of multiplicative separability in the solution (the wave function). Observe that the non-linear functional given by the Rayleigh quotient $\mathcal{E}(\varphi) = \langle ( -\Delta + V) \varphi, \varphi \rangle/\lVert \varphi \rVert^2$ is stationary, \ie\ $\delta \mathcal{E} = 0$, if and only if  $(-\Delta  + V) \varphi = E\varphi$, \tie, its Euler--Lagrange equation gives the Schr\"odinger equation and its Lagrange multiplier gives the eigenvalue $E$. The Hartree--Fock approximation seeks stationarity under the additional condition of multiplicative separability, \ie\ $\varphi = \varphi_1 \otimes \dots  \otimes \varphi_d$. If the non-linear functional
\[
\mathcal{E}(\varphi_1, \ldots, \varphi_d) =
\langle ( -\Delta + V) \varphi_1\otimes \dots \otimes \varphi_d, \varphi_1 \otimes \dots \otimes \varphi_d \rangle
\]
with constraints $\lVert \varphi_1 \rVert^2 =\dots = \lVert \varphi_d \rVert^2 = 1$ is stationary with respect to variations in $\varphi_1,\ldots,\varphi_d$, \ie\ $\delta \mathcal{L} = 0$, where
\[
\mathcal{L}(\varphi_1, \ldots, \varphi_d,  \lambda_1,\ldots,\lambda_d) =
\mathcal{E}(\varphi_1, \ldots, \varphi_d)  - \lambda_1 \lVert \varphi_1 \rVert^2 - \dots - \lambda_d \lVert \varphi_d \rVert^2
\]
is the Lagrangian with Lagrange multipliers $\lambda_1,\ldots,\lambda_d$, then the Euler--Lagrange equations give us
\begin{equation}\label{eq:HF}
\biggl[ -\Delta_i + \sum_{j \ne i} \int_{\mathbb{R}^3} \lvert \varphi_j(y) \rvert^2 \,V(x,y) \D y \biggr] \varphi_i = \lambda_i \varphi_i, \quad i =1,\ldots,d.
\end{equation}
These equations make physical sense. For a fixed $i$, \eqref{eq:HF} is the Schr\"odinger equation for particle $i$ in a potential field due to the charge of  particle $j$; this charge is spread over space with density $\lvert \varphi_j \rvert^2$, and we have to sum over the potential fields created by all particles $j=1,\ldots, i-1, i+1,\ldots,d$.
While the system of coupled integro-differential equations \eqref{eq:HF} generally would not have an analytic solution like the one in Example~\ref{eg:ide}, it readily lends itself to numerical solution via a combination of quadrature and finite difference \cite[Chapter~6, Section~7]{Fischer}.
\end{example}

We have ignored the spin variables in the last and the next examples, as spin has already been addressed in Examples~\ref{eg:spin} and would only be an unnecessary distraction here. So, strictly speaking, these examples are about Hartree approximation \cite{Hartree}, \ie\ no spin, and not Hartree--Fock approximation \cite{Fock}, \ie\ with spin.

\begin{example}[multiconfiguration Hartree--Fock]\label{eg:MCHF}
We will continue our discussion in the last example but with some simplifications. We let $d =2$ and $V(x,y) = V(y,x)$. We emulate the Hartree--Fock approximation for the time-independent Schr\"odinger equation $(-\Delta + V) f = E f$ but now with an ansatz of the form
\begin{equation}\label{eq:twoterms}
f = a_1 \varphi_1 \otimes \varphi_1  + a_2 \varphi_2 \otimes \varphi_2,
\end{equation}
where $\varphi_1,\varphi_2 \in L^2(\mathbb{R}^3)$ are orthonormal and $a=(a_1,a_2) \in \mathbb{R}^2$ is a unit vector, \tie, we consider the non-linear functional
\begin{align*}
\mathcal{E}(\varphi_1, \varphi_2,a_1,a_2) &=
a_1^2\langle ( -\Delta + V) \varphi_1\otimes \varphi_1, \varphi_1 \otimes \varphi_1 \rangle \\*
&\quad+2a_1 a_2\langle ( -\Delta + V) \varphi_1\otimes \varphi_1, \varphi_2 \otimes \varphi_2 \rangle \\*
&\quad +a_2^2 \langle ( -\Delta + V) \varphi_2\otimes \varphi_2, \varphi_2 \otimes \varphi_2 \rangle
\end{align*}
with constraints $\lVert \varphi_1 \rVert^2 =\lVert \varphi_2 \rVert^2 =1$, $\langle \varphi_1, \varphi_2 \rangle = 0$, $\lVert a \rVert^2 =1$. Introducing Lagrange multipliers $\lambda_{11}, \lambda_{12}, \lambda_{22}$ and $\lambda$  for the constraints, we obtain the Lagrangian
\begin{align*}
&  \mathcal{L}(\varphi_1, \varphi_2,a_1, a_2, \lambda_{11},\lambda_{12},\lambda_{22},\lambda) \\*
  &\quad = \mathcal{E}(\varphi_1, \varphi_2,a_1,a_2)  
+ \lambda_{11}\lVert \varphi_1 \rVert^2 + \lambda_{12} \langle \varphi_1, \varphi_2 \rangle + \lambda_{22}\lVert \varphi_2 \rVert^2  - \lambda \lVert a \rVert^2.
\end{align*}
For $i,j \in \{1,2\}$ we write
\[
b_{ij} =\langle ( -\Delta + V) \varphi_i\otimes \varphi_i, \varphi_j \otimes \varphi_j \rangle , \quad
c_{ij}(x) = \int_{\mathbb{R}^3}  \varphi_i(x) \varphi_j(y) V(x,y)  \D y.
\]
The stationarity conditions $\nabla_a \mathcal{L} = 0$ and $\delta \mathcal{L}$ = 0 give us
\[
\begin{bmatrix}
b_{11} -\lambda & b_{12} \\
b_{12} & b_{22} -\lambda
\end{bmatrix}
\begin{bmatrix}
a_1\\
a_2
\end{bmatrix}
=0
\]
and
\[
\begin{bmatrix}
-\Delta_1 \varphi_1 \\
-\Delta_2 \varphi_2
\end{bmatrix}
=
\begin{bmatrix}
c_{11} (x) -\lambda_{11} & (a_2/a_1)( c_{12} (x) - \lambda_{12}) \\
(a_1/a_2) ( c_{12} (x) - \lambda_{12}) & c_{22} (x) - \lambda_{22}
\end{bmatrix}
\begin{bmatrix}
\varphi_1 \\
\varphi_2
\end{bmatrix}
\]
respectively. What we have described here is a simplified version of the multiconfiguration Hartree--Fock approximation \cite[Chapters~3 and 4]{Fischer}, and it may also be viewed as the starting point for a variety of other heuristics, one of which will be our next example.
\end{example}

We will describe an extension of the ansatz in \eqref{eq:twoterms} to more summands and to higher order.

\begin{example}[tensor networks]\label{eg:tennet}
  Instead of restricting ourselves to $L^2(\mathbb{R}^3)$ as in the last two examples, we will now assume arbitrary separable Hilbert spaces $\mathbb{H}_1,\ldots,\mathbb{H}_d$ to allow for spin and other properties. We will seek a solution $f \in \mathbb{H}_1 \otimes \dots \otimes \mathbb{H}_d$ to the Schr\"odinger equation for $d$ particles.
  Readers are reminded of the discussion at the end of Example~\ref{eg:mrank}. Note that by definition of $\otimes$, $f$ has finite rank regardless of whether $\mathbb{H}_1,\dots,\mathbb{H}_d$ are finite- or infinite-dimensional.
  Let $\mrank(f) =(r_1,\ldots,r_d)$. Then there is a decomposition
\begin{equation}\label{eq:onansatz}
f = \sum_{i=1}^{r_1} \sum_{j=1}^{r_2} \cdots \sum_{k=1}^{r_d} c_{ij\cdots k} \, \varphi_i \otimes \psi_j  \otimes \dots \otimes \theta_k
\end{equation}
with orthogonal factors \cite[equation~13]{DeLathauwer}
\[
\langle \varphi_i, \varphi_j \rangle = \langle \psi_i, \psi_j \rangle = \dots = \langle \theta_i, \theta_j \rangle = \begin{cases} 0 & i \ne j, \\ 1 & i = j. \end{cases}
\]
Note that when our spaces have inner products, any multilinear rank decomposition \eqref{eq:multdecomp} may have its factors orthogonalized
\cite[Theorem~2]{DeLathauwer},
\tie, the orthogonality constraints do not limit the range of possibilities in \eqref{eq:onansatz}.
By definition of multilinear rank, there exist subspaces $\mathbb{H}_1', \ldots, \mathbb{H}_d'$ of dimensions $r_1,\ldots,r_d$ with $f \in \mathbb{H}_1' \otimes  \dots \otimes \mathbb{H}_d'$ that attain the minimum in \eqref{eq:mrank}. As such, we may replace $\mathbb{H}_i$ with $\mathbb{H}_i'$ at the outset, and to simplify our discussion we may as well assume that $\mathbb{H}_1, \dots , \mathbb{H}_d$ are of dimensions $r_1,\ldots, r_d$.

The issue with \eqref{eq:onansatz} is that there is an exponential number of rank-one terms as $d$ increases. Suppose $r_1=\dots =r_d=r$; then there are $r^d$ summands in \eqref{eq:onansatz}. This is not unexpected because  \eqref{eq:onansatz} is the most general form a finite-rank tensor can take. Here the ansatz in \eqref{eq:onansatz} describes the whole space $\mathbb{H}_1 \otimes \dots \otimes \mathbb{H}_d$ and does not quite serve its purpose;  an ansatz is supposed to be an educated guess, typically based on physical insights, that captures a small region of the space where the solution likely lies. The goal of \emph{tensor networks} is to provide such an ansatz by limiting the coefficients $[c_{ij\cdots k}] \in \mathbb{R}^{r_1 \times \dots \times r_d}$ to a much smaller set. The first and best-known example is the \emph{matrix product states} tensor network \cite{Anderson,White,WH}, which imposes on the coefficients the structure
\[
c_{ij\cdots k} = \tr(A_i B_j \cdots C_k), \quad A_i \in \mathbb{R}^{n_1 \times n_2}, B_j \in \mathbb{R}^{n_2 \times n_3},\ldots, C_k \in \mathbb{R}^{n_d \times n_1}
\]
for $i=1,\ldots,r_1$, $j =1,\ldots,r_2, \ldots, k =1,\ldots, r_d$. An ansatz of the form
\begin{equation}\label{eq:mps}
f = \sum_{i=1}^{r_1} \sum_{j=1}^{r_2} \cdots \sum_{k=1}^{r_d} \tr(A_i B_j \cdots C_k) \, \varphi_i \otimes \psi_j  \otimes \dots \otimes \theta_k
\end{equation}
is called a matrix product state or MPS \cite{ITEH}. Note that the coefficients are now parametrized by $r_1 + r_2 +\dots + r_d$ matrices of various sizes. For easy comparison, if $r_1=\dots =r_d=r$ and $n_1=\dots =n_d=n$, then the coefficients in \eqref{eq:onansatz} have $r^d$ degrees of freedom whereas those in \eqref{eq:mps} only have $rdn^2$. When $n_1=1$, the first and last matrices in \eqref{eq:mps} are a row and a column vector respectively; as the trace of a $1\times 1$ matrix is itself, we may drop the `$\tr$' in \eqref{eq:mps}.  This special case with $n_1 = 1$ is sometimes called MPS with open boundary conditions \cite{Anderson} and the more general case is called MPS with periodic conditions.

The above discussion of matrix product states conceals an important structure. Take $d= 3$ and denote the entries of the matrices as
$A_i = [a^{(i)}_{\alpha\beta} ]$,
$B_j = [b^{(j)}_{\beta\gamma} ]$
and
$C_k = [c^{(k)}_{\gamma\alpha} ]$.
Then
\begin{align*}
f &= \sum_{i,j,k=1}^{r_1,r_2,r_3} \tr(A_i B_j  C_k) \, \varphi_i \otimes \psi_j  \otimes \theta_k \\*
&=\sum_{i,j,k=1}^{r_1,r_2,r_3} \biggl[ \sum_{\alpha,\beta,\gamma=1}^{n_1, n_2, n_3} a^{(i)}_{\alpha\beta} b^{(j)}_{\beta\gamma} c^{(k)}_{\gamma\alpha} \, \varphi_i \otimes \psi_j  \otimes \theta_k\biggr] \\
&= \sum_{\alpha,\beta,\gamma=1}^{n_1, n_2, n_3}
\biggl[\sum_{i=1}^{r_1}  a^{(i)}_{\alpha\beta} \, \varphi_i \biggr]\otimes \biggl[ \sum_{j=1}^{r_2} b^{(j)}_{\beta\gamma} \, \psi_j \biggr] \otimes \biggl[  \sum_{k=1}^{r_3} c^{(k)}_{\gamma\alpha} \, \theta_k\biggr]\\*
&=  \sum_{\alpha,\beta,\gamma=1}^{n_1, n_2, n_3} \varphi_{\alpha\beta} \otimes \psi_{\beta\gamma} \otimes \theta_{\gamma \alpha},
\end{align*}
where
\[
\varphi_{\alpha\beta} \coloneqq \sum_{i=1}^{r_1}  a^{(i)}_{\alpha\beta} \, \varphi_i,\quad
\psi_{\beta\gamma} \coloneqq \sum_{j=1}^{r_2} b^{(j)}_{\beta\gamma} \, \psi_j \quad\text{and}\quad
\theta_{\gamma \alpha} \coloneqq \sum_{k=1}^{r_3} c^{(k)}_{\gamma\alpha} \, \theta_k.
\]
In other words, the indices have the incidence structure of an undirected graph, in this case a triangle. This was first observed in \citet{LQY} and later generalized in \citet{tnr}, the bottom line being that any tensor network state is  a sum of separable functions indexed by a graph.  In the following, we show some of the most common tensor network states, written in this simplified form, together with the graphs they correspond to in Figure~\ref{fig:networks}.
\medskip

\noindent{Periodic matrix product states:}
\[
f(x,y,z) = \sum_{i,j,k=1}^{n_1,n_2,n_3} \varphi_{ij}(x) \psi_{jk}(y) \theta_{ki}(z).
\]

\noindent{Tree tensor network states:}
\[
f(x,y,z,w) = \sum_{i,j,k=1}^{n_1,n_2,n_3} \varphi_{ijk}(x) \psi_{i}(y) \theta_{j}(z) \pi_{k}(w).
\]
\noindent{Open matrix product states:}
\[
f(x,y,z,u,v) = \sum_{i,j,k,l=1}^{n_1,n_2,n_3,n_4}\varphi_{i}(x) \psi_{ij}(y) \theta_{jk}(z)\pi_{kl}(u)\rho_{l}(v).
\]
\noindent{Projected entangled pair states:}
\[
f(x,y,z,u,v,w) =\sum_{i,j,k,l,m,n,o=1}^{n_1,n_2,n_3,n_4,n_5,n_6,n_7} \varphi_{ij}(x) \psi_{jkl}(y) \theta_{lm}(z) \pi_{mn}(u) \rho_{nko}(v) \sigma_{oi}(w).
\]
The second and the last, often abbreviated to TTNS and PEPS, were proposed by \citet{SDV} and \citet{VC} respectively.  More generally, any periodic MPS corresponds to a cycle graph and any open MPS corresponds to a line graph. We have deliberately written them without the $\otimes$ symbol to emphasize that all these {ans\"atze} 
are just sums of separable functions, differing only in terms of how their factors are indexed.
\begin{figure}[ht]
\centering
\begin{tikzpicture}[scale=0.8]
\filldraw
     (4,1) node[align=center, above] {\textsc{mps} (open)}
     (0,0) circle (2pt) node[align=center, below] {$\varphi$}
 -- (1,0) circle (0pt) node[align=center, above] {$i$}
 -- (2,0) circle (2pt)  node[align=center, below] {$\psi$}  
  -- (3,0) circle (0pt) node[align=center, above] {$j$}
 -- (4,0) circle (2pt)  node[align=center, below] {$\theta$}
  -- (5,0) circle (0pt) node[align=center, above] {$k$}
   -- (6,0) circle (2pt) node[align=center, below] {$\pi$}
    -- (7,0) circle (0pt) node[align=center, above] {$l$}
 -- (8,0) circle (2pt) node[align=center, below] {$\rho$} 
-- cycle;
\end{tikzpicture}
\begin{tikzpicture}[scale=0.8]
\filldraw
     (0.75,1.5) node[align=center, below] {\textsc{ttns}}
     (0,0) circle (2pt) node[align=center, below] {$\psi$}
      -- (1,0) circle (0pt) node[align=center, above] {$i$}
 -- (2,0) circle (2pt)  node[align=center, below] {$\varphi$}  
 -- (4,0) circle (2pt)  node[align=center, below] {$\theta$}  
  -- (3,0) circle (0pt) node[align=center, above] {$j$}
 -- (2,0) circle (2pt) node[align=center, below]{} %{$\varphi$}
  -- (2,1) circle (0pt) node[align=center, right] {$k$}
 -- (2,2) circle (2pt) node[align=center, above] {$\pi$} 
-- cycle;
\end{tikzpicture}
\begin{tikzpicture}[scale=0.8]
\filldraw
    (-1,0.75) node[align=center, left] {\textsc{mps} (periodic)}
    (0,0) circle (2pt) node[align=center, below] {$\varphi$}
     -- (0.7,0.7) circle (0pt) node[align=center, above] {$j$}
 -- (2,2) circle (2pt)  node[align=center, below] {$\psi$}  
  -- (0,2) circle (0pt) node[align=center, below] {$k$}
 -- (-2,2) circle (2pt)  node[align=center, below] {$\theta$}  
  -- (-0.7,0.7) circle (0pt) node[align=center, above] {$i$}
 -- (0,0) circle (2pt) node[align=center, below]{}% {$\varphi$}
-- cycle;
\end{tikzpicture}
\begin{tikzpicture}[scale=0.8]
\filldraw
    (2,2.5) node[align=center, above] {\textsc{peps}}
    (0,0) circle (2pt) node[align=center, below] {$\varphi$}
     -- (1,0) circle (0pt) node[align=center, above] {$j$}
 --(2,0) circle (2pt)  node[align=center, below] {$\psi$}  
  -- (3,0) circle (0pt) node[align=center, above] {$l$}
 --(4,0) circle (2pt) node[align=center, below] {$\theta$}
  -- (4,1) circle (0pt) node[align=center, left] {$m$}
 --(4,2) circle (2pt)  node[align=center, above] {$\pi$}
   -- (3,2) circle (0pt) node[align=center, below] {$n$}
 -- (2,2) circle (2pt)  node[align=center, above] {$\rho$}
    -- (2,1) circle (0pt) node[align=center, right] {$k$}
  --(2,0) circle (2pt) node[align=center, below]{}% {$\psi$}  
   -- (2,2) circle (2pt) node[align=center, above]{}% {$v$}
       -- (1,2) circle (0pt) node[align=center, below] {$o$}
 -- (0,2) circle (2pt)  node[align=center, above] {$\sigma$}
     -- (0,1) circle (0pt) node[align=center, right] {$i$}
--   (0,0) circle (2pt) node[align=center, below]{}% {$\varphi$}
 -- cycle;
\end{tikzpicture}
\parbox{270pt}{\caption{Graphs associated with common tensor networks.}}
\label{fig:networks}
\end{figure}
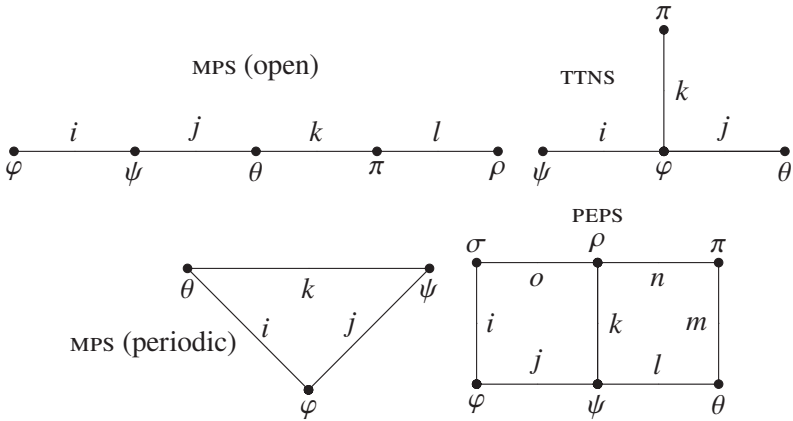

One difference between the Hartree--Fock approximations in Examples~\ref{eg:HF} and \ref{eg:MCHF} and the tensor network methods here is that the factors $\varphi_i, \psi_j, \ldots, \theta_k$ in \eqref{eq:mps} are often fixed in advance as some standard bases of $\mathbb{H}_1,\ldots,\mathbb{H}_d$, called a \emph{local basis} in this context. The computational effort is solely in determining the coefficients $c_{ij\cdots k}$;  in the case of MPS this can be done via several singular value decompositions \cite{Orus}. Indeed, because of the way it is computed, the coefficients of the MPS ansatz are sometimes represented as
\[
\tr(Q_1 \Sigma_1 Q_2 \Sigma_2 \dots \Sigma_d Q_{d+1}), \quad Q_i \in \Or(n_i), \; \Sigma_i \in \mathbb{R}^{n_i \times n_{i+1}}
\]
for $i =1,\dots d$, with $\Sigma_i$ a diagonal matrix. That any $\tr(A_1 A_2 \cdots A_d)$ with $A_i \in \mathbb{R}^{n_i \times n_{i+1}}$ may be written in this form follows from the singular value decompositions $A_i = U_i \Sigma_i V_i^\tp$ and setting $Q_{i+1} = V_i^\tp U_{i+1}$ for $i=1,\ldots,d-1$, and $Q_1 = U_1$, $Q_{d+1} = V_d^\tp$.

We end with a brief word about
our initial assumption that $f \in \mathbb{H}_1 \otimes \cdots \otimes \mathbb{H}_d$.
For a solution to a PDE, we expect that $f \in \mathbb{H}_1 \hatotimes_\F \cdots \hatotimes_\F \mathbb{H}_d$, where $\hatotimes_\F$ is the Hilbert--Schmidt tensor product discussed in Examples~\ref{eg:Hilbert} and \ref{eg:compact}.
Since any tensor product basis  \eqref{eq:tenprodbasis2} formed from orthonormal bases of  $\mathbb{H}_1,\ldots,\mathbb{H}_d$ is  an orthonormal basis of
$\mathbb{H}_1 \hatotimes \dots \hatotimes \mathbb{H}_d$ with respect to the inner product in \eqref{eq:innprod1}, it is always possible to express
$f \in \mathbb{H}_1 \hatotimes \dots \hatotimes \mathbb{H}_d$
in the form \eqref{eq:onansatz} as long as we allow $r_1,\ldots, r_d$ to be infinite. The ones with finite $r_1,\ldots, r_d$ are precisely the finite-rank tensors in $\mathbb{H}_1 \otimes \cdots \otimes \mathbb{H}_d$, but any other
$f \in \mathbb{H}_1 \hatotimes_\F \cdots \hatotimes_\F \mathbb{H}_d$
can be approximated to arbitrary accuracy by
elements of $\mathbb{H}_1 \otimes \dots \otimes \mathbb{H}_d$, \tie, $f$
is a limit of finite-rank tensors with respect to the Hilbert--Schmidt norm.
\end{example}

We will now embark on a different line of applications of definition~\ref{st:tensor3} tied to our discussions of integral kernels in Examples~\ref{eg:Calderon}, \ref{eg:distributions} and \ref{eg:Mercer}. As we saw in Example~\ref{eg:Calderon}, some of the most important kernels take the form $K(v,v') = f(\lVert v - v' \rVert)$ for some real-valued function $f$ on a norm space $\mathbb{V}$. In this case we usually write $\kappa(v-v') = f(\lVert v - v' \rVert)$ and call $\kappa$ a radial basis kernel \cite{Buhmann}. These are usually not separable like those in Examples~\ref{eg:sepIE} or \ref{eg:septransf}. The fast multipole method \cite{FMM} may be viewed as an algorithm for approximating the convolution of certain radial basis kernels with certain functions. We will discuss how tensors arise in this algorithm.

\begin{example}[multipole tensors]\label{eg:multpole}
Let $\mathbb{V}$ be an inner product space. As usual we will denote the inner product and its induced norm by $\langle \, \cdot\,, \cdot\,  \rangle$ and $\lVert \, \cdot \, \rVert$ respectively. As we saw in Examples~\ref{eg:norms} and \ref{eg:tenalg}, this induces an inner product and norm on $\mathbb{V}^{\otimes d}$ that we may denote with the same notation. For any $s \ge 0$, the function
\[
\kappa_s \colon  \mathbb{V}\setminus \{0\} \to \mathbb{R}, \quad \kappa_s(v) =
\begin{cases}
\log \, \lVert v \rVert  & s = 0,\\
1/\lVert v \rVert^s & s > 0,
\end{cases}
\]
is called a Newtonian kernel. Henceforth we will drop the subscript $s$ from $\kappa_s$ since this is normally fixed throughout. If $\mathbb{V}$ is equipped with some measure, then the convolution
\[
\kappa \ast f(v) = \int_{v \in \mathbb{V}} \kappa(v - w )f(w)  \D w
\]
is called the Newtonian potential\footnote{Strictly speaking, the names are used for the case when $s = \dim \mathbb{V}-2$, and for more general $s$ they are called the Riesz kernel and Riesz potentials. Because of the singularity in $\kappa$, the integral is interpreted in a principal value sense as in \eqref{eq:biHT} but we ignore such details \cite[Chapter~1, Section~8.18]{Stein}. We have also dropped some multiplicative constants.}  of $f$.  We will assume $s > 0$ in the following; the $s=0$ case is similar. Our interpretation of higher-order derivatives as multilinear maps in Example~\ref{eg:hod} makes their calculation an effortless undertaking:
\begin{align*}
[D\kappa(v)](h) &= -\dfrac{s }{\lVert v \rVert^{s + 2}}\langle v, h \rangle, \\*
[D^2\kappa(v)](h_1, h_2) &= -\dfrac{s }{\lVert v \rVert^{s + 2}}\langle h_2, h_1 \rangle + 
 \dfrac{s(s+2)}{\lVert v \rVert^{s + 4}}\langle v, h_1 \rangle\langle v, h_2 \rangle, \\
       [D^3\kappa(v)](h_1, h_2, h_3) &\\*
&\hspace{-36pt} = \dfrac{s(s+2)}{\lVert v \rVert^{s + 4}}
[
\langle v, h_1 \rangle\langle h_2, h_3 \rangle 
+\langle v, h_2 \rangle\langle h_1, h_3 \rangle 
+\langle v, h_3 \rangle\langle h_1, h_2 \rangle
]& \\*
&\hspace{-24pt} -\dfrac{s(s+2)(s+4)}{\lVert v \rVert^{s + 6}}\langle v, h_1 \rangle\langle v, h_2 \rangle \langle v, h_3\rangle,\\
[D^4\kappa(v)](h_1, h_2, h_3,h_4) &\\*
&\hspace{-36pt} = \dfrac{s(s+2)}{\lVert v \rVert^{s + 4}}
[
\langle h_1, h_2 \rangle\langle h_3, h_4 \rangle
+\langle h_1, h_3 \rangle \langle h_2, h_4 \rangle
+\langle h_1, h_4 \rangle\langle h_2, h_3 \rangle
]\\*
&\hspace{-24pt}  -\dfrac{s(s+2)(s+4)}{\lVert v \rVert^{s + 6}}[
\langle v, h_1 \rangle \langle v, h_2 \rangle \langle h_3, h_4\rangle
+\langle v, h_1 \rangle\langle v, h_3 \rangle \langle h_2, h_4\rangle 
\\
&\hspace{-24pt}
+\langle v, h_1 \rangle \langle v, h_4 \rangle \langle h_2, h_3\rangle 
+\langle v, h_2 \rangle \langle v, h_3 \rangle \langle h_1, h_4\rangle
\\*
&\hspace{-24pt}
+\langle v, h_2 \rangle \langle v, h_4 \rangle \langle h_1, h_3\rangle
+\langle v, h_3 \rangle \langle v, h_4 \rangle \langle h_1, h_2\rangle 
]\\*
&\hspace{-24pt}
+\dfrac{s(s+2)(s+4)(s+6)}{\lVert v \rVert^{s + 8}}\langle v, h_1 \rangle
\langle v, h_2 \rangle \langle v, h_3\rangle\langle v, h_4\rangle.
\end{align*}
To obtain these, all we need is the expression for $[D\kappa(v)](h)$, which follows from binomial expanding to linear terms in $h$, \tie,
\[
\kappa(v + h) = \dfrac{1}{\lVert v \rVert^s}\biggl[ 1 + \dfrac{2\langle v, h\rangle}{\lVert v \rVert^2} + \dfrac{\lVert h\rVert^2}{\lVert v \rVert^2} \biggr]^{-s/2} = \kappa(v)  -\dfrac{s }{\lVert v \rVert^{s + 2}}\langle v, h \rangle + O(\lVert h \rVert^2),
\]
along with the observation that $[D\langle v, h'\rangle ](h) = \langle h, h' \rangle$ and the product rule $D f \cdot g(v) =  f(v) \cdot Dg(v) +g(v) \cdot Df(v) $. Applying these repeatedly gives us $D^k\kappa(v)$ in any coordinate system without having to calculate a single partial derivative. 

As we saw in Example~\ref{eg:linhod}, these derivatives may be linearized via the universal factorization property. With an inner product, this is particularly simple. We will take $s =1$, as this gives us the Coulomb potential, the most common case.  We may rewrite the above expressions as
\begin{align*}
[D\kappa(v)](h) &= \biggl\langle -\dfrac{v}{\lVert v \rVert^3}, h \biggr\rangle, \\*
[D^2\kappa(v)](h_1,h_2) &= \biggl\langle -\dfrac{I}{\lVert v \rVert^3} + \dfrac{3\,v \otimes v }{\lVert v \rVert^5}, h_1 \otimes h_2 \biggr\rangle,\\*
[D^3\kappa(v)](h_1, h_2, h_3) &= \biggl\langle \dfrac{3}{\lVert v \rVert^5} \biggl(v \otimes I + \stackrel[{\displaystyle v}]{{\displaystyle I}}{\otimes} + I \otimes v \biggr) - \dfrac{15\,  v \otimes v \otimes v}{\lVert v \rVert^7} , h_1 \otimes h_2 \otimes h_3 \biggr\rangle.
\end{align*}
The symbol where $I$ and $v$ appear above and below $\otimes$ is intended to mean the following. Take any orthonormal basis $e_1,\ldots,e_n$ of $\mathbb{V}$; then $I =\sum_{i=1}^n e_i \otimes e_i$, and we have
\[
v \otimes I = \sum_{i=1}^n v \otimes e_i \otimes e_i, \quad  {\stackrel[{\displaystyle v}]{{\displaystyle I}}{\otimes} } = \sum_{i=1}^n  e_i \otimes v \otimes e_i, \quad I \otimes v = \sum_{i=1}^n  e_i \otimes e_i \otimes v.
\]
The tensors appearing in the first argument of the inner products are precisely the linearized derivatives and they carry important physical meanings:
\begin{alignat*}{3}
&\text{monopole} & \kappa(v) &= \dfrac{1}{\lVert v \rVert},\\*
&\text{dipole} & \partial\kappa(v) &= \dfrac{-1}{\lVert v \rVert^2}\, \hat{v},\\
&\text{quadrupole}\quad & \partial^2\kappa(v) &= \dfrac{1}{\lVert v \rVert^3}(3\,\hat{v} \otimes \hat{v} - I),\\*
&\text{octupole} & \partial^3\kappa(v) &= \dfrac{-1}{\lVert v \rVert^4} \biggl[ 15\, \hat{v} \otimes \hat{v} \otimes \hat{v} - 3\biggl(\hat{v} \otimes I + \stackrel[{\displaystyle \hat{v}}]{{\displaystyle I}}{\otimes} + I \otimes \hat{v} \biggr) \biggr],
\end{alignat*}
where $\hat{v} \coloneqq v/\lVert v \rVert$. The norm-scaled $k$th linearized derivative
\[
M^k \kappa(v) \coloneqq \lVert v \rVert^{k+1} \, \partial^k \kappa(v) \in \mathbb{V}^{\otimes k}
\]
is called a \emph{multipole tensor} or $k$-pole tensor.  The Taylor expansion \eqref{eq:taylor0}  may then be written in terms of the multipole tensors,
\[
\kappa(v)= \kappa( v_0) \sum_{k=0}^d\dfrac{1}{k!} \biggl\langle M^k \kappa(v_0), \dfrac{(v-v_0)^{\otimes k}}{\lVert v_0 \rVert^k} \biggr\rangle+ R(v-v_0),
\]
and this is called a \emph{multipole expansion} of $\kappa(v) =1/\lVert v \rVert$. We may rewrite this in terms of a tensor Taylor series and a tensor geometric series as in \eqref{eq:taylor},
\[
\partial \kappa (v_0)  = \sum_{k=0}^\infty \dfrac{1}{k!}  M^k \kappa(v_0)\in \widehat{\Ten}(\mathbb{V}),
\quad S\biggl(\dfrac{v-v_0}{\lVert v_0 \rVert}\biggr)=  \sum_{k=0}^\infty \dfrac{(v-v_0)^{\otimes k}}{\lVert v_0 \rVert^k}\in \widehat{\Ten}(\mathbb{V}).
\]
Note that
\[
\partial \kappa (v) = 
1 - \dfrac{1}{2!} \hat{v} + \dfrac{1}{3!}(3\,\hat{v} \otimes \hat{v} - I) - \dfrac{1}{4!} \biggl[ 15\, \hat{v} \otimes \hat{v} \otimes \hat{v} - 3\biggl(\hat{v} \otimes I + \stackrel[{\displaystyle \hat{v}}]{{\displaystyle I}}{\otimes} + I \otimes \hat{v} \biggr) \biggr] + \cdots
\]
is an expansion in the tensor algebra $\widehat{\Ten}(\mathbb{V})$.

Let $f \colon  \mathbb{V} \to \mathbb{R}$ be a compactly supported function that usually represents a charge distribution confined to a region within $\mathbb{V}$. The tensor-valued integral
\[
(M^k \kappa) \ast f (v)= \int_{v \in \mathbb{V}} M^k \kappa(v - w )f(w)  \D w
\]
is called a \emph{multipole moment} or $k$-pole moment, and the Taylor expansion
\[
\kappa \ast f (v)= \kappa (v_0) \sum_{k=0}^d\dfrac{1}{k!} \biggl\langle (M^k \kappa) \ast f (v_0), \dfrac{(v-v_0)^{\otimes k}}{\lVert  v_0 \rVert^k} \biggr\rangle+ R(v-v_0),
\]
is called a \emph{multipole expansion} of $f$. One of the most common scenarios is when we have a finite collection of $n$ point charges at positions $v_1,\ldots,v_n \in \mathbb{V}$ with charges $q_1,\ldots,q_n \in \mathbb{R}$. Then $f$ is given by
\begin{equation}\label{eq:chargedist}
f(v) = q_1 \delta (v - v_1) + \dots + q_n \delta(v - v_n).
\end{equation}
Our goal is to estimate the potential at a point $v$ some distance away from $v_1,\ldots,v_n$. Now consider a distinct point $v_0$ that is much nearer to $v_1,\ldots,v_n$ than to $v$ in the sense that $\lVert v_i - v_0 \rVert \le c \lVert v - v_0 \rVert$, $i =1,\ldots,n$, for some $c < 1$. Suppose there is just a single point charge with $f(v)  = q_i \delta (v - v_i)$; then $(M^k \kappa) \ast f (v-v_0) = q_i M^k \kappa (v_i-v_0)$ and the multipole expansion of $f$ at the point $v - v_i = (v - v_0) + (v_0 - v_i)$ about the point $v - v_0$ is
\begin{align*}
\kappa \ast f (v - v_i) &= \kappa  (v- v_0) \sum_{k=0}^d\dfrac{q_i}{k!} \biggl\langle M^k \kappa (v_i-v_0), \dfrac{(v_0-v_i)^{\otimes k}}{\lVert v - v_0 \rVert^k} \biggr\rangle +R(v_0 - v_i)  \\*
&= \sum_{k=0}^d  \varphi_k(v_i-v_0)\psi_k(v - v_0) + O(c^{d+1}),
\end{align*}
where we have written
\[
\varphi_k(v) \coloneqq \dfrac{(-1)^k q_i}{k!} \langle M^k \kappa(v) , v^{\otimes k} \rangle, \quad \psi_k(v) \coloneqq \dfrac{\kappa(v)}{\lVert v \rVert^k},
\]
and since $R(v_0 - v_i)/\lVert v_0 - v_i \rVert^d \to 0$, we assume $R(v_0 - v_i) = O(c^{d+1})$ for simplicity. For the general case in \eqref{eq:chargedist} with $n$ point charges, the potential at $v$ is then approximated by summing over $i$:
\begin{equation}\label{eq:approxsepker}
\kappa \ast f(v) \approx \sum_{k=0}^d  \biggl[\sum_{i=1}^n \varphi_k(v_i-v_0) \biggr] \psi_k(v - v_0).
\end{equation}
This sum can be computed in $O(nd)$ complexity or $O(n \log (1/\varepsilon) )$ for an $\varepsilon$-accurate algorithm. While the high level idea in \eqref{eq:approxsepker} is still one of approximating a function by a sum of $d$ separable functions, the fast multipole method involves a host of other clever ideas, not least among which are the techniques for performing such sums in more general situations \cite{Demaine}, for subdividing the region containing $v_1,\ldots,v_n$ into cubic cells (when  $\mathbb{V} = \mathbb{R}^3$) and thereby organizing these computations into a tree-like multilevel algorithm \cite{Barnes}. Clearly the approximation is good only when $v$ is far from $v_1,\ldots,v_n$ relative to $v_0$ but the algorithm allows one to circumvent this requirement.  We refer to \citet{FMM} for further information.
\end{example}

Multipole moments, multipole tensors and multipole expansions are usually discussed in terms of coordinates \cite[Chapter~4]{Jackson}. The coordinate-free approach, the multipole expansion as an element of the tensor algebra, {\em etc.}, are results of our working with definitions~\ref{st:tensor2} and \ref{st:tensor3}.

Although non-separable kernels make for more interesting examples, separable kernels arise as some of the most common multidimensional integral transforms, as we saw in Example~\ref{eg:septransf}. Also, they warrant a mention if only to illustrate why separability in a kernel is computationally desirable.

\begin{example}[discrete multidimensional transforms]\label{eg:DMT}\hspace{-8pt}%
  Three of the best-known discrete transforms are the discrete Fourier transform we
encountered in Example~\ref{eg:int}, the discrete Z-transform and the discrete cosine transform:
\begin{align*}
F(x_1,\dots,x_d) &= \sum_{k_1=-\infty}^\infty \cdots \sum_{k_d=-\infty}^\infty f(k_1,\dots,k_d) \, \rme^{-\rmi  k_1 x_1 - \cdots -\rmi  k_d x_d},\\
F(z_1,\ldots,z_d)&= \sum_{k_1=-\infty}^{\infty} \cdots \sum_{k_d=-\infty}^{\infty} f(k_1,\ldots,k_d) \, z_1^{-k_1}  \cdots z_d^{-k_d},\\
F(j_1,\ldots,j_d )&=\\*
&\hspace{-24pt} \sum_{k_1=0}^{n_1-1} \cdots \sum_{k_d=0}^{n_d-1} f(k_1,\ldots,k_d) \cos\biggl( \frac{ \pi (2j_1+1) j_1}{2n_1} \biggr)\cdots \cos \biggl(\frac{ \pi (2j_d+1) j_d}{2n_d}\biggr),
\end{align*}
where the $x_i$ are real, the $z_i$ are complex, and $j_i$ and $k_i$ are integer variables. We refer the reader to \citet[Sections~2.2 and 4.2]{MSP} for the first two and \citet[Chapter~5]{DCT} for the last. While we have stated them for general $d$, in practice $d =2,3$ are the most useful. The separability of these kernels is exploited in the \emph{row--column decomposition} for their evaluation \cite[Section~2.3.2]{MSP}. Assuming, for simplicity, that we have $d=2$, a kernel $K(j,k) = \varphi(j)\psi(k)$ and all integer variables, then 
\[
F(x,y)
=\sum_{j=-\infty}^{\infty} \sum_{k=-\infty}^{\infty} \varphi(j)\psi(k)f(x-j,y-k)
=\sum_{j=-\infty}^{\infty}\varphi(j)\Bigg[ \sum_{k=-\infty}^{\infty} \psi(k)f(x-j,y-k)\Bigg].
\]
We store the sum in the bracket, which we then re-use when evaluating $F$ at other points $(x',y)$ where only the first argument $x' = x +\delta$ is changed:
\begin{align*}
F(x+\delta,y)&=\sum_{j=-\infty}^{\infty}\varphi(j)\Bigg[ \sum_{k=-\infty}^{\infty} \psi(k)f(x-j+\delta,y-k)\Bigg]\\
&=\sum_{j=-\infty}^{\infty}\varphi(j+\delta)\Bigg[ \sum_{k=-\infty}^{\infty} \psi(k)f(x-j,y-k)\Bigg].
\end{align*}
We have assumed that the indices run over all integers to avoid having to deal with boundary complications. In reality, when we have a finite sum as in the discrete cosine transform, evaluating $F$ in the direct manner would have taken $n_1^2 n_2^2 \cdots n_d^2$ additions and multiplications, whereas the  row--column decomposition would just require
\[
n_1 n_2 \cdots n_d(n_1 + n_2 + \cdots + n_d)
\]
additions and multiplications. In cases where there are fast algorithms available for the one-dimensional transform, say, if we employ one-dimensional FFT in an evaluation of the $d$-dimensional DFT via row--column decomposition, the number of additions and multiplications could be further reduced to
\[
n_1 n_2 \cdots n_d \log_2(n_1 + n_2 + \cdots + n_d) \quad\text{and}\quad \frac{1}{2} n_1 n_2 \cdots n_d \log_2(n_1 + n_2 + \cdots + n_d)
\]
respectively. Supposing $d =2$ and $n_1 = n_2 = 2^{10}$, the approximate number of multiplications required to evaluate a two-dimensional DFT using the direct method, the row--column decomposition method and the row--column decomposition with FFT method are $10^{12}$, $2 \times 10^9$, and $10^7$ respectively \cite[Section~2.3.2]{MSP}.
\end{example} 

The fact that we may construct a tensor product basis $\mathscr{B}_1 \otimes \dots \otimes \mathscr{B}_d$ out of given bases $\mathscr{B}_1,\ldots,\mathscr{B}_d$ allows us to build multivariate bases out of univariate ones.  While this is not necessarily the best way to construct multivariate bases, it is the simplest and often serves as a basic standard that more sophisticated multivariate bases are compared against. This discussion is most interesting with Hilbert spaces, and this will be the backdrop for our next example, where we will see the construction for various notions of bases and beyond.

\begin{example}[tensor product wavelets and splines]\label{eg:wavelet}
  We let $\mathbb{H}$ be any separable Banach space
  with norm $\lVert \, \cdot \, \rVert$.
  Then $\mathscr{B} = \{\varphi_i \in \mathbb{H} \colon  i \in  \mathbb{N} \}$ is said to be a \emph{Schauder basis} of $\mathbb{H}$ if, for every $f \in \mathbb{H}$, there is a unique sequence $(a_i)_{i=1}^\infty$ such that
\[
f = \sum_{i=1}^\infty  a_i \varphi_i,
\]
where, as usual, this means the series on the right converges to $f$ in $\lVert \, \cdot \, \rVert$. A Banach space may not have a Schauder basis but a Hilbert space always does. If $\mathbb{H}$ is a Hilbert space with inner product $\langle \,\cdot\,,\cdot \, \rangle$ and $\lVert \, \cdot \, \rVert$ its induced norm, then its Schauder basis has specific names when it satisfies additional conditions;  two of the best-known ones are the orthonormal basis and the \emph{Riesz basis}. Obviously the elements of a Schauder basis must be linearly independent, but this can be unnecessarily restrictive since overcompleteness
can
be a desirable feature \cite{Mallat}, leading us to the notion of \emph{frames}. For easy reference, we define them as follows \cite{Heil}:
\begin{alignat*}{5}
&\text{orthonormal basis}\quad &\lVert f \rVert^2 &= \sum_{i=1}^\infty \lvert \langle f, \varphi_i \rangle \rvert^2, & & f\in \mathbb{H},\\*
&\text{Riesz basis} & \alpha \sum_{i=1}^\infty \lvert a_i \rvert^2  &\le \biggl\lVert \sum_{i=1}^\infty  a_i \varphi_i \biggr\rVert^2 \le \beta \sum_{i=1}^\infty \lvert a_i \rvert^2, \quad &&  (a_i)_{i=1}^\infty \in l^2(\mathbb{N}),\\*
&\text{frame} &  \alpha \lVert f \rVert^2 &\le  \sum_{i=1}^\infty \lvert \langle f, \varphi_i \rangle \rvert^2 \le \beta \lVert f \rVert^2, && f \in \mathbb{H},
\end{alignat*}
where the constants $0 < \alpha < \beta$ are called \emph{frame constants} and if $\alpha = \beta$, then the frame is \emph{tight}. Clearly every orthonormal basis is a Riesz basis and every Riesz basis is a frame. 

Let $\mathbb{H}_1,\ldots,\mathbb{H}_d$ be separable Hilbert spaces and let $\mathscr{B}_1,\ldots,\mathscr{B}_d$ be countable dense spanning sets. Let
$\mathbb{H}_1\hatotimes \dots \hatotimes \mathbb{H}_d$ be their Hilbert--Schmidt tensor product as discussed in Examples~\ref{eg:Hilbert} and \ref{eg:compact} and let $\mathscr{B}_1 \otimes \dots \otimes \mathscr{B}_d$ be as defined in \eqref{eq:tenprodbasis2}. Then
\[
\mathscr{B}_1,\ldots,\mathscr{B}_d \text{ are}
\begin{cases}
\text{orthonormal bases},\\
\text{Riesz bases},\\
\text{frames},
\end{cases}\hspace{-8pt}
\Leftrightarrow \ \ 
\mathscr{B}_1 \otimes \dots \otimes \mathscr{B}_d \text{ is}
\begin{cases}
\text{an orthonormal basis},\\
\text{a Riesz basis},\\
\text{a frame}.
\end{cases}
\]
The forward implication is straightforward, and if the frame constants of $\mathscr{B}_i$ are $\alpha_i$ and $\beta_i$, $i=1,\ldots,d$, then the frame constants of $\mathscr{B}_1 \otimes \dots \otimes \mathscr{B}_d$ are $\prod_{i=1}^d \alpha_i$ and $\prod_{i=1}^d \beta_i$ \cite[Lemma~8.18]{Feichtinger}. The converse is, however, more surprising \cite{Bourouihiya}.

If $\psi \in L^2(\mathbb{R})$ is a wavelet, \tie, $\mathscr{B}_\psi \coloneqq \{\psi_{m,n} \colon  (m,n)\in \mathbb{Z} \times \mathbb{Z} \}$ with $\psi_{m,n}(x) \coloneqq 2^{m/2}\psi(2^m x - n)$ is an orthonormal basis, Riesz basis or frame of $L^2(\mathbb{R})$, then
\begin{equation}\label{eq:tenprodbasis}
\mathscr{B}_\psi \otimes \mathscr{B}_{\psi} = \{\psi_{m,n} \otimes \psi_{p,q} \colon  (m,n,p,q) \in \mathbb{Z} \times \mathbb{Z} \times \mathbb{Z} \times \mathbb{Z}\}
\end{equation}
is an orthonormal basis, Riesz basis or frame of $L^2(\mathbb{R}^2)$ by our general discussion above. More generally, $\mathscr{B}_\psi^{\otimes d}$ gives a tensor product wavelet basis/frame for $L^2(\mathbb{R}^d)$ \cite[Section~7.7.4]{Mallat}. While this is certainly the most straightforward way to obtain a wavelet basis for multivariate functions out of univariate wavelets, there are often better options.

Some orthonormal wavelets have a \emph{multiresolution analysis}, \tie, a sequence of nested subspaces whose intersection is $\{0\}$ and union is dense in $L^2(\mathbb{R})$,
\[
\{0\} \subseteq \dots \subseteq \mathbb{V}_1 \subseteq \mathbb{V}_{0}\subseteq \mathbb{V}_{-1}\subseteq \mathbb{V}_{-2}\subseteq \dots \subseteq L^2(\mathbb{R}),
\]
defined by a \emph{scaling function} $\varphi \in L^2(\mathbb{R})$ so that $\{ \varphi_{m,n} \colon   n \in \mathbb{Z} \}$ is an orthonormal basis of $\mathbb{V}_m$ for any $m \in \mathbb{Z}$. If we write the orthogonal complement of $\mathbb{V}_m$ in $\mathbb{V}_{m-1}$ as $\mathbb{W}_m$, then  $\{\psi_{m,n} \colon  n\in \mathbb{Z} \}$ is an orthonormal basis of $\mathbb{W}_m$ \cite[Chapter~2]{Meyer}. Multiresolution analysis interacts with tensor products via the rules in Example~\ref{eg:calc}:
\[
\mathbb{V}_{m-1} \otimes \mathbb{V}_{m-1} = ( \mathbb{V}_m \otimes \mathbb{V}_m ) \oplus ( \mathbb{V}_m \otimes \mathbb{W}_m ) \oplus ( \mathbb{W}_m \otimes \mathbb{V}_m ) \oplus ( \mathbb{W}_m \otimes \mathbb{W}_m ).
\]
Here $\oplus$ denotes orthogonal direct sum. In image processing lingo, the subspaces
\begin{equation}\label{eq:multres}
\begin{aligned}
\mathbb{V}_m \otimes \mathbb{V}_m &=\overline{\spn}\{\varphi_{m,n} \otimes \varphi_{p,q} \colon  (n,q) \in \mathbb{Z} \times \mathbb{Z}\}, \\*
\mathbb{V}_m \otimes \mathbb{W}_m &=\overline{\spn}\{\varphi_{m,n} \otimes \psi_{p,q} \colon  (n,q) \in \mathbb{Z} \times \mathbb{Z}\}, \\
\mathbb{W}_m \otimes \mathbb{V}_m &=\overline{\spn}\{\psi_{m,n} \otimes \varphi_{p,q} \colon  (n,q) \in \mathbb{Z} \times \mathbb{Z} \}, \\
\mathbb{W}_m \otimes \mathbb{W}_m &=\overline{\spn}\{\psi_{m,n} \otimes \psi_{p,q} \colon  (n,q) \in \mathbb{Z} \times \mathbb{Z} \}
\end{aligned}
\end{equation}
are called the LL-, LH-, HL- and HH-subbands respectively, with L for low-pass and H for high-pass. For an image $f \in \mathbb{V}_{m-1} \otimes \mathbb{V}_{m-1}$, its component in the LL-subband $\mathbb{V}_m \otimes \mathbb{V}_m$ represents a low-resolution approximation to $f$; the high-resolution details are contained in its components in the LH-, HL- and HH-subbands, whose direct sum is sometimes called the detailed space.

The orthonormal basis obtained via \eqref{eq:multres} can be a better option than simply taking the tensor product \eqref{eq:tenprodbasis}. For instance, in evaluating an $n \times n$ dense matrix--vector product
  representing an integral transform like the ones in Example~\ref{eg:multpole},
  the \citet{Beylkin} wavelet-variant of the fast multipole algorithm mentioned in Example~\ref{eg:Krylov} runs in time $O(n \log n)$ using the basis in \eqref{eq:tenprodbasis} and in time $O(n)$ using that in \eqref{eq:multres}, both impressive compared to the usual $O(n^2)$, but the latter is clearly superior when $n$ is very large.

The main advantage of \eqref{eq:tenprodbasis} is its generality; with other types of bases or frames we also have similar constructions. For instance, the simplest types of multivariate B-splines \cite{Hollig} are constructed out of tensor products of univariate B-splines; expressed in a multilinear rank decomposition \eqref{eq:multdecomp}, we have
\begin{equation}\label{eq:spline}
B_{l,m,n}(x,y,z) = \sum_{i=1}^p \sum_{j=1}^q \sum_{k=1}^r a_{ijk} B_{i,l}(x) B_{j,m}(y) B_{k,n}(z),
\end{equation}
where $B_{i,l}, B_{j,m}, B_{k,n}$ are univariate B-splines of degrees $l,m,n$. Nevertheless,
the main drawback of
straightforward tensor product constructions
is that they attach
undue importance to the directions of the coordinate axes \cite{Cohen}. There are often better alternatives such as box splines or beamlets, curvelets, ridgelets, shearlets, wedgelets, {\em etc.}, that exploit the geometry of $\mathbb{R}^2$ or $\mathbb{R}^3$.
\end{example}

We next discuss a covariant counterpart to the contravariant example above: a univariate quadrature on $[-1,1]$ is a covariant $1$-tensor, and a multivariate quadrature on $[-1,1]^d$ is a covariant $d$-tensor.

\begin{example}[quadrature]\label{eg:quadrature}
Let $f \colon  [-1,1] \to \mathbb{R}$ be a univariate polynomial function of degree not more than $n$. Without knowing anything else about $f$, we know that there exist $n$ distinct points $x_0,x_1,\ldots, x_n \in [-1,1]$, called \emph{nodes}, and coefficients $w_0,w_1,\ldots,w_n \in \mathbb{R}$, called \emph{weights}, so that
\[
\int_{-1}^1 f(x)  \D x = w_0 f(x_0) + w_1 f(x_1) + \dots + w_n f(x_n).
\]
In fact, since $f$ is arbitrary, the formula holds with the same nodes and weights for all polynomials of degree $n$ or less. This is called a \emph{quadrature} formula and its existence is simply a consequence of the following observations.
A
definite integral is a linear functional,
\begin{equation*}
  I \colon  C([-1,1]) \to \mathbb{R}, \quad I(f) = \int_{-1}^1 f(x)  \D x,
  \end{equation*}
as is point evaluation, introduced in Example~\ref{eg:multmult},
\begin{equation*}
  \varepsilon_x \colon  C([-1,1]) \to \mathbb{R},\quad \varepsilon_x(f) = f(x).
\end{equation*}
Since $\mathbb{V} = \{ f \in C([-1,1]) \colon  f(x) = c_0 + c_1 x + \dots + c_n x^n \}$ is a vector space of dimension $n+1$, 
its dual space $\mathbb{V}^*$ also has dimension $n+1$, and the $n+1$ linear functionals $\varepsilon_{x_0},\varepsilon_{x_1},\ldots, \varepsilon_{x_n}$, obviously linearly independent, form a basis of $\mathbb{V}^*$. Thus it must be possible to write $I$ as a unique linear combination
\begin{equation}\label{eq:quadrature2}
I = w_0 \varepsilon_{x_0} + w_1 \varepsilon_{x_1} + \dots + w_n \varepsilon_{x_n}.
\end{equation}
This observation is presented in \citet[Chapter~2, Theorem~7]{Lax} as an instance in linear algebra where `even trivial-looking material has interesting consequences'.

More generally, a quadrature formula is simply a linear functional in $\mathbb{V}^*$ that is a linear combination of point evaluation functionals, as on the right-hand side of \eqref{eq:quadrature2}. If $\mathbb{V} \subseteq C(X)$ is any $(n+1)$-dimensional subspace and $X$  is compact Hausdorff with a Borel measure, then the same argument in the previous paragraph shows that there is a quadrature for $I = \int_X$ that gives the exact value of the integral for all $ f \in \mathbb{V}$ in terms of the values of $f$ at $n+1$ nodes $x_0,\ldots,x_n \in X$. However, this abstract existential argument  does not make it any easier to find $x_0,\ldots,x_n \in X$ for a specific type of function space $\mathbb{V}$ and a specific domain $X$, nor does it tell us about the error if we apply the formula to some $f \notin \mathbb{V}$.  In fact, for specific  $\mathbb{V}$ and $X$ we can often do better. For instance, the three-node Gauss quadrature on $X=[-1,1]$ has $x_0 =0$, $x_1 =\sqrt{3/5} = -x_2$ and $w_0 = 8/9$, $w_1 = 5/9 = w_2$, with
\[
\int_{-1}^1 f(x)  \D x = \frac{8}{9} f(0) + \frac{5}{9} f\biggl(\sqrt{\frac{3}{5}}\biggr) +  \frac{5}{9} f\biggl(-\sqrt{\frac{3}{5}}\biggr)
\]
giving exact answers for all \emph{quintic} polynomials $f$. So $\dim\mathbb{V}=6$ but we only need three point evaluations to determine $I$. There is no contradiction here: we are not saying that every linear functional in $\mathbb{V}^*$ can be determined with these three point evaluations, just this one linear functional $I$. This extra efficacy comes from exploiting the structure of $\mathbb{V}$ as a space of polynomials;  while the moment functionals $f \mapsto \int_{-1}^1 x^n f(x)  \D x$ are linearly independent, they satisfy non-linear relations \cite{Schmudgen} that can be exploited to obtain various $n$-node quadratures with exact answers for all polynomials up to degree $2n-1$ \cite{Golub}. Generally speaking, the more nodes there are in a quadrature, the more accurate it is but also the more expensive it is. As such we will take the number of nodes as a crude measure of accuracy and computational cost.

As univariate quadratures are well studied in theory \cite{Brass} and practice \cite{Davis}, it is natural to build multivariate quadratures from univariate ones. We assume without too much loss of generality that $X = [-1,1]^d$, noting that any $[a_1,b_1] \times \dots \times [a_d,b_d]$ can be transformed into $[-1,1]^d$ with an affine change of variables. Since a univariate quadrature is a linear functional $\varphi \colon  C([-1,1]) \to \mathbb{R}$, it seems natural to define a multivariate quadrature as a $d$-linear functional
\begin{equation}\label{eq:multquad0}
\Phi \colon  C([-1,1]) \times \dots \times C([-1,1]) \to \mathbb{R},
\end{equation}
but this does not appear to work as we want multivariate quadratures to be defined on $C([-1,1]^d)$. Fortunately the universal factorization property \eqref{eq:commdiag} gives us a unique linear map
\[
F_\Phi \colon  C([-1,1]) \otimes \dots \otimes C([-1,1]) \to \mathbb{R}
\]
with $\Phi = F_\Phi \circ \seg$, and as $C([-1,1]^d)=C([-1,1])^{\otimes d}$ by \eqref{eq:ctsint}, we get a linear functional
\begin{equation}\label{eq:multquad3}
F_\Phi \colon  C([-1,1]^d) \to \mathbb{R}.
\end{equation}
A multivariate quadrature is a linear functional $F \colon  C([-1,1]^d) \to \mathbb{R}$ that is a linear combination of point evaluation functionals on $[-1,1]^d$. Note that given $F$ we may also work backwards to obtain a multilinear functional $\Phi_F = F \circ \seg$ as in \eqref{eq:multquad0}. So \eqref{eq:multquad0} and \eqref{eq:multquad3} are equivalent.

The simplest multivariate quadratures are simply separable products of univariate quadratures. Let
\[
\varphi_1 = \sum_{i=1}^{n_1} a_i \varepsilon_{x_i},\;
\varphi_2 = \sum_{j=1}^{n_2} b_j \varepsilon_{y_j},\ldots,
\varphi_d = \sum_{k=1}^{n_d} c_k \varepsilon_{z_k}
\]
be $d$ univariate quadratures with nodes $N_i \subseteq [-1,1]$ of size $n_i$, $i =1,\ldots,d$. Then the covariant $d$-tensor
\[
\varphi_1 \otimes \varphi_2 \otimes \dots \otimes \varphi_d = \sum_{i=1}^{n_1} \sum_{j=1}^{n_2} \cdots \sum_{k=1}^{n_d}  a_i b_j \cdots c_k \, \varepsilon_{x_i} \otimes \varepsilon_{y_j} \otimes \dots \otimes \varepsilon_{z_k}
\]
is a multivariate quadrature on $[-1,1]^d$ usually called the tensor product quadrature.
The downside of such a quadrature is immediate: for an $f \in C([-1,1]^d)$, evaluating $\varphi_1 \otimes \dots \otimes \varphi_d(f)$ requires $f$ to be evaluated on $n_1n_2\cdots n_d$ nodes, but for $d=2,3$ this is still viable. For instance, the nine-node tensor product Gauss quadrature on $[-1,1] \times [-1,1]$ is the covariant $2$-tensor
\[
\biggl( \frac{5}{9} \varepsilon_{\scriptscriptstyle-\sqrt{3/5}} + \frac{8}{9} \varepsilon_0 + \frac{5}{9} \varepsilon_{\scriptscriptstyle\sqrt{3/5}} \biggr)^{\otimes 2},
\]
with nodes
\[
\sqrt{\frac{3}{5}} \biggl\{
\begin{bmatrix}
-1\\
-1
\end{bmatrix}\!,
\begin{bmatrix}
0\\
-1
\end{bmatrix}\!,
\begin{bmatrix}
1\\
-1
\end{bmatrix}\!,
\begin{bmatrix}
-1\\
0
\end{bmatrix}\!,
\begin{bmatrix}
0\\
0
\end{bmatrix}\!,
\begin{bmatrix}
1\\
0
\end{bmatrix}\!,
\begin{bmatrix}
-1\\
1
\end{bmatrix}\!,
\begin{bmatrix}
0\\
1
\end{bmatrix}\!,
\begin{bmatrix}
1\\
1
\end{bmatrix}
\biggr\},
\]
and corresponding weights
\[
\biggl\{\frac{25}{81},\frac{40}{81},\frac{25}{81},\frac{40}{81},\frac{64}{81},\frac{40}{81},\frac{25}{81},\frac{40}{81},\frac{25}{81} \biggr\}.
\]
One observation is that the weights in the tensor product quadrature are always positive whenever those in its univariate factors $\varphi_1,\ldots,\varphi_d$ are. So tensor product quadrature is essentially just a form of weighted average. Taking a leaf from Strassen's algorithm in Example~\ref{eg:Strassen}, allowing for negative coefficients in a tensor decomposition can be a starting point for superior algorithms; we will examine this in the context of quadrature.

The Smolyak quadrature \cite{Smolyak} is a more sophisticated multivariate quadrature that is also based on tensor products. Our description is adapted from \citet{Kaarnioja} and \citet{Kirby}. As before, we are given univariate quadratures $\varphi_i$ with nodes $N_i \subseteq [-1,1]$ of cardinality $n_i$, $i =1,\ldots,d$. An important difference is that we now order $\varphi_1,\varphi_2,\ldots,\varphi_d$ so that $n_1 \le n_2 \le \dots \le n_d$, \tie, these univariate quadratures get increasingly accurate but also increasingly expensive as $i$ increases. A common scenario would have $n_i = 2^{i-1}$. The goal is to define a tensor product quadrature that uses the cheaper quadratures on the left end of the list $\varphi_1, \varphi_2,\ldots$ liberally and the expensive ones on the right end $\ldots, \varphi_{d-1},\varphi_d$ sparingly. We set $\varphi_0 \coloneqq 0$,
\[
\Delta_i \coloneqq \varphi_i - \varphi_{i-1}, \quad i =1,\ldots,d,
\]
and observe that
\[
\varphi_i = \Delta_1 + \Delta_2 + \dots + \Delta_i, \quad i=1,\ldots,d.
\]
The separable multivariate quadrature we defined earlier is
\begin{align*}
\varphi_1 \otimes \varphi_2 \otimes \dots \otimes \varphi_d &= \Delta_1 \otimes (\Delta_1 + \Delta_2) \otimes \dots \otimes (\Delta_1 + \Delta_2 + \dots + \Delta_d);
\intertext{it could be made less expensive and accurate by replacing the $\varphi_i$ with $\varphi_1$,}
\varphi_1 \otimes \varphi_1 \otimes \dots \otimes \varphi_1 &= \Delta_1 \otimes \Delta_1 \otimes \dots \otimes \Delta_1,
\intertext{or more expensive and accurate by replacing the $\varphi_i$ with $\varphi_d$,}
\varphi_d \otimes \varphi_d \otimes \dots \otimes \varphi_d &= (\Delta_1 + \Delta_2 + \dots + \Delta_d)\otimes \dots \otimes (\Delta_1 + \Delta_2 + \dots + \Delta_d).
\end{align*}
The Smolyak quadrature allows one to adjust the level of accuracy and computational cost between these two extremes. 
For any $r \in \mathbb{N}$, the level-$r$ Smolyak quadrature is the covariant $d$-tensor
\begin{align*}
S_{r,d} &= \sum_{k=0}^{r-1} \sum_{i_1 + \dots + i_d = d + k} \Delta_{i_1} \otimes \dots \otimes \Delta_{i_d} \\*
&=\sum_{k = r -d}^{r-1} (-1)^{r-1-k} \binom{d-1}{r-1-k} \sum_{i_1 + \dots + i_d = d + k} \varphi_{i_1} \otimes \dots \otimes \varphi_{i_d}.
\end{align*}
Note that when $r=1$ and $(d-1)d$, we get
\[
S_{1,d} = \Delta_1^{\otimes d}=\varphi_1^{\otimes d}, \quad S_{(d-1)d,d} = (\Delta_1 + \Delta_2 + \dots + \Delta_d)^{\otimes d} =\varphi_d^{\otimes d}
\]
respectively. Take the smallest non-trivial case with $d =r =2$; we have $\Delta_1 = \varphi_1$, $\Delta_2 = \varphi_2 - \varphi_1$ and
\begin{align*}
S_{2,2} &= \Delta_1 \otimes \Delta_1 + \Delta_1 \otimes \Delta_2 + \Delta_2 \otimes \Delta_1 \\*
&
= \varphi_1 \otimes \varphi_2 + \varphi_2 \otimes \varphi_1 - \varphi_1 \otimes \varphi_1.
\end{align*}
Observe that the weights can now be negative. For the one-point and two-point Gauss quadratures $\varphi_1 = 2 \varepsilon_0$ and $\varphi_2 =  \varepsilon_\mrt + \varepsilon_\prt$, we have
\begin{equation}\label{eq:S22}
S_{2,2} (f) = 2 f\biggl(0,\frac{-1}{\sqrt{3}}\biggr) + 2 f\biggl(0,\frac{1}{\sqrt{3}}\biggr) + 2 f\biggl(\frac{-1}{\sqrt{3}},0\biggr)  + 2f\biggl(\frac{1}{\sqrt{3}},0\biggr) - 4f(0,0)
\end{equation}
for any $f \in C([-1,1]^2)$. We will mention just one feature of Smolyak quadrature: if the univariate quadrature $\varphi_i$ is exact for univariate polynomials of degree $2i-1$, $i=1,\ldots,d$, then $S_{r,d}$ is exact for $d$-variate polynomials of degree $2r-1$ \cite{Heiss}. Thus \eqref{eq:S22} gives the exact answers when integrating bivariate cubic polynomials $f$.

For any $f \in C([-1,1]^d)$, $S_{r,d}(f)$ has to be evaluated at the nodes in
\[
\bigcup_{r \le i_1 + \dots + i_d \le r +d-1} N_{i_1} \times \dots \times N_{i_d},
\]
which is called a \emph{sparse grid} \cite{Gerstner}. Smolyak quadrature has led to further developments under the heading of sparse grids; we refer interested readers to \citet{Garcke} for more information. Nevertheless, many multivariate quadratures tend to be `grid-free' and involve elements of randomness or pseudorandomness, essentially \eqref{eq:quadrature2} with the nodes $x_0,\ldots,x_n$ chosen (pseudo)randomly in $[-1,1]^d$ to have low discrepancy, as defined in \citet{Niederreiter}. As a result, modern study of multivariate quadratures tends to be quite different from the classical approach described above.
\end{example}

We revisit a point that we made in Example~\ref{eg:hypmatrep}: it is generally not a good idea to identify tensors in $\mathbb{V}_1\otimes \dots \otimes \mathbb{V}_d$ with hypermatrices in $\mathbb{R}^{n_1 \times \dots \times n_d}$, as the bases often carry crucial information, \tie, in a decomposition such as
\[
T = \sum_{i=1}^p \sum_{j=1}^q \cdots \sum_{k=1}^r a_{ij\cdots k} \, e_i \otimes f_j \otimes \dots \otimes g_k,
\]
the basis elements $e_i, f_j,\ldots,g_k$ can be far more important than the coefficients $a_{ij\cdots k}$. Quadrature provides a fitting example. The basis elements are the point evaluation functionals at the nodes and these are key to the quadrature: once they are chosen, the coefficients, \ie\ the weights, can be determined almost as an afterthought. Furthermore, while these basis elements are by definition \emph{linearly} independent,
in the context of quadrature
there are \emph{non-linear} relations between them that will be lost if one just looks at the coefficients.

All our examples in this article have focused on computations in the classical sense; we will end with a quantum computing example, adapted from \citet[Chapter~7]{Nakahara}, \citet[Chapter~6]{Nielsen} and \citet[Section~2.3]{Wallach}. While tensors appear at the beginning, ultimately Grover's quantum search algorithm \cite{Grover} reduces to the good old power method.

\begin{example}[Grover's quantum search]\label{eg:grover}
Let $\mathbb{C}^2$ be equipped with its usual Hermitian inner product $\langle x,y\rangle = x^* y$ and let $e_0, e_1$ be any pair of orthonormal vectors in $\mathbb{C}^2$. Let $d \in \mathbb{N}$ and  $n = 2^d$. Suppose we have a function
\begin{equation}\label{eq:groverf}
f \colon  \{0,1,\ldots,n-1\} \to \{-1,+1\}, \quad
f(i) =\begin{cases}
-1 & \text{if } i = j,\\
+1 & \text{if } i \ne j,
\end{cases}
\end{equation}
\tie, $f$ takes the value $-1$ for exactly one $j \in \{0,1,\ldots,n-1\}$ but we do not know which $j$. Also,  $f$ is only accessible as a black box;  we may evaluate $f(i)$ for any given $i$ and so to find $j$ in this manner would require $O(n)$ evaluations in the worst case. Grover's algorithm finds $j$ in $O(\sqrt{n})$ complexity with a quantum computer.

We  begin by observing that the $d$-tensor $u$ below may be expanded as
\begin{equation}\label{eq:u}
u \coloneqq \biggl[\dfrac{1}{\sqrt{2}} (e_0 + e_1) \biggr]^{\otimes d} = \dfrac{1}{\sqrt{n}} \sum_{i=0}^{n-1} u_i \in (\mathbb{C}^2)^{\otimes d} ,
\end{equation}
with
\[
u_i \coloneqq e_{i_1} \otimes e_{i_2}  \otimes \dots \otimes e_{i_d},
\]
where $i_1,\ldots,i_d \in \{0,1\}$ are given by expressing the integer $i$ in binary:
\[
[i]_2 = i_d i _{d-1} \cdots i_2 i_1.
\]
Furthermore, with respect to the inner product on $(\mathbb{C}^2)^{\otimes d}$ given by \eqref{eq:innprod1} and its induced norm, $u$ is of unit norm and $\{u_0,\ldots,u_{n-1}\}$ is an orthonormal basis. Recall the notion of a Householder reflector in Example~\ref{eg:GHG}. We define two Householder reflectors $H_f$ and $H_u \colon  (\mathbb{C}^2)^{\otimes d} \to (\mathbb{C}^2)^{\otimes d}$,
\[
H_f (v) = v - 2\langle u_j, v \rangle u_j, \quad H_u(v) = v - 2 \langle u, v\rangle u,
\]
reflecting $v$ about the hyperplane orthogonal to $u_j$ and $u$ respectively. Like the function $f$, the Householder reflectors $H_f$ and $H_u$ are only accessible as black boxes. Given any $v \in (\mathbb{C}^2)^{\otimes d}$, a quantum computer allows us to evaluate $H_f(v)$ and $H_u(v)$ in logarithmic complexity $O(\log n) = O(d)$ \cite[p.~251]{Nielsen}.

Grover's algorithm is essentially the power method with $H_u H_f$ applied to $u$ as the initial vector. Note that $H_u H_f$ is unitary and thus norm-preserving and we may skip the normalization step. If we write
\[
(H_u H_f)^k u = a_k u_j + b_k \sum_{i \ne j} u_i,
\]
then we may show \cite[Propositions~7.2--7.4]{Nakahara} that
\[
a_0 = b_0 = \dfrac{1}{\sqrt{n}}, \quad
\begin{cases}
a_k = \dfrac{2-n}{n}a_{k-1} + \dfrac{2(1-n)}{n} b_{k-1},\\[7pt]
b_k = \dfrac{2}{n}a_{k-1} + \dfrac{2-n}{n} b_{k-1},
\end{cases}
\]
for $k \in \mathbb{N}$, with the solution given by
\[
a_k = -\sin(2k+1)\theta, \quad b_k = \dfrac{-1}{\sqrt{n-1}} \cos (2k+1)\theta, \quad \theta \coloneqq \sin^{-1} (1/\sqrt{n}).
\]
The probability of a quantum measurement giving the correct answer $j$ is given by
\[
\lvert a_k \rvert^2 = \sin^2 (2k+1)\theta,
\]
and for $n \gg 1$ this is greater than $1-1/n$ when $k \approx \pi\sqrt{n}/4$. This algorithm can be extended to search an $f$ in \eqref{eq:groverf} with $m$ copies of $-1$s \cite[Section~7.2]{Nakahara}.
\end{example}

Observe that after setting up the stage with tensors, Grover's algorithm deals only with \emph{linear} operators on tensor spaces. It is the same with Shor's quantum algorithm for integer factorization and discrete log \cite{Shor}, which also relies on the tensor $u$ in \eqref{eq:u} but uses DFT instead of Householder reflectors, although the details \cite[Chapter~8]{Nakahara} are too involved to go into here.  In either case there is no computation, classical or quantum, involving higher-order tensors.

\section{Odds and ends}

We conclude our article by recording some parting thoughts. Section~\ref{sec:omit} is about the major topics that had to be omitted. Sections~\ref{sec:comp} and \ref{sec:hypana} may be viewed respectively as explanations as to why our article is not written like a higher-order analogue of \citet{GVL} or \citet{Horn1}.

\subsection{Omissions}\label{sec:omit}

There are three major topics in our original plans for this article that we did not manage to cover: (i) symmetric and alternating tensors, (ii)~algorithms based on tensor contractions, and (iii)~provably correct algorithms for higher-order tensor problems.

In (i) we would have included a discussion of the following: polynomials and differential forms as coordinate representations of symmetric and alternating tensors; polynomial kernels as Veronese embedding of symmetric tensors, Slater determinants as Pl\"ucker embedding of alternating tensors; moments and sums-of-squares theory for symmetric tensors; matrix--matrix multiplication as a special case of the Cauchy--Binet formula for alternating tensors.
We also planned to discuss the symmetric and alternating analogues of the three definitions of a tensor, the three constructions of tensor products and the tensor algebra. The last in particular gives us the symmetric and alternating Fock spaces with their connections to bosons and fermions. We would also have discussed various notions of symmetric and alternating tensor ranks.

In (ii) we had intended to demonstrate how the Cooley--Tukey FFT,
the multidimensional DFT,
the Walsh transform, wavelet packet transform, Yates' method in factorial designs,
even FFT on finite non-Abelian groups,
{\em etc.}, are all tensor contractions. We would also show how the Strassen tensor in Example~\ref{eg:Strassen},
the tensors corresponding to the multiplication operations in Lie, Jordan and Clifford algebras,
and any tensor network in Example~\ref{eg:tennet}  may each be realized as the self-contraction of a  rank-one tensor. There would also be an explanation of matchgate tensors and holographic algorithms tailored for a numerical linear algebra readership.

In (iii) we had planned to provide a reasonably complete overview of the handful of provably correct algorithms for various NP-hard problems involving higher-order tensors
from three different communities: polynomial optimization, symbolic computing/computer algebra and theoretical computer science.
We will have more to say about these below.

\subsection{Numerical multilinear algebra?}\label{sec:comp}

This article is mostly about how tensors arise in computations and not so much about how one might compute tensorial quantities. Could one perhaps develop algorithms analogous to those in numerical linear algebra \cite{Demmel,GVL,HighamBook,Trefethen} but for higher-order tensors, or more accurately, for hypermatrices since most models of computations  require that we work in terms of coordinates? To be more precise, given a hypermatrix $A \in \mathbb{R}^{n_1 \times \dots \times n_d}$ with $d > 2$, can we compute its ranks, norms, determinants, eigenvalues and eigenvectors, various decompositions and approximations,
solutions to multilinear systems of equations, multilinear programming and multilinear regression with coefficients given by $A$, {\em etc.}?
The answer we think is no. One reason is that most of these problems are NP-hard \cite{HL}. Another related reason -- but more of an observation -- is that it seems difficult to construct \emph{provably correct} algorithms when $d > 2$.

If we sort a list of numbers in Excel or compute the QR factorization of a matrix in MATLAB, we count on the program to give us the correct answer: a sorted list or a QR factorization up to some rounding error. Such is the nature of algorithms: producing the correct answer or an approximation with a controllable error to the problem under consideration is an integral part of its definition \cite[Chapter~1]{Skiena}. An algorithm needs to have a proof of correctness or, better, a certificate of correctness that the user may employ to check the output. In numerical computing, an example of the former would be the Eckart--Young theorem, which guarantees that SVD produces a best rank-$r$ approximation, or the Lax equivalence theorem, which guarantees convergence of consistent finite difference schemes; an example of the latter would be relative backward error for a matrix decomposition or duality gap for a convex optimization problem having strong duality.

Note that we are not asking for efficient algorithms:  exponential running time or numerical instability are acceptable. The basic requirement is just that for some non-trivial set of inputs, the output
in exact arithmetic
can be guaranteed to be the answer we seek up to some quantifiable error. While such algorithms for tensor problems do exist in the literature -- notably the line of work in \citet{BCMT} and \citet{Nie} that extended \citet{Sylvester} and \citet{Reznick} to give a provably correct algorithm for symmetric tensor decomposition -- they tend to be the exception rather than the rule. It is more common to find `algorithms' for tensor problems that are demonstrably wrong.

The scarcity of provably correct algorithms goes hand in hand with the NP-hardness of tensor problems. For illustration, we consider the best rank-one approximation problem for a hypermatrix $A \in \mathbb{R}^{n \times n \times n}$:
\begin{equation}\label{eq:approx2}
\min_{x,y,z\in \mathbb{R}^n}\; \lVert A - x \otimes y \otimes z \rVert_\F.
\end{equation}
This is particularly pertinent as every tensor approximation problem contains \eqref{eq:approx2} as a special case; all notions of rank  --  tensor rank, multilinear rank, tensor network rank, slice rank, {\em etc.}\   --  agree on the non-zero lowest-ranked tensor and \eqref{eq:approx2} represents the
lowest order beyond $d=2$.
If we are given a purported solution $(x,y,z)$, there is no straightforward way to ascertain that it is indeed a solution unless $A$ happens to be rank-one and equal to $x \otimes y \otimes z$. One way to put this is that \eqref{eq:approx2} is NP-hard but not NP-complete:  given a candidate solution to an NP-complete problem, one can check whether it is indeed a solution in polynomial time; an NP-hard problem does not have this feature. In fact, any graph $G$ can be encoded as a hypermatrix $A_G \in \mathbb{R}^{n \times n \times n}$ in such a way that finding a best rank-one approximation to $A_G$ gives us the chromatic number of the graph \cite{HL}, which is NP-hard even to \emph{approximate} \cite{KLS}.

In case it is not clear, applying general-purpose methods such as gradient descent or alternating least-squares to an NP-hard optimization problem does not count as an algorithm. With few exceptions, an algorithm for a problem, say, the Golub--Kahan bidiagonalization for SVD, must exploit the special structure of the problem. Even though one would not expect gradient descent to yield SVD, it is not uncommon to find `algorithms' for tensor problems that do  no more than apply off-the-shelf non-linear optimization methods (\eg\  steepest descent, Newton's method) or coordinate cycling methods (\eg\  non-linear Jacobi or Gauss--Seidel). 

Furthermore, it is a misconception to think that it is only NP-hard to find global optimizers, and that if one just needs a local optimizer then NP-hardness is no obstacle. In fact, finding local optimizers \cite[Theorem~5.1]{Vavasis} or even just stationary points \cite{Khachiyan2}
may well be NP-hard problems too. It is also a fallacy to assume that since NP-hardness is an asymptotic notion, one should be fine with `real-world' instances of moderate dimensions.
As David S.\ Johnson, co-author of the standard compendium \cite{GJ}, puts it \cite{Johnson}:
\begin{quote}
What does it mean for practitioners? I believe that the years have taught us to take the warnings of NP-completeness seriously. If an optimization problem is NP-hard, it is rare that we find algorithms that, even when restricted to `real-world' instances, always seem to find optimal solutions, and do so in empirical polynomial time.
\end{quote} 
A pertinent example for us is the Strassen tensor $\Mu_{3,3,3}$ in Example~\ref{eg:Strassen} with $m=n=p=3$. With respect to the standard basis on $\mathbb{R}^{3 \times 3}$, it is a simple $9 \times 9 \times 9$ hypermatrix with mostly zeros and a few ones, but its tensor rank is still unknown today.

\subsection{Hypermatrix analysis?}\label{sec:hypana}

Readers conversant with matrix analysis \cite{Bellman,Bhatia,Horn1,Horn2}  or applied linear algebra \cite{Boyd2,Dym,Lax,Strang} would certainly have noticed that \emph{abstract linear algebra}  --  vector spaces, dual spaces, bases, dual bases, change of basis, {\em etc.}  --  can be largely or, with some effort, even completely avoided. Essentially, one just needs to know how to work with matrices. In fact, pared down to an absolute minimum, there is just one prerequisite: how to add and multiply matrices.

There is a high\-brow partial explanation for this. The category of matrices $\mathsf{Mat}$ and the category of vector spaces $\mathsf{Vec}$ are equivalent:  abstract linear algebra and concrete linear algebra are by and large the same thing \cite[Corollary~1.5.11]{Riehl}. The analogous statement is false for hypermatrices and tensor spaces except when interpreted in ways that reduce it to variations of the $\mathsf{Mat}$--$\mathsf{Vec}$ equivalence
as in Example~\ref{eg:homult2}. A main reason is that it is impossible to multiply hypermatrices in the manner one multiplies matrices, as we have discussed in Example~\ref{eg:homult}. The reader may notice that if we remove everything that involves matrix multiplication from \citet{Bellman,Bhatia,Horn1,Horn2}  or \citet{Boyd2,Dym,Lax,Strang}, there would hardly be any content left.

Another big difference in going to order $d \ge 3$ is the lack of canonical forms that we mentioned in Section~\ref{sec:trans}. The $2$-tensor transformation rules, \ie\ equivalence, similarity and congruence, allow us to transform matrices to various canonical forms: rational, Jordan, Weyr, Smith, Turnbull--Aitken, Hodge--Pedoe, {\em etc.}, or Schur and singular value decompositions with orthogonal/unitary similarity/congruence and equivalence. These play an outsize role in matrix analysis and applied linear algebra. For $d \ge 3$, while there is the Kronecker--Weierstrass form for $m \times n \times 2$ hypermatrices, this is the only exception; there is no canonical form for $m \times n \times p$ hypermatrices, even if we fix $p=3$, and therefore none for $n_1 \times \dots \times n_d$ hypermatrices since these include $m \times n \times 3$ as a special case by setting $n_4 = \dots = n_d =1$. Note that we did not say `no \emph{known}': the `no'  above is in the sense of `mathematically proven to be non-existent'.

Tensor rank appears to give us a handle -- after all, writing a matrix $A$ as a sum of $\rank(A)$ matrices, each of rank one, is essentially the Smith normal form -- but this is an illusion. The $\Mu_{3,3,3}$ example mentioned at the end of Section~\ref{sec:comp} is one among the hundreds of things we do not know about tensor rank. In fact, $4 \times 4 \times 4$ hypermatrices appear to be more or less the boundary of our current knowledge. Two of the biggest breakthroughs in recent times -- finding the border rank of $\Mu_{2,2,2} \in \mathbb{C}^{4 \times 4 \times 4}$ \cite{LandsbergJAMS} and finding the equations that define border rank-$4$ tensors $\{ A \in \mathbb{C}^{4 \times 4 \times 4} \colon \overline{\rank}(A) \le 4\}$ \cite{Friedland1} -- were both about $4 \times 4 \times 4$ hypermatrices and both required Herculean efforts, notwithstanding the fact that border rank is a simpler notion than rank.

\section*{Acknowledgement}

Most of what I learned about tensors in the last ten years I learned from my two former postdocs Yang Qi and Ke Ye; I gratefully acknowledge the education they have given me. Aside from them, I would like to thank Keith Conrad, Zhen Dai, Shmuel Friedland, Edinah Gnang, Shenglong Hu, Risi Kondor, and J.~M. Landsberg, Visu Makam, Peter McCullagh, Emily Riehl and Thomas Schultz for answering my questions and/or helpful discussions. I would like to express heartfelt gratitude to Arieh Iserles for his kind encouragement and Glennis Starling for her excellent copy-editing. Figure~\ref{fig:sepfun} is adapted from Walmes~M. Zeviani's gallery of TikZ art and Figure~\ref{fig:svm} is reproduced from Yifan Peng's blog.
  This work is partially supported by DARPA HR00112190040, NSF DMS-11854831, and the Eckhardt Faculty Fund.

This article was written while taking refuge from Covid-19 at my parents' home in Singapore. After two decades in the US, it is good to be reminded of how fortunate I am to have such wonderful parents. I dedicate this article to them: my father Lim Pok Beng and my mother Ong Aik Kuan.

\label{lastpage}

\end{document}